\documentclass[11pt,reqno]{amsart}
\usepackage[T1]{fontenc}
\usepackage[utf8]{inputenc}
\usepackage{lmodern}
\usepackage[margin=1.1in]{geometry}
\usepackage[expansion=false]{microtype}
\usepackage{amsmath,amssymb,mathtools}
\usepackage{enumitem,booktabs,array}
\usepackage{placeins}
\usepackage{xcolor}
\definecolor{DiagramPrimary}{gray}{0.15}
\definecolor{DiagramSecondary}{gray}{0.45}
\usepackage{tikz}
\usetikzlibrary{arrows.meta,calc,positioning}
\usepackage{hyperref}
\usepackage{bookmark}
\hypersetup{hidelinks,pdftitle={Canonical Cauchy bases on Ahlfors regular sets},pdfauthor={Mahishanka Withanachchi},pdfsubject={Harmonic and complex analysis}}
\setlist[enumerate]{leftmargin=1.8em,itemsep=3pt,topsep=4pt}
\newtheorem{theorem}{Theorem}
\newtheorem{lemma}{Lemma}[section]
\newtheorem{proposition}[lemma]{Proposition}
\newtheorem{corollary}[lemma]{Corollary}
\theoremstyle{remark}
\newtheorem{remark}[lemma]{Remark}
\numberwithin{equation}{section}
\newcommand{\R}{\mathbb R}
\newcommand{\C}{\mathbb C}

\newcommand{\G}{\mathcal G}
\newcommand{\supp}{\operatorname{supp}}
\newcommand{\dist}{\operatorname{dist}}

\newcommand{\diam}{\operatorname{diam}}

\DeclareMathOperator*{\esssup}{ess\,sup}

\title[Canonical Cauchy bases]{Canonical Cauchy bases on Ahlfors regular sets}
\author{Mahishanka Withanachchi}
\address{Department of Mathematics \& Statistics, University of Calgary, Calgary, Canada}
\date{}
\allowdisplaybreaks[1]
\begin{document}
\lefthyphenmin=3
\righthyphenmin=3
\brokenpenalty=10000
\raggedbottom
\begin{abstract}
Let $\mu$ be an $s$-Ahlfors regular probability measure on the real
line, with $0<s<1$. The zeros of its Cauchy transform determine a
canonical rational basis. We characterize bounded and compact pole
corrections by Carleson packing and obtain upper and lower estimates
for the operator norm and essential norm, including weighted estimates. Local moment
tests give a Schatten criterion for every positive exponent, with
coefficients that measure both pole displacement and tangent change.
On finitely many endpoint chains in a fixed thin Cantor set, an exact
Toeplitz reduction determines the basis and frame regimes, the Fredholm
index, and the critical case of dense nonclosed range. For $C^1$
deformations with a lower chord bound, square Dini continuity of the derivative gives
a compact transport correction on every interior $L^p$ space. Cantor constructions
show that this modulus condition is sharp and distinguish it from the
Dini condition at $L^1$. At positive regularity, the correction has the
sharp approximation rate $n^{-(m+\alpha)/s}$ for
$C^{m+1,\alpha}$ deformations, including both endpoints. Ordering by
gap length gives endpoint Schauder bases, almost everywhere
reconstruction, variation estimates, and a sharp H\"older
approximation transition.
\end{abstract}
\maketitle

\section{Introduction}

Let $\mu$ be an $s$-Ahlfors regular probability measure on a compact
subset $E$ of the real line, where $0<s<1$. Its Cauchy transform
\[
 F(z)=\int_E\frac{d\mu(t)}{t-z}
\]
has a unique zero $a_g$ in each bounded complementary gap $g$.
The functions
\[
 e_g(x)=\frac{\sqrt{w_g}}{x-a_g},\qquad
 w_g=F'(a_g)^{-1},
\]
together with the constant function form an orthonormal basis of
$L^2(\mu)$. This paper studies the effect of moving the poles and
deforming the supporting set. The main results give geometric
criteria for operator ideals and compactness, and an exact basis and
frame classification for a family of nested pole perturbations.

Proposition~\ref{thm:basis} supplies the quantitative foundation.
The canonical family is unconditional in every $L^p(\mu)$,
$1<p<\infty$, with constants depending only on $p$ and the
regularity data. The identity for the reciprocal Cauchy transform
\[
 -\frac1{F(z)}=z-m_1+\int\frac{d\nu(a)}{a-z},
 \qquad \nu=\sum_gw_g\delta_{a_g},\qquad m_1=\int t\,d\mu(t),
\]
and the estimate $\nu(B(x,r))\lesssim r^{2-s}$ for $x\in E$
give Calder\'on--Zygmund bounds for all bounded coefficient
multipliers. The uniform kernel smoothness exponent $1-s$ is sharp.

\subsection*{Packing and spectral criteria}
Write $\ell_g=|g|$ and consider complex poles $b_g$ with
$|b_g-a_g|\le c_0\ell_g$, where $c_0$ is a fixed admissible
geometric constant. Put
\[
 q_g(x)=\frac{\sqrt{w_g}}{x-b_g},\qquad
 \eta_g=\frac{|b_g-a_g|}{\ell_g},\qquad
 \lambda_g=\eta_g^2\ell_g^s.
\]
Theorem~\ref{thm:stability} shows that $(q_g)$ is Bessel exactly
when $(\lambda_g)$ has finite Carleson packing. If $T1=0$ and
$Te_g=q_g-e_g$, then, for every $1<p<\infty$,
\[
 \|T\|_{p\to p}\asymp_p\mathcal C(\lambda)^{1/2},\qquad
 \|T\|_{\mathrm{ess},L^p\to L^p}
       \asymp_p\gamma(\lambda)^{1/2},
\]
where $\gamma$ is the limiting packing of the small gaps.
Vanishing tail packing gives an unconditional basis without requiring
the total packing norm to be small. The weighted version holds on
$L^p(v\,d\mu)$ for every $v\in A_p(\mu)$;
see Proposition~\ref{prop:weighted-packing}.

Theorem~\ref{thm:geometric-schatten} gives a spectral criterion that
also includes tangent changes. For complex coefficients
$\alpha_g,\beta_g$, let
\[
 R1=0,\qquad
 Re_g(x)=
 \frac{\sqrt{w_g}\,[\alpha_g(x-a_g)+\ell_g\beta_g]}
      {(x-a_g)(x-b_g)}.
\]
For every $r>0$,
\[
 R\in S^r
 \quad\Longleftrightarrow\quad
 \sum_g(|\alpha_g|^2+|\beta_g|^2)^{r/2}<\infty,
\]
with equivalence of the corresponding norms or quasinorms.
In particular, pole displacements characterize $S^r$ for every
positive exponent. The upper estimate follows from Carleson
truncation on the gap filtration. For $r\le1$, the lower estimate
uses local tests obtained by subtracting polynomial projections.
Their vanishing moments control interactions with larger gaps.
A finite separation of scales and locations makes the remaining
interactions small enough to recover the diagonal coefficients.

\subsection*{An exact transition on Cantor chains}
Finite packing does not imply a lower basis bound.
Proposition~\ref{prop:dyadic-model-transfer} realizes classical
paraproduct models by admissible Cauchy poles on a fixed thin Cantor
set. It yields frames with every finite excess and with infinite
excess, with uniform packing and frame bounds.
The comparison has an error of order $\rho$, where $\rho$ is the
Cantor contraction ratio.

Theorem~\ref{prop:exact-chain-index} removes that error for finitely
many endpoint chains. It applies at every fixed
$0<\rho\le1/256$, to complex relative shifts with a limit on each
chain. Self-similarity produces a bilateral Cauchy transform with a
central zero $a_\rho$. If the limiting relative shift on the $j$th
chain is $\eta_j$, set
\[
 \beta_j=a_\rho+(1-2\rho)\eta_j.
\]
The synthesis operator is Fredholm exactly when
$|\beta_j|\ne\sqrt2a_\rho$ for every chain. In that case it is
surjective and its kernel dimension is
\[
 \#\{j:|\beta_j|>\sqrt2a_\rho\}.
\]
All these arrays are complete and Bessel. Thus the strict interior
of each critical circle gives a Riesz basis, crossing a circle adds
one redundant coordinate, and a point on a circle gives dense
nonclosed synthesis range. The proof combines a rational product
argument for completeness with a compact reduction to Toeplitz
operators. A classical theta identity reduces each limiting block
to $I-\delta S^*$, up to invertible factors, where
$\delta=\beta/(\sqrt2a_\rho)$.
Finitely many additional admissible replacements are allowed.
More general compact changes preserve the Fredholm criterion and
index; completeness for those changes is not asserted.

\subsection*{Curve transport and its regularity thresholds}
Let the parameter poles $b_g$ be real, and let
$\phi:[0,1]\to\C$ be $C^1$ with a lower chord bound. The tangent
normalization
\[
 Q_g(z)=\frac{\sqrt{w_g}\phi'(b_g)}{z-\phi(b_g)}
\]
defines a correction $H_\phi e_g=Q_g\circ\phi-q_g$.
Theorem~\ref{thm:curve-transport} proves compactness on every
$L^p(\mu)$, $1<p<\infty$, under the square Dini condition
\[
 \int_0^1\omega_\phi(t)^2\,\frac{dt}{t}<\infty,
 \qquad
 \omega_\phi(t)=\sup_{|x-y|\le t}|\phi'(x)-\phi'(y)|.
\]
The deformation need not fix $E$. A more local sufficient condition
uses the derivative oscillation around each parameter pole, summed
over annuli and tested by Carleson packing; see
Proposition~\ref{prop:curve-local-oscillation}.
The argument first proves a Hilbert space estimate and then combines
uniform singular integral bounds with interpolation. Compactness preserves the essential norm of the pole correction and the Fredholm index of the corresponding synthesis operator.
It does not by itself imply injectivity of the transported synthesis.

The square Dini threshold is sharp as a condition on a prescribed
concave modulus. Proposition~\ref{prop:sharp-modulus-transport}
constructs arbitrarily small $C^1$ deformations for which interior
$L^p$ boundedness holds exactly when the square Dini integral is
finite. In the same examples, $L^1$ boundedness holds exactly when
the stronger Dini integral is finite. M\"obius postcomposition
preserves these thresholds and gives examples on circular supports
that are not contained in a line. The Dini condition also gives
compactness on both endpoint spaces and on $C(E)$.

For maps sufficiently close to the identity in $C^1$ that fix $E$
pointwise, Proposition~\ref{prop:c1-support-fixing}
gives an exact packing criterion. Theorem~\ref{thm:geometric-schatten}
then yields the joint spectral criterion
\[
 H_\phi\in S^r
 \quad\Longleftrightarrow\quad
 \sum_g\left(|\phi'(a_g)-1|^2+
       \frac{|\phi(a_g)-a_g|^2}{\ell_g^2}\right)^{r/2}<\infty,
 \qquad r>0.
\]
The dependence on the extension through the gaps is essential.
Gap bump constructions vary the pole positions and tangent factors.
They give compact, bounded noncompact, and unbounded corrections
with $\phi$ arbitrarily close to the identity in $C^1$.

At positive regularity, $\phi\in C^{m+1,\alpha}$ gives
\[
 a_n^{(p)}(H_\phi)\lesssim n^{-(m+\alpha)/s},
 \qquad 1\le p\le\infty,
\]
with the same estimate on $C(E)$ and a corresponding local
Lipschitz bound in the derivative. The rate is attained by deformations parametrized by arc length whose support is noncollinear. The approximation removes
a polynomial term of finite rank, retains the coefficients associated with large gaps,
and projects the remaining output onto polynomials on the gap cells.
Smooth deformations give decay faster than every power, and
holomorphically univalent deformations give exponential decay.
M\"obius postcomposition changes the correction by an explicit operator of rank one and preserves invertibility of the bounded transport
operators under the realizations specified below.

\subsection*{Ordered reconstruction}
Proposition~\ref{thm:approximation-norming} supplies the endpoint
reconstruction used in the transport results. Every ordering by nonincreasing gap length, including every ordering among ties, gives
Schauder bases of $C(E)$ and $L^1(\mu)$. The partial sums converge
almost everywhere for $L^1$ data and satisfy variation estimates of
weak type $(1,1)$ and strong type $(p,p)$ for variation exponent
greater than two. The same projections give the sharp H\"older
approximation transition at exponent $1-s$.

\subsection*{Relation to earlier work}
The qualitative spectral resolution is part of Clark theory and
the theory of perturbations of rank one
\cite{Clark1972,GesztesySimon1995,Simon1995}.
We give a direct proof because the subsequent estimates require
the location of the poles and the growth of the reciprocal measure.
Rational bases have a much broader history.
Petrushev~\cite[Corollary~7.1]{Petrushev2000} constructs bases of
uniformly bounded rational degree in classical function spaces;
Bultheel and Gonz\'alez-Vera~\cite{BultheelGonzalezVera2000} develop
rational multiscale decompositions for probability measures.
Takenaka--Malmquist systems are Schauder bases of their closed spans
in $H^p$, $1<p<\infty$ \cite{QianChenTan2014}.
The present system consists of rational functions of degree one with poles selected by the measure, and its estimates
are uniform over singular $s$-regular supports with fixed data.

For real poles the Clark representation identifies the perturbed
functions with normalized kernels in model spaces, up to uniformly
bounded scalar factors. Baranov's stability theory
\cite{Baranov2005} and its extension in
\cite{BaranovFricainMashreghi2010} are therefore direct comparisons.
The scale given by the reciprocal derivative in \cite[Theorem~1.4]{Baranov2005} is
$w_g\asymp\ell_g^{2-s}$. Our examples on sparse families of gaps permit
displacements of order $\ell_g$, and the ratio of the corresponding
unnormalized kernel norms tends to zero. The comparability assumptions on
kernel norms in \cite[Theorem~1.1]{Baranov2005} and
\cite[Theorem~7.2]{BaranovFricainMashreghi2010} consequently do not
hold for those examples. Appendix~\ref{app:modelspaces} gives the
normalizations and the precise calculation.

Finite replacement of kernel bases in spaces with the division
property is an established principle; see \cite{Baranov2023}.
Its use here is combined with the geometric tail norm and a direct
$L^p$ construction of the dual functionals. General embedding theory for model spaces provides another comparison
\cite{Baranov2009,LaceySawyerShenUriarteWick2018}.
In particular, \cite[Theorems~1.2, 1.9, and 1.10]{LaceySawyerShenUriarteWick2018}
give criteria with two weights for the Cauchy transform, together with consequences for model spaces and compactness. Our packing criterion explicitly
evaluates the norm for the particular measures indexed by gaps;
the essential norm concerns the difference from canonical synthesis.
The full perturbed synthesis is never compact, as explained in
Appendix~\ref{app:modelspaces}.

The $L^2$ Clark identification alone gives no assertion on the
$L^p(\mu)$ spaces. For a finite positive singular measure
on the circle, boundedness of the normalized Cauchy transform from
$L^p(\mu)$ to $H^p$ for some $p>2$ forces the measure to be
discrete; see Aleksandrov's theorem recalled in
\cite[Theorem~1]{Poltoratski2003}. We prove the estimates for the Cauchy basis directly from the gap geometry.

Belov, Mengestie, and Seip
\cite[Theorem~1.3 and Section~5]{BelovMengestieSeip2011} give general
invertibility criteria for weighted discrete Hilbert transforms.
Their introductory Example~1 exhibits
bounded discrete Cauchy transforms that become invertible after
one point is removed from the target sequence. Taking adjoints
gives a surjective synthesis with one redundant coordinate.
The paraproduct obstruction and its reverse H\"older criterion are
also classical \cite{Pereyra1994,PereyraWard1998}.
The additional assertion here is a uniform comparison, for a
whole class of real Carleson arrays, with actual Cauchy poles on
a nonatomic regular Cantor set. It produces every finite excess
and infinite excess while keeping the poles admissible in their
original gaps and the packing norm bounded.
A fixed change of Clark parameter instead produces relative
shifts tending to zero; Appendix~\ref{app:modelspaces} gives the
scale calculation. For an inverted endpoint chain, the weights in
the discrete transform satisfy $v_n\asymp2^n$. Hence the additional
condition $v_n=o(\sum_{k<n}v_k)$ in
\cite[Theorem~5.1]{BelovMengestieSeip2011} fails. Our exact chain
theorem evaluates the limiting matrix for this geometry, including
its critical circle and index, within the general framework for discrete
transforms in their Theorem~1.3. The failure of the additional
weight hypothesis excludes only the specialized Theorem~5.1.
The bilateral identity used to evaluate its symbol is
classical \cite[\S\,17.8]{DLMF}.
Proposition~\ref{prop:clark-projection-obstruction}
also excludes a full Clark projection realizing positive
finite excess for every inner dilating factor, when the coordinates
are exactly the original gap nodes and infinity and the original
Clark identification is used. A realization of this form with infinite excess requires a factor that does not extend holomorphically
through all of $E$.

Our construction using operators of finite rank belongs to the piecewise polynomial
approximation tradition of Birman and Solomyak
\cite{BirmanSolomyak1977}. Metric entropy estimates for integral
operators on $C(X)$ provide a further comparison through
\cite{CarlHeinrichKuhn1988}; spectral estimates for operators
with singular measures and potential kernels are developed in
\cite{RozenblumTashchiyanJFA2022,RozenblumTashchiyan2022}.
The localized expansion method of
Rochberg and Semmes~\cite{RochbergSemmes1989} treats singular values
of nonlinear singular integral operators and already gives
estimates for approximation numbers on interior $L^p$ spaces in
Section~VII.D. Here the operator is specified by
Cauchy coefficients indexed by complementary gaps. The two Schur
estimates use the input gaps and output cells together and yield
the same bound on every $L^p$, including both endpoints, and on
$C(E)$. They cover every positive finite H\"older regularity and
include the integer endpoints without a logarithmic loss.
The matching examples retain both unit speed and noncollinear
support. Exponential decay from a uniformly analytic kernel is
classical; its $L^2$ form is covered by
\cite[Theorem~11.3]{BirmanSolomyak1977}.
Integrable envelopes for the projection kernels and stability under changes
of summation order have wavelet precedents in
\cite{KellyKonRaphael1994}. The cancellation with the omitted tail used for our projections has the same form as the wavelet argument
in \cite[proof of Theorem~1]{Tao1996}. The required estimates here
are for the prescribed Cauchy family of rational degree one ordered by gap length. The variation argument uses L\'epingle's inequality
\cite{Lepingle1976} and the martingale comparison of
Jones, Seeger, and Wright~\cite{JonesSeegerWright2008}.
Variation estimates for approximate identities and spatial
truncations are treated in
\cite{Liu2019,ZorinKranich2020}; our cutoffs instead select rational
coefficients by gap length.

The spectral and weighted arguments use established methods from harmonic analysis. The nearly weakly orthogonal sequence method
of Rochberg and Semmes~\cite{RochbergSemmes1989} and its metric formulation
in \cite[Section~7]{Hytonen2026} supply the Carleson truncation
principle. Schatten characterizations at every positive exponent
are already known for dyadic paraproducts
\cite[Theorem~2.1]{PottSmith2004} and uniformly adapted wavelet
families \cite[preprint version, Theorem~4.1 and Lemma~4.15]{LaceyLiWick2026}.
These arguments use compressions on separated scales; the wavelet
treatment also separates locations.
Our rational columns vary with the poles and numerator coefficients,
and their endpoint rescalings involve varying singular measures.
Lemma~\ref{lem:local-moment-schatten} supplies a uniform nonzero
polynomial residual for this family, which yields the local moment
tests needed for the lower estimate. The proof adapts these methods
to a uniform geometric criterion for the rational gap family.
The sparse theorem of Lorist~\cite[Theorem~6.1]{Lorist2021}
gives the weighted upper bounds; interpolation with change of
measure \cite{SteinWeiss1958} transfers the geometric lower bounds.

Square Dini regularity also
occurs in the sharp estimates of Dong, Jung, and Kim
\cite{DongJungKim2026} for nonlocal operators with variable
coefficients. Their operators and function spaces differ from the transport indexed by gaps considered here. Our conclusion is compactness
of this transport, with a sharp uniform modulus condition and a local
Carleson sufficient condition; no new PDE solvability theorem is
claimed. To realize displacements by maps that fix the
support, we use scaled smooth bumps inside the gaps and verify the
required regularity directly.

Finally, polynomial Remez inequalities on Ahlfors regular sets are
established results
\cite{Brudnyi2002,BrudnyiBrudnyi2006,Brudnyi2010}.
The Remez consequence proved below combines such a norming estimate
with a bound for the product clearing the retained poles.
We make no optimality claim for its constants.

Section~\ref{sec:prelim} develops the gap geometry and the analytic
estimates shared by the proofs. Section~\ref{sec:results} states
the four main theorems and the reconstruction proposition.
Sections~\ref{subsec:basis-proof} through~\ref{subsec:reconstruction-proof}
give the proofs and their consequences.
Section~\ref{sec:future} records the remaining questions, and
Appendix~\ref{app:modelspaces} gives the comparison with model spaces.

\section{Preliminaries}
\label{sec:prelim}

\subsection{Regular sets and gap geometry}
Let $E\subset[0,1]$ be compact, with $\min E=0$ and $\max E=1$.
Let $\mu$ be a probability measure supported on $E$. Fix $0<s<1$ and
$A\ge1$ such that
\begin{equation}
 A^{-1}r^s\le\mu(B(x,r))\le Ar^s,
 \qquad x\in E,\quad 0<r\le1.
 \label{sp:regularity}
\end{equation}
Write $B_E(x,r)=B(x,r)\cap E$. The implicit constants in $\lesssim$ and
$\asymp$ may depend on $A,s$. Any additional dependence is stated explicitly.
The upper bound extends to $r>1$, since
\[
 \mu(B(x,r))\le1\le r^s.
\]

If a ball $B(a,r)$ meets $E$, choose $x\in B(a,r)\cap E$. Then
\[
 B(a,r)\cap E\subset B_E(x,2r),\qquad
 \mu(B(a,r))\le C r^s.
\]
For each fixed $L\ge1$, the two regularity bounds give
\begin{equation}
 \mu(B_E(x,Lr))\le C_L\mu(B_E(x,r)),
 \qquad x\in E,\quad r>0.
 \label{eq:prelim-doubling}
\end{equation}
For $r\le1$, the quotient is at most $A^2L^s$; for $r>1$, both
balls have full measure. Thus $(E,|\cdot|,\mu)$ is a doubling metric
measure space.

The measure has no atoms because
\[
 \mu(\{x\})\le\inf_{r>0}Ar^s=0.
\]
Every ball centered on $E$ has positive mass, so $E$ has no isolated
points. A maximal $r$-separated family in $E$ has at most $Cr^{-s}$
points: its disjoint $r/3$-balls have mass at least $cr^s$ each.
If $N(E,r)$ is the least number of radius-$r$ balls covering $E$, then
\begin{equation}
 N(E,r)\le Cr^{-s},\qquad 0<r\le1.
 \label{sp:covering}
\end{equation}
Consequently,
\[
 |E|\le C r^{-s}r=C r^{1-s}\longrightarrow0.
\]

Let $\G$ denote the collection of bounded complementary intervals $g=(u_g,v_g)$,
and put $\ell_g=v_g-u_g$. For $0<r\le R\le2$ and $x\in E$,
\begin{equation}
 \#\{g:r/2<\ell_g\le r,\ \dist(x,g)\le R\}
 \le C(R/r)^s,
 \qquad \#\{g:\ell_g\ge r\}\le Cr^{-s}.
 \label{eq:prelim-gapcount}
\end{equation}
Indeed, the left endpoints of gaps of length greater than $r/2$
are $r/2$-separated. Their $r/8$-balls are disjoint, and
\[
 c r^s\,\#\{g:r/2<\ell_g\le r,\ \dist(x,g)\le R\}
 \le \mu(B_E(x,R+2r))\le C R^s.
\]
The global count follows by applying the same argument to all gaps
of length at least $r$, since the total mass is one.

\subsection{Cauchy transforms and basis terminology}
For $z\in\C\setminus E$, define
\[
 F(z)=C_\mu(z)=\int_E\frac{d\mu(t)}{t-z},
 \qquad m_j=\int_Et^j\,d\mu(t).
\]
Differentiation under the integral is valid on compact subsets of
$\C\setminus E$. For $\Im z>0$,
\[
 \Im F(z)=\int_E\frac{\Im z}{|t-z|^2}\,d\mu(t)>0.
\]
A holomorphic function with nonnegative imaginary part in the upper
half-plane is called a Herglotz function. In particular, $-1/F$ is
Herglotz, since
\[
 \Im\left(-\frac1F\right)=\frac{\Im F}{|F|^2}>0.
\]

We use the inner product
\[
 \langle f,h\rangle=\int_E f\overline h\,d\mu,
\]
linear in its first argument. For a bounded self-adjoint operator $T$
with spectrum in $[0,1]$ and a vector $h$, the spectral theorem gives
a finite positive measure $\nu_h$ such that
\[
\langle(T-z)^{-1}h,h\rangle
 =\int_{[0,1]}\frac{d\nu_h(a)}{a-z},
 \qquad \nu_h([0,1])=\|h\|^2.
\]
If $h$ is cyclic, the map $p\mapsto p(T)h$ extends from polynomials
to a unitary map from $L^2(\nu_h)$ onto the underlying Hilbert space.
The vector is cyclic when the span of $\{T^nh:n\ge0\}$ is dense.

The measure associated with $G(z)=\int(a-z)^{-1}d\nu(a)$ is recovered by
\begin{align*}
 \nu(\{a_0\})
 &=\lim_{y\downarrow0}y\Im G(a_0+iy),\\
 \nu((a,b))
 &=\lim_{y\downarrow0}\frac1\pi\int_a^b\Im G(x+iy)\,dx
 \quad\text{if }\nu(\{a,b\})=0.
\end{align*}
The first identity follows by dominated convergence applied to
$y^2/((a-a_0)^2+y^2)$. The second follows from Tonelli's theorem and
the integral of the Poisson kernel. We use these formulas to identify
the atoms in the gaps and to show that no other part of the measure remains.

A family $(v_j)$ in a Hilbert space is Bessel if
\[
 \sum_j|\langle f,v_j\rangle|^2\le B\|f\|^2.
\]
By Hilbert space duality, this is equivalent to
\[
 \left\|\sum_jc_jv_j\right\|^2\le B\sum_j|c_j|^2
\]
for finite sequences. A complete family is a Riesz basis if the last
inequality also has a positive lower bound. For an unconditional basis
in $L^p$, every function has a unique expansion that converges under every
rearrangement. Below, we use the criterion of uniform boundedness of the
finite coordinate multipliers.
The synthesis map of a Bessel family is
$\mathcal A(c)=\sum_jc_jv_j$ from $\ell^2$ to the Hilbert space.
The family is a frame if its analysis inequality also has a
positive lower bound, or equivalently if $\mathcal A$ is surjective.
Its excess is $\dim\ker\mathcal A$, which may be infinite.

For $0<\alpha\le1$, write
\[
 [f]_{C^\alpha(E)}
 =\sup_{x\ne y\in E}\frac{|f(x)-f(y)|}{|x-y|^\alpha}.
\]
Restrictions of polynomials are dense in $C(E)$ and in $L^p(\mu)$,
$1\le p<\infty$. This follows from Stone--Weierstrass and the density
of continuous functions in $L^p$ for a finite Borel measure on a
compact metric space.

\subsection{The filtration determined by the gaps}
We construct a filtration used both for curve approximation and
for the variation estimates in reconstruction. Put $r_k=2^{-k}$, $k\ge0$.
There is $c_{\rm por}>0$ such that every interval $[a,b]$ with
$a,b\in E$ and $a<b$ contains a complementary gap of length at
least $c_{\rm por}(b-a)$.
To verify this, put $L=b-a$ and choose a maximal set whose points are
separated by at least $\eta L$ in $E\cap[a,b]$.

If this set has $N$ points, disjointness of its
balls of radius $\eta L/3$ and regularity give
\[
 cN(\eta L)^s
 \le \mu(B_E(a,2L))
 \le CL^s,\qquad N\le C\eta^{-s}.
\]
Its balls of radius $\eta L$ cover $E\cap[a,b]$.
If every gap inside $(a,b)$ had length less than $2\eta L$,
each point of $[a,b]$ would be within $\eta L$ of $E\cap[a,b]$.
The balls of radius $2\eta L$ centered at the chosen points would then cover $[a,b]$.

Their total Lebesgue length is at most
\[
 4\eta LN\le C\eta^{1-s}L.
\]
Choose $\eta>0$ so small that $C\eta^{1-s}<1$.
The resulting contradiction proves the claim with
$c_{\rm por}=2\eta$.

For $r=r_k$, cut $E$ at all gaps of length at least $r$.
Let $\mathcal Q_k$ be the resulting finite partition into compact
subsets of $E$. Distinct cells are separated by a gap of length
at least $r$, and the partitions are nested.

For $Q\in\mathcal Q_k$, the interval between $\min Q$ and $\max Q$
contains no gap of length at least $r_k$.
The porosity estimate therefore gives
$\operatorname{diam}Q\le c_{\rm por}^{-1}r_k$.
Moreover, for every $x\in Q$,
\[
 B_E(x,r_k)\subset Q\subset B_E(x,Cr_k).
\]
Thus, for a uniform constant $C_0\ge1$,
\begin{equation}
 \begin{gathered}
 \operatorname{diam}Q\le C_0r_k,\qquad
 c r_k^s\le\mu(Q)\le C r_k^s,\\
 \operatorname{dist}(Q,Q')\ge r_k
 \quad(Q,Q'\in\mathcal Q_k,\ Q\ne Q').
 \end{gathered}
 \label{var:cell-geometry}
\end{equation}

\begin{figure}[htbp]
\centering
\begin{tikzpicture}[x=11.2cm,y=1.2cm,font=\small,>=Stealth]
 \foreach \y/\lab in {2.6/{r=3^{-1}},1.5/{r=3^{-2}},.4/{r=3^{-3}}}
   \node[anchor=east] at (-.035,\y+.06) {$\lab$};
 \foreach \a/\b in {0/.333333,.666667/1}{
   \fill[DiagramPrimary!8] (\a,2.43) rectangle (\b,2.78);
   \draw[DiagramPrimary] (\a,2.43) rectangle (\b,2.78);
 }
 \foreach \a/\b in {0/.111111,.222222/.333333,.666667/.777778,.888889/1}{
   \fill[DiagramPrimary!8] (\a,1.33) rectangle (\b,1.68);
   \draw[DiagramPrimary] (\a,1.33) rectangle (\b,1.68);
 }
 \foreach \a/\b in {0/.037037,.074074/.111111,.222222/.259259,.296296/.333333,.666667/.703704,.740741/.777778,.888889/.925926,.962963/1}{
   \fill[DiagramPrimary!8] (\a,.23) rectangle (\b,.58);
   \draw[DiagramPrimary] (\a,.23) rectangle (\b,.58);
   \foreach \y in {2.6,1.5,.4}
     \draw[black,line width=1.6pt] (\a,\y)--(\b,\y);
 }
 \foreach \a/\b in {.166667/.055556,.166667/.277778,.833333/.722222,.833333/.944444}
   \draw[->,black!35] (\a,2.4)--(\b,1.72);
 \foreach \a/\b in {.055556/.018519,.055556/.092593,.277778/.240741,.277778/.314815,.722222/.685185,.722222/.759259,.944444/.907407,.944444/.981481}
   \draw[->,black!35] (\a,1.3)--(\b,.62);
 \draw[<->,DiagramSecondary] (.333333,2.96)--(.666667,2.96);
 \node[DiagramSecondary,anchor=south] at (.5,2.99) {retained gap};
 \node[anchor=north,align=center] at (.5,-.04)
 {$\operatorname{diam}Q\le C r,\quad \mu(Q)\asymp r^s,\quad
   \operatorname{dist}(Q,Q')\ge r$\\[4pt]
  $\displaystyle \mathbb E_r f|_Q=\frac{1}{\mu(Q)}\int_Qf\,d\mu$};
\end{tikzpicture}
\caption{The partition obtained by cutting at every gap of length at least
$r$, illustrated on the middle-thirds Cantor set at three triadic scales.
Each rectangle represents the part of the support in one cell; the black
segments indicate only a finite stage of that support. Smaller cutoffs
refine the partition. At dyadic scales, these cells define the
averages used in compact curve approximation and in the comparison
of $P_{2^{-k}}$ with $\mathbb E_k$.}
\label{fig:gap-filtration}
\end{figure}
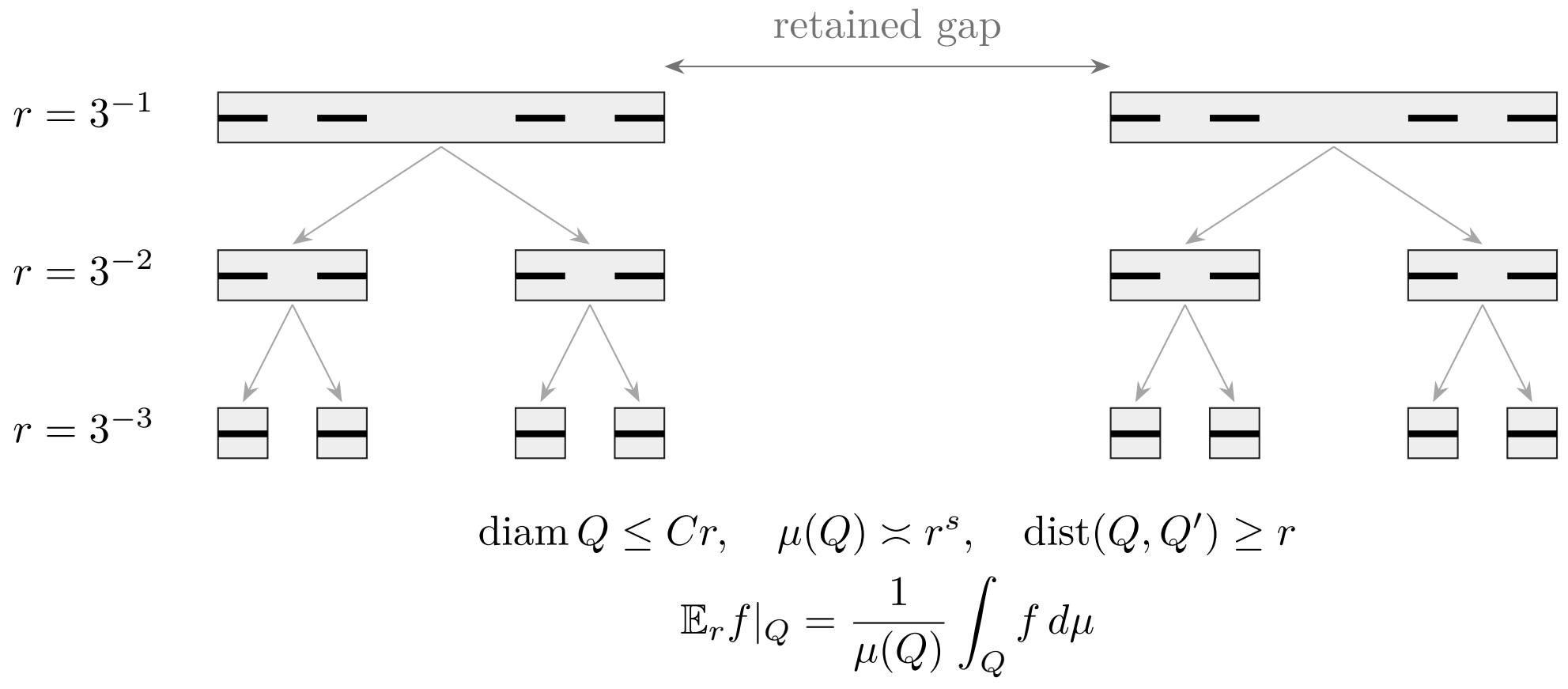

Figure~\ref{fig:gap-filtration} illustrates the nested cells.

Let $\mathbb E_k$ be conditional expectation for this partition:
\[
 \mathbb E_k f(x)
 =\frac1{\mu(Q_k(x))}\int_{Q_k(x)}f\,d\mu,
 \qquad Q_k(x)\in\mathcal Q_k,\quad x\in Q_k(x).
\]
The partition at level zero consists of $E$, so
$\mathbb E_0f=\langle f,1\rangle1$.
For $f\in C(E)$, \eqref{var:cell-geometry} implies
\[
 |\mathbb E_kf(x)-f(x)|
 \le\sup_{\substack{y\in E\\|x-y|\le C_0r_k}}
                         |f(y)-f(x)|\longrightarrow0
\]
uniformly on $E$. Conditional expectations are contractions on
every $L^p$, $1\le p\le\infty$.
By density, $\mathbb E_kf\to f$ in $L^p$ for
$1\le p<\infty$.
\subsection{Maximal, singular integral, and martingale estimates}
The centered maximal function is
\[
 M_\mu f(x)=\sup_{r>0}
 \frac1{\mu(B_E(x,r))}\int_{B_E(x,r)}|f|\,d\mu.
\]
For a fixed radius, the quotient is continuous in the center, since
the boundaries of intervals have zero $\mu$ and $|f|\mu$ measure.
Thus $M_\mu f$ is lower semicontinuous.
The $5r$ covering lemma and \eqref{eq:prelim-doubling} give
\begin{equation}
 \mu\{M_\mu f>t\}\le \frac Ct\|f\|_1,\qquad t>0.
 \label{eq:prelim-maximal-weak}
\end{equation}
Indeed, for $t>\|f\|_1$, select disjoint balls $B_j$ with
average greater than $t$ whose fivefold enlargements cover the level set.
Then
\[
 \mu\{M_\mu f>t\}
 \le\sum_j\mu(5B_j)
 \le C\sum_j\mu(B_j)
 \le\frac Ct\sum_j\int_{B_j}|f|\,d\mu
 \le\frac Ct\|f\|_1.
\]
For smaller $t$, the estimate follows from $\mu(E)=1$.
Splitting $f$ at height $t/2$ and integrating the resulting distribution estimate
also gives
\begin{equation}
 \|M_\mu f\|_p\le C_p\|f\|_p,\qquad 1<p<\infty.
 \label{eq:prelim-maximal-lp}
\end{equation}
For example,
\[
 \|M_\mu f\|_p^p
 \le 2Cp\int_E |f(x)|
       \int_0^{2|f(x)|}t^{p-2}\,dt\,d\mu(x)
 =\frac{2^pCp}{p-1}\|f\|_p^p.
\]
Applying \eqref{eq:prelim-maximal-weak} to approximations by continuous functions
proves the differentiation theorem for $L^1(\mu)$ functions with respect to balls.
The cell geometry likewise gives
$\mathbb E_k|f|\le C M_\mu f$.
Since $\mathbb E_k h\to h$ uniformly for continuous $h$, the same
density argument proves $\mathbb E_kf\to f$ almost everywhere for
every $f\in L^1(\mu)$.

We use the following form of the Calder\'on--Zygmund estimate.
We include the proof to specify the vector-valued version used
in the variation estimates. The underlying principle is classical
\cite{CoifmanWeiss1977,NTV1998}.

\begin{lemma}
\label{lem:vector-cz}
Suppose that each of $X,Y$ is either a finite-dimensional Banach space or a Hilbert space.
Suppose the linear operator $S:L^2(\mu;X)\to L^2(\mu;Y)$ has norm at most $B_2$
and an off-diagonal kernel $K(x,y)\in\mathcal L(X,Y)$ satisfying
\[
 \|K(x,y)\|\le H|x-y|^{-s}
\]
and
\[
 \|K(x,y)-K(x',y)\|+\|K(y,x)-K(y,x')\|
 \le H\frac{|x-x'|^\delta}{|x-y|^{s+\delta}}
\]
when $|x-x'|\le |x-y|/2$, where $0<\delta\le1$.
The off-diagonal representation means that, for bounded $f$ supported
in a closed set $F\subset E$,
\[
 Sf(x)=\int_F K(x,y)f(y)\,d\mu(y)
 \quad\text{for almost every }x\in E\setminus F.
\]
Then
\begin{equation}
 \begin{gathered}
 \sup_{t>0}t\,\mu\{\|Sf\|_Y>t\}
       \le C(B_2+H)\|f\|_{L^1(X)},\\
 \|Sf\|_{L^p(Y)}
       \le C_p(B_2+H)\|f\|_{L^p(X)},\qquad 1<p<\infty.
 \end{gathered}
 \label{eq:prelim-cz-conclusion}
\end{equation}
The constants depend on $A,s,\delta,p$, but not on the dimensions
of $X,Y$.
\end{lemma}

\begin{proof}
We may normalize $B_2+H=1$, since the zero case is immediate.
Let $f\in L^2(\mu;X)$ and $t>\|f\|_{L^1(X)}$.
In the gap filtration, select the first cell $Q$ along each branch for which
\[
 \frac1{\mu(Q)}\int_Q\|f\|_X\,d\mu>t.
\]
The stopping cells are disjoint. Their parents have average at most $t$,
and \eqref{var:cell-geometry} bounds the ratio between the measures of a parent and its child.
For the Bochner averages $f_Q=\mu(Q)^{-1}\int_Qf\,d\mu$, set
\[
 b_Q=(f-f_Q)\mathbf1_Q,\qquad
 g=f-\sum_Qb_Q,\qquad \Omega=\bigcup_Q Q.
\]
Differentiation by the filtration and the stopping construction give
\[
 \begin{gathered}
 \mu(\Omega)\le t^{-1}\|f\|_{L^1(X)},\qquad
 \|g\|_{L^\infty(X)}\le Ct,\qquad
 \|g\|_{L^1(X)}\le\|f\|_{L^1(X)},\\
 \|g\|_{L^2(X)}^2\le Ct\|f\|_{L^1(X)},\qquad
 \int b_Q\,d\mu=0,\qquad
 \sum_Q\|b_Q\|_{L^1(X)}\le2\|f\|_{L^1(X)}.
 \end{gathered}
\]
The sum of the bad functions converges in $L^2(\mu;X)$ by
disjointness and domination by $\|f\|_X+Ct\mathbf1_\Omega$.

If $Q$ occurs at scale $r_Q$, choose $z_Q\in Q$ and put
$Q^*=B_E(z_Q,4C_0r_Q)$. Regularity gives
\[
 \mu\left(\bigcup_QQ^*\right)
 \le C\sum_Q\mu(Q)
 \le Ct^{-1}\|f\|_{L^1(X)}.
\]
Bounded truncations of mean zero supported on $Q$ converge to $b_Q$
in both $L^2$ and $L^1$. The $L^2$ continuity of $S$ and the kernel
size bound therefore extend the off-diagonal representation to $b_Q$
outside $Q^*$. For $x\notin Q^*$, cancellation gives
\[
 \begin{aligned}
 \|Sb_Q(x)\|_Y
 &\le\int_Q\|K(x,y)-K(x,z_Q)\|\,\|b_Q(y)\|_X\,d\mu(y)\\
 &\le C r_Q^\delta |x-z_Q|^{-s-\delta}
                  \|b_Q\|_{L^1(X)}.
 \end{aligned}
\]
The exterior integral is bounded by
\[
 r_Q^\delta
 \int_{|x-z_Q|\ge4C_0r_Q}|x-z_Q|^{-s-\delta}\,d\mu(x)
 \le C\sum_{j\ge0}2^{-j\delta}\le C.
\]
Apply this first to finite sums of the bad functions and then use their $L^2$
convergence and Fatou's lemma. Chebyshev's inequality yields
\[
 \begin{aligned}
 \mu\{\|Sf\|_Y>t\}
 &\le \mu\left(\bigcup_QQ^*\right)
       +\frac4{t^2}\|Sg\|_{L^2(Y)}^2
       +\frac2t\int_{E\setminus\cup_QQ^*}
                     \left\|S\sum_Qb_Q\right\|_Y\,d\mu\\
 &\le \frac Ct\|f\|_{L^1(X)}.
 \end{aligned}
\]
The remaining range of $t$ follows from $\mu(E)=1$.
Approximation by bounded functions defines the weak $L^1$ extension.

For $1<p<2$, split $f$ at height $t$ according to its $X$ norm.
The weak estimate on the part above the threshold and the $L^2$ estimate on
the part below it give, after integrating in $t$,
\[
 \|Sf\|_{L^p(Y)}^p
 \le C_p\left(\frac1{p-1}+\frac1{2-p}\right)
                  \|f\|_{L^p(X)}^p.
\]
For $p>2$, the adjoint maps $L^2(\mu;Y^*)$ to
$L^2(\mu;X^*)$ with the same bound and has kernel
$K(y,x)^*$. The preceding argument applies to it, and duality
gives the desired estimate. The finite-dimensional or Hilbert
space hypothesis guarantees these standard duality identifications for
Bochner spaces.
The adjoint kernel identity follows by Fubini on separated supports,
where the size estimate makes the integral absolutely convergent.
Restoring $B_2+H$ proves \eqref{eq:prelim-cz-conclusion}.
\end{proof}

We also use L\'epingle's martingale variation inequality
\cite{Lepingle1976}:
\[
 \left\|\sup_{k_0<\cdots<k_M}
       \left(\sum_{j=1}^M
       |\mathbb E_{k_j}f-\mathbb E_{k_{j-1}}f|^\rho\right)^{1/\rho}
       \right\|_p
 \le C_{p,\rho}\|f\|_p,
 \qquad 1<p<\infty,\quad \rho>2.
\]
The supremum is over finite increasing sequences. Only its $L^2$
case is needed for our comparison argument; the other exponents
for the rational projections follow from Lemma~\ref{lem:vector-cz}.
The comparison is in the spirit of
Jones, Seeger, and Wright \cite{JonesSeegerWright2008}.
Once the zeros in the gaps have been constructed, we write
\[
 w_g=F'(a_g)^{-1},\qquad
 e_g(x)=\frac{\sqrt{w_g}}{x-a_g},\qquad e_0=1.
\]

The following proposition records the canonical basis and the uniform
estimates used throughout. Its proof is given in
Section~\ref{subsec:basis-proof}.

\begin{proposition}[Canonical Cauchy basis]
\label{thm:basis}
Each gap $g$ contains a unique zero $a_g$ of $F$. Put
\[
 w_g=F'(a_g)^{-1},\qquad e_g(x)=\frac{\sqrt{w_g}}{x-a_g}.
\]
The family $\{1\}\cup\{e_g:g\in\mathcal G\}$ is an orthonormal
basis of $L^2(\mu)$ and an unconditional basis of $L^p(\mu)$,
$1<p<\infty$. For $\|\beta\|_\infty\le1$, the multipliers
\[
 T_\beta f=\beta_0\langle f,1\rangle1+
             \sum_g\beta_g\langle f,e_g\rangle e_g
\]
are uniformly bounded on each $L^p$, $1<p<\infty$, and of weak
type $(1,1)$. After removing the constant component, their kernels satisfy
\[
 |K_\beta(x,y)|\le C|x-y|^{-s},\qquad
 |K_\beta(x,y)-K_\beta(x',y)|
 \le C\frac{|x-x'|^{1-s}}{|x-y|}
\]
for $|x-x'|\le|x-y|/2$, and likewise in the second variable.
The smoothness exponent $1-s$ is optimal uniformly over the multipliers.
\end{proposition}

\section{Main results}
\label{sec:results}
Throughout, $E,\mu,s$ satisfy \eqref{sp:regularity}, and
$\mathcal G$ denotes the bounded complementary gaps.

The gap estimates give $w_g\asymp\ell_g^{2-s}$ and a uniform lower
bound for the relative distance of each zero from the gap endpoints. Choose any constant $c_*>0$ such that
$\dist(a_g,E)\ge c_*\ell_g$ for every gap.
Fix $0<c_0<c_*/2$ and choose complex poles satisfying
$|b_g-a_g|\le c_0\ell_g$. Write
\[
 q_g(x)=\frac{\sqrt{w_g}}{x-b_g},\qquad
 d_g=q_g-e_g,\qquad
 \eta_g=\frac{|b_g-a_g|}{\ell_g},\qquad
 \lambda_g=\eta_g^2\ell_g^s.
\]
Choose an endpoint $\xi_g$ of each gap and set
\begin{equation}
 \mathcal C(\lambda)=
 \sup_{\substack{x\in E\\0<R\le2}}R^{-s}
 \sum_{\substack{\xi_g\in B(x,R)\\\ell_g\le R}}\lambda_g,
 \qquad
 \gamma(\lambda)=\lim_{r\downarrow0}\mathcal C_r(\lambda),
 \label{eq:packing-tail-definition}
\end{equation}
where
\[
 \mathcal C_r(\lambda)
 =\mathcal C\bigl((\lambda_g\mathbf1_{\{\ell_g<r\}})_g\bigr).
\]
The tail norm is tested over all boxes, including those with $R>r$.
The limit exists by monotonicity. Different choices of endpoints give
equivalent norms after the box radius is increased by a fixed factor.
Define
\[
 Dc=\sum_gc_gd_g,\qquad
 Tf=\sum_g\langle f,e_g\rangle d_g,
\]
initially for finite sums. For a bounded operator $S:X\to Y$, let
$\|S\|_{\mathrm{ess},X\to Y}$ be its distance in operator norm
from the compact operators. On $L^2(\mu)$, $S^r$ denotes the
Schatten class whose singular values belong to $\ell^r$, $r>0$.

\begin{theorem}[Packing, pole stability, and its obstruction]
\label{thm:stability}
The family $(q_g)$ is Bessel in $L^2(\mu)$ if and only if
$\mathcal C(\lambda)<\infty$. In that case
\[
 \|D\|_{\ell^2\to L^2}\asymp\mathcal C(\lambda)^{1/2},
 \qquad
 \|D\|_{\mathrm{ess},\ell^2\to L^2}
 =\|T\|_{\mathrm{ess},L^2\to L^2}
 \asymp\gamma(\lambda)^{1/2}.
\]
For every $1<p<\infty$, there are constants $c_p,C_p>0$, depending
only on $p$ and the fixed geometric data, such that
\[
 \begin{aligned}
 c_p\mathcal C(\lambda)^{1/2}
 &\le \|T\|_{p\to p}
 \le C_p\mathcal C(\lambda)^{1/2},\\
 c_p\gamma(\lambda)^{1/2}
 &\le \|T\|_{\mathrm{ess},L^p\to L^p}
 \le C_p\gamma(\lambda)^{1/2}.
 \end{aligned}
\]
The upper constant $C_p$ can be chosen uniformly when any subset of
the displacements is replaced by zero.  In particular, $T$ is compact
on one, and hence every, space $L^p(\mu)$ in the stated range if and only if
$\gamma(\lambda)=0$.
If
\[
 C_p\gamma(\lambda)^{1/2}<1,
\]
then $I+T$ is invertible on $L^p(\mu)$ and
$\{1\}\cup\{q_g:g\in\mathcal G\}$ is an unconditional basis there.
In particular, vanishing tail packing gives an unconditional basis
for every $1<p<\infty$ and a Riesz basis in $L^2(\mu)$, without
requiring the total packing norm to be small. The basis constants may depend
on the finitely many poles outside a sufficiently small tail.

Finiteness of $\mathcal C(\lambda)$ alone does not imply the
basis property. There is a fixed $s$-regular Cantor probability measure,
for some $0<s<1$, such that for every positive integer $N$, one
can choose admissible real poles for which $\{1\}\cup(q_g)$ is
a frame with synthesis kernel of dimension $N$. There is also
such a frame with infinite-dimensional synthesis kernel. The
packing norms and both frame bounds can be chosen uniformly.
In the finite case, deleting the constant and $N-1$ specified
unshifted gap functions gives a Riesz basis. In particular, for
$N=1$, the gap family $(q_g)$ alone is a Riesz basis of $L^2(\mu)$.
\end{theorem}

\begin{theorem}[Trace ideals for rational gap families]
\label{thm:geometric-schatten}
Let $(b_g)$ be an admissible array of complex poles, and let
$(\alpha_g)$ and $(\beta_g)$ be complex sequences. Define
\[
 r_g(x)=
 \frac{\sqrt{w_g}\,[\alpha_g(x-a_g)+\ell_g\beta_g]}
      {(x-a_g)(x-b_g)},
 \qquad
 \zeta_g=(|\alpha_g|^2+|\beta_g|^2)^{1/2},
\]
and set $R1=0$, $Re_g=r_g$ on the finite canonical span.
For every $0<r<\infty$, this operator extends to $S^r(L^2(\mu))$
if and only if $(\zeta_g)\in\ell^r(\mathcal G)$. Moreover,
\begin{equation}
 \|R\|_{S^r}\asymp_r
       \left(\sum_g\zeta_g^r\right)^{1/r},
 \qquad 0<r<\infty.
 \label{st:rational-schatten-criterion}
\end{equation}
The constants depend only on $r$ and the fixed geometric and
admissibility parameters.
\end{theorem}

For curve transport, the parameter poles $b_g$ are real. Assume that
$\phi\in C^1([0,1];\C)$ satisfies
\[
 |\phi(x)-\phi(y)|\ge m_0|x-y|\quad(m_0>0).
\]
For $\omega_\phi(t)=\sup_{|x-y|\le t}|\phi'(x)-\phi'(y)|$, put
\[
 \mathfrak D_{2,\phi}
   =\left(\sum_{j\ge0}\omega_\phi(2^{-j})^2\right)^{1/2},
 \qquad
 \mathfrak D_\phi=\sum_{j\ge0}\omega_\phi(2^{-j}).
\]
These quantities are allowed to be infinite. Their finiteness is
equivalent to the square Dini and Dini integral conditions,
respectively.
On $\phi(E)$, put $\mu_\phi=\phi_\#\mu$ and define
\[
 Q_g(z)=\frac{\sqrt{w_g}\phi'(b_g)}{z-\phi(b_g)},\qquad
 H_\phi f=\sum_g\langle f,e_g\rangle(Q_g\circ\phi-q_g).
\]
We write $H=H_\phi$ when the parametrization is fixed.
For a bounded operator on a Banach space $X$, its approximation
numbers are
\[
 a_n^X(S)=\inf_{\operatorname{rank}R<n}\|S-R\|_{X\to X}.
\]
We abbreviate $a_n^{L^p(\mu)}$ to $a_n^{(p)}$; on $L^2$ these
are the singular values of a compact operator.

\begin{theorem}[Curve transport and the sharp regularity scale]
\label{thm:curve-transport}
If $\mathfrak D_{2,\phi}<\infty$, the finite sums defining $H$
converge in operator norm to a compact operator on every
$L^p(\mu)$, $1<p<\infty$, and
\[
 \|H\|_{2\to2}\le C\mathfrak D_{2,\phi}/m_0.
\]
The family $(Q_g)$ is then Bessel in $L^2(\mu_\phi)$ exactly when
$\mathcal C(\lambda)<\infty$. Under this packing condition,
\[
 \|T+H\|_{\mathrm{ess},L^p\to L^p}
       \asymp_p\gamma(\lambda)^{1/2},\qquad 1<p<\infty.
\]
On each such space, $I+T+H$ is Fredholm if and only if $I+T$ is
Fredholm, and their indices agree. In particular, the transport for the canonical poles $I+H$ is Fredholm of index zero.

If $\mathfrak D_\phi<\infty$, the convergence and compactness also
hold on $L^1(\mu)$, $L^\infty(\mu)$, and $C(E)$, with
\[
 \|H\|_{p\to p}\le C\mathfrak D_\phi/m_0
 \quad(1\le p\le\infty),\qquad
 \|H\|_{C(E)\to C(E)}\le C\mathfrak D_\phi/m_0.
\]

Fix an integer $m\ge0$, $0<\alpha\le1$, and put $\sigma=m+\alpha$.
If $\phi\in C^{m+1,\alpha}$, then
\[
 a_n^{(p)}(H)\le Cn^{-\sigma/s},
 \qquad n\ge1,\quad1\le p\le\infty,
\]
and the same bound holds on $C(E)$. For the same fixed parameter
poles, under a common lower chord bound and common $C^{m,\alpha}$ bounds
for the derivatives, one also has
\[
 a_n^{(p)}(H_\phi-H_\psi)
 \le C\|\phi'-\psi'\|_{C^{m,\alpha}([0,1])}n^{-\sigma/s},
\]
with the same conclusion on $C(E)$. The constants are independent of $p$.

For every $0<s<1$, integer $m\ge0$, and $0<\alpha\le1$, there
are an $s$-regular Cantor measure, $\varepsilon_0>0$, and,
for each $0<\varepsilon\le\varepsilon_0$, a map
$\phi_\varepsilon\in C^{m+1,\alpha}$ such that
\[
 |\phi_\varepsilon'|=1,\qquad
 |\phi_\varepsilon(x)-\phi_\varepsilon(y)|\ge\tfrac12|x-y|,
 \qquad \|\phi_\varepsilon'-1\|_{C^{m,\alpha}}\le C\varepsilon,
\]
the image $\phi_\varepsilon(E)$ is not contained in any affine line,
and, for the canonical parameter poles $b_g=a_g$,
\[
 a_n^{(p)}(H_{\phi_\varepsilon})
 \asymp\varepsilon n^{-\sigma/s}
 \quad(n\ge1,\ 1\le p\le\infty).
\]
These upper and lower bounds also hold on $C(E)$, with constants and
$\varepsilon_0$ independent of $p$.

For $1<p<\infty$, whenever $H$ is bounded, if $A_b=I+T$ is invertible on $L^p(\mu)$ and
$\|A_b^{-1}H\|_{p\to p}<1$, then $\{1\}\cup(Q_g)$ is an
unconditional basis of $L^p(\mu_\phi)$. Under the Dini condition,
a sufficient condition is
\[
 C_p\mathcal C(\lambda)^{1/2}+C\mathfrak D_\phi/m_0<1.
\]
When $b_g=a_g$ and $C\mathfrak D_\phi/m_0<1$, the family ordered by
gap length is also a Schauder basis of $C(\phi(E))$ and $L^1(\mu_\phi)$.
\end{theorem}

The approximation constants may depend on the regularity data,
$c_0,m_0,m,\alpha$, and the $C^{m,\alpha}$ norm of the derivative.
For $q>0$, the theorem gives
\[
 H\in S^q\quad(q>s/\sigma),
\]
and hence trace class when $\sigma>s$. For the examples above,
$H_{\phi_\varepsilon}\in S^q$ exactly when $q>s/\sigma$.
These are sharp uniform estimates for the stated regularity class;
an individual parametrization can give faster decay, including $H=0$.
Corollary~\ref{cor:smooth-analytic-curves} gives decay faster than
every power for smooth deformations and exponential decay for
holomorphically univalent deformations. Proposition~\ref{prop:mobius-covariance}
shows that M\"obius postcomposition changes $H$ by an operator of
rank at most one and preserves the invertibility of the transport.

For deformations sufficiently close to the identity in $C^1$ that
fix $E$ pointwise, the exact criterion
and the spectral criterion involving both positions and tangents are given
in Proposition~\ref{prop:c1-support-fixing} and
Corollary~\ref{cor:c1-joint-schatten}.
Corollary~\ref{cor:c1-transport-trichotomy} shows that $C^1$
regularity alone does not ensure boundedness or compactness,
and proves a strict separation between interior $L^p$ and $L^1$.
Proposition~\ref{prop:sharp-modulus-transport} makes this
separation sharp for every prescribed concave modulus: the
square Dini integral is the exact interior threshold in its examples,
and the Dini integral is the exact $L^1$ threshold.

We next give the exact classification for finite endpoint chains.
The terms \emph{frame} and \emph{excess} refer to the Hilbert space
synthesis described in Section~\ref{sec:prelim}.

Fix $0<\rho\le1/256$, and let $\mu_\rho$ be the equally weighted measure
on the symmetric Cantor set $E_\rho$ generated by the maps
$t\mapsto\rho t$ and $t\mapsto1-\rho+\rho t$. Put $u(t)=1-\rho+\rho t$ and define
\begin{equation}
 \mathfrak F_\rho(z)=
 \sum_{k\in\mathbb Z}2^{-k-1}
       \int_{E_\rho}\frac{d\mu_\rho(t)}{\rho^k u(t)-z}.
 \label{idx:blowup-transform}
\end{equation}
The series converges locally uniformly away from
$\bigcup_{k\in\mathbb Z}\rho^k u(E_\rho)$ and $0$.
It has a unique zero $a_\rho$ in $(\rho,1-\rho)$, and
\[
 |a_\rho-\tfrac12|\le2\rho.
\]
These facts are proved below.

Let $J_1,\ldots,J_N$ be pairwise disjoint Cantor cylinders. The
central gaps of the successive left descendants $J_j0^n$, $n\ge0$,
form the endpoint chain based at
$J_j$. Write their lengths and canonical poles as
$\ell_{j,n}$ and $a_{j,n}$. Consider complex pole shifts
\begin{equation}
 b_{j,n}=a_{j,n}+\eta_{j,n}\ell_{j,n},\qquad
 |\eta_{j,n}|\le\eta_*\le\frac6{25},\qquad
 \eta_{j,n}\longrightarrow\eta_j,
 \qquad
 \beta_j=a_\rho+(1-2\rho)\eta_j.
 \label{idx:chain-array}
\end{equation}
All other poles are canonical, apart from finitely many additional
admissible replacements. The bound on $\rho$ and
\eqref{dyadic:canonical-estimates} imply
$\dist(a_I,E_\rho)\ge(49/100)\ell_I$.
For this Cantor result we take $c_*=49/100$ and $c_0=6/25$;
then $c_0<c_*/2$, and the shifts in \eqref{idx:chain-array}
satisfy the admissibility condition.

Write
\[
 \mathcal A_b:\C\oplus\ell^2(\mathcal G)\longrightarrow L^2(\mu_\rho),
 \qquad \mathcal A_b(u,c)=u1+\sum_gc_gq_g.
\]
Recall that a bounded operator is Fredholm if its range is closed
and its kernel and cokernel are finite dimensional. Its index is
the dimension of its kernel minus the dimension of its cokernel.

\begin{theorem}[Exact basis, frame, and index classification]
\label{prop:exact-chain-index}
For the pole array \eqref{idx:chain-array}, the family
$\{1\}\cup(q_g)_g$ is complete and Bessel in $L^2(\mu_\rho)$.
The operator $\mathcal A_b$ is Fredholm exactly when
\begin{equation}
 |\beta_j|\ne\sqrt2\,a_\rho
 \qquad(1\le j\le N).
 \label{idx:critical-circles}
\end{equation}
When \eqref{idx:critical-circles} holds, the family is a frame and
\begin{equation}
 \operatorname{ind}\mathcal A_b
 =\dim\ker\mathcal A_b
 =\#\{j:|\beta_j|>\sqrt2\,a_\rho\}.
 \label{idx:index-count}
\end{equation}
In particular, it is a Riesz basis precisely when
$|\beta_j|<\sqrt2\,a_\rho$ for every $j$.
If equality holds for any $j$, its synthesis range is dense and
not closed, so the complete Bessel family is not a frame.
For real constant shifts on one chain, the positive transition is
\begin{equation}
 \eta_{\rm c}(\rho)
 =\frac{(\sqrt2-1)a_\rho}{1-2\rho},
 \qquad
 \eta_{\rm c}(\rho)\longrightarrow\frac{\sqrt2-1}{2}
 \quad\hbox{as }\rho\downarrow0.
 \label{idx:real-threshold}
\end{equation}
\end{theorem}

For reconstruction, define
\[
 \mathcal V_r=\operatorname{span}\bigl(\{1\}\cup
                       \{e_g:\ell_g\ge r\}\bigr),\qquad
 P_rf=\langle f,1\rangle1+
             \sum_{\ell_g\ge r}\langle f,e_g\rangle e_g,
 \quad0<r\le1,
\]
and let $E_r(f)=\inf_{v\in\mathcal V_r}\|f-v\|_{C(E)}$.
Write
\[
 V_\rho(Pf)=
 \sup_{\substack{M\ge1\\1\ge r_0>\cdots>r_M>0}}
 \left(\sum_{j=1}^M|P_{r_j}f-P_{r_{j-1}}f|^\rho\right)^{1/\rho}.
\]
For $0<\alpha\le1$, put
\[
 \Psi_{\alpha,s}(r)=
 \begin{cases}
 r^\alpha,&\alpha<1-s,\\
 r^{1-s}\log(e/r),&\alpha=1-s,\\
 r^{1-s},&\alpha>1-s.
 \end{cases}
\]

\begin{proposition}[Ordered reconstruction and sharp approximation]
\label{thm:approximation-norming}
The projections \(P_r\) are uniformly bounded on every \(L^p(\mu)\),
\(1\le p\le\infty\). As \(r\downarrow0\), they converge to the
identity strongly on \(L^p(\mu)\), \(1\le p<\infty\), uniformly
on \(E\) for every \(f\in C(E)\), and almost everywhere for every
\(f\in L^1(\mu)\).

Let \(g_1,g_2,\ldots\) be any enumeration by nonincreasing gap
length, with arbitrary ordering among equal lengths, and set
\[
 Q_0f=\langle f,1\rangle1,\qquad
 Q_Nf=Q_0f+\sum_{j=1}^N\langle f,e_{g_j}\rangle e_{g_j}.
\]
The ordered family \(\{1,e_{g_1},e_{g_2},\ldots\}\) is a
Schauder basis of \(C(E)\) and \(L^1(\mu)\).
The sums \(Q_Nf\) also converge almost everywhere for every
\(f\in L^1(\mu)\).
For every \(\rho>2\), the variation over the full sequence satisfies
\[
 \sup_{\tau>0}\tau\mu\{V_\rho(Qf)>\tau\}
 \le C_\rho\|f\|_1,\qquad
 \|V_\rho(Qf)\|_p\le C_{p,\rho}\|f\|_p
 \quad(1<p<\infty).
\]
Here the variation is over all finite increasing sequences of
indices, and the constants are independent of the ordering among
equal lengths. The same inequalities hold for \(V_\rho(Pf)\).

Moreover, \(\dim\mathcal V_r\le1+Cr^{-s}\) and
\[
 E_r(f)\asymp\|f-P_rf\|_{C(E)}
 \quad(f\in C(E)),\qquad
 \sup_{[f]_{C^\alpha(E)}\le1}E_r(f)
 \asymp_\alpha\Psi_{\alpha,s}(r)
 \quad(r\downarrow0).
\]
\end{proposition}

These approximation estimates concern the fixed spaces $\mathcal V_r$;
spaces of comparable dimension can admit faster approximation.
The proof of the reconstruction proposition also gives the following Remez inequality:
if $S\subset E$ is Borel and $\mu(S)\ge\theta>0$, then
\[
 \|v\|_{C(E)}
 \le \exp\!\left(Cr^{-s}\log\frac{C}{\theta}\right)
       \esssup_{x\in S}|v(x)|,
 \qquad v\in\mathcal V_r.
\]
The constant depends only on the regularity data.

\section{The canonical basis and its multipliers}
\label{subsec:basis-proof}
\begin{proof}[Proof of Proposition~\ref{thm:basis}]
Write $F=C_\mu$ and $m_j=\int_Et^j\,d\mu(t)$. Every constant in
this proof depends only on $A,s$, unless stated otherwise.
The preliminary estimates show that $\mu$ has no atoms, that $E$ is
perfect and has Lebesgue measure zero, and that the covering bound
\eqref{sp:covering} holds.

Fix $g=(b,c)$ and $a\in g$. Put
\[
 p=a-b,\qquad q=c-a,
 \qquad p+q=\ell_g,
\]
and write the Cauchy integral in terms of
\[
 L(a)=\int_{E\cap(-\infty,b]}\frac{d\mu(t)}{a-t},
 \qquad
 R(a)=\int_{E\cap[c,\infty)}\frac{d\mu(t)}{t-a}.
\]
The ball $B(b,p)$ meets $E$ only to the left of $b$. On this ball,
$a-t\le2p$, so
\[
 L(a)\ge\frac{\mu(B(b,p))}{2p}
 \ge\frac1{2A}p^{s-1}.
\]
For $k\ge0$, the set
\[
 A_k=\{t\in E:t\le b,\ 2^kp\le a-t<2^{k+1}p\}
\]
lies in $B(b,2^{k+1}p)$. Therefore
\[
 \begin{aligned}
 L(a)
 &=\sum_{k\ge0}\int_{A_k}\frac{d\mu(t)}{a-t}\\
 &\le\sum_{k\ge0}(2^kp)^{-1}\mu(B(b,2^{k+1}p))\\
 &\le A2^s p^{s-1}\sum_{k\ge0}2^{-k(1-s)}
 =\frac{A2^s}{1-2^{-(1-s)}}p^{s-1}.
 \end{aligned}
\]
The upper mass bound remains valid when the radius exceeds one,
since $\mu(E)=1\le Ar^s$ for $r\ge1$.
The same argument applies at $c$. With fixed positive constants
$c_1,C_1$, we obtain
\begin{equation}
 c_1p^{s-1}\le L(a)\le C_1p^{s-1},
 \qquad
 c_1q^{s-1}\le R(a)\le C_1q^{s-1}.
 \label{sp:one-sided}
\end{equation}
On the gap, differentiation under the integral gives
\[
 F(a)=R(a)-L(a),
 \qquad
 F'(a)=\int_E\frac{d\mu(t)}{(t-a)^2}>0.
\]
Because $s-1<0$, the endpoint limits are
\[
 \lim_{a\downarrow b}F(a)=-\infty,
 \qquad
 \lim_{a\uparrow c}F(a)=+\infty.
\]
Continuity and strict monotonicity show that there is exactly one zero $a_g$.
At this zero, $L(a_g)=R(a_g)$. Writing $p=a_g-b$ and $q=c-a_g$,
\eqref{sp:one-sided} gives
\[
 \frac{c_1}{C_1}\le\frac{p^{s-1}}{q^{s-1}}
 \le\frac{C_1}{c_1}.
\]
Thus, with $B=(C_1/c_1)^{1/(1-s)}$,
\[
 B^{-1}\le\frac pq\le B,
 \qquad
 \min(p,q)\ge\frac{\ell_g}{1+B}.
\]
Restricting the derivative integral to $B(b,p)$ gives
\[
 F'(a_g)\ge\frac{\mu(B(b,p))}{4p^2}
 \ge\frac1{4A}p^{s-2}
 \ge c\ell_g^{s-2}.
\]
For the reverse estimate, set $\rho_g=\dist(a_g,E)\asymp\ell_g$.
The upper mass bound for balls centered outside $E$ gives
\[
 \begin{aligned}
 F'(a_g)
 &\le\sum_{k\ge0}(2^k\rho_g)^{-2}
                \mu(B(a_g,2^{k+1}\rho_g))\\
 &\le C\rho_g^{s-2}\sum_{k\ge0}2^{-k(2-s)}
 \le C\ell_g^{s-2}.
 \end{aligned}
\]
Consequently,
\begin{equation}
 c\ell_g\le\min(a_g-u_g,v_g-a_g),
 \qquad
 c\ell_g^{2-s}\le w_g:=F'(a_g)^{-1}\le C\ell_g^{2-s}.
 \label{sp:gap-bounds}
\end{equation}

We construct the reciprocal measure by compressing the coordinate
multiplication operator. Let
\[
 Mf(t)=tf(t),\qquad
 Pf=f-\langle f,1\rangle1,
 \qquad \mathcal K=P L^2(\mu)=\{1\}^{\perp}.
\]
Set
\[
 T=PM|_{\mathcal K},\qquad h=PM1=t-m_1.
\]
For $u,v\in\mathcal K$, self-adjointness follows from
\[
 \langle Tu,v\rangle
 =\langle Mu,v\rangle
 =\langle u,Mv\rangle
 =\langle u,Tv\rangle.
\]
Moreover,
\[
 0\le\langle Tu,u\rangle
 =\int_Et|u(t)|^2\,d\mu(t)
 \le\|u\|_2^2.
\]
Thus $0\le T\le I$ and $\sigma(T)\subset[0,1]$.
For $c\in\C$ and $u\in\mathcal K$,
\[
 \begin{aligned}
 M(c1+u)
 &=\bigl(m_1c+\langle u,h\rangle\bigr)1+ch+Tu,\\
 M&=
 \begin{pmatrix}
 m_1&\langle\,\cdot\,,h\rangle\\
 h&T
 \end{pmatrix}_{\C1\oplus\mathcal K}.
 \end{aligned}
\]
Take $z\in\C\setminus[0,1]$ and solve $(M-z)(c1+u)=1$.
The two component equations are
\[
 (m_1-z)c+\langle u,h\rangle=1,
 \qquad ch+(T-z)u=0.
\]
Substitution of $u=-c(T-z)^{-1}h$ gives
\[
 c\bigl(m_1-z-\langle(T-z)^{-1}h,h\rangle\bigr)=1.
\]
Since $(M-z)^{-1}1=(t-z)^{-1}$, the constant component is
\[
 c=\langle(M-z)^{-1}1,1\rangle=F(z).
\]
It follows that
\begin{equation}
 -\frac1{F(z)}=z-m_1+\langle(T-z)^{-1}h,h\rangle.
 \label{sp:schur}
\end{equation}
Let $\mathsf E_T$ be the spectral resolution of $T$, and define
$\nu(J)=\langle\mathsf E_T(J)h,h\rangle$. Then
\begin{equation}
 \langle(T-z)^{-1}h,h\rangle
 =\int_{[0,1]}\frac{d\nu(a)}{a-z},
 \qquad
 \nu([0,1])=\|h\|_2^2=m_2-m_1^2.
 \label{sp:spectral-measure}
\end{equation}
In particular, $\nu$ is finite and positive.
For $x\in E$ and $0<y\le1$, positivity gives
\[
 \begin{aligned}
 \Im F(x+iy)
 &=\int_E\frac{y\,d\mu(t)}{(t-x)^2+y^2}\\
 &\ge\frac{\mu(B(x,y))}{2y}
 \ge\frac1{2A}y^{s-1}.
 \end{aligned}
\]
For the corresponding upper estimate, use the ball of radius $y$
and the annuli of radii $2^ky$:
\[
 \begin{aligned}
 |F(x+iy)|
 &\le\frac{\mu(B(x,y))}{y}
   +\sum_{k\ge0}(2^ky)^{-1}\mu(B(x,2^{k+1}y))\\
 &\le Ay^{s-1}+A2^sy^{s-1}\sum_{k\ge0}2^{-k(1-s)}
 \le Cy^{s-1}.
 \end{aligned}
\]
Hence
\[
cy^{s-1}\le\Im F(x+iy)\le|F(x+iy)|\le Cy^{s-1}.
\]
Put $H=-1/F$. The elementary identities
\[
 |H|=|F|^{-1},\qquad \Im H=\frac{\Im F}{|F|^2}
\]
give
\[
|H(x+iy)|\le Cy^{1-s},\qquad
 cy^{1-s}\le\Im H(x+iy)\le Cy^{1-s}.
\]
For example, the bounds for the imaginary part follow from
\[
 \frac{cy^{s-1}}{C^2y^{2s-2}}
 \le\frac{\Im F}{|F|^2}
 \le\frac1{\Im F}
 \le Cy^{1-s}.
\]
The representations \eqref{sp:schur} and \eqref{sp:spectral-measure}
also give
\[
 \Im H(x+iy)
 =y+\int\frac{y\,d\nu(a)}{(a-x)^2+y^2}
 \ge\frac{\nu(B(x,y))}{2y}.
\]
Therefore
\[
 \nu(B(x,y))\le Cy^{2-s},\qquad x\in E,\quad0<y\le1.
\]
Cover $E$ by the balls in \eqref{sp:covering}. We obtain
\[
 0\le\nu(E)
 \le\sum_{j=1}^{N(E,r)}\nu(B(x_j,r))
 \le Cr^{-s}r^{2-s}=Cr^{2-2s}\longrightarrow0.
\]
Thus the reciprocal measure gives no mass to $E$.
Fix a closed interval $[c,d]\subset g\setminus\{a_g\}$. The
function $H$ extends analytically across this interval and is real
there. Taylor expansion on a slightly larger compact interval gives
\[
 \sup_{c\le x\le d}|\Im H(x+iy)|=O(y).
\]
On the other hand, Tonelli's theorem gives
\[
 \begin{aligned}
 \int_c^d\bigl(\Im H(x+iy)-y\bigr)\,dx
 &=\int\left(\int_c^d
       \frac{y\,dx}{(a-x)^2+y^2}\right)d\nu(a)\\
 &=\int\left[
       \arctan\frac{d-a}{y}-\arctan\frac{c-a}{y}
          \right]d\nu(a).
 \end{aligned}
\]
The expression in brackets lies in $[0,\pi]$ and tends to $\pi$ on $(c,d)$.
Dominated convergence therefore yields
\[
 0=\pi\nu((c,d))+\frac\pi2\nu(\{c,d\}).
\]
Countably many such intervals cover $g\setminus\{a_g\}$, so
$\nu(g\setminus\{a_g\})=0$.
We can also compute the mass at the zero directly. From the
spectral representation,
\[
 y\Im H(a_g+iy)
 =y^2+\int\frac{y^2}{(a-a_g)^2+y^2}\,d\nu(a)
 \longrightarrow\nu(\{a_g\}).
\]
Since the zero of $F$ is simple,
\[
 F(a_g+iy)=iyF'(a_g)+O(y^2),
 \qquad
 H(a_g+iy)=\frac{i}{yF'(a_g)}+O(1).
\]
Consequently,
\[
 \nu(\{a_g\})=F'(a_g)^{-1}=w_g.
\]
The interval $[0,1]$ is the disjoint union of $E$ and its bounded
gaps. We have proved
\begin{equation}
 -\frac1{F(z)}
 =z-m_1+\sum_{g\in\G}\frac{w_g}{a_g-z},
 \qquad z\in\C\setminus[0,1].
 \label{sp:reciprocal}
\end{equation}
The mass and growth identities are
\begin{equation}
 \sum_gw_g=m_2-m_1^2,
 \qquad
 \nu(B(x,r))\le Cr^{2-s}\quad(x\in E, r>0).
 \label{eq:dual-growth}
\end{equation}
For $r\ge1$, the last bound follows from $\nu([0,1])\le1$.
The series in \eqref{sp:reciprocal} converges locally uniformly.
Indeed, if $K\subset\C\setminus[0,1]$ is compact and
$\delta_K=\dist(K,[0,1])>0$, then
\[
 \sup_{z\in K}\sum_{g\notin\mathcal F}
       \left|\frac{w_g}{a_g-z}\right|
 \le\delta_K^{-1}\sum_{g\notin\mathcal F}w_g
 \longrightarrow0
\]
as the finite families $\mathcal F$ increase to all the gaps.
Define $e_0=1$ and $e_g(t)=\sqrt{w_g}/(t-a_g)$. We have
\[
 \langle e_g,1\rangle=\sqrt{w_g}F(a_g)=0,
 \qquad
 \|e_g\|_2^2=w_gF'(a_g)=1.
\]
For $g\ne g'$, partial fractions give
\[
 \begin{aligned}
 \langle e_g,e_{g'}\rangle
 &=\sqrt{w_gw_{g'}}
    \int_E\frac{d\mu(t)}{(t-a_g)(t-a_{g'})}\\
 &=\frac{\sqrt{w_gw_{g'}}}{a_g-a_{g'}}
      \bigl(F(a_g)-F(a_{g'})\bigr)=0.
 \end{aligned}
\]
Thus the family is orthonormal.
To establish completeness, let $v_n=PM^n1$ for $n\ge1$. Then
\[
 M^n1=m_n1+v_n,
 \qquad v_1=h,
 \qquad v_{n+1}=Tv_n+m_nh.
\]
Induction gives
\[
 v_n=T^{n-1}h+m_1T^{n-2}h+\cdots+m_{n-1}h.
\]
In particular,
\[
 \operatorname{span}(v_1,\ldots,v_n)
 =\operatorname{span}(h,Th,\ldots,T^{n-1}h).
\]
If $u\in\mathcal K$ is orthogonal to every $v_n$, then
\[
 \langle u,t^n\rangle
 =\langle u,v_n+m_n1\rangle=0
 \qquad(n\ge0).
\]
Polynomials are dense in $L^2(\mu)$, so $u=0$. Hence $h$ is
cyclic for $T$.

For a polynomial $p$, define $\mathcal U p=p(T)h$. The spectral
theorem gives
\[
 \|\mathcal U p\|_2^2=\int|p(a)|^2\,d\nu(a),
 \qquad
 T\mathcal U p=\mathcal U(ap).
\]
The density of polynomials in $L^2(\nu)$, together with cyclicity, extends
$\mathcal U$ to a unitary map from $L^2(\nu)$ onto $\mathcal K$.
Because $\nu=\sum_gw_g\delta_{a_g}$, the functions
\[
 \varphi_g=w_g^{-1/2}\mathbf1_{\{a_g\}}
\]
form an orthonormal basis of $L^2(\nu)$. They satisfy
\[
 T\mathcal U\varphi_g=a_g\mathcal U\varphi_g.
\]
If $Tu=a_gu$ with $u\in\mathcal K$, then
\[
 P((t-a_g)u)=0,
 \qquad (t-a_g)u=\text{constant},
\]
so this eigenspace is spanned by $e_g$.
The phase is fixed by
\[
 \begin{aligned}
 \langle h,e_g\rangle
 &=\sqrt{w_g}\int_E\frac{t-m_1}{t-a_g}\,d\mu(t)\\
 &=\sqrt{w_g}\bigl(1+(a_g-m_1)F(a_g)\bigr)=\sqrt{w_g},\\
 \langle h,\mathcal U\varphi_g\rangle
 &=\langle1,\varphi_g\rangle_{L^2(\nu)}=\sqrt{w_g}.
 \end{aligned}
\]
Thus $\mathcal U\varphi_g=e_g$. This proves completeness and Parseval's identity:
\[
\|f\|_2^2
 =\left|\int_Ef\,d\mu\right|^2
  +\sum_gw_g\left|\int_E\frac{f(t)}{t-a_g}\,d\mu(t)\right|^2.
\]
The closed support of $\nu$ may have accumulation points in
$E$. We have proved a complete pure point spectral resolution,
which does not require the closed spectrum to be discrete.

Choose either endpoint $\xi_g$ of $g$. The lower bounds for the
distances from $a_g$ to the gap endpoints give
\[
 |x-a_g|\asymp\ell_g+|x-\xi_g|,
 \qquad\sqrt{w_g}\asymp\ell_g^{1-s/2}.
\]
Consequently,
\begin{equation}
 |e_g(x)|\le C\ell_g^{-s/2}
   \left(1+\frac{|x-\xi_g|}{\ell_g}\right)^{-1}.
 \label{sp:size}
\end{equation}
The exact difference formula is
\[
 e_g(x)-e_g(x')
 =\frac{\sqrt{w_g}(x'-x)}{(x-a_g)(x'-a_g)}.
\]
For $1\le p<\infty$, we have $p>s$. Annular summation gives
\[
 \begin{aligned}
 \int_E|x-a_g|^{-p}\,d\mu(x)
 &\le C\ell_g^{s-p}\sum_{k\ge0}2^{-k(p-s)}\\
 &=\frac{C}{1-2^{-(p-s)}}\ell_g^{s-p}.
 \end{aligned}
\]
On $B_E(\xi_g,\ell_g)$, the distance $|x-a_g|$ is at most
$2\ell_g$. Hence
\[
 \int_E|x-a_g|^{-p}\,d\mu(x)
 \ge(2\ell_g)^{-p}\mu(B_E(\xi_g,\ell_g))
 \ge c_p\ell_g^{s-p}.
\]
Multiplying by $w_g^{p/2}$ and taking the $p$th root gives
\[
 \|e_g\|_p\asymp_p
 \ell_g^{1-s/2}\ell_g^{s/p-1}
 =\ell_g^{s(1/p-1/2)}.
\]
The pointwise bound \eqref{sp:size} and evaluation at $\xi_g$
similarly give $\|e_g\|_\infty\asymp\ell_g^{-s/2}$.
Continuity on $E$ and the full support of $\mu$ show that the supremum
equals the essential supremum.

We next derive three integral estimates from \eqref{eq:dual-growth}.
For $x\in E$ and $R>0$,
\begin{align}
 \int_{|a-x|<R}|a-x|^{-1}\,d\nu(a)&\le CR^{1-s},
 \label{eq:dual-near-one}\\
 \int_{|a-x|\ge R}|a-x|^{-2}\,d\nu(a)&\le CR^{-s},
 \label{eq:dual-far-two}\\
 \int_{|a-x|\ge R}|a-x|^{-3}\,d\nu(a)&\le CR^{-1-s}.
 \label{eq:dual-far-three}
\end{align}
For the first estimate, sum over annuli approaching the center:
\[
 \begin{aligned}
 \int_{0<|a-x|<R}|a-x|^{-1}\,d\nu(a)
 &\le\sum_{j\ge0}\frac{2^{j+1}}R\nu(B(x,2^{-j}R))\\
 &\le CR^{1-s}\sum_{j\ge0}2^{-j(1-s)}
 \le CR^{1-s}.
 \end{aligned}
\]
There is no atom at $x$, since $\nu(E)=0$.
For $k=2,3$, summing over annuli away from the center gives
\[
 \begin{aligned}
 \int_{|a-x|\ge R}|a-x|^{-k}\,d\nu(a)
 &\le\sum_{j\ge0}(2^jR)^{-k}\nu(B(x,2^{j+1}R))\\
 &\le CR^{2-s-k}\sum_{j\ge0}2^{-j(k+s-2)}.
 \end{aligned}
\]
The exponents $k+s-2$ are $s$ and $1+s$, both positive.
This proves the other two estimates.

Fix $|\beta_0|\le1$ and $|\beta_g|\le1$. For a finite family
$\mathcal F\subset\G$, the kernel of the part with zero mean is
\[
 K_{\beta,\mathcal F}(x,y)
 =\sum_{g\in\mathcal F}\frac{\beta_gw_g}{(x-a_g)(y-a_g)}.
\]
The absolute value of every such sum is bounded by the corresponding
integral against $\nu$. All estimates below are therefore uniform in
$\mathcal F$ and the coefficients.
Let $r=|x-y|>0$. On $|a-x|<r/3$, the distance to $y$ is at least
$2r/3$. Thus
\[
 \begin{aligned}
 \int_{|a-x|<r/3}\frac{d\nu(a)}{|x-a||y-a|}
 &\le\frac{3}{2r}\int_{|a-x|<r/3}|x-a|^{-1}\,d\nu(a)\\
 &\le Cr^{-1}r^{1-s}=Cr^{-s}.
 \end{aligned}
\]
Interchanging $x,y$ proves the same estimate on $|a-y|<r/3$.
On the remaining region, both distances are at least $r/3$, and
\[
 |a-x|\le|a-y|+r\le4|a-y|,
 \qquad |a-y|\le4|a-x|.
\]
Therefore
\[
 \int_{\substack{|a-x|\ge r/3\\|a-y|\ge r/3}}
       \frac{d\nu(a)}{|x-a||y-a|}
 \le4\int_{|a-x|\ge r/3}|x-a|^{-2}\,d\nu(a)
 \le Cr^{-s}.
\]
The full series converges absolutely off the diagonal and satisfies
\[
|K_\beta(x,y)|\le C|x-y|^{-s}.
\]

For the continuity estimate, put $h=|x-x'|$. Assume first that
$0<h\le r/8$, and let $N=\{a:|a-x|<2h\}$. On $N$,
\[
 |a-y|\ge r-2h\ge3r/4,
 \qquad |a-x'|\le|a-x|+h<3h.
\]
Consequently,
\[
 \begin{aligned}
 &\int_N\left|
 \frac1{(x-a)(y-a)}-\frac1{(x'-a)(y-a)}\right|d\nu(a)\\
 &\quad\le\frac C r\left(
    \int_{|a-x|<2h}|x-a|^{-1}\,d\nu(a)
   +\int_{|a-x'|<3h}|x'-a|^{-1}\,d\nu(a)\right)\\
 &\quad\le C\frac{h^{1-s}}r.
 \end{aligned}
\]
Outside $N$, the distances from $a$ to $x,x'$ are comparable:
\[
 \tfrac12|a-x|\le|a-x'|\le\tfrac32|a-x|.
\]
The exact rational difference therefore gives
\[
 \left|
 \frac1{(x-a)(y-a)}-\frac1{(x'-a)(y-a)}\right|
 \le\frac{2h}{|a-x|^2|a-y|}.
\]
Partition this region into the disjoint sets
\[
 \begin{aligned}
 U&=\{a:2h\le|a-x|<r/3\},\\
 V&=\{a:|a-y|<r/3\},\\
 W&=\{a:|a-x|\ge r/3, |a-y|\ge r/3\}.
 \end{aligned}
\]
On $U$, we have $|a-y|\ge2r/3$. Using
\eqref{eq:dual-far-two},
\[
 \begin{aligned}
 \int_U\frac{2h\,d\nu(a)}{|a-x|^2|a-y|}
 &\le C\frac h r\int_{|a-x|\ge2h}|a-x|^{-2}\,d\nu(a)\\
 &\le C\frac h r h^{-s}=C\frac{h^{1-s}}r.
 \end{aligned}
\]
On $V$, the inequality $|a-x|\ge2r/3$ and
\eqref{eq:dual-near-one} give
\[
 \begin{aligned}
 \int_V\frac{2h\,d\nu(a)}{|a-x|^2|a-y|}
 &\le C\frac h{r^2}\int_{|a-y|<r/3}|a-y|^{-1}\,d\nu(a)\\
 &\le C\frac h{r^2}r^{1-s}=Chr^{-1-s}.
 \end{aligned}
\]
On $W$, comparability of the distances and
\eqref{eq:dual-far-three} give
\[
 \int_W\frac{2h\,d\nu(a)}{|a-x|^2|a-y|}
 \le Ch\int_{|a-x|\ge r/3}|a-x|^{-3}\,d\nu(a)
 \le Chr^{-1-s}.
\]
The last two bounds imply the required estimate because
\[
 hr^{-1-s}
 =\left(\frac hr\right)^s\frac{h^{1-s}}r
 \le\frac{h^{1-s}}r.
\]
If $r/8<h\le r/2$, then $|x'-y|\ge r/2$. The size estimate
alone gives
\[
 |K_\beta(x,y)-K_\beta(x',y)|
 \le Cr^{-s}
 \le C8^{1-s}\frac{h^{1-s}}r.
\]
The case $h=0$ is immediate. We have proved
\[
|K_\beta(x,y)-K_\beta(x',y)|
 \le C\frac{|x-x'|^{1-s}}{|x-y|},
 \qquad |x-x'|\le\tfrac12|x-y|.
\]
The kernel is symmetric in $x,y$, so the same estimate holds in
the second variable. Both estimates also hold for finite truncations.
The preceding estimates were obtained after taking absolute values
term by term. They therefore give the stronger bounds
\begin{align}
 \sum_{g\in\mathcal G}|e_g(x)e_g(y)|
 &\le C|x-y|^{-s},
 \label{eq:absolute-kernel-size}\\
 \sum_{g\in\mathcal G}
       |e_g(x)-e_g(x')|\,|e_g(y)|
 &\le C\frac{|x-x'|^{1-s}}{|x-y|},
 \label{eq:absolute-kernel-smoothness}
\end{align}
for $x,y,x'\in E$, $x\ne y$, and
$|x-x'|\le\tfrac12|x-y|$. By interchanging the variables, the
same estimate holds for differences in the second variable.
Both inequalities remain valid when the sum is restricted to an
arbitrary subfamily of gaps.

For the finite multiplier including the constant coordinate, Parseval gives
\[
 \begin{aligned}
 \|T_{\beta,\mathcal F}f\|_2^2
 &=|\beta_0|^2|\langle f,1\rangle|^2
   +\sum_{g\in\mathcal F}|\beta_g|^2|\langle f,e_g\rangle|^2\\
 &\le\|f\|_2^2.
 \end{aligned}
\]
Its kernel is $\beta_0+K_{\beta,\mathcal F}$. Since
$|x-y|\le1$ on $E$, the constant term satisfies the same size estimate
and does not affect the difference estimates.
Apply \eqref{eq:prelim-cz-conclusion} with dimension parameter $s$
and H\"older exponent $1-s$. The finite multipliers satisfy
\[
 \|T_{\beta,\mathcal F}f\|_p\le C_p\|f\|_p,
 \qquad
 \mu\{|T_{\beta,\mathcal F}f|>\tau\}
 \le\frac C\tau\|f\|_1.
\]
The constants are independent of $\mathcal F$ and $\beta$.

Let $\mathcal F_N\uparrow\G$. On $L^2$, orthogonality implies
\[
 \|T_{\beta,\mathcal F_N}f-T_\beta f\|_2^2
 =\sum_{g\notin\mathcal F_N}|\beta_g|^2
                           |\langle f,e_g\rangle|^2
 \longrightarrow0.
\]
Choose a subsequence that converges almost everywhere and retain the
index $N$ for it. For $f\in L^p\cap L^2$, Fatou's lemma gives
\[
\|T_\beta f\|_p\le C_p\|f\|_p,
 \qquad1<p<\infty.
\]
Density gives the bounded extension to all of $L^p$.
For the weak estimate, fix $0<\varepsilon<\tau$. Almost everywhere
convergence gives
\[
 \begin{aligned}
 \mu\{|T_\beta f|>\tau\}
 &\le\liminf_N
       \mu\{|T_{\beta,\mathcal F_N}f|>\tau-\varepsilon\}\\
 &\le\frac C{\tau-\varepsilon}\|f\|_1.
 \end{aligned}
\]
Letting $\varepsilon\downarrow0$ proves the weak bound on
$L^1\cap L^2$.
For general $f\in L^1$, choose bounded $f_j\to f$ in $L^1$.
The inequality
\[
 \mu\{|T_\beta f_j-T_\beta f_k|>\tau\}
 \le\frac C\tau\|f_j-f_k\|_1
\]
shows that the images are Cauchy in measure. Their limit defines
$T_\beta f$, independently of the approximating sequence.
Passing to an almost everywhere convergent subsequence gives
\[
 \mu\{|T_\beta f|>\tau\}\le\frac C\tau\|f\|_1.
\]
Linearity and agreement with the $L^2$ operator follow by taking
limits in measure. This establishes the weak $(1,1)$ extension.

For $f\in L^2$ and $x$ outside its closed essential support,
let $\delta=\dist(x,\operatorname{ess\,supp}f)>0$. The uniform bound
\[
 |K_{\beta,\mathcal F_N}(x,y)f(y)|
 \le C\delta^{-s}|f(y)|,
 \qquad y\in\operatorname{ess\,supp}f,
\]
permits dominated convergence in the finite integral formulas.
Comparison with the $L^2$ limit identifies $K_\beta$ as the
kernel of the part with zero mean, away from the diagonal.

Fix an enumeration of the basis beginning with $e_0=1$, and let
$P_N$ project onto its first $N$ coordinates. The multiplier estimate
gives $\sup_N\|P_N\|_{p\to p}<\infty$.
For bounded $f$ and $1<p\le2$, the probability normalization gives
\[
 \|P_Nf-f\|_p\le\|P_Nf-f\|_2\longrightarrow0.
\]
For $p>2$, choose $q>p$ and
$\vartheta\in(0,1)$ with $1/p=\vartheta/2+(1-\vartheta)/q$.
Then
\[
 \begin{aligned}
 \|P_Nf-f\|_p
 &\le\|P_Nf-f\|_2^\vartheta
       \|P_Nf-f\|_q^{1-\vartheta}\\
 &\le\|P_Nf-f\|_2^\vartheta
       \bigl((C_q+1)\|f\|_q\bigr)^{1-\vartheta}
 \longrightarrow0.
 \end{aligned}
\]
For arbitrary $f\in L^p$, approximate by a bounded function $b$.
The estimate
\[
 \|P_Nf-f\|_p
 \le(C_p+1)\|f-b\|_p+\|P_Nb-b\|_p
\]
proves convergence by first letting $N\to\infty$ and then
$\|f-b\|_p\to0$.
The convergence is unconditional. If a finite set $J$ contains the
first $N$ coordinates, then
\[
 f-P_Jf=(I-P_J)(f-P_Nf),
 \qquad
 \|f-P_Jf\|_p\le(1+C_p)\|f-P_Nf\|_p.
\]
Thus the net of all finite coordinate sums converges to $f$.
Each coefficient functional is continuous, since
\[
 |\langle f,e_g\rangle|
 \le\|f\|_p\|e_g\|_{p'}<\infty.
\]
This proves the unconditional basis assertion in every enumeration.

For a finite expansion, write $c_0=\langle f,1\rangle$,
$c_g=\langle f,e_g\rangle$, and take independent Rademacher signs
$\varepsilon_0,\varepsilon_g$. Khintchine's inequality gives
\[
 \mathbb E\left|\varepsilon_0c_0+
                   \sum_g\varepsilon_gc_ge_g(x)\right|^p
 \asymp_p\left(|c_0|^2+\sum_g|c_g|^2|e_g(x)|^2\right)^{p/2}.
\]
Integrating in $x$ yields
\[
 \mathbb E\|T_\varepsilon f\|_p^p\asymp_p\|Sf\|_p^p,
 \qquad
 Sf=\left(|c_0|^2+\sum_g|c_g|^2|e_g|^2\right)^{1/2}.
\]
The upper multiplier bound and the identity $T_\varepsilon^2f=f$ give
\[
 C_p^{-1}\|f\|_p
 \le\|T_\varepsilon f\|_p
 \le C_p\|f\|_p.
\]
Averaging proves the upper and lower square function estimates for finite
expansions. For general $f$, use
\[
 P_Nf\to f\text{ in }L^p,
 \qquad S(P_Nf)\uparrow Sf,
\]
and monotone convergence to conclude
\[
 c_p\|f\|_p\le\|Sf\|_p\le C_p\|f\|_p.
\]

Finally, the exponent $1-s$ is sharp uniformly over the multipliers.
There must be infinitely many bounded gaps. Otherwise $E$ would
be a finite union of closed intervals and points; its zero Lebesgue
measure and absence of isolated points rule this out.
Choose distinct gaps $g_n=(u_n,v_n)$ with $\ell_{g_n}\to0$.
After passing to a subsequence,
\[
 u_n,v_n,a_{g_n}\longrightarrow x_*\in E.
\]
Fix $y\in E\setminus\{x_*\}$. For the multiplier selecting only
the coordinate $g_n$, the exact endpoint difference is
\[
 \begin{aligned}
 |K_{g_n}(u_n,y)-K_{g_n}(v_n,y)|
 &=\frac{w_{g_n}\ell_{g_n}}
 {|u_n-a_{g_n}|\,|v_n-a_{g_n}|\,|y-a_{g_n}|}\\
 &\asymp\ell_{g_n}^{1-s}.
 \end{aligned}
\]
Here $|u_n-y|$ is bounded above and bounded away from zero, and
$\ell_{g_n}\le|u_n-y|/2$ for all sufficiently large $n$.
A uniform continuity estimate with exponent $\eta>1-s$ would imply
\[
 c\ell_{g_n}^{1-s}
 \le C\frac{\ell_{g_n}^\eta}{|u_n-y|^{s+\eta}},
 \qquad
 \ell_{g_n}^{1-s-\eta}\le C',
\]
which is impossible as $n\to\infty$. This proves the asserted
sharpness and completes the proof.
\end{proof}
\section{Pole perturbations}
\label{subsec:stability-proof}
\subsection{Packing and operator bounds}
\begin{proof}[Norm estimates in Theorem~\ref{thm:stability}]
All constants in this proof depend only on $s$, the regularity data in
\eqref{sp:regularity}, and the fixed relative displacement bound $c_0$,
unless otherwise indicated. We write
$B_E(x,r)=E\cap B(x,r)$.  The gap estimates \eqref{sp:gap-bounds}
give
\[
 \operatorname{dist}(a_g,E)\geq c_*\ell_g,
 \qquad w_g\asymp\ell_g^{2-s}.
\]
The restriction on the displacement gives
\[
 \operatorname{dist}(b_g,E)
 \geq\operatorname{dist}(a_g,E)-|b_g-a_g|
 \geq(c_*-c_0)\ell_g\geq\tfrac12c_*\ell_g.
\]
For either endpoint $\xi_g$, the gap estimates for $a_g$ and
the displacement bound give
\[
 |b_g-\xi_g|\leq |a_g-\xi_g|+|b_g-a_g|
 \leq(1-c_*+c_0)\ell_g\leq\ell_g.
\]
These distance estimates also hold for complex $b_g$.
For $p\in\{a_g,b_g\}$ we therefore have
\[
 |x-p|\leq |x-\xi_g|+\ell_g,
 \qquad
 |x-\xi_g|+\ell_g\leq |x-p|+2\ell_g
                      \leq(1+4/c_*)|x-p|.
\]
Hence
\begin{equation}
 |x-a_g|\asymp |x-b_g|\asymp \ell_g+|x-\xi_g|,
 \qquad x\in E.
 \label{st:pole-distance}
\end{equation}
The difference is
\[
d_g(x):=q_g(x)-e_g(x)
 =\frac{\sqrt{w_g}(b_g-a_g)}{(x-b_g)(x-a_g)}.
\]
Consequently,
\begin{align}
 |d_g(x)|
 &\leq C\eta_g\frac{\ell_g^{2-s/2}}
                         {(\ell_g+|x-\xi_g|)^2},
 \label{st:d-envelope}\\
 |e_g(x)|
 &\leq C\frac{\ell_g^{1-s/2}}{\ell_g+|x-\xi_g|}.
 \label{st:e-envelope}
\end{align}
On the corresponding ball, we also have the lower bound
\begin{equation}
 |d_g(x)|\geq c\eta_g\ell_g^{-s/2}
 \qquad\bigl(x\in B_E(\xi_g,\ell_g)\bigr).
 \label{st:d-lower}
\end{equation}
Indeed, on this ball,
\[
 \sqrt{w_g}\,|b_g-a_g|\geq c\eta_g\ell_g^{2-s/2},
 \qquad
 |x-b_g|\,|x-a_g|\leq C\ell_g^2.
\]
For real displacements, each difference has a fixed sign. For
complex displacements, we use a sector estimate. Put
\[
 \delta_g=b_g-a_g,\qquad \rho_0=\frac{c_0}{c_*}<\frac12,
 \qquad \zeta_g(x)=\frac{\delta_g}{x-a_g}.
\]
The lower bound for the distance from $a_g$ to $E$ gives
$|\zeta_g(x)|\leq\rho_0$ for every $x\in E$.
If $\delta_g\ne0$, define
$\tau_g=\overline{\delta_g}/|\delta_g|$, which has modulus one. Then
\[
 \tau_g d_g(x)=
 \frac{|\delta_g|\sqrt{w_g}}{(x-a_g)^2}
                 \frac1{1-\zeta_g(x)}.
\]
For $|z|\leq\rho_0$, one has
\[
 \Re\frac1{1-z}=\frac{1-\Re z}{|1-z|^2}
                 \geq\frac1{1+\rho_0}.
\]
Indeed, after multiplication by the positive denominator, the
required difference is
\[
 (1+\rho_0)(1-\Re z)-|1-z|^2
 =(1-\rho_0)(\rho_0+\Re z)+\rho_0^2-|z|^2\geq0.
\]
Since $|1-z|^{-1}\leq(1-\rho_0)^{-1}$, it follows that
\begin{equation}
 \Re(\tau_gd_g(x))
 \geq\frac{1-\rho_0}{1+\rho_0}|d_g(x)|,
 \qquad x\in E.
 \label{st:sector-lower}
\end{equation}
Thus each rotated difference takes values in a fixed sector of
the right half-plane, preventing cancellation in its integral.
If $\delta_g=0$, put
$\tau_g=1$; then $d_g=0$, and all these inequalities remain valid.
Let
\[
 M_\mu f(z)=\sup_{r>0}\frac{1}{\mu(B_E(z,r))}
                         \int_{B_E(z,r)}|f|\,d\mu.
\]
Fix $z\in B_E(\xi_g,\ell_g)$. The annuli of radii
$2^k\ell_g$, together with \eqref{st:d-envelope}, give
\begin{align*}
 |\langle f,d_g\rangle|
 &\leq C\eta_g\ell_g^{2-s/2}
       \sum_{k\geq0}(2^k\ell_g)^{-2}
       \int_{B_E(\xi_g,2^{k+1}\ell_g)}|f|\,d\mu\\
 &\leq C\eta_g\ell_g^{2-s/2}
       \sum_{k\geq0}(2^k\ell_g)^{s-2}M_\mu f(z)\\
 &=C\eta_g\ell_g^{s/2}
       \sum_{k\geq0}2^{k(s-2)}M_\mu f(z)\\
 &\leq C\lambda_g^{1/2}M_\mu f(z).
\end{align*}
For the second inequality, we used
\[
 B_E(\xi_g,2^{k+1}\ell_g)\subset B_E(z,2^{k+2}\ell_g),
 \qquad
 \int_{B_E(z,2^{k+2}\ell_g)}|f|\,d\mu
 \leq C(2^k\ell_g)^sM_\mu f(z).
\]
For radii above one, the upper mass estimate follows from
$\mu(E)=1$. The remaining factors satisfy
\[
 \ell_g^{2-s/2}(2^k\ell_g)^{s-2}
   =\ell_g^{s/2}2^{-k(2-s)},
 \qquad
 \sum_{k\geq0}2^{-k(2-s)}=\frac1{1-2^{s-2}}.
\]
Since $\lambda_g^{1/2}=\eta_g\ell_g^{s/2}$, taking the infimum in
$z$ yields
\begin{equation}
 |\langle f,d_g\rangle|
 \leq C\lambda_g^{1/2}
          \inf_{B_E(\xi_g,\ell_g)}M_\mu f.
 \label{st:maximal-pairing}
\end{equation}
To sum \eqref{st:maximal-pairing}, we prove the Carleson embedding
\begin{equation}
 \sum_{g\in\mathcal G}\lambda_g
       \left(\inf_{B_E(\xi_g,\ell_g)}F\right)^2
 \leq C\mathcal C(\lambda)\int_E F^2\,d\mu
 \label{st:carleson-embedding}
\end{equation}
for every nonnegative lower semicontinuous function $F$ on $E$.
All sums of nonnegative terms are limits of finite subsums.
Fix $t>0$ and write
\[
 \Omega_t=\{z\in E:F(z)>t\},
 \qquad
 \mathcal G_t=
 \{g:\inf_{B_E(\xi_g,\ell_g)}F>t\}.
\]
By lower semicontinuity, $\Omega_t$ is relatively open. If
$\Omega_t=\varnothing$, then $\mathcal G_t=\varnothing$. If
$\Omega_t=E$, the Carleson box of radius two gives
\[
 \sum_g\lambda_g\leq 2^s\mathcal C(\lambda)
                  =2^s\mathcal C(\lambda)\mu(\Omega_t).
\]
Suppose, therefore, that both $\Omega_t$ and its complement in $E$
are nonempty.  Apply the $5r$ covering lemma to the balls
\[
 B_E(z,r(z)),\qquad
 r(z)=\frac1{10}\operatorname{dist}(z,E\setminus\Omega_t),
 \qquad z\in\Omega_t.
\]
This gives a countable pairwise disjoint subfamily
$B_i=B_E(z_i,r_i)$ such that
\[
 B_i\subset\Omega_t,
 \qquad \Omega_t\subset\bigcup_i B_E(z_i,5r_i),
 \qquad 0<r_i\leq\frac1{10}.
\]
The bound $r_i\leq1/10$ follows from $\operatorname{diam}E=1$.
For $g\in\mathcal G_t$,
\[
 \xi_g\in\Omega_t,
 \qquad B_E(\xi_g,\ell_g)\subset\Omega_t.
\]
Assign $g$ to one index $i$ with $\xi_g\in5B_i$. Every point of
$E\setminus\Omega_t$ has distance at least $\ell_g$ from $\xi_g$.
The triangle inequality therefore gives
\[
\ell_g
 \leq\operatorname{dist}(\xi_g,E\setminus\Omega_t)
 \leq |\xi_g-z_i|+
        \operatorname{dist}(z_i,E\setminus\Omega_t)
 \leq 15r_i.
\]
Thus every gap assigned to $i$ is counted in the Carleson box with
center $z_i$ and radius $15r_i\leq3/2$.  The definition of
$\mathcal C(\lambda)$, followed by lower regularity and disjointness,
yields
\begin{align*}
 \sum_{\inf_{B_E(\xi_g,\ell_g)}F>t}\lambda_g
 &\leq \sum_i
       \sum_{\substack{\xi_g\in B_E(z_i,15r_i)\\
                       \ell_g\leq15r_i}}\lambda_g\\
 &\leq 15^s\mathcal C(\lambda)\sum_i r_i^s\\
 &\leq C\mathcal C(\lambda)\sum_i\mu(B_i)\\
 &\leq C\mathcal C(\lambda)\mu(\Omega_t).
\end{align*}
The estimate also holds in the empty and full cases.
Tonelli's theorem and integration of the distribution function now give
\begin{align*}
 \sum_g\lambda_g
       \left(\inf_{B_E(\xi_g,\ell_g)}F\right)^2
 &=\int_0^\infty 2t
       \sum_{\inf_{B_E(\xi_g,\ell_g)}F>t}\lambda_g\,dt\\
 &\leq C\mathcal C(\lambda)
       \int_0^\infty 2t\mu(\Omega_t)\,dt\\
 &=C\mathcal C(\lambda)\int_EF^2\,d\mu.
\end{align*}
This proves \eqref{st:carleson-embedding}.
Figure~\ref{fig:carleson} illustrates the assignment and the packing estimate.

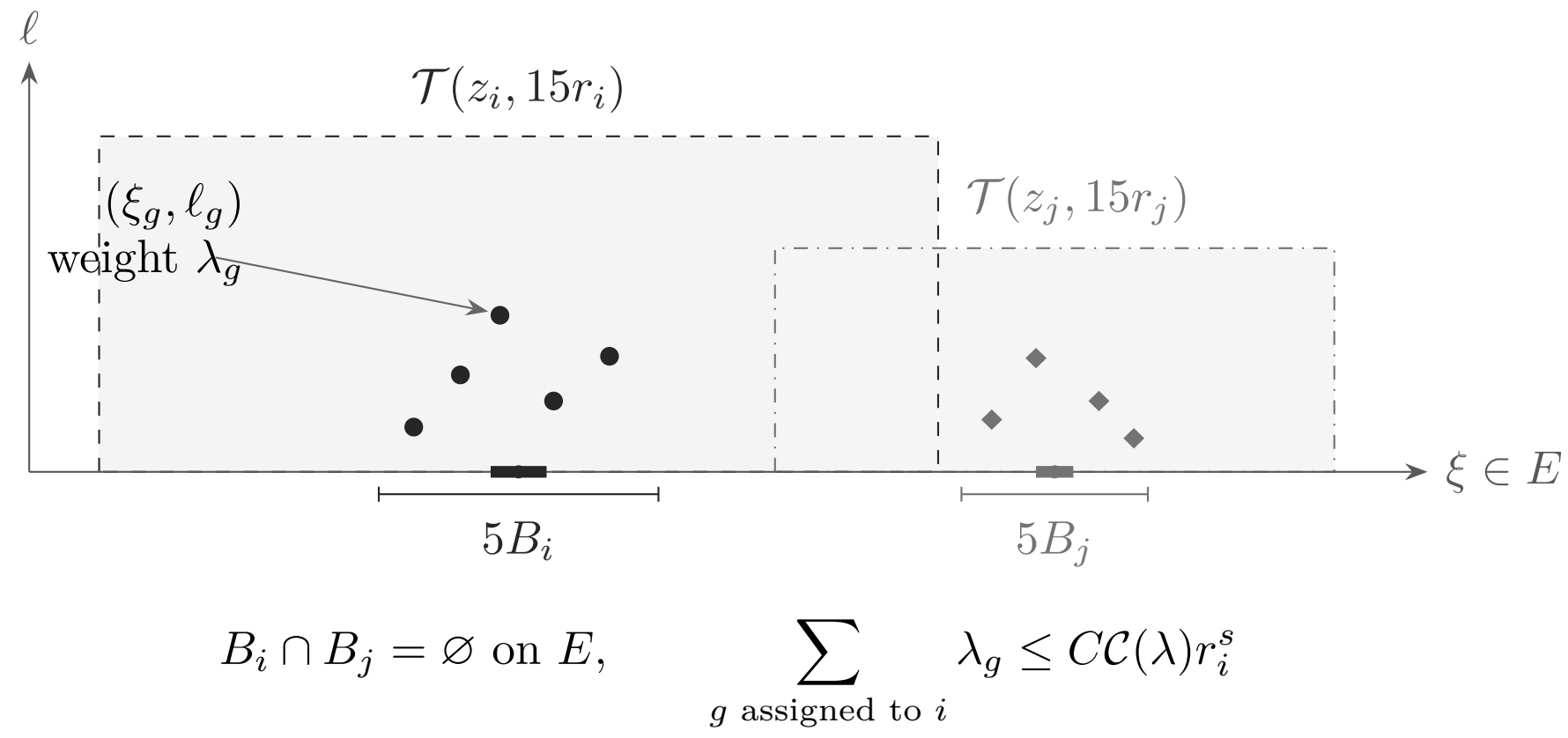
\begin{figure}[htbp]
\centering
\begin{tikzpicture}[x=1.75cm,y=1.4cm,>=Stealth,font=\small]
 \fill[DiagramPrimary!5] (1,0) rectangle (4.6,1.8);
 \fill[DiagramSecondary!7] (3.9,0) rectangle (6.3,1.2);
 \draw[DiagramPrimary,dashed] (1,0) rectangle (4.6,1.8);
 \draw[DiagramSecondary,dash dot] (3.9,0) rectangle (6.3,1.2);
 \draw[->,black!65] (.7,0)--(6.7,0) node[right] {$\xi\in E$};
 \draw[->,black!65] (.7,0)--(.7,2.2) node[above] {$\ell$};
 \node[DiagramPrimary,anchor=south] at (2.8,1.83) {$\mathcal T(z_i,15r_i)$};
 \node[DiagramSecondary,anchor=south] at (5.2,1.23) {$\mathcal T(z_j,15r_j)$};
 \draw[DiagramPrimary,line width=2.5pt] (2.68,0)--(2.92,0);
 \draw[DiagramSecondary,line width=2.5pt] (5.02,0)--(5.18,0);
 \fill[DiagramPrimary] (2.8,0) circle (1.4pt);
 \fill[DiagramSecondary] (5.1,0) circle (1.4pt);
 \draw[DiagramPrimary,thin] (2.2,-.12)--(3.4,-.12);
 \draw[DiagramPrimary,thin] (2.2,-.08)--(2.2,-.16);
 \draw[DiagramPrimary,thin] (3.4,-.08)--(3.4,-.16);
 \node[DiagramPrimary,anchor=north] at (2.8,-.17) {$5B_i$};
 \draw[DiagramSecondary,thin] (4.7,-.12)--(5.5,-.12);
 \draw[DiagramSecondary,thin] (4.7,-.08)--(4.7,-.16);
 \draw[DiagramSecondary,thin] (5.5,-.08)--(5.5,-.16);
 \node[DiagramSecondary,anchor=north] at (5.1,-.17) {$5B_j$};
 \foreach \x/\y in {2.35/.24,2.55/.52,2.72/.84,2.95/.38,3.19/.62}
   \fill[DiagramPrimary] (\x,\y) circle (2pt);
 \foreach \x/\y in {4.83/.28,5.02/.61,5.29/.38,5.44/.18}
   \fill[DiagramSecondary] (\x,\y) ++(-.045,0)--++(.045,.055)
       --++(.045,-.055)--++(-.045,-.055)--cycle;
 \draw[->,black!60] (1.5,1.15)--(2.67,.86);
 \node[anchor=east,align=right] at (1.7,1.28)
   {$(\xi_g,\ell_g)$\\weight $\lambda_g$};
 \node[anchor=north,align=center] at (3.7,-.67)
   {$B_i\cap B_j=\varnothing\text{ on }E,\qquad
     \displaystyle\sum_{g\text{ assigned to }i}\lambda_g
       \le C\mathcal C(\lambda)r_i^s$};
\end{tikzpicture}
\caption{The decomposition of a level set used in the stability proof.
Here $\mathcal T(z,R)=\{(\xi,\ell):\xi\in B_E(z,R),\ 0<\ell\le R\}$.
The disjoint Whitney balls $B_i$ lie in the relatively open set
$\Omega_t\subset E$, and their fivefold enlargements cover it.
A gap with $B_E(\xi_g,\ell_g)\subset\Omega_t$ is assigned to one index
with $\xi_g\in5B_i$; then $\ell_g\le15r_i$.
The enlarged boxes may overlap. Only the base balls are disjoint;
this gives a bound for the sum of $r_i^s$ in terms of
$\mu(\Omega_t)$. The support and the assignment are schematic.}
\label{fig:carleson}
\end{figure}

We verify that $F=M_\mu f$ is lower semicontinuous.
Upper regularity implies $\mu(\{x\})=0$. Consequently,
\[
 \mu(\{z-r,z+r\})=0,
 \qquad
 \int_{\{z-r,z+r\}}|f|\,d\mu=0.
\]
For fixed $r>0$ and $z_j\to z$ in $E$, dominated convergence gives
\[
 \int_{B_E(z_j,r)}|f|\,d\mu
   \longrightarrow\int_{B_E(z,r)}|f|\,d\mu,
 \qquad
 \mu(B_E(z_j,r))\longrightarrow\mu(B_E(z,r))>0.
\]
Each average is continuous in its center. Its supremum over
$r>0$ is lower semicontinuous. Applying
\eqref{eq:prelim-maximal-lp} and \eqref{st:carleson-embedding} gives
\[
\sum_g|\langle f,d_g\rangle|^2
 \leq C\mathcal C(\lambda)\|M_\mu f\|_2^2
 \leq C\mathcal C(\lambda)\|f\|_2^2.
\]
For a finitely supported sequence $c$, Cauchy--Schwarz gives
\[
 \begin{aligned}
 \left\|\sum_g c_gd_g\right\|_2
 &=\sup_{\|f\|_2=1}
       \left|\sum_g c_g\langle d_g,f\rangle\right|\\
 &\leq\left(\sum_g|c_g|^2\right)^{1/2}
       \sup_{\|f\|_2=1}
       \left(\sum_g|\langle f,d_g\rangle|^2\right)^{1/2}\\
 &\leq C\mathcal C(\lambda)^{1/2}\|c\|_{\ell^2}.
 \end{aligned}
\]
By completion, $D$ extends to $\ell^2(\mathcal G)$, with
\begin{equation}
 \|D\|_{\ell^2\to L^2}\leq C\mathcal C(\lambda)^{1/2}.
 \label{st:D-upper}
\end{equation}
For necessity, suppose that $D$ is bounded and fix $x\in E$ and
$0<R\leq2$.  Set $f=\mathbf1_{B_E(x,4R)}$.  If
$\xi_g\in B_E(x,R)$ and $\ell_g\leq R$, then
$B_E(\xi_g,\ell_g)\subset B_E(x,4R)$.  The sector estimate
\eqref{st:sector-lower}, the pointwise lower bound
\eqref{st:d-lower}, and lower regularity give
\begin{align*}
 |\langle f,d_g\rangle|
 &=\left|\int_{B_E(x,4R)}d_g\,d\mu\right|\\
 &\geq\Re\int_{B_E(x,4R)}\tau_gd_g\,d\mu\\
 &\geq c\int_{B_E(x,4R)}|d_g|\,d\mu\\
 &\geq c\eta_g\ell_g^{-s/2}
             \mu\bigl(B_E(\xi_g,\ell_g)\bigr)\\
 &\geq c\eta_g\ell_g^{s/2}
  =c\lambda_g^{1/2}.
\end{align*}
The analysis operator $D^*$ has the same norm as $D$.  Therefore
\begin{align*}
 \sum_{\substack{\xi_g\in B_E(x,R)\\\ell_g\leq R}}\lambda_g
 &\leq C\sum_g|\langle f,d_g\rangle|^2\\
 &\leq C\|D\|^2\mu\bigl(B_E(x,4R)\bigr)\\
 &\leq C\|D\|^2R^s.
\end{align*}
In the last line, radii above one are again handled by $\mu(E)=1$.
Taking the supremum proves
\[
c\mathcal C(\lambda)^{1/2}
 \leq\|D\|_{\ell^2\to L^2}
 \leq C\mathcal C(\lambda)^{1/2}.
\]
Let $U(c)=\sum_gc_ge_g$ and $V(c)=\sum_gc_gq_g$, initially on
finite sequences. Proposition~\ref{thm:basis} gives
\[
 \|U\|_{\ell^2\to L^2}=1,
 \qquad V=U+D,
 \qquad D=V-U.
\]
Thus
\[
 \|V\|\leq1+\|D\|,
 \qquad \|D\|\leq1+\|V\|.
\]
The boundedness of $V$ is equivalent, by Hilbert space duality,
to the Bessel inequality for $(q_g)$. This proves the asserted
necessary and sufficient condition.

When the total packing norm is small, this identification also gives
the Riesz basis conclusion.
Define $T=DU^*$ on $L^2(\mu)$.
The orthonormal basis in Proposition~\ref{thm:basis} gives
\[
 T1=0,\qquad Te_g=d_g,\qquad
 \|T\|_{2\to2}\leq C\mathcal C(\lambda)^{1/2}.
\]
If $\varepsilon=C\mathcal C(\lambda)^{1/2}<1$, then
\[
 \sum_{k=0}^{\infty}\|(-T)^k\|_{2\to2}
 \leq\sum_{k=0}^{\infty}\varepsilon^k
 =\frac1{1-\varepsilon}.
\]
The series $S=\sum_{k\geq0}(-T)^k$ converges in operator norm and satisfies
\[
 (I+T)S=S(I+T)=I,
 \qquad (I+T)1=1,
 \qquad (I+T)e_g=q_g.
\]
Hence the perturbed system is complete. For a finite combination
$v=c_0'1+\sum_gc_ge_g$, the triangle inequality gives
\[
 (1-\varepsilon)\|v\|_2
 \leq\|(I+T)v\|_2
 \leq(1+\varepsilon)\|v\|_2,
 \qquad
 \|v\|_2^2=|c_0'|^2+\sum_g|c_g|^2.
\]
Squaring gives the Riesz bounds $(1-\varepsilon)^2$ and
$(1+\varepsilon)^2$, and continuity extends them
to all square summable coefficient sequences.

We now estimate the kernel of the perturbation operator. The sums
over scales below also make the dependence on $s<1$ explicit.
For a finite family $\mathcal F\subset\mathcal G$, let
\[
 T_{\mathcal F}f(x)
 =\sum_{g\in\mathcal F}\langle f,e_g\rangle d_g(x)
 =\int_EK_{\mathcal F}(x,y)f(y)\,d\mu(y),
\]
where
\[
 K_{\mathcal F}(x,y)
 =\sum_{g\in\mathcal F}d_g(x)\overline{e_g(y)}
 =\sum_{g\in\mathcal F}d_g(x)e_g(y).
\]
The last equality holds because $e_g$ is real, even when $d_g$
is complex. The adjoint kernel is
\[
 \overline{K_{\mathcal F}(y,x)}
 =\sum_{g\in\mathcal F}e_g(x)\overline{d_g(y)}.
\]
Its size and continuity estimates follow from those of $K_{\mathcal F}$
by conjugation and interchange of variables.
Set $a=\sup_g\eta_g$.  The Carleson box centered at $\xi_g$ with
radius $\ell_g$ contains the gap $g$ itself, so
\[
a\leq\mathcal C(\lambda)^{1/2}.
\]
We prove the uniform size estimate
\begin{equation}
 |K_{\mathcal F}(x,y)|
 \leq Ca|x-y|^{-s},
 \label{st:kernel-size}
\end{equation}
and, when $|x-x'|\leq\tfrac12|x-y|$, the continuity estimate
\begin{equation}
 \begin{aligned}
 &|K_{\mathcal F}(x,y)-K_{\mathcal F}(x',y)|\\
 &\qquad+|K_{\mathcal F}(y,x)-K_{\mathcal F}(y,x')|
 \leq Ca\frac{|x-x'|^{1-s}}{|x-y|}.
 \end{aligned}
 \label{st:kernel-smoothness}
\end{equation}
The right side of the second estimate equals
$Ca(|x-x'|/|x-y|)^{1-s}|x-y|^{-s}$, the standard
Calder\'on--Zygmund continuity bound of order $1-s$ for an
$s$-dimensional growth condition.

We first count gaps at a single scale. If $0<t\leq R\leq2$
and $z\in E$, then
\begin{equation}
 \#\{g:t\leq\ell_g<2t,\ \xi_g\in B_E(z,R)\}
 \leq C(R/t)^s.
 \label{st:band-count}
\end{equation}
Let $I$ be the family counted in \eqref{st:band-count}. Replace
$\xi_g$ by the left endpoint $u_g$. For distinct gaps in this band,
\[
 |u_g-u_{g'}|\geq t,
 \qquad |u_g-\xi_g|<2t.
\]
The balls $B_E(u_g,t/4)$ are pairwise disjoint, and
\[
 B_E(u_g,t/4)\subset B_E(z,R+3t).
\]
Regularity therefore gives
\[
 c\,\#I\,t^s
 \leq\sum_{g\in I}\mu(B_E(u_g,t/4))
 \leq\mu(B_E(z,R+3t))
 \leq C(R+3t)^s\leq C'R^s.
\]
Dividing by $ct^s$ proves the count. For $R>2$, use the global
bound $\#\{g:t\leq\ell_g<2t\}\leq Ct^{-s}$, which follows
from the same calculation with disjoint balls.
Two annular consequences will be used repeatedly.  For every
$\beta>s$,
\begin{equation}
 \sum_{t\leq\ell_g<2t}(t+|z-\xi_g|)^{-\beta}
 \leq C_\beta t^{-\beta}.
 \label{st:band-potential}
\end{equation}
To verify this estimate, separate the centers with
$|z-\xi_g|<t$ from the successive annuli. The count just proved gives
\[
 \begin{aligned}
 \sum_{t\leq\ell_g<2t}(t+|z-\xi_g|)^{-\beta}
 &\leq Ct^{-\beta}
   +C\sum_{k\geq0}2^{(k+1)s}(2^kt)^{-\beta}\\
 &\leq Ct^{-\beta}
       \left(1+2^s\sum_{k\geq0}2^{-k(\beta-s)}\right)\\
 &=Ct^{-\beta}\left(1+\frac{2^s}{1-2^{s-\beta}}\right).
 \end{aligned}
\]
For $R\geq t$, the corresponding tail estimate is
\[
\sum_{\substack{t\leq\ell_g<2t\\|z-\xi_g|\geq R}}
       (t+|z-\xi_g|)^{-\beta}
 \leq C_\beta t^{-s}R^{s-\beta}.
\]
Indeed,
\[
 \begin{aligned}
 \sum_{\substack{t\leq\ell_g<2t\\|z-\xi_g|\geq R}}
       (t+|z-\xi_g|)^{-\beta}
 &\leq C\sum_{k\geq0}
       \left(\frac{2^{k+1}R}{t}\right)^s(2^kR)^{-\beta}\\
 &\leq\frac{C2^s}{1-2^{s-\beta}}\,t^{-s}R^{s-\beta}.
 \end{aligned}
\]
We use these estimates with $\beta=1,2,3,4$. The case
$\beta=1$ requires $s<1$.

Fix distinct $x,y\in E$ and write $r=|x-y|$.  Split the gaps into
dyadic bands $t\leq\ell_g<2t$, with $t=2^{-j-1}$ and $j\geq0$;
all bounded gaps have length less than one. For one band, write
\[
 L_x(g)=t+|x-\xi_g|,\qquad L_y(g)=t+|y-\xi_g|.
\]
By \eqref{st:d-envelope} and \eqref{st:e-envelope}, the sum of
absolute values of its kernel terms is at most
\[
Ca t^{3-s}
 \sum_{t\leq\ell_g<2t}L_x(g)^{-2}L_y(g)^{-1}.
\]
If $t\geq r$, then $L_x(g)$ and $L_y(g)$ are comparable, and
\[
 L_y(g)\leq L_x(g)+r\leq2L_x(g),
 \qquad L_x(g)\leq2L_y(g).
\]
Using \eqref{st:band-potential} with $\beta=3$ gives
\[
 Ca t^{3-s}\sum_gL_x(g)^{-2}L_y(g)^{-1}
 \leq Ca t^{3-s}\sum_gL_x(g)^{-3}
 \leq Ca t^{3-s}t^{-3}=Ca t^{-s}.
\]
If $t<r$, partition the centers into the three disjoint sets
\[
 \begin{aligned}
 \mathcal N_x&=\{g:|x-\xi_g|<r/2\},\\
 \mathcal N_y&=\{g:|y-\xi_g|<r/2\},\\
 \mathcal R&=\{g:|x-\xi_g|\geq r/2,\ |y-\xi_g|\geq r/2\},
 \end{aligned}
\]
always restricted to the band under consideration.  On
$\mathcal N_x$, $L_y(g)\geq r/2$, so the contribution is at most
\[
 Ca t^{3-s}r^{-1}\sum_gL_x(g)^{-2}
 \leq Ca\frac{t^{1-s}}r.
\]
On $\mathcal N_y$, $L_x(g)\geq r/2$, and the corresponding bound is
\[
 Ca t^{3-s}r^{-2}\sum_gL_y(g)^{-1}
 \leq Ca\frac{t^{2-s}}{r^2}.
\]
On $\mathcal R$ the two distances are comparable within a fixed
factor.  The tail estimate with $\beta=3$, at radius comparable to
$r$ (or the full estimate if $r/2<t$), gives
\[
 Ca t^{3-s}\sum_{\mathcal R}L_x(g)^{-3}
 \leq Ca t^{3-2s}r^{s-3}.
\]
The last two terms are bounded by the first, as follows from
\[
 \frac{t^{2-s}/r^2}{t^{1-s}/r}=\frac tr\leq1,
 \qquad
 \frac{t^{3-2s}r^{s-3}}{t^{1-s}/r}
       =\left(\frac tr\right)^{2-s}\leq1.
\]
Consequently, the bound for the band is
\begin{equation}
 \sum_{t\leq\ell_g<2t}|d_g(x)e_g(y)|
 \leq Ca
 \begin{cases}
 t^{-s},&t\geq r,\\
 t^{1-s}/r,&t<r.
 \end{cases}
 \label{st:band-size-bound}
\end{equation}
The two geometric sums satisfy
\[
 \begin{aligned}
 \sum_{t\geq r}t^{-s}
 &\leq Cr^{-s}\sum_{k\geq0}2^{-ks}
 =\frac{Cr^{-s}}{1-2^{-s}},\\
 \sum_{t<r}t^{1-s}
 &\leq Cr^{1-s}\sum_{k\geq0}2^{-k(1-s)}
 =\frac{Cr^{1-s}}{1-2^{s-1}}.
 \end{aligned}
\]
Here and below $t$ runs over the dyadic band parameters. Hence
\[
 Ca\sum_{t\geq r}t^{-s}
 +\frac{Ca}{r}\sum_{t<r}t^{1-s}
 \leq Ca r^{-s}.
\]
This proves \eqref{st:kernel-size}, uniformly for finite subfamilies,
and also proves absolute convergence of the full kernel away from
the diagonal. Convergence at small scales uses $1-s>0$.

To prove the continuity estimate, first suppose that
$0<h=|x-x'|\leq r/4$.  For the bands with $t\leq h$, apply
\eqref{st:band-size-bound} at both pairs $(x,y)$ and $(x',y)$.
Their separations are comparable to $r$, and the same statement
holds for the transposed pairs.  These bands contribute at most
\begin{equation}
 \frac{Ca}{r}\sum_{t\leq h}t^{1-s}
 \leq Ca\frac{h^{1-s}}r
 \label{st:small-band-difference}
\end{equation}
to the sum of the two differences.
For $t>h$,
\[
 t+|x'-\xi_g|\leq t+|x-\xi_g|+h
               \leq2\bigl(t+|x-\xi_g|\bigr),
\]
and the reverse comparison follows by interchanging $x$ and $x'$.
The rational identities give
\begin{align*}
 e_g(x)-e_g(x')
 &=\frac{\sqrt{w_g}(x'-x)}{(x-a_g)(x'-a_g)},\\
 d_g(x)-d_g(x')
 &=\frac{\sqrt{w_g}(b_g-a_g)(x'-x)
                  (x+x'-a_g-b_g)}
          {(x-a_g)(x-b_g)(x'-a_g)(x'-b_g)}.
\end{align*}
By \eqref{st:pole-distance}, every denominator factor at $x$ or
$x'$ is comparable in absolute value to $L_x(g)$. The numerators
satisfy
\[
 \begin{aligned}
 \sqrt{w_g}&\leq Ct^{1-s/2},\\
 |b_g-a_g|&\leq2at,\\
 |x+x'-a_g-b_g|
 &\leq |x-a_g|+|x'-b_g|\leq CL_x(g).
 \end{aligned}
\]
It follows that
\begin{align}
 |e_g(x)-e_g(x')|
 &\leq C h t^{1-s/2}L_x(g)^{-2},
 \label{st:e-difference}\\
 |d_g(x)-d_g(x')|
 &\leq Ca h t^{2-s/2}L_x(g)^{-3}. \notag
\end{align}
These identities use values of the rational functions only at points of
$E$ and require no differentiation across a gap.
For the kernel difference in the first variable, the resulting band sum
is bounded by
\[
 Ca h t^{3-s}\sum_g L_x(g)^{-3}L_y(g)^{-1}.
\]
If $t\geq r$, then
\[
 Ca h t^{3-s}\sum_gL_x(g)^{-3}L_y(g)^{-1}
 \leq Ca h t^{3-s}\sum_gL_x(g)^{-4}
 \leq Ca h t^{-1-s}.
\]
For $h<t<r$, the portion near $x$ gives
\[
 Ca h t^{3-s}
 \sum_{\mathcal N_x}L_x(g)^{-3}L_y(g)^{-1}
 \leq Ca h t^{3-s}r^{-1}t^{-3}
 =Ca h\frac{t^{-s}}r.
\]
The part near $y$ gives
\[
 Ca h t^{3-s}
 \sum_{\mathcal N_y}L_x(g)^{-3}L_y(g)^{-1}
 \leq Ca h t^{3-s}r^{-3}t^{-1}
 =Ca h\frac{t^{2-s}}{r^3}.
\]
On the remaining set, the distance comparison and the tail
estimate with $\beta=4$ give
\[
 Ca h t^{3-s}\sum_{\mathcal R}L_x(g)^{-3}L_y(g)^{-1}
 \leq Ca h t^{3-s}t^{-s}r^{s-4}
 =Ca h t^{3-2s}r^{s-4}.
\]
The ratios of the last two bounds to the first are
\[
 \frac{t^{2-s}/r^3}{t^{-s}/r}=\left(\frac tr\right)^2,
 \qquad
 \frac{t^{3-2s}r^{s-4}}{t^{-s}/r}
       =\left(\frac tr\right)^{3-s}.
\]
Both are at most one. The middle and large bands satisfy
\begin{align*}
 \sum_{h<t<r}Ca h\frac{t^{-s}}r
 +\sum_{t\geq r}Ca h t^{-1-s}
 &\leq\frac{Ca h}{r}
       h^{-s}\sum_{k\geq0}2^{-ks}
       +Ca h r^{-1-s}\sum_{k\geq0}2^{-k(1+s)}\\
 &\leq Ca\frac{h^{1-s}}r+Ca h r^{-1-s}\\
 &\leq Ca\frac{h^{1-s}}r.
\end{align*}
The last comparison uses
\[
 \frac{h r^{-1-s}}{h^{1-s}/r}
   =\left(\frac hr\right)^s\leq1.
\]
Together with \eqref{st:small-band-difference}, this is the required
bound for $K_{\mathcal F}(x,y)-K_{\mathcal F}(x',y)$.

For the other difference, use \eqref{st:e-difference} in the second
factor.  Its band sum is bounded by
\[
 Ca h t^{3-s}\sum_gL_y(g)^{-2}L_x(g)^{-2}.
\]
For $t\geq r$, the same potential estimate with power four gives
$Ca h t^{-1-s}$. For $h<t<r$, the contributions near the two points give
\[
 \begin{aligned}
 Ca h t^{3-s}\sum_{\mathcal N_x}L_y(g)^{-2}L_x(g)^{-2}
 &\leq Ca h t^{3-s}r^{-2}t^{-2}
   =Ca h\frac{t^{1-s}}{r^2},\\
 Ca h t^{3-s}\sum_{\mathcal N_y}L_y(g)^{-2}L_x(g)^{-2}
 &\leq Ca h t^{3-s}t^{-2}r^{-2}
   =Ca h\frac{t^{1-s}}{r^2}.
 \end{aligned}
\]
The remaining set contributes at most $Ca h t^{3-2s}r^{s-4}$,
and
\[
 \frac{t^{3-2s}r^{s-4}}{t^{1-s}/r^2}
       =\left(\frac tr\right)^{2-s}\leq1.
\]
Consequently,
\[
 \sum_{h<t<r}Ca h\frac{t^{1-s}}{r^2}
 +\sum_{t\geq r}Ca h t^{-1-s}
 \leq Ca h r^{-1-s}
 \leq Ca\frac{h^{1-s}}r.
\]
Adding the small bands proves the estimate in the second variable for
$h\leq r/4$. If $r/4<h\leq r/2$, then
\[
 |x'-y|\geq r-h\geq r/2,
 \qquad
 r^{-s}\leq4^{1-s}\frac{h^{1-s}}r.
\]
The size estimate bounds each difference by $Ca r^{-s}$ and hence
gives the required bound in this range. At $h=0$ both differences
vanish. This proves \eqref{st:kernel-smoothness}.

Suppose that $\mathcal C(\lambda)<\infty$ and put
$\kappa=\mathcal C(\lambda)^{1/2}$. If $\kappa=0$, then
$\eta_g=0$ for every gap, and $T=0$. Assume from now on that
$\kappa>0$.
For finite $\mathcal F$, \eqref{st:D-upper} gives
\[
 \|\kappa^{-1}T_{\mathcal F}\|_{2\to2}
 \leq\kappa^{-1}\|D\|_{\ell^2\to L^2}\|U^*\|_{L^2\to\ell^2}
 \leq C.
\]
The kernel bounds give
\[
 |\kappa^{-1}K_{\mathcal F}(x,y)|\leq C|x-y|^{-s},
 \qquad
 |\kappa^{-1}\Delta_xK_{\mathcal F}(x,y)|
 \leq C\frac{|x-x'|^{1-s}}{|x-y|},
\]
because $a/\kappa\leq1$. The same estimate holds in the second
variable. Here $\Delta_xK(x,y)=K(x,y)-K(x',y)$.
Apply \eqref{eq:prelim-cz-conclusion} with dimension parameter $s$
and smoothness exponent $1-s$. The normalized $L^2$ and kernel
constants are uniform in $\mathcal F$. Restoring the factor
$\kappa$ proves
\[
\|T_{\mathcal F}\|_{L^p\to L^p}
 \leq C_p\mathcal C(\lambda)^{1/2},
 \qquad 1<p<\infty,
\]
uniformly in $\mathcal F$.
Choose finite families $\mathcal F_j\uparrow\mathcal G$, and let
$P_{\mathcal F_j}$ denote the corresponding coordinate projections
on $\ell^2(\mathcal G)$. For $f\in L^2(\mu)$,
\[
 \begin{aligned}
 \|T_{\mathcal F_j}f-Tf\|_2
 &=\|D(P_{\mathcal F_j}-I)U^*f\|_2\\
 &\leq C\kappa
   \left(\sum_{g\notin\mathcal F_j}
                  |\langle f,e_g\rangle|^2\right)^{1/2}
 \longrightarrow0.
 \end{aligned}
\]
For $f\in L^2(\mu)\cap L^p(\mu)$, choose an almost everywhere
convergent subsequence. Fatou's lemma yields
\[
 \|Tf\|_p^p
 \leq\liminf_j\|T_{\mathcal F_j}f\|_p^p
 \leq C_p^p\kappa^p\|f\|_p^p.
\]
Bounded functions are dense in $L^p(\mu)$, so this defines the
unique bounded extension with
\[
\|T\|_{L^p\to L^p}
 \leq C_p\mathcal C(\lambda)^{1/2}.
\]

The series defining $Tf$ converges unconditionally in $L^p$.
Indeed, finite linear combinations of the original basis are dense
by Proposition~\ref{thm:basis}. If $v$ is such a combination and
$\mathcal F$ contains all its gap indices, then
\[
 T_{\mathcal F}v=Tv.
\]
For every $f\in L^p(\mu)$,
\[
 \begin{aligned}
 \|T_{\mathcal F}f-Tf\|_p
 &\leq\|T_{\mathcal F}(f-v)\|_p+\|T(f-v)\|_p\\
 &\leq2C_p\kappa\|f-v\|_p.
 \end{aligned}
\]
Approximate $f$ by $v$, and then choose a finite family containing
its indices. This proves convergence along the directed family of
all finite subsets, independently of their enumeration.
If $C_p\kappa<1$, the inverse on $L^p$ is
\[
 (I+T)^{-1}=\sum_{k=0}^{\infty}(-T)^k,
 \qquad
 \|(I+T)^{-1}\|_{p\to p}\leq\frac1{1-C_p\kappa}.
\]
For $h=c+\sum_g c_ge_g$ in $L^p$, unconditional convergence and
boundedness give
\[
 (I+T)h=c+\sum_g c_gq_g.
\]
Surjectivity of $I+T$ gives existence of this expansion for every
element of $L^p$; injectivity and the original basis give uniqueness.
The same bounded isomorphism preserves unconditional convergence.
\end{proof}

\subsection{Essential norm and finite pole replacements}
\label{subsec:tail-stability}
We use the tail packing norm from \eqref{eq:packing-tail-definition}.
Finite replacement of kernel bases in spaces with the division property
is classical; see \cite[proof of Theorem~4.3]{Baranov2023}.
The following proof gives an explicit
construction on $L^p(\mu)$ in which the constant coordinate remains fixed.

\begin{proof}[Essential norm and basis assertions of Theorem~\ref{thm:stability}]
Let \(\pi_r\) be the coordinate projection in \(\ell^2(\mathcal G)\)
onto the gaps with \(\ell_g\ge r\), and let \(P_r\) be the canonical
projection onto these gaps and the constant.  Both have finite rank.
Write
\[
 T_r=T(I-P_r).
\]
The estimates already proved, applied to the omitted gaps, give
\[
 c\mathcal C_r(\lambda)^{1/2}
 \le \|D(I-\pi_r)\|_{\ell^2\to L^2}
 =\|T_r\|_{2\to2},
 \qquad
 \|T_r\|_{q\to q}\le C_q\mathcal C_r(\lambda)^{1/2}
 \quad(1<q<\infty).
\]
These imply upper and lower norm estimates on every interior \(L^p\).
For \(p\ne2\), choose \(q\in(1,\infty)\) on the opposite side
of \(2\) from \(p\), and choose \(0<\vartheta<1\) so that
\(1/2=\vartheta/p+(1-\vartheta)/q\).  If
\(\mathcal C_r(\lambda)>0\), interpolation gives
\[
 c\mathcal C_r(\lambda)^{1/2}
 \le \|T_r\|_{p\to p}^{\vartheta}
       \bigl(C_q\mathcal C_r(\lambda)^{1/2}\bigr)^{1-\vartheta},
\]
and hence \(\|T_r\|_{p\to p}\ge c_p\mathcal C_r(\lambda)^{1/2}\).
If \(\mathcal C_r(\lambda)=0\), the upper estimate gives \(T_r=0\),
so the same conclusion holds.  Together with the \(p=2\) case,
\begin{equation}
 c_p\mathcal C_r(\lambda)^{1/2}
 \le \|T_r\|_{p\to p}
 \le C_p\mathcal C_r(\lambda)^{1/2},
 \qquad 1<p<\infty.
 \label{eq:tail-lp-two-sided}
\end{equation}
The identical argument for \(T\), with \(\mathcal C_r\) replaced by
\(\mathcal C\), proves the full upper and lower estimates for the operator norm.

For every compact \(K:\ell^2\to L^2\), one has
\(\|K(I-\pi_r)\|\to0\).  Indeed, the projections converge strongly
and uniformly on the compact closure of the image of the unit ball
under \(K^*\).  Since \(D\pi_r\) has finite rank,
\[
 \|D\|_{\mathrm{ess},\ell^2\to L^2}
 =\lim_{r\downarrow0}\|D(I-\pi_r)\|_{\ell^2\to L^2}
 \asymp\gamma(\lambda)^{1/2}.
\]
If \(U:\ell^2(\mathcal G)\to L^2(\mu)\) is the canonical isometric
synthesis map, the identities \(T=DU^*\) and \(D=TU\) give
\[
 \|D\|_{\mathrm{ess},\ell^2\to L^2}
 =\|T\|_{\mathrm{ess},L^2\to L^2}.
\]

Fix \(1<p<\infty\), and put
\(B_p=\sup_r\|I-P_r\|_{p\to p}<\infty\).
The operators \(TP_r\) have finite rank, so
\eqref{eq:tail-lp-two-sided} gives
\[
 \|T\|_{\mathrm{ess},L^p\to L^p}
 \le C_p\gamma(\lambda)^{1/2}.
\]
For the reverse estimate, let \(K\) be any compact operator on
\(L^p(\mu)\).  The canonical projections are self-adjoint under
the integral pairing and converge strongly on \(L^{p'}\).
Consequently,
\[
 \|K(I-P_r)\|_{p\to p}
 =\|(I-P_r)^*K^*\|_{p'\to p'}\longrightarrow0,
\]
because the projections converge uniformly on the compact closure of
the image of the unit ball under \(K^*\). Thus
\[
 B_p\|T-K\|_{p\to p}
 \ge \liminf_{r\downarrow0}
       \|(T-K)(I-P_r)\|_{p\to p}
 \ge c_p\gamma(\lambda)^{1/2}.
\]
Taking the infimum over \(K\) and absorbing \(B_p\) into the
lower constant proves
\[
 \|T\|_{\mathrm{ess},L^p\to L^p}
 \asymp_p\gamma(\lambda)^{1/2}.
\]
In particular, \(T\) is compact on one, and hence every, interior
\(L^p(\mu)\) exactly when \(\gamma(\lambda)=0\).

We now prove the basis assertion. Fix \(p\) with
\(C_p\gamma(\lambda)^{1/2}<1\), and choose \(r>0\) so that
\[
 \kappa=C_p\mathcal C_r(\lambda)^{1/2}<1.
\]
Let
\[
 \mathcal T=\{g:\ell_g<r\},\qquad
 \mathcal J=\{g:\ell_g\ge r,\ b_g\ne a_g\}.
\]
The set \(\mathcal J\) is finite. Define
\[
 A_{\mathcal T}=I+T_r.
\]
The Neumann series shows that \(A_{\mathcal T}\) is bounded and
invertible on \(L^p\). It maps \(e_g\) to \(q_g\) for
\(g\in\mathcal T\), fixes every other \(e_g\), and fixes \(1\).
The resulting family is an unconditional basis. It remains to
replace the finitely many vectors indexed by \(\mathcal J\).

For a finite \(\mathcal S\subset\mathcal T\), put
\[
 R_{\mathcal S}(z)=
 \prod_{g\in\mathcal S}\frac{z-b_g}{z-a_g},
 \qquad
 A_{\mathcal S}1=1,\qquad
 A_{\mathcal S}e_g=
 \begin{cases}
  q_g,&g\in\mathcal S,\\
  e_g,&g\notin\mathcal S.
 \end{cases}
\]
The same perturbation estimate gives
\[
 \|A_{\mathcal S}^{-1}\|_{p\to p}\le(1-\kappa)^{-1}.
\]
Partial fractions and \(F(a_g)=0\) imply
\[
 \int_E R_{\mathcal S}(x)\,d\mu(x)=1,\qquad
 \int_E\frac{R_{\mathcal S}(x)}{x-z}\,d\mu(x)
       =F(z)R_{\mathcal S}(z).
\]
To check the second identity, write
\(R_{\mathcal S}(z)=1+\sum_{g\in\mathcal S}c_g/(z-a_g)\)
and use
\[
 \int_E\frac{d\mu(x)}{(x-a_g)(x-z)}
 =\frac{F(z)-F(a_g)}{z-a_g}
 =\frac{F(z)}{z-a_g}.
\]
Both sides extend across each removable point \(a_g\).

The bilinear functional
\[
 \Lambda_{\mathcal S}(f)
      =\int_E f(x)R_{\mathcal S}(x)\,d\mu(x)
\]
annihilates all the nonconstant vectors of the
\(A_{\mathcal S}\)-basis and takes the value \(1\) on the
constant. Consequently,
\[
 \Lambda_{\mathcal S}(f)
 =\int_E A_{\mathcal S}^{-1}f\,d\mu,\qquad
 \|R_{\mathcal S}\|_{p'}
 =\|\Lambda_{\mathcal S}\|
 \le(1-\kappa)^{-1}.
\]
Here the pairing is bilinear; no complex conjugation is inserted
into \(R_{\mathcal S}\). Its norm is nevertheless the usual
\(L^{p'}\) norm.

Since the gaps are disjoint,
\[
 \sum_g|b_g-a_g|\le c_0\sum_g\ell_g\le c_0.
\]
Therefore the product
\[
 P(z)=\prod_{g\in\mathcal T}\frac{z-b_g}{z-a_g}
\]
converges locally uniformly away from \(E\) and the
poles \(a_g\). On a compact set at positive distance from \(E\),
the denominators corresponding to small gaps are uniformly bounded
below, so the sum of the absolute perturbations of the factors is finite.
The product is nonzero wherever none of its factors vanishes.

Take an increasing finite exhaustion of \(\mathcal T\).
Reflexivity of \(L^{p'}\) gives a subsequence for which
\(R_{\mathcal S}\) converges weakly to a function \(R\in L^{p'}\).
The preceding identity for the Cauchy transform yields
\[
 G(z):=\int_E\frac{R(x)}{x-z}\,d\mu(x)=F(z)P(z),
 \qquad \int_E R\,d\mu=1.
\]
The identity initially holds away from these poles and extends
across them by analyticity of the left side.
For every \(g\in\mathcal T\), one has \(G(b_g)=0\): a moved
pole is a zero of its own product factor, whereas an unchanged
pole remains the zero \(a_g\) of \(F\). For every
\(g\notin\mathcal T\), the product has a finite nonzero value at
\(a_g\), and hence \(G(a_g)=0\).

For \(j\in\mathcal J\), none of the factors of \(P\) has a
zero or pole at \(a_j\) or \(b_j\). Indeed, the admissible disk
around each \(a_g\) lies in the vertical strip over its own
gap, and these strips are disjoint. Thus
\[
 G(a_j)=0,\qquad
 G'(a_j)=F'(a_j)P(a_j)\ne0.
\]
Also \(G(b_j)\ne0\): if \(b_j\) is nonreal then
\(\operatorname{Im}F(b_j)\ne0\), while if \(b_j\) is real
then its own gap contains the unique zero \(a_j\ne b_j\).
The product factor \(P(b_j)\) is nonzero in both cases.

Define, for \(j\in\mathcal J\), the bounded linear functional
\[
 L_j(f)=
 \frac{1}{\sqrt{w_j}\,G'(a_j)}
 \int_E\frac{f(x)R(x)}{x-a_j}\,d\mu(x).
\]
It annihilates the constant, takes the value \(1\) on \(e_j\),
and vanishes on every other vector of the
\(A_{\mathcal T}\)-basis. These assertions follow from the
zeros of \(G\) and
\[
 \int_E\frac{R(x)}{(x-a_j)(x-b)}\,d\mu(x)
 =\frac{G(b)-G(a_j)}{b-a_j}.
\]
Hence the \(L_j\) are precisely the corresponding coordinate
functionals.

The finite replacement matrix is
\[
 M_{ji}=L_j(q_i)
 =\frac{\sqrt{w_i}\,G(b_i)}
        {\sqrt{w_j}\,G'(a_j)(b_i-a_j)},
 \qquad i,j\in\mathcal J.
\]
All row and column factors are nonzero. The Cauchy matrix
\((1/(b_i-a_j))_{j,i}\) is invertible: the entries in each list
are distinct, and the two lists are disjoint. Explicitly,
\[
 \det\left(\frac1{b_i-a_j}\right)_{j,i=1}^N
 =
 \frac{\displaystyle
       \prod_{j<k}(a_j-a_k)\prod_{i<k}(b_k-b_i)}
      {\displaystyle\prod_{j,i}(b_i-a_j)}\ne0.
\]
Thus \(M\) is invertible.

For completeness, set
\[
 Lf=(L_j(f))_{j\in\mathcal J},\qquad
 Sc=\sum_{i\in\mathcal J}c_i(q_i-e_i).
\]
The finite replacement operator is \(B=I+SL\), and
\(I+LS=M\). Its inverse is the bounded operator
\[
 B^{-1}=I-SM^{-1}L.
\]
The operator \(B\) fixes the constant and every vector of the tail
basis outside \(\mathcal J\), while \(Be_j=q_j\) for
\(j\in\mathcal J\). Therefore \(BA_{\mathcal T}\) maps the
canonical unconditional basis onto the complete perturbed
family and is boundedly invertible on \(L^p\).
The bounded operators $BA_{\mathcal T}$ and $I+T$ agree on the
constant and every canonical basis vector, hence on their dense span.
Thus $BA_{\mathcal T}=I+T$, which proves the stated invertibility.
At \(p=2\), this is a Riesz basis.
\end{proof}

\begin{remark}
The conclusion is unchanged by finitely many admissible pole
displacements. In particular, large perturbations on finitely many
scales do not obstruct the basis conclusion when the tail packing
vanishes. The basis constants can depend on the finite replacement
matrix and need not be controlled by the tail norm alone.
The final assertion of Theorem~\ref{thm:stability} follows from
the dyadic comparison in Proposition~\ref{prop:dyadic-model-transfer}
and the constructions in Corollary~\ref{cor:dyadic-prescribed-excess}.
\end{remark}

\subsection{Summable displacements and endpoint stability}
\label{subsec:endpoint-pole-stability}
\begin{proposition}[Pole stability at the endpoints]
\label{prop:endpoint-pole-stability}
For the admissible complex poles of Theorem~\ref{thm:stability}, put
\[
 \eta_j=
 \sup\left\{
 \frac{|b_g-a_g|}{\ell_g}:
       2^{-j-1}\le\ell_g<2^{-j}
 \right\},
 \qquad j\ge0,
\]
with the supremum of an empty set defined as zero.
If
\begin{equation}
 \sum_{j\ge0}\eta_j<\infty,
 \label{eq:endpoint-displacement-sum}
\end{equation}
then the correction series defining \(T\) converges in operator
norm to a compact operator on every \(L^p(\mu)\),
\(1\le p\le\infty\), and on \(C(E)\). On all these spaces,
\[
 \|T\|\le C\sum_{j\ge0}\eta_j,
\]
where \(C\) depends only on the regularity data and the bound on
admissible pole displacements. The operator \(I+T\) is invertible on
every one of these spaces, without a smallness assumption on the
sum in \eqref{eq:endpoint-displacement-sum}.

Consequently, \(\{1\}\cup\{q_g\}\) is an unconditional basis in
every \(L^p(\mu)\), \(1<p<\infty\), and a Riesz basis in
\(L^2(\mu)\). Ordered by nonincreasing gap length, with arbitrary
ordering among equal lengths, it is a Schauder basis in
\(L^1(\mu)\) and \(C(E)\).
\end{proposition}

\begin{proof}
Write \(\mathcal G_j=\{g:2^{-j-1}\le\ell_g<2^{-j}\}\).
The envelopes \eqref{st:d-envelope} and~\eqref{st:e-envelope}
and upper regularity give
\[
 \|e_g\|_1\le C\ell_g^{s/2},\qquad
 \|d_g\|_1\le C\eta_g\ell_g^{s/2}.
\]
Indeed, for denominators of orders one and two, integration over
annuli gives convergent series with respective ratios
\(2^{s-1}\) and \(2^{s-2}\). Thus, for \(g\in\mathcal G_j\),
\[
 \begin{aligned}
 |d_g(x)|\,\|e_g\|_1
 &\le C\eta_j
       \frac{\ell_g^2}{(\ell_g+|x-\xi_g|)^2},\\
 |e_g(y)|\,\|d_g\|_1
 &\le C\eta_j
       \frac{\ell_g}{\ell_g+|y-\xi_g|}.
 \end{aligned}
\]
At the fixed scale \(t=2^{-j}\), the localized gap count bounds
the number of endpoints in a ball of radius \(2^kt\) by
\(C2^{ks}\). Summing the last two estimates over distance
annuli therefore gives, for every subfamily
\(\mathcal F\subset\mathcal G_j\),
\[
 \begin{aligned}
 \sup_x\int_E
   \sum_{g\in\mathcal F}|d_g(x)e_g(y)|\,d\mu(y)
 &\le C\eta_j\sum_{k\ge0}2^{-k(2-s)}
 \le C\eta_j,\\
 \sup_y\int_E
   \sum_{g\in\mathcal F}|d_g(x)e_g(y)|\,d\mu(x)
 &\le C\eta_j\sum_{k\ge0}2^{-k(1-s)}
 \le C\eta_j.
 \end{aligned}
\]
Both sums converge because \(s<1\).
The Schur estimates bound the norm by \(C\eta_j\) on every
\(L^p\), with one constant for all \(p\), and the row estimate
gives the same bound on \(C(E)\). Each band is finite and its
correction has continuous finite-dimensional range.
Summability in \eqref{eq:endpoint-displacement-sum} now proves
convergence in operator norm over arbitrary finite subsets and
compactness on every stated space.

Choose an integer \(J\) so large that
\[
 \kappa=C\sum_{j\ge J}\eta_j<1,
 \qquad
 \mathcal T=\{g:\ell_g<2^{-J}\},
 \qquad
 \mathcal J=\{g:\ell_g\ge2^{-J},\ b_g\ne a_g\}.
\]
The set \(\mathcal J\) is finite. If
\[
 T_{\mathcal T}f=\sum_{g\in\mathcal T}
                    \langle f,e_g\rangle d_g,
 \qquad A_{\mathcal T}=I+T_{\mathcal T},
\]
the Neumann series makes \(A_{\mathcal T}\) invertible on all
the stated spaces. Its inverse has norm at most
\((1-\kappa)^{-1}\). The same holds for
\(A_{\mathcal S}=I+\sum_{g\in\mathcal S}
\langle\,\cdot\,,e_g\rangle d_g\), for every finite
\(\mathcal S\subset\mathcal T\). All these operators and their
Neumann inverses are consistent on intersections of the spaces,
since their finite sums are the same integral operators.

If \(\mathcal J\) is empty, then \(I+T=A_{\mathcal T}\)
and invertibility is already proved. Otherwise, we adapt the
finite replacement argument above to the endpoints.
For a finite \(\mathcal S\subset\mathcal T\), set
\[
 R_{\mathcal S}(z)=
       \prod_{g\in\mathcal S}\frac{z-b_g}{z-a_g}.
\]
The preceding partial fraction identities give
\[
 \int_E R_{\mathcal S}\,d\mu=1,\qquad
 \int_E\frac{R_{\mathcal S}(x)}{x-z}\,d\mu(x)
       =F(z)R_{\mathcal S}(z).
\]
Hence the bilinear functional
\(\Lambda_{\mathcal S}(f)=\int_E fR_{\mathcal S}\,d\mu\)
annihilates all the nonconstant images of the canonical basis
under \(A_{\mathcal S}\) and takes the value one on the constant.
Completeness of that basis in \(L^2\) therefore gives
\[
 \Lambda_{\mathcal S}(f)
       =\int_E A_{\mathcal S}^{-1}f\,d\mu.
\]
The identity extends to \(L^1\): \(L^2\) is dense there,
\(R_{\mathcal S}\) is bounded on \(E\), and the \(L^1\) and
\(L^2\) Neumann inverses agree on \(L^2\).
Taking the \(L^1\) functional norm gives the endpoint bound
\begin{equation}
 \|R_{\mathcal S}\|_\infty
       \le(1-\kappa)^{-1}.
 \label{eq:endpoint-product-bound}
\end{equation}
No complex conjugation is inserted in this pairing.

Take an increasing finite exhaustion of \(\mathcal T\).
By \eqref{eq:endpoint-product-bound} and weak-* compactness,
a subsequence converges weak-* in \(L^\infty(\mu)\) to a
function \(R\). Here \(L^1(\mu)\) is separable, so a
subsequence suffices. The same locally uniform product argument
as above gives
\[
 P(z)=\prod_{g\in\mathcal T}\frac{z-b_g}{z-a_g},
 \qquad
 G(z):=\int_E\frac{R(x)}{x-z}\,d\mu(x)=F(z)P(z),
 \qquad
 \int_E R\,d\mu=1.
\]
In passing to the weak-* limit, we use the bounded test functions
\(1\) and \((x-z)^{-1}\), which belong to \(L^1(\mu)\) for
\(z\notin E\). The transform identity extends across the
removable product poles. For \(g\in\mathcal T\) it gives
\(G(b_g)=0\), and for \(g\notin\mathcal T\) it gives
\(G(a_g)=0\).
The admissible disks lie in disjoint vertical strips over the gaps.
Consequently,
for every \(j\in\mathcal J\),
\[
 G'(a_j)=F'(a_j)P(a_j)\ne0,\qquad G(b_j)\ne0,
\]
by the same argument excluding zeros as in the finite replacement
proof.

Define
\[
 L_j(f)=\frac{1}{\sqrt{w_j}\,G'(a_j)}
          \int_E\frac{f(x)R(x)}{x-a_j}\,d\mu(x),
 \qquad j\in\mathcal J.
\]
These functionals are bounded on \(L^1\), since \(R\in L^\infty\)
and \(a_j\) has positive distance from \(E\). They are therefore
bounded on every \(L^p\) and on \(C(E)\).
The same identity for divided differences shows that they annihilate
the constant and all the \(A_{\mathcal T}\)-images of the
canonical vectors except the corresponding \(e_j\), on which
they take the value one. Their finite replacement matrix is
\[
 M_{ji}=L_j(q_i)
       =\frac{\sqrt{w_i}\,G(b_i)}
       {\sqrt{w_j}\,G'(a_j)(b_i-a_j)},
 \qquad i,j\in\mathcal J.
\]
As in the preceding proof, its nonzero row and column factors
reduce invertibility to the Cauchy determinant for the distinct,
disjoint lists \((a_j)\) and \((b_i)\).

Let
\[
 Lf=(L_j(f))_{j\in\mathcal J},\qquad
 Sc=\sum_{i\in\mathcal J}c_i(q_i-e_i),\qquad B=I+SL.
\]
The range of \(S\) is contained in \(C(E)\), and
\(I+LS=M\). Thus both
\[
 B=I+SL,\qquad B^{-1}=I-SM^{-1}L
\]
are bounded on every stated space. On the constant and every
canonical basis vector, \(BA_{\mathcal T}\) agrees with \(I+T\).
The canonical span is dense in \(L^2\), so
\[
 BA_{\mathcal T}=I+T\quad\hbox{on }L^2.
\]
All these operators are also bounded on \(L^1\), and density
of \(L^2\) in \(L^1\) proves the identity there.
Their consistency then gives the identity on each \(L^p\),
including \(L^\infty\subset L^1\). On \(C(E)\), the
\(L^1\) identity is an equality of continuous functions and
therefore holds everywhere, since \(\supp\mu=E\).
No density of finite canonical expansions in \(L^\infty\)
is needed in this argument.
The factorization proves invertibility on every space, with
inverse \(A_{\mathcal T}^{-1}B^{-1}\).

Finally, \(I+T\) fixes the constant and sends \(e_g\) to \(q_g\).
It therefore maps the unconditional bases of
Proposition~\ref{thm:basis} and the ordered Schauder bases of
Proposition~\ref{thm:approximation-norming} to the asserted
perturbed bases. The inverse norms may depend on the finite
replacement matrix; only the tail was required to be small.
\end{proof}

\subsection{Cantor examples}
Let $E$ be the middle-thirds Cantor set, let $\mu$ be its natural
probability measure, and set $s=\log2/\log3$.
A basic interval of generation $n$ has length $3^{-n}$ and measure
$2^{-n}=(3^{-n})^s$. If $3^{-n}\leq r<3^{1-n}$, the basic
interval containing $x$ is contained in $B(x,r)$ up to endpoints
of zero measure. Hence
\[
 \mu(B(x,r))\geq2^{-n}
 \geq(r/3)^s=\tfrac12r^s.
\]
Distinct basic intervals of generation $n$ have separating gaps of
length at least $3^{-n}$. A ball of radius $r<3^{1-n}$ meets at
most eight such intervals, so
\[
 \mu(B(x,r))\leq8\cdot2^{-n}\leq8r^s.
\]
Thus the measure satisfies the regularity hypotheses with fixed constants.
There are $2^{n-1}$ complementary gaps of generation $n$, each
of length $3^{-n}$. Choose $0\leq\eta_n\leq c_0$, arbitrary
complex numbers $\sigma_g$ of modulus one, and set
\[
 b_g-a_g=\sigma_g\eta_n\ell_g.
\]
The criterion below is independent of these complex directions.
The total mass of the packing sequence is exactly
\[
\sum_g\lambda_g
 =\sum_{n\geq1}2^{n-1}\eta_n^2(3^{-n})^s
 =\frac12\sum_{n\geq1}\eta_n^2.
\]
The ball of radius two gives
\[
 \mathcal C(\lambda)\geq2^{-s}\sum_g\lambda_g
              =2^{-s-1}\sum_n\eta_n^2.
\]
For the upper bound, fix $B_E(x,R)$. If $3^{-n}\leq R$, the
endpoint packing argument in \eqref{st:band-count} gives
\[
 \#\{g:\ell_g=3^{-n},\ \xi_g\in B_E(x,R)\}
 \leq C(R/3^{-n})^s.
\]
Multiplying by the common weight and summing in $n$ yields
\[
 \begin{aligned}
 \sum_{\substack{\xi_g\in B_E(x,R)\\\ell_g\leq R}}\lambda_g
 &\leq C\sum_{3^{-n}\leq R}
       (R/3^{-n})^s\eta_n^2(3^{-n})^s\\
 &\leq CR^s\sum_{n\geq1}\eta_n^2.
 \end{aligned}
\]
Taking the supremum proves, with infinite values allowed,
\[
\mathcal C(\lambda)\asymp\sum_{n\geq1}\eta_n^2.
\]
\begin{corollary}[A Cantor classification]
\label{cor:generation-classification}
Let \(E\) be the middle-thirds Cantor set with its natural
probability measure, and let \(s=\log2/\log3\).
For each gap \(g\) of generation \(n\), set
\[
 b_g=a_g+\sigma_g\eta_n\ell_g,
 \qquad |\sigma_g|=1,\qquad 0\le\eta_n\le c_0,
\]
where the complex directions \(\sigma_g\) are arbitrary.
The following conditions are equivalent:
\begin{enumerate}
\item \(\sum_{n\ge1}\eta_n^2<\infty\).
\item \((q_g)_g\) is a Bessel sequence in \(L^2(\mu)\).
\item \(\{1\}\cup\{q_g\}_g\) is a Riesz basis of \(L^2(\mu)\).
\item \(\{1\}\cup\{q_g\}_g\) is an unconditional basis of
      \(L^p(\mu)\) for every \(1<p<\infty\).
\end{enumerate}
Whenever these conditions hold, \(T\) is compact on every
\(L^p(\mu)\), \(1<p<\infty\). No smallness of
\(\sum_n\eta_n^2\) is required.
\end{corollary}

\begin{proof}
Applying the preceding generation count after discarding the
first \(N\) generations gives
\[
 \mathcal C(\lambda)\asymp\sum_{n\ge1}\eta_n^2,
 \qquad
 \mathcal C_{3^{-N}}(\lambda)
       \asymp\sum_{n>N}\eta_n^2,
\]
with constants independent of the complex directions and of
\(N\). The lower bound uses the ball of radius two and the
identity \(\sum_{g\text{ of generation }n}\lambda_g
=\eta_n^2/2\); the upper bound is the same endpoint count as
before. Thus (1) is equivalent to (2), and (1) gives vanishing
tail packing, so Theorem~\ref{thm:stability} gives (3), (4),
and compactness. Condition (3) implies (2). Finally,
\(\|q_g\|_2\asymp1\) by the gap estimates; hence an
unconditional \(L^2\) basis is a Riesz basis, which proves
that (4) implies (3).
\end{proof}

For example, if \(\eta_n=\eta n^{-\beta}\), where
\(0<\eta\le c_0\) and \(\beta\ge0\), all four conditions
hold exactly when \(\beta>1/2\). At \(\beta=1/2\), the
system fails to be Bessel for every \(\eta>0\).

\begin{remark}[Comparison with quadratic perturbations]
The envelope estimate and the lower bounds on balls for \(d_g=q_g-e_g\)
give \(\|d_g\|_2^2\asymp\eta_g^2\). Therefore, in this
model,
\[
 \sum_g\|q_g-e_g\|_2^2
 \asymp\sum_{n\ge1}2^{n-1}\eta_n^2.
\]
The criterion in Corollary~\ref{cor:generation-classification}
is thus strictly weaker than quadratic closeness to the
canonical orthonormal basis. For \(\eta_n=\eta/n\),
\(0<\eta\le c_0\), the perturbation is compact and all the
interior \(L^p\) basis conclusions hold, although
\(D:\ell^2\to L^2(\mu)\) is not Hilbert--Schmidt.
\end{remark}

Taking $\eta_n=\eta>0$ at every generation gives
\[
 \sum_g\lambda_g=\frac{\eta^2}{2}\sum_{n\geq1}1=\infty.
\]
Hence a fixed, arbitrarily small relative displacement of every gap
pole need not produce a Bessel sequence. The failure can also be
verified directly.
Proposition~\ref{thm:basis} gives $\langle1,e_g\rangle=0$, and the
sector bound \eqref{st:sector-lower} together with
\eqref{st:d-lower} implies
\[
 \begin{aligned}
 |\langle1,q_g\rangle|^2
 &=|\langle1,d_g\rangle|^2\\
 &\geq\left(\Re\int_E\tau_gd_g\,d\mu\right)^2\\
 &\geq c\left(\int_E|d_g|\,d\mu\right)^2
 \geq c\eta^2\ell_g^s.
 \end{aligned}
\]
The sum of these squares is infinite by the preceding generation
count, contradicting the Bessel inequality tested at the constant
function.
In contrast, consider only the gaps
\[
 g_n=(3^{-n},2\cdot3^{-n}),\qquad n\geq1,
\]
move their poles vertically by $b_{g_n}-a_{g_n}=i\eta\ell_{g_n}$, and leave
all other poles fixed.  For every $0<R\leq2$, independently of
the center of a Carleson box,
\begin{align*}
 \sum_{\ell_{g_n}\leq R}\lambda_{g_n}
 &=\eta^2\sum_{3^{-n}\leq R}2^{-n}\\
 &\leq2\eta^2R^s.
\end{align*}
Indeed, let $N\geq1$ be least with $3^{-N}\leq R$. Then
\[
 \sum_{3^{-n}\leq R}2^{-n}
 =\sum_{n\geq N}2^{-n}
 =2^{1-N}
 =2(3^{-N})^s\leq2R^s.
\]
This includes $N=1$ and proves
$\mathcal C(\lambda)\leq2\eta^2$.
The sparse shifts preserve a Riesz basis when $\eta$ is sufficiently
small. For each fixed $p$, a sufficiently small $\eta$ in terms of
$p$ preserves an unconditional $L^p$ basis.
Their size is of order $\ell_{g_n}$, whereas
$w_{g_n}\asymp\ell_{g_n}^{2-s}$. Thus
\[
 \frac{|b_{g_n}-a_{g_n}|}{w_{g_n}}
 \asymp\eta\ell_{g_n}^{s-1}\longrightarrow\infty.
\]
Stability therefore depends on collective packing across gaps
and scales. The uniform and sparse constructions exhibit this
distinction for the same Cantor measure and the same bound on
relative displacements.

\subsection{Uniform realization of dyadic paraproduct models}
\label{subsec:dyadic-model-transfer}

To study the lower basis bound, we compare the Cauchy synthesis
operator with a dyadic model in operator norm.
We use the equally weighted Cantor measure with contraction ratio
$\rho$, and identify its cylinder tree with the binary tree $\mathcal T$.
For a cylinder $I$ of generation $n$ put
\[
 L_I=\rho^n,\qquad m_I=2^{-n},\qquad
 h_I=m_I^{-1/2}(\mathbf1_{I_{\rm right}}-\mathbf1_{I_{\rm left}}),
 \qquad u_I=m_I^{-1/2}\mathbf1_I.
\]
The root cylinder is included in $\mathcal T$.  Its central gap is
$g_I$, with $\ell_I=(1-2\rho)L_I$.  We abbreviate the canonical
functions and their poles by $e_I$ and $a_I$.
We identify all these $L^2$ spaces with the same Bernoulli space,
so their Haar coordinates are independent of $\rho$.

For a real array $\eta=(\eta_I)_{I\in\mathcal T}$ define
\begin{equation}
 M(\eta)^2=\sup_{J\in\mathcal T}\frac1{m_J}
                  \sum_{I\subset J}\eta_I^2m_I.
 \label{dyadic:carleson}
\end{equation}
Here and below $I\subset J$ allows equality.  The associated dyadic
paraproduct is
\[
 \Pi_\eta f=\sum_I\eta_I\sqrt{m_I}\,\langle f\rangle_I h_I,
 \qquad \langle f\rangle_I=m_I^{-1}\int_I f\,d\mu_\rho.
\]
The Carleson embedding theorem gives
$\|\Pi_\eta\|_{2\to2}\le2M(\eta)$.

\begin{proposition}[Uniform dyadic model transfer]
\label{prop:dyadic-model-transfer}
Fix $0<\eta_*<1/4$. There exist an absolute constant $\rho_0>0$ and a constant
$C_{\eta_*}$ with the following property.  Suppose
$0<\rho\le\rho_0$, $\eta_I\in\mathbb R$,
$|\eta_I|\le\eta_*$, and $M=M(\eta)<\infty$.  Set
\[
 b_I=a_I+\eta_I\ell_I,\qquad
 q_I(x)=\frac{\sqrt{w_I}}{x-b_I},\qquad
 \alpha_I=\frac1{1-4\eta_I^2},\qquad
 \beta_I=\frac{2\eta_I}{1-4\eta_I^2}.
\]
Let $W_\rho$ be the unitary defined by $W_\rho1=1$ and
$W_\rho h_I=e_I$.  Define, initially on finite Haar expansions,
\[
 A_\rho1=1,\quad A_\rho h_I=q_I,
 \qquad
 B_\eta1=1,\quad B_\eta h_I=\alpha_Ih_I+\beta_Iu_I.
\]
Then $A_\rho$ and $B_\eta$ extend to bounded operators and
\begin{equation}
 \|W_\rho-I\|_{2\to2}\le C\rho,
 \qquad
 \|A_\rho-W_\rho B_\eta\|_{2\to2}
       \le C_{\eta_*}\rho M.
 \label{dyadic:transfer}
\end{equation}
If $D_\alpha1=1$ and $D_\alpha h_I=\alpha_Ih_I$, then
\begin{equation}
 B_\eta=(I+2\Pi_\eta^*)D_\alpha.
 \label{dyadic:paraproduct-factorization}
\end{equation}
Moreover, the geometric packing norm of
$\lambda_I=\eta_I^2\ell_I^{s_\rho}$, where
$s_\rho=\log2/\log(1/\rho)$, satisfies
\begin{equation}
 \mathcal C(\lambda)\asymp M^2
 \label{dyadic:packing-equivalence}
\end{equation}
with absolute constants for $0<\rho\le\rho_0$.
If $\eta_*\le c_0<1/4$ and $2c_0<c_*<1/2$ are fixed, a further
decrease of $\rho$, depending on these constants, ensures
$|b_I-a_I|\le c_0\ell_I$ and
$\operatorname{dist}(a_I,E_\rho)\ge c_*\ell_I$ for every $I$.
Thus the original admissibility condition within each gap is retained.

In particular, suppose $B_\eta$ maps a closed complemented subspace
$X$ isomorphically onto $L^2$, and write
$R=(B_\eta|_X)^{-1}:L^2\to X$.
If $C_{\eta_*}\rho M\|R\|<1$, then
$A_\rho|_X:X\to L^2$ is also an isomorphism.
For any fixed closed complement $Z$ of $X$, $\ker A_\rho$ is a graph
over $Z$.  Consequently, whenever $B_\eta$ is surjective, sufficiently
small $\rho$ preserves its synthesis excess, whether finite or infinite.
Whenever $B_\eta$ is invertible, sufficiently small $\rho$ gives a
Riesz basis $\{1\}\cup\{q_I:I\in\mathcal T\}$.
\end{proposition}

\begin{proof}
We may take $\rho_0\le1/256$ and decrease it further as needed.
Every constant in the proof is uniform for this range of $\rho$.

\emph{Canonical estimates and comparison with the Haar basis.}
If $z_I$ is the midpoint of $I$, symmetry of $\mu_\rho|_I$ and the
decomposition into sibling cylinders along the path to the root give
\[
 \left|\int_{E_\rho\setminus I}\frac{d\mu_\rho(t)}{t-z}\right|
 \le \frac{m_I}{L_I}\frac{\rho}{(1-2\rho)^2},\qquad z\in g_I.
\]
The local derivative is at least $m_I/L_I^2$. Evaluating at
$z_I\pm2\rho L_I$ therefore locates the canonical zero between these
two points. On $I$,
the measure lies in the two end cylinders of relative length $\rho$;
the external derivative contributes only $O(\rho^2)m_I/L_I^2$.
It follows that
\begin{equation}
 |a_I-z_I|\le2\rho L_I,
 \qquad
 \sqrt{w_I}=\frac{L_I}{2\sqrt{m_I}}(1+O(\rho)),
 \qquad |e_I-h_I|\le C\rho u_I\quad\hbox{on }I.
 \label{dyadic:canonical-estimates}
\end{equation}
If $K$ is the ancestor of $I$ at distance $a\ge1$, then on the
sibling shell in $K\setminus I$ corresponding to that ancestor,
\[
 |e_I(x)|\le C\frac{\rho^a}{\sqrt{m_I}}.
\]
The shell has mass $2^{a-1}m_I$.  Thus
$\|e_I-h_I\|_2\le C\rho$, and, since both vectors are real and
normalized,
\begin{equation}
 |\langle e_I,h_I\rangle-1|\le C\rho^2.
 \label{dyadic:diagonal}
\end{equation}

We next estimate the off-diagonal entries, including the multiplicities
needed to bound the operator norm. If $J$ is a strict descendant of
$I$ at distance $d$, cancellation of $h_J$ and differentiation of
$e_I$ on $J$ give
\[
 |\langle e_I,h_J\rangle|
 \le C\frac{L_J}{L_I}\sqrt{\frac{m_J}{m_I}}
 =C(\rho/\sqrt2)^d.
\]
If $J$ is a strict ancestor of $I$ at distance $d$, let $c_J$ be
the constant value of $h_J$ on the child containing $I$.
The canonical cancellation $\int e_I\,d\mu_\rho=0$ gives
\[
 \langle e_I,h_J\rangle
 =\int e_I(h_J-c_J)\,d\mu_\rho.
\]
The other child of $J$ contributes at most
$C\rho^d\sqrt{m_J/m_I}$.  The ancestor shells outside $J$
contribute at most this quantity times
$\sum_{a\ge1}(2\rho)^a$.  Consequently,
\[
 |\langle e_I,h_J\rangle|\le C(\sqrt2\rho)^d.
\]
Finally, if $I,J$ are disjoint and have least common ancestor $K$,
write $u$ and $v$ for their respective distances below $K$.
The cylinders lie in opposite children of $K$, so cancellation on
$J$ gives
\[
 |\langle e_I,h_J\rangle|
 \le C\frac{L_IL_J}{L_K^2}\sqrt{\frac{m_J}{m_I}}
 =C\rho^{u+v}2^{(u-v)/2}.
\]
For a fixed input $I$, the descendant contribution to the column
sum is bounded by $C\sum_{d\ge1}(\sqrt2\rho)^d$, the ancestor
contribution by the same sum, and the disjoint contribution by
\[
 C\sum_{u,v\ge1}(\sqrt2\rho)^{u+v}.
\]
For a fixed output $J$, the contribution to the row sum from ancestor
inputs is bounded by $C\sum_{d\ge1}(\rho/\sqrt2)^d$, the
contribution from descendant inputs by $C\sum_{d\ge1}(2\sqrt2\rho)^d$,
and the disjoint contribution by
\[
 C\sum_{u,v\ge1}(2\sqrt2\rho)^u(\rho/\sqrt2)^v.
\]
The factors $2^d$ and $2^{u-1}$ count the relevant descendants.
Together with \eqref{dyadic:diagonal}, these estimates show that the absolute
row and column sums of the matrix of $W_\rho-I$ are $O(\rho)$.
The constant row and column vanish.  The Schur test proves the
first estimate in \eqref{dyadic:transfer}.

\emph{Two Carleson embedding estimates.}
For completeness, suppose that nonnegative numbers $c_I$ satisfy
$\sum_{I\subset J}c_I\le M^2m_J$. Selecting the maximal cylinders
for which $|\langle f\rangle_I|>t$ and integrating the distribution function give
\[
 \sum_Ic_I|\langle f\rangle_I|^2
 \le M^2\|\mathcal M_{\rm dyad}f\|_2^2
 \le4M^2\|f\|_2^2.
\]
Taking $c_I=\eta_I^2m_I$ and dualizing proves
\begin{equation}
 \left\|\sum_I c_I'|\eta_I|u_I\right\|_2
 \le2M\|(c_I')\|_{\ell^2}.
 \label{dyadic:indicator-embedding}
\end{equation}
This estimate also holds with $|c_I'|$ in place of $c_I'$.
For the ancestor $I^{(a)}$ at distance $a$, define
\[
 T_ac=\sum_{I:\,\operatorname{gen}(I)\ge a}
               c_I|\eta_I|u_{I^{(a)}}.
\]
For its adjoint, we have
\[
 \|T_a^*f\|_{\ell^2}^2
 =\sum_J m_J\left(\sum_{I^{(a)}=J}\eta_I^2\right)
                         |\langle f\rangle_J|^2.
\]
The coefficient measure on the right has Carleson norm at most
$2^aM^2$, because for every cylinder $K$,
\[
 \sum_{J\subset K}m_J\sum_{I^{(a)}=J}\eta_I^2
 =2^a\sum_{I:\,I^{(a)}\subset K}\eta_I^2m_I
 \le2^aM^2m_K.
\]
Thus
\begin{equation}
 \|T_a\|_{\ell^2\to L^2}\le2\,2^{a/2}M.
 \label{dyadic:ancestor-embedding}
\end{equation}

\emph{Comparison of the perturbation columns.}
Put
\[
 d_I=q_I-e_I,\qquad
 \widehat d_I=(\alpha_I-1)h_I+\beta_Iu_I.
\]
On the left and right children of $I$, respectively, the limiting
values of $q_I$ in units of $m_I^{-1/2}$ are
$-(1+2\eta_I)^{-1}$ and $(1-2\eta_I)^{-1}$.
They equal the values of $\alpha_Ih_I+\beta_Iu_I$.
After division by $L_I$, the pole is
$1/2+O(\rho)+\eta_I(1-2\rho)$, the weight factor is
$(1+O(\rho))/2$, and the support lies in
$[0,\rho]\cup[1-\rho,1]$.
The derivatives with respect to $\eta_I$ of these rational formulas
differ from their limiting derivatives by at most
$C_{\eta_*}\rho$.  All denominators stay uniformly away from zero.
At $\eta_I=0$ both perturbations vanish; integrating the derivative
comparison from $0$ to $\eta_I$ gives
\[
 |d_I-\widehat d_I|\le C_{\eta_*}\rho|\eta_I|u_I
                   \quad\hbox{on }I.
\]
On an ancestor shell corresponding to $K=I^{(a)}$, the identity
\[
 d_I(x)=\frac{\sqrt{w_I}\eta_I\ell_I}
                  {(x-b_I)(x-a_I)}
\]
gives
\[
 |d_I(x)|\le C_{\eta_*}|\eta_I|
            \frac{\rho^{2a}}{\sqrt{m_I}}
 =C_{\eta_*}|\eta_I|(\sqrt2\rho^2)^a u_{I^{(a)}}(x).
\]
Since $\widehat d_I$ is supported on $I$, the global bound is
\[
 |d_I-\widehat d_I|
 \le C_{\eta_*}|\eta_I|\left(
       \rho u_I+\sum_{a=1}^{\operatorname{gen}(I)}
                     (\sqrt2\rho^2)^a u_{I^{(a)}}\right).
\]
Apply \eqref{dyadic:indicator-embedding} to the first term and
\eqref{dyadic:ancestor-embedding} to the others.  For finite sequences,
\begin{align*}
 \left\|\sum_I c_I(d_I-\widehat d_I)\right\|_2
 &\le C_{\eta_*}M
       \left(\rho+\sum_{a\ge1}(2\rho^2)^a\right)\|c\|_{\ell^2}\\
 &\le C_{\eta_*}\rho M\|c\|_{\ell^2}.
\end{align*}
Because $\sup_I|\eta_I|\le M$,
$|\alpha_I-1|\le C_{\eta_*}|\eta_I|$, and
$|\beta_I|\le C_{\eta_*}|\eta_I|$, the same embedding also gives
$\|B_\eta-I\|\le C_{\eta_*}M$.
We have proved
\[
 \|(A_\rho-W_\rho)-(B_\eta-I)\|\le C_{\eta_*}\rho M.
\]
Subtract $(W_\rho-I)(B_\eta-I)$ and use the first part of
\eqref{dyadic:transfer} to obtain its second part.  This also proves
boundedness of $A_\rho$.  The factorization
\eqref{dyadic:paraproduct-factorization} follows directly from
$\Pi_\eta^*h_I=\eta_Iu_I$ and $\Pi_\eta^*1=0$.

\emph{Packing and preservation of the range.}
For $\rho<1/4$, the measure satisfies
\[
 \tfrac12 R^{s_\rho}\le\mu_\rho(B(x,R))\le2R^{s_\rho},
 \qquad x\in E_\rho,\quad0<R\le1.
\]
Indeed, if $\rho^n\le R<\rho^{n-1}$, the level-$n$ cylinder
containing $x$ lies in the ball, whereas the ball can meet at most
one level-$(n-1)$ cylinder: distinct such cylinders have separation
at least $(1-2\rho)\rho^{n-2}>2\rho^{n-1}$.
Endpoint equalities concern sets of measure zero.
The upper bound extends to $R>1$.

Since $\lambda_I=(1-2\rho)^{s_\rho}\eta_I^2m_I$, the box centered
at the left endpoint of a cylinder $J$ with radius $L_J$ contains
all coefficients below $J$.  This gives
$M^2\le(1-2\rho)^{-s_\rho}\mathcal C(\lambda)$.
Conversely, every cylinder counted by a box $(x,R)$ lies in
$B(x,3R)$, since its gap length is at most $R$ and its length is
at most $R/(1-2\rho)$.
Take the maximal cylinders among those counted by the box.
They are disjoint, so \eqref{dyadic:carleson} bounds the coefficient
sum by $M^2\mu_\rho(B(x,3R))\le6M^2R^{s_\rho}$.
This proves \eqref{dyadic:packing-equivalence}.
The centrality estimate in \eqref{dyadic:canonical-estimates} is
$\operatorname{dist}(a_I,E_\rho)/\ell_I\ge
(1/2-3\rho)/(1-2\rho)$, which tends to $1/2$.
It gives the additional admissibility assertion.

Finally set $\widetilde A=W_\rho^*A_\rho$.
Then $\|\widetilde A-B_\eta\|\le C_{\eta_*}\rho M$, and
\[
 \widetilde A R=I+(\widetilde A-B_\eta)R
\]
is invertible under the asserted smallness condition.  It follows
that $\widetilde A|_X$, and therefore $A_\rho|_X$, is an isomorphism.
For $L^2=X\oplus Z$, its kernel is exactly
\[
 \{z-(A_\rho|_X)^{-1}A_\rho z:z\in Z\}.
\]
If $B_\eta$ is surjective, choose $X=(\ker B_\eta)^\perp$;
the open mapping theorem supplies the bounded inverse $R$.
The graph formula proves the assertion about excess, including
infinite-dimensional kernels.
\end{proof}

\begin{corollary}[A classical reverse H\"older sufficient condition]
\label{cor:dyadic-rh-positive}
Under the hypotheses of Proposition~\ref{prop:dyadic-model-transfer},
assume that the paraexponential
\[
 R_\eta=\prod_{I\in\mathcal T}
              (1-2\eta_I\sqrt{m_I}h_I)
\]
exists as a positive weight of mean one and belongs to the dyadic
reverse H\"older class $RH_2^{\rm d}$; that is,
\[
 \sup_{J\in\mathcal T}
 \frac{\bigl(m_J^{-1}\int_J R_\eta^2\,d\mu\bigr)^{1/2}}
      {m_J^{-1}\int_J R_\eta\,d\mu}<\infty.
\]
The product is understood as the almost everywhere limit of products
over the first finitely many generations; the required dyadic doubling
follows from the uniform bound on the factors, as verified below.
Then $B_\eta$ is invertible, and the Cauchy family $\{1\}\cup(q_I)_I$
is a Riesz basis whenever $\rho$
is sufficiently small.  This is a sufficient condition for the
model, obtained from the classical paraproduct resolvent theorem;
no converse for a fixed positive $\rho$ is asserted.
\end{corollary}

\begin{proof}
The symbol $b=\sum_I\eta_I\sqrt{m_I}h_I$ is in dyadic BMO by
\eqref{dyadic:carleson}.  The pointwise factors satisfy
$|2\eta_I\sqrt{m_I}h_I|\le2\eta_*<1$, uniformly in $I$.
Write $R_n$ for the product through generation $n-1$ and
$\mathcal F_n$ for the level-$n$ cylinder sigma-algebra.
The positive martingale $R_n$ has mean one; its assumed almost everywhere
limit $R_\eta$ has the same mean.  Scheff\'e's lemma gives
$R_n\to R_\eta$ in $L^1$, and hence
$\mathbb E(R_\eta\mid\mathcal F_n)=R_n$.
For either child $I'$ of $I$ this implies
\[
 \frac{\int_{I'}R_\eta\,d\mu}{\int_I R_\eta\,d\mu}
 =\frac{1\mp2\eta_I}{2}\ge\frac{1-2\eta_*}{2}>0.
\]
Thus $R_\eta\,d\mu$ is dyadically doubling, as required in the
class $RH_2^{\rm d}$ used by the cited theorem.
Identifying the Bernoulli filtration with the ordinary dyadic
filtration, Pereyra's resolvent theorem, in the form restated in
\cite[Theorem~1.3]{PereyraWard1998}, applies with parameter $-2$.
It gives invertibility of $I+2\Pi_\eta$ on the subspace of functions with mean zero
$L^2_0$. On $\mathbb C1\oplus L^2_0$ its matrix is
\[
 \begin{pmatrix}
  1&0\\
  2\Pi_\eta1&(I+2\Pi_\eta)|_{L^2_0}
 \end{pmatrix}.
\]
This matrix is invertible, so its adjoint $I+2\Pi_\eta^*$ is
invertible as well.  The diagonal $D_\alpha$ is bounded and
invertible.  The factorization
\eqref{dyadic:paraproduct-factorization} and the transfer theorem
complete the proof.
\end{proof}

\begin{corollary}[Prescribed finite or infinite excess]
\label{cor:dyadic-prescribed-excess}
There is a sufficiently small contraction ratio $\rho>0$ such that,
for every positive integer $N$, the same Cantor measure $\mu_\rho$
admits real shifts within the original gaps satisfying
\[
 |b_I-a_I|\le\frac6{25}\ell_I,
 \qquad \operatorname{dist}(a_I,E_\rho)\ge\frac{49}{100}\ell_I,
 \qquad \mathcal C(\lambda)\le C
\]
for which the augmented family $\{1\}\cup(q_I)_I$ is a frame with
synthesis kernel of dimension exactly $N$.
One may delete the constant and $N-1$ specified unshifted gap
functions to obtain a Riesz basis.  There is also such a frame whose
synthesis kernel is infinite-dimensional.  Both frame bounds and
the packing bounds can be taken independent of $N$, including the
example with infinite excess, and every construction has
$\gamma(\lambda)>0$.
\end{corollary}

\begin{proof}
Fix $\eta=9/40$.  Choose a finite full binary subtree with $N$
leaves, whose cylinder sets $J_1,\ldots,J_N$ partition $E_\rho$.
On each leaf put $\eta_I=\eta$ exactly for the chain
$I=J_j0^n$, $n\ge0$, and put $\eta_I=0$ elsewhere.
For every cylinder $K$, the sum of $\eta_I^2m_I$ below $K$ is at
most $2\eta^2m_K$: below a selected root the chain masses form a
geometric series, and distinct roots are disjoint.
Thus $M^2=2\eta^2$, independently of $N$.

Let $h_{j,n}=h_{J_j0^n}$ and
$Y_j=\overline{\operatorname{span}}\{u_{J_j},h_{j,n}:n\ge0\}$.
These subspaces are mutually orthogonal, are invariant under
$B_\eta$, and $B_\eta$ is the identity on their joint orthogonal
complement.  Also $B_\eta u_{J_j}=u_{J_j}$: each $u_{J_j}$ is a
finite combination of the constant and the Haar functions at the
unshifted internal nodes of the finite subtree.
The recurrence $u_{J0^{n+1}}=(u_{J0^n}-h_{J0^n})/\sqrt2$ gives
\[
 u_{J0^n}=r^nu_J-\sum_{k=0}^{n-1}r^{n-k}h_{J0^k},
 \qquad r=2^{-1/2}.
\]
In the coordinates $\mathbb C u_{J_j}\oplus\ell^2(\mathbb N_0)$,
each block is therefore identical to
\[
 \begin{pmatrix}
  1&\beta(r^n)_{n\ge0}\\
  0&\alpha(I-\sigma S^*)(I-rS^*)^{-1}
 \end{pmatrix},\qquad
 r=2^{-1/2},\quad
 \alpha=(1-4\eta^2)^{-1},\quad\beta=2\eta\alpha,
 \quad\sigma=(1+2\eta)r>1.
\]
Here $S^*$ is the backward shift.  If $S$ is the forward shift,
then $-\sigma^{-1}S(I-\sigma^{-1}S)^{-1}$ is a bounded right
inverse of $I-\sigma S^*$, whose kernel is spanned by
$(\sigma^{-n})_{n\ge0}$.  Thus the lower-right block is
surjective with this same one-dimensional kernel.
The kernel vector of the full block has nonzero $u_{J_j}$
coordinate, namely $-\beta/(1-r/\sigma)$ when its chain
coordinates are $(\sigma^{-n})_{n\ge0}$.
Thus the restriction of each block to the closed span of the
Haar functions along its chain is an isomorphism onto $Y_j$, with inverse
norm depending only on $\eta$.

Put $F=\operatorname{span}\{u_{J_j}:1\le j\le N\}$ and
$X=F^\perp$.  The preceding orthogonal decomposition proves that
$B_\eta|_X:X\to L^2$ is an isomorphism with inverse norm bounded
independently of $N$.  The transfer theorem therefore makes
$A_\rho|_X$ an isomorphism for one sufficiently small $\rho$
independent of $N$.  The space $F$ is spanned by the constant and
the $N-1$ Haar functions at the internal nodes of the finite
subtree.  Hence $A_\rho|_X$ is precisely the synthesis map obtained
by deleting those $N$ coordinates, and the graph formula gives
$\dim\ker A_\rho=N$.

For infinite excess use the prefix-free roots $J_j=1^j0$,
$j\ge0$.  They are pairwise disjoint and cover $E_\rho$ up to the
single null endpoint corresponding to $111\cdots$.
The same Carleson bound and the same mutually orthogonal model
blocks apply.  Their uniform inverse bound implies that
$B_\eta$ restricts to an isomorphism on
$X=(\overline{\operatorname{span}}\{u_{J_j}:j\ge0\})^\perp$.
The same $\rho$ gives the corresponding isomorphism for $A_\rho$,
and its synthesis kernel is a graph over the infinite-dimensional
span of the $u_{J_j}$.

The canonical estimate \eqref{dyadic:canonical-estimates} gives
\[
 \operatorname{dist}(a_I,E_\rho)
 \ge\frac{1/2-3\rho}{1-2\rho}\ell_I.
\]
For sufficiently small $\rho$ this is at least
$(49/100)\ell_I$, whereas $\eta<6/25<(49/100)/2$.
Thus the shifts satisfy the original geometric admissibility
condition.  The uniform packing bound follows from
\eqref{dyadic:packing-equivalence}.
Every selected chain has arbitrarily small gaps with relative
shift $\eta$.  The box centered at an endpoint of one such gap,
with radius its length, gives
$\gamma(\lambda)\ge\eta^2>0$.
Finally, the bound $\|B_\eta-I\|\le C M$ and the transfer estimate
give a uniform upper bound for $\|A_\rho\|$.
The proof by a Neumann series also gives a uniform bound for
$\|(A_\rho|_X)^{-1}\|$, since $M$ and the model inverse bound are
independent of the forest.  If $R_A=(A_\rho|_X)^{-1}$, then
$R_A^*A_\rho^*=I$, so
\[
 \|R_A\|^{-2}\|f\|_2^2
 \le |\langle f,1\rangle|^2+\sum_I|\langle f,q_I\rangle|^2
 \le\|A_\rho\|^2\|f\|_2^2.
\]
These are the asserted uniform frame bounds for both finite and
infinite excess.
\end{proof}

\begin{remark}[Classical model and geometric realization]
\label{rem:dyadic-prior-work}
The paraproduct model and its resolvent obstructions are classical;
see Pereyra~\cite{Pereyra1994} and
Pereyra and Ward~\cite[Theorem~1.3]{PereyraWard1998}.
For a single chain, with $c=2\eta$, the corresponding
paraexponential takes the value $(1-c)(1+c)^n$ on the $n$th
sibling shell.  Its square integral changes from finite to
infinite at $(1+c)^2=2$, exactly the threshold $\sigma=1$
in the displayed model.
Finite excess for discrete Cauchy transforms also appears in
\cite[Introduction, Example~1]{BelovMengestieSeip2011}.
Proposition~\ref{prop:dyadic-model-transfer} gives a realization,
uniform in operator norm, of all the displayed Carleson paraproduct
models with $\sup_I|\eta_I|<1/4$ by Cauchy poles on a
nonatomic regular Cantor set, with the poles remaining in their
original gaps.
\end{remark}

\subsection{An exact index criterion on endpoint chains}
\label{subsec:exact-chain-index}

The comparison with a dyadic model leaves an error of order $\rho$.
For finitely many endpoint chains, self-similarity gives an exact
criterion at each fixed contraction ratio. The limiting matrix is
Toeplitz, and its symbol is a quotient of theta functions. Its only
possible zero in the admissible range gives the transition between a
basis and a frame with positive excess.

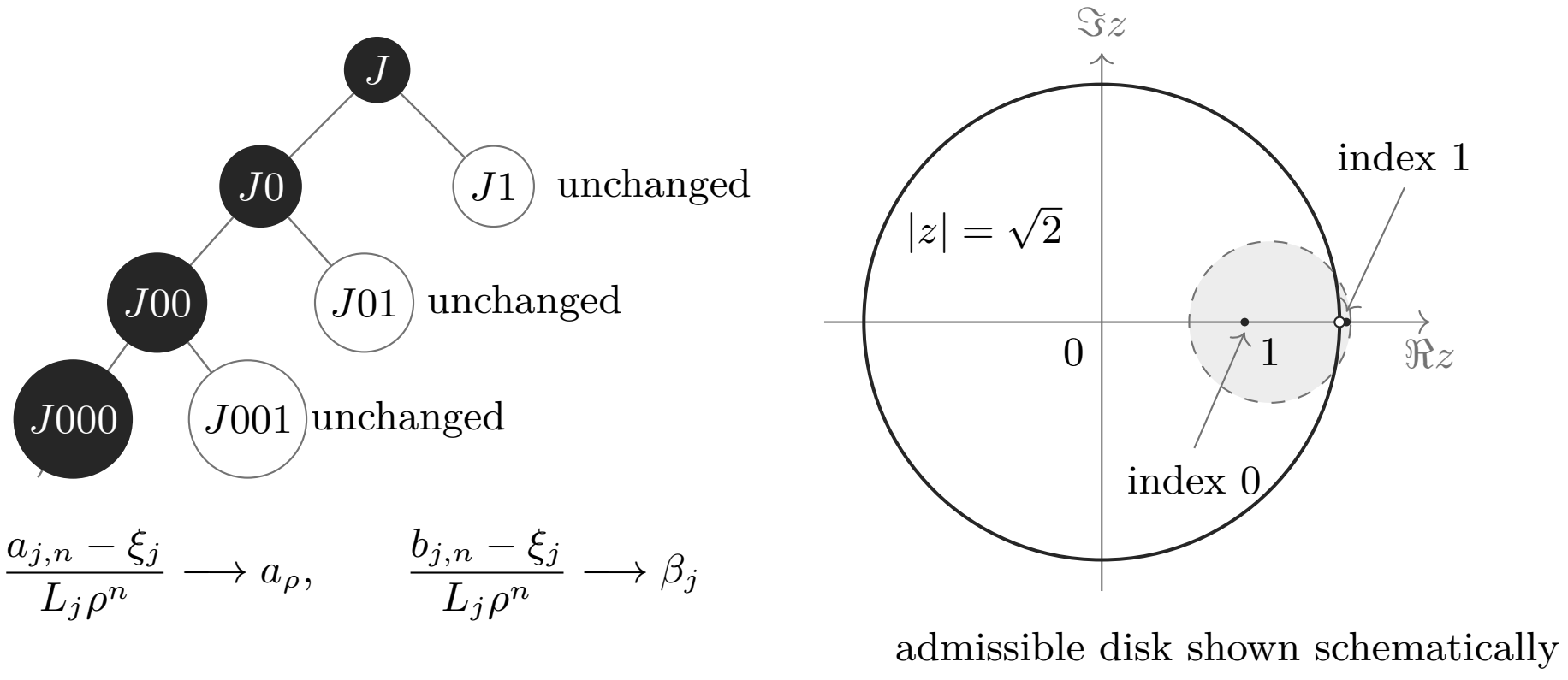
\begin{figure}[htbp]
\centering
\begin{tikzpicture}[x=1cm,y=1cm,font=\footnotesize,
 branch/.style={draw=DiagramSecondary,line width=.45pt},
 active/.style={circle,draw=DiagramPrimary,fill=DiagramPrimary,
                 text=white,inner sep=2pt},
 passive/.style={circle,draw=DiagramSecondary,fill=white,inner sep=2pt}]
 \node at (2.5,4.6) {(a) An endpoint chain};
 \coordinate (r) at (2.7,3.7);
 \coordinate (a) at (1.8,2.8);
 \coordinate (b) at (1.0,1.9);
 \coordinate (c) at (.35,1.0);
 \draw[branch] (r)--(a) (r)--(3.6,2.8)
               (a)--(b) (a)--(2.6,1.9)
               (b)--(c) (b)--(1.7,1.0);
 \draw[branch,dashed] (c)--(.05,.5);
 \node[active] at (r) {$J$};
 \node[active] at (a) {$J0$};
 \node[active] at (b) {$J00$};
 \node[active] at (c) {$J000$};
 \node[passive] at (3.6,2.8) {$J1$};
 \node[passive] at (2.6,1.9) {$J01$};
 \node[passive] at (1.7,1.0) {$J001$};
 \node[anchor=west] at (3.95,2.8) {unchanged};
 \node[anchor=west] at (2.95,1.9) {unchanged};
 \node[anchor=west] at (2.05,1.0) {unchanged};
 \node[align=center] at (2.5,-.2)
  {$\displaystyle
    \frac{a_{j,n}-\xi_j}{L_j\rho^n}\longrightarrow a_\rho,
    \qquad
    \frac{b_{j,n}-\xi_j}{L_j\rho^n}\longrightarrow\beta_j$};

 \begin{scope}[shift={(8.3,1.75)},scale=1.3]
  \fill[DiagramSecondary!13] (1,0) circle (.48);
  \draw[->,DiagramSecondary] (-1.65,0)--(1.95,0)
       node[below] {$\Re z$};
  \draw[->,DiagramSecondary] (0,-1.6)--(0,1.6)
       node[above] {$\Im z$};
  \draw[DiagramPrimary,line width=.75pt] (0,0) circle ({sqrt(2)});
  \draw[DiagramSecondary,dashed] (1,0) circle (.48);
  \fill[DiagramPrimary] (.85,0) circle (.025);
  \fill[DiagramPrimary] (1.455,0) circle (.025);
  \draw[DiagramPrimary,fill=white] ({sqrt(2)},0) circle (.03);
  \node[below left] at (0,0) {$0$};
  \node[below] at (1,0) {$1$};
  \node at (-.7,.55)
       {$|z|=\sqrt2$};
  \draw[->,DiagramSecondary] (.55,-.75)--(.85,-.06);
  \node[below,align=center] at (.55,-.75) {index $0$};
  \draw[->,DiagramSecondary] (1.8,.8)--(1.455,.06);
  \node[above,align=center] at (1.8,.8) {index $1$};
  \node[align=center] at (.75,-1.95) {admissible disk shown schematically};
 \end{scope}
 \node at (8.65,4.6) {(b) The normalized pole $z=\beta/a_\rho$};
\end{tikzpicture}
\caption{Endpoint chains and the index transition. Scaling each gap by
its cylinder length gives limiting poles $a_\rho$ and $\beta_j$.
Within the admissible disk, each chain contributes zero to the index
inside the critical circle and one outside. On the circle, the family is complete but has no
frame lower bound. The contributions of finitely many disjoint
chains add. The tree spacing and admissible disk are schematic.}
\label{fig:chain-index}
\end{figure}

The notation and hypotheses are those of
Theorem~\ref{prop:exact-chain-index}. Figure~\ref{fig:chain-index}
shows the endpoint geometry and the normalized critical circle.

\begin{proof}[Proof of Theorem~\ref{prop:exact-chain-index}]
We give the completeness argument and the operator reduction separately.
Neither requires passage to the limit $\rho\downarrow0$.

Let $\xi_j$ be the left endpoint of $J_j$, and let
$L_j$ and $m_j$ be its length and mass. Write
\[
 a_{j,n}=\xi_j+L_j\rho^n\alpha_{j,n},\qquad
 b_{j,n}=\xi_j+L_j\rho^n\beta_{j,n},\qquad
 \tau=\sqrt2\rho.
\]
The canonical estimates give, uniformly in $j,n$,
\begin{equation}
 |\alpha_{j,n}-\tfrac12|\le2\rho,
 \qquad
 \beta_{j,n}=\alpha_{j,n}+(1-2\rho)\eta_{j,n},
 \qquad
 \sqrt{w_{j,n}}\asymp\frac{L_j}{\sqrt{m_j}}\tau^n.
 \label{idx:chain-centrality}
\end{equation}
The coefficient array has finite packing: within a selected chain,
its descendant masses sum to at most twice the mass of its largest
member, and the roots are disjoint. The finite additional changes
also have finite packing. Theorem~\ref{thm:stability} therefore gives
bounded synthesis.

Suppose $f$ is orthogonal to the family. The canonical vectors are
real, and the unchanged ones, including $1$, remain in the family.
Consequently $h=\overline f$ belongs to the closed span of the
canonical vectors at the changed gaps. Its expansion has square
summable coefficients, say $h=\sum c_g e_g$, supported on the
selected chains and a finite set. For $z$ outside $E_\rho$ and the
canonical poles,
\[
 \int\frac{h(x)}{x-z}\,d\mu_\rho(x)=F(z)R(z),
 \qquad
 R(z)=\sum_g\frac{c_g\sqrt{w_g}}{z-a_g}.
\]
Indeed, the identity holds term by term because $F(a_g)=0$;
local uniform convergence follows from Cauchy--Schwarz and
\eqref{idx:chain-centrality}. Orthogonality gives $R(b_g)=0$
whenever $b_g\ne a_g$. Here $F(b_g)\ne0$: nonreal points have a
nonzero imaginary part under the Cauchy transform, and the canonical
pole is its only zero in the corresponding real gap.
If $b_g=a_g$, orthogonality to $e_g$ gives $c_g=0$.
The selected canonical poles accumulate only at the finitely many
chain endpoints. On a compact set avoiding those endpoints and
poles, the denominators in $R$ are uniformly bounded below.
Since the selected weights have finite sum, Cauchy--Schwarz gives
local uniform convergence there. Thus $R$ extends meromorphically
through the rest of $E_\rho$.

The product
\[
 P(z)=\prod_g\frac{z-b_g}{z-a_g}
\]
converges locally uniformly away from the finitely many endpoints
$\xi_j$ and its poles. This follows from
$\sum_g|b_g-a_g|<\infty$. Its only zeros and poles are the indicated
ones, with factors equal to one omitted. Thus $G=R/P$ is holomorphic
on $\mathbb C\setminus\{\xi_1,\ldots,\xi_N\}$.

Fix one endpoint $\xi_j$ and set
$r_n=L_j\rho^{n+1/2}$. For all sufficiently large $n$, the circle
$|z-\xi_j|=r_n$ avoids every other chain and all finite exceptional
poles. From \eqref{idx:chain-centrality},
\[
 \begin{split}
 \sum_{k\ge0}\frac{w_{j,k}}{|z-a_{j,k}|^2}
 &\le C_j\left(\sum_{k\le n}2^k+
       r_n^{-2}L_j^2\sum_{k>n}(2\rho^2)^k\right)\\
 &\le C_j2^n.
 \end{split}
\]
The other terms in $R$ stay bounded there, so
$|R(z)|\le C_f2^{n/2}$. On the same circle, factor the chain product
at level $n$:
\[
 P_j(z)=
 \prod_{k\le n}\frac{b_{j,k}-\xi_j}{a_{j,k}-\xi_j}
 \prod_{k\le n}
 \frac{1-(z-\xi_j)/(b_{j,k}-\xi_j)}
      {1-(z-\xi_j)/(a_{j,k}-\xi_j)}
 \prod_{k>n}\frac{1-(b_{j,k}-\xi_j)/(z-\xi_j)}
                   {1-(a_{j,k}-\xi_j)/(z-\xi_j)}.
\]
The last two products and their reciprocals are uniformly bounded.
To see this, each argument subtracted from one has modulus at most a fixed
constant times $\rho^{|n-k|+1/2}$, and the resulting geometric
series is bounded. All these arguments have modulus less than $1/2$ by
$\rho\le1/256$ and \eqref{idx:chain-centrality}. Moreover,
\[
 \left|\frac{b_{j,k}-\xi_j}{a_{j,k}-\xi_j}\right|
 =\frac{|\beta_{j,k}|}{\alpha_{j,k}}>\frac14.
\]
The factors from the other chains and the finite exceptional set
are bounded above and below near $\xi_j$. Hence
\[
 |P(z)|\ge c_j4^{-n},\qquad
 |G(z)|\le C_f(4\sqrt2)^n
 \quad\hbox{on }|z-\xi_j|=r_n.
\]
For each integer $k\ge1$, the coefficient of
$(z-\xi_j)^{-k}$ in the Laurent expansion of $G$ is bounded by
\[
 C_f r_n^k(4\sqrt2)^n
 =C_fL_j^k\rho^{k/2}(4\sqrt2\rho^k)^n\longrightarrow0.
\]
All negative Laurent coefficients vanish. Every endpoint is
therefore removable. Finally $R(z)=O(1/z)$ and $P(z)\to1$ at
infinity. Liouville's theorem gives $G=0$, then $R=0$ and $f=0$.
This proves completeness for every admissible array on these finitely
many chains, even without the convergence assumption in
\eqref{idx:chain-array}.

We now identify the limiting matrix. The shell decomposition
$E_\rho\setminus\{0\}=\bigcup_{k\ge0}\rho^k u(E_\rho)$ gives
\[
 (2\rho)^nF(\rho^nz)=
 \sum_{k\ge-n}2^{-k-1}
       \int\frac{d\mu_\rho(t)}{\rho^k u(t)-z}.
\]
The tail over positive indices converges geometrically with ratio $1/2$;
the tail over negative indices converges geometrically with ratio $2\rho$.
The same local convergence holds after differentiation. This proves
the assertions of convergence in \eqref{idx:blowup-transform}.
For a chain based at $J_j$, put
\[
 F_{j,n}(z)=\frac{L_j}{m_j}(2\rho)^n
                  F(\xi_j+L_j\rho^nz).
\]
Its contribution from $J_j$ is the preceding truncated series.
The contribution from $E_\rho\setminus J_j$ and its derivatives
tend locally uniformly to zero. Indeed, this support has positive
distance from $\xi_j$, whereas $2\rho<1$.
Thus $F_{j,n}\to\mathfrak F_\rho$ locally uniformly with derivatives.
The zeros $\alpha_{j,n}$ lie in the fixed central interval in
\eqref{idx:chain-centrality}. Since $\mathfrak F_\rho'>0$ in the
base gap, their limits agree with its unique zero $a_\rho$.
Writing
\[
 d_{j,n}=F_{j,n}'(\alpha_{j,n})^{-1/2},
 \qquad d=\mathfrak F_\rho'(a_\rho)^{-1/2},
\]
we obtain
\begin{equation}
 \alpha_{j,n}\longrightarrow a_\rho,
 \quad \beta_{j,n}\longrightarrow\beta_j,
 \quad d_{j,n}\longrightarrow d,
 \quad
 \sqrt{w_{j,n}}=\frac{L_j}{\sqrt{m_j}}\tau^nd_{j,n}.
 \label{idx:blowup-limits}
\end{equation}
The bound on $a_\rho$ follows by taking limits in
\eqref{idx:chain-centrality}. Positivity of the imaginary part shows
that $\mathfrak F_\rho$ has no nonreal zero.

Identify the source of synthesis with $L^2(\mu_\rho)$ by the canonical
orthonormal basis. Let $\mathcal H$ be the closed span of the
canonical vectors on the selected chains. Finite exceptional
replacements change the operator by finite rank and can be omitted
for the Fredholm calculation. Relative to
$\mathcal H\oplus\mathcal H^\perp$, synthesis has the form
\[
 A=\begin{pmatrix}C&0\\D&I\end{pmatrix},\qquad
 A\begin{pmatrix}I&0\\-D&I\end{pmatrix}
 =\begin{pmatrix}C&0\\0&I\end{pmatrix}.
\]
The right factor is bounded and invertible. Hence the Fredholm
problem reduces to the compression $C$.

Within the $j$th chain, the Cauchy identity gives
\[
 \begin{split}
 C^{(j)}_{mn}
 &=\langle q_{j,n},e_{j,m}\rangle
 =\sqrt{w_{j,n}w_{j,m}}
       \frac{F(b_{j,n})}{b_{j,n}-a_{j,m}}\\
 &=d_{j,m}d_{j,n}F_{j,n}(\beta_{j,n})
       \frac{\tau^{m-n}}
            {\beta_{j,n}-\rho^{m-n}\alpha_{j,m}},
 \end{split}
\]
At $b_{j,n}=a_{j,n}$ the diagonal quotient is interpreted by
continuity. The centrality bounds imply
\[
 |C^{(j)}_{mn}|\le
 \begin{cases}
 C\tau^{m-n},&m>n,\\
 C,&m=n,\\
 C2^{(m-n)/2},&m<n.
 \end{cases}
\]
The same bound holds for the limiting diagonals
\[
 c_k(\beta)=
 \frac{\mathfrak F_\rho(\beta)}{\mathfrak F_\rho'(a_\rho)}
       \frac{\tau^k}{\beta-a_\rho\rho^k},\qquad k\in\mathbb Z.
\]
When $\beta=a_\rho$, these are $c_0=1$ and $c_k=0$ for $k\ne0$.
For every fixed $k$, \eqref{idx:blowup-limits} gives
$C^{(j)}_{n+k,n}\to c_k(\beta_j)$.
A finite band whose entries tend to zero on every diagonal is
compact; the displayed geometric bounds make the omitted bands
small in operator norm by the Schur test. It follows that
\[
 C^{(j)}-T_{s_{\beta_j}}\text{ is compact},\qquad
 s_\beta(z)=\sum_{k\in\mathbb Z}c_k(\beta)z^k,
 \quad |z|=1.
\]
Here $T_s$ is the Toeplitz matrix $(c_{m-n})_{m,n\ge0}$.
The Fourier series is absolutely convergent.
For different chains, separation of their roots and
\eqref{idx:blowup-limits} give
\[
 |\langle q_{j,n},e_{i,m}\rangle|
 \le C_{i,j}\tau^m2^{-n/2}\qquad(i\ne j).
\]
Each off-diagonal block is Hilbert--Schmidt. Thus $C$ is a compact
perturbation of $\bigoplus_{j=1}^N T_{s_{\beta_j}}$.

It remains to factor the symbol. Put
\[
 \Theta_\rho(t)=\prod_{k\ge0}(1-t\rho^k)(1-\rho^{k+1}/t),
 \qquad Q_\rho=\prod_{k\ge1}(1-\rho^k).
\]
The classical bilateral identity
\[
 \sum_{k\in\mathbb Z}\frac{t^k}{1-x\rho^k}
 =Q_\rho^2\frac{\Theta_\rho(xt)}
                   {\Theta_\rho(x)\Theta_\rho(t)},
 \qquad \rho<|t|<1,
\]
is the specialization $b=\rho a$ of Ramanujan's summation
\cite[Eq.~17.8.2]{DLMF}. One can verify it directly: both sides,
as meromorphic functions of $x$, obey $H(\rho x)=t^{-1}H(x)$
and have the same principal part at $x=1$ and hence at every
$x=\rho^k$. Their difference is holomorphic on $\mathbb C^*$.
Its Laurent coefficients $h_m$ satisfy
$(\rho^m-t^{-1})h_m=0$; the stated annulus makes every coefficient
zero.

For $\beta\ne a_\rho$, substitution of
$x=a_\rho/\beta$ and $t=\tau z$ yields
\begin{equation}
 s_\beta(z)=
 \frac{\mathfrak F_\rho(\beta)Q_\rho^2}
      {\beta\mathfrak F_\rho'(a_\rho)\Theta_\rho(a_\rho/\beta)}
 \frac{\Theta_\rho(a_\rho\tau z/\beta)}
      {\Theta_\rho(\tau z)}.
 \label{idx:theta-symbol}
\end{equation}
The prefactor is nonzero. The admissible disk gives
$1/2<|\beta/a_\rho|<3/2$. Therefore, with
$c=a_\rho\tau/\beta$,
\[
 \rho^2<|c|<1,
 \qquad \rho<\tau<1.
\]
All factors of \eqref{idx:theta-symbol}, and their reciprocals,
are bounded analytic or coanalytic factors on the unit circle,
except possibly
\[
 1-\frac{\rho}{c}z^{-1}
 =1-\frac{\beta}{\sqrt2a_\rho}z^{-1}.
\]
To display the factors, denote the prefactor in
\eqref{idx:theta-symbol} by $\kappa_\beta$ and put
\[
 \begin{aligned}
 A_+(z)&=\prod_{k\ge0}
           \frac{1-c\rho^kz}{1-\tau\rho^kz},\\
 A_-(z)&=\frac{\displaystyle\prod_{k\ge1}
                         (1-\rho^{k+1}/(cz))}
                    {\displaystyle\prod_{k\ge0}
                         (1-\rho^{k+1}/(\tau z))}.
 \end{aligned}
\]
Each parameter multiplying $z$ or $z^{-1}$ in these linear factors
has modulus less than one.
The functions $A_+$ and $A_+^{-1}$ are bounded analytic factors;
$A_-$ and $A_-^{-1}$ are bounded coanalytic factors. The products
converge uniformly on the unit circle, and
\[
 s_\beta(z)=\kappa_\beta A_-(z)
       (1-\delta_\beta z^{-1})A_+(z),\qquad
 \delta_\beta=\frac{\beta}{\sqrt2a_\rho}.
\]
The exact Toeplitz product rule for a coanalytic factor followed
by an analytic factor now gives
\begin{equation}
 T_{s_\beta}=L_\beta(I-\delta_\beta S^*)R_\beta,
 \qquad
 \delta_\beta=\frac{\beta}{\sqrt2a_\rho},
 \label{idx:shift-factorization}
\end{equation}
with the bounded invertible factors
$L_\beta=\kappa_\beta T_{A_-}$ and $R_\beta=T_{A_+}$.
The case $\beta=a_\rho$ was already seen to give the identity.

For $|\delta|<1$, $I-\delta S^*$ is invertible by its Neumann
series. For $|\delta|>1$, it is surjective with one-dimensional
kernel spanned by $(\delta^{-n})_{n\ge0}$ and bounded right inverse
\[
 -\delta^{-1}S(I-\delta^{-1}S)^{-1}.
\]
At $|\delta|=1$, normalized finite segments of
$(\delta^{-n})$ give weakly null approximate kernels for both
$I-\delta S^*$ and its adjoint; their residual norms tend to zero.
The invertible factors in \eqref{idx:shift-factorization} transport
these sequences to the Toeplitz operator and its adjoint, and
compact errors send the resulting weakly null sequences to zero.
For a critical block, take the resulting normalized weakly null
sequence $v_k$ in the chain space, so that $C^*v_k\to0$.
The displayed block form of $A$ gives
\[
 A^*(v_k,0)=(C^*v_k,0)\longrightarrow0.
\]
Thus the adjoint of synthesis is not bounded below. The finite
extra replacements have compact adjoints and do not change this
conclusion.

Compact invariance and additivity of the Fredholm index give
\eqref{idx:critical-circles} and \eqref{idx:index-count} away from
the critical circles. Completeness makes every closed synthesis
range equal to $L^2$; hence Fredholm synthesis is surjective, and
its index is precisely its excess. At a critical circle the
adjoint has no lower bound, so synthesis is not surjective. Its
range is dense by completeness and is therefore not closed.
Finally \eqref{idx:real-threshold} follows from
$\beta=a_\rho+(1-2\rho)\eta$ and $a_\rho\to1/2$.
\end{proof}

\begin{corollary}[Stability of the index under a compact change]
\label{cor:chain-index-compact}
Start with an array in Theorem~\ref{prop:exact-chain-index}.
Allow further admissible replacements on the remaining gaps and
suppose that their canonical relative displacement array has
vanishing tail packing. The resulting synthesis is Fredholm
precisely under \eqref{idx:critical-circles}, and then has index
\eqref{idx:index-count}.
\end{corollary}

\begin{proof}
The extra synthesis correction is compact by
Theorem~\ref{thm:stability}. Compact perturbations preserve both
Fredholmness and the Fredholm index.
\end{proof}

This corollary asserts an index law. The completeness argument above
uses finitely many accumulation points and does not apply to arbitrary
additional changes on infinitely many gaps. Thus its basis and frame
classification is asserted only for the chain arrays and finite extra
replacements in Theorem~\ref{prop:exact-chain-index}.

The limiting matrix belongs to the theory of discrete Cauchy
transforms. Belov, Mengestie, and Seip
\cite[Theorem~1.3 and Section~5]{BelovMengestieSeip2011} give general
invertibility criteria for weighted discrete Cauchy transforms.
After translating the endpoint to zero and inverting, the present chain has
$\gamma_n\asymp\rho^{-n}$ and weights
$v_n=w_n/a_n^2\asymp2^n$; in particular, the additional hypothesis
$v_n=o(\sum_{k<n}v_k)$ in their more explicit Theorem~5.1 does not
hold here. The result above identifies the exact limiting matrix
for the prescribed Cantor geometry and evaluates its symbol,
critical circle, and index. The bilateral identity and the operator
factorization are classical; the stated conclusion concerns this
geometric realization.

\subsection{Geometric Schatten criteria}
\label{subsec:schatten-displacements}

The upper estimate in Theorem~\ref{thm:geometric-schatten} uses
the Carleson truncation method of Rochberg and Semmes
\cite{RochbergSemmes1989}; see also the formulation on spaces of
homogeneous type in \cite[Section~7]{Hytonen2026}.
For lower estimates below exponent one, scale separation and
compression already appear in the paraproduct arguments of
\cite[Proposition~3.10]{PottSmith2004} and
\cite[preprint version, Section~4.2]{LaceyLiWick2026}.
We first construct local tests with enough vanishing moments for
the present rational columns. The key point is uniform
nondegeneracy as both the rational profile and the rescaled
regular measure vary.

\begin{lemma}[A spectral lower bound from local moments]
\label{lem:local-moment-schatten}
Let $(b_g)$ be an admissible array of complex poles. For complex
numbers $\alpha_g,\beta_g$, put
\[
 r_g(x)=\frac{\sqrt{w_g}\,[\alpha_g(x-a_g)+\ell_g\beta_g]}
                 {(x-a_g)(x-b_g)},
 \qquad
 z_g=(|\alpha_g|^2+|\beta_g|^2)^{1/2}.
\]
Suppose a bounded operator $R$ on $L^2(\mu)$ satisfies
$Re_g=r_g$ for every gap. For every $0<p\le1$,
\begin{equation}
 \sum_g z_g^p\le C_p\|R\|_{S^p}^p.
 \label{eq:local-moment-schatten-lower}
\end{equation}
The right side is interpreted as infinity if $R\notin S^p$.
The constants depend only on $p$ and the geometric and admissibility
parameters.
\end{lemma}

\begin{proof}
We may assume $R\in S^p$. Write $\xi_g=u_g$ for the left endpoint
of $g$. Choose an integer $N\ge0$ such that
\begin{equation}
 p(N+1+s/2)>s.
 \label{eq:moment-order-choice}
\end{equation}
Fix $0<\kappa<c_*/8$ and a real function
$\chi\in C_c^\infty((-\kappa,\kappa))$ such that
$0\le\chi\le1$ and $\chi=1$ on $[-\kappa/2,\kappa/2]$.

For $z_g>0$, define $\psi_g=r_g/z_g$. If $z_g=0$, make the same
definition using $\alpha_g=1$, $\beta_g=0$ in the formula for $r_g$.
The normalized functions have the form
\begin{equation}
 \psi_g(\xi_g+\ell_gt)
 =\ell_g^{-s/2}c_g
       \frac{U_g(t-A_g)+V_g}{(t-A_g)(t-B_g)},
 \quad
 |U_g|^2+|V_g|^2=1,
 \label{eq:normalized-local-rational}
\end{equation}
where
\[
 c\le c_g=\sqrt{w_g}\,\ell_g^{s/2-1}\le C,
 \qquad c_*\le A_g\le1-c_*,
 \qquad |B_g-A_g|\le c_0<c_*/2.
\]
We first construct local test functions. On
$J=[-\kappa,\kappa]$, let
\begin{align*}
 \nu_g(S)&=\ell_g^{-s}\mu\bigl((\xi_g+\ell_g S)\cap E\bigr)
 \quad(S\subset J\text{ Borel}),\\
 f_g(t)&=c_g\frac{U_g(t-A_g)+V_g}{(t-A_g)(t-B_g)}.
\end{align*}
Let $P_g$ be the orthogonal projection of $f_g$ onto the
polynomials of degree at most $N$ in $L^2(\chi\,d\nu_g)$.
The following bounds are uniform in $g$:
\begin{equation}
 \|P_g\|_{C^{N+1}(J)}\le C_N,
 \qquad
 \int_J\chi|f_g-P_g|^2\,d\nu_g\ge c_N>0.
 \label{eq:uniform-polynomial-residual}
\end{equation}
Here the $(N+1)$st derivative of $P_g$ is zero.

To verify \eqref{eq:uniform-polynomial-residual}, note that the
measures $\nu_g$ have uniformly bounded mass, satisfy
$\nu_g((t-r,t+r))\le Cr^s$, and obey
\[
 \int_J\chi\,d\nu_g\ge A^{-1}(\kappa/2)^s.
\]
Every weak limit therefore has no atoms and has positive
$\chi$-weighted mass. The moment matrices of the polynomials
$1,t,\ldots,t^N$ are uniformly positive definite. Otherwise,
a convergent sequence of polynomials whose coefficient vectors have
norm one would vanish almost everywhere for the nonzero nonatomic
limit measure $\chi\,d\nu$. A nonzero polynomial has only finitely many real
zeros, which is a contradiction. Hence the coefficients of $P_g$
are uniformly bounded.

The parameters in \eqref{eq:normalized-local-rational} range over a
compact set, and the denominators stay uniformly away from zero on
$J$. If the second bound in \eqref{eq:uniform-polynomial-residual}
failed, a subsequence would yield a polynomial $P$ and a nonzero
rational function
\[
 f(t)=c\frac{u(t-A)+v}{(t-A)(t-B)},\qquad |u|^2+|v|^2=1,
\]
that agree on a set of positive nonatomic measure. The identity
theorem gives
\[
 c[u(t-A)+v]=P(t)(t-A)(t-B).
\]
The left side has degree at most one. The right side has degree at
least two unless $P=0$; that case would force $u=v=0$. This
contradiction proves the lower bound.

Define the smooth function on the real line
\[
 \varphi_g(x)=\ell_g^{-s/2}
 \chi\left(\frac{x-\xi_g}{\ell_g}\right)
 \left[
 f_g\left(\frac{x-\xi_g}{\ell_g}\right)
 -P_g\left(\frac{x-\xi_g}{\ell_g}\right)
 \right],
\]
extended by zero outside $\xi_g+\ell_gJ$. The poles lie outside a
fixed neighborhood of this interval. Projection and
\eqref{eq:uniform-polynomial-residual} imply
\begin{equation}
 \begin{gathered}
 \int_E x^j\varphi_g(x)\,d\mu(x)=0\quad(0\le j\le N),\\
 \langle\psi_g,\varphi_g\rangle
     =\int_J\chi|f_g-P_g|^2\,d\nu_g\ge c_N,\\
 \operatorname{supp}\varphi_g\subset
       [\xi_g-\kappa\ell_g,\xi_g+\kappa\ell_g],\qquad
 \|\varphi_g^{(j)}\|_\infty
       \le C_{N,j}\ell_g^{-s/2-j}\quad(j\ge0).
 \end{gathered}
 \label{eq:local-moment-tests}
\end{equation}
In particular, $\|\varphi_g\|_1\le C_N\ell_g^{s/2}$.

The family $(\varphi_g)$ is Bessel. If $\ell_g\le\ell_h$, Taylor
expansion of $\varphi_h$ at $\xi_g$ and the vanishing moments give
\[
 |\langle\varphi_g,\varphi_h\rangle|
 \le C_N(\ell_g/\ell_h)^{N+1+s/2}.
\]
The pairing is zero unless the support intervals meet. For a fixed
$h$, the number of intersecting intervals at scale
$\ell_g\asymp2^{-j}\ell_h$ is at most $C2^{js}$, by
\eqref{eq:prelim-gapcount}. At scale
$\ell_g\asymp2^j\ell_h$ there are at most $C$ such intervals.
Since $N+1-s/2>0$, the row and column sums of the Gram matrix are
bounded. The Schur test therefore gives
\[
 \left\|\bigl(\langle f,\varphi_g\rangle\bigr)_g\right\|_{\ell^2}
 \le C_N\|f\|_2.
\]

We next record the decay of the mixed pairings. Set
$d_{gh}=|\xi_g-\xi_h|$. The rational formula and gap distances give
\[
 |\psi_h(x)|\le C\frac{\ell_h^{1-s/2}}
                           {\ell_h+|x-\xi_h|},\qquad x\in E,
 \qquad
 \|\psi_h\|_1\le C\ell_h^{s/2}.
\]
The last estimate follows by annular summation, since $s<1$.
For $\ell_g\le\ell_h$, the Taylor segments inside
$[\xi_g-\kappa\ell_g,\xi_g+\kappa\ell_g]$ stay a fixed multiple
of $\ell_h$ away from both poles of $\psi_h$. Differentiating the
rational formula and using the moments in
\eqref{eq:local-moment-tests} gives
\begin{equation}
 |\langle\psi_h,\varphi_g\rangle|
 \le C_N
 \left(\frac{\ell_g}{\ell_h}\right)^{N+1+s/2}
 \left(1+\frac{d_{gh}}{\ell_h}\right)^{-N-2}.
 \label{eq:mixed-fine-moments}
\end{equation}
Indeed, the derivative of order $N+1$ is at most
$C_N\ell_h^{1-s/2}(\ell_h+d_{gh})^{-N-2}$ on those segments;
multiplication by $\ell_g^{N+1}\|\varphi_g\|_1$ gives the result.

For $\ell_g>\ell_h$, the $L^1$ estimate gives
$|\langle\psi_h,\varphi_g\rangle|
 \le C_N(\ell_h/\ell_g)^{s/2}$.
If $d_{gh}>2\ell_g$, Taylor expansion is again valid and improves
this bound to
\[
 C_N\ell_g^{N+1+s/2}\ell_h^{1-s/2}d_{gh}^{-N-2}
 =C_N\left(\frac{\ell_h}{\ell_g}\right)^{1-s/2}
      \left(\frac{\ell_g}{d_{gh}}\right)^{N+2}.
\]
Since $1-s/2>s/2$, both cases give
\begin{equation}
 |\langle\psi_h,\varphi_g\rangle|
 \le C_N
 \left(\frac{\ell_h}{\ell_g}\right)^{s/2}
 \left(1+\frac{d_{gh}}{\ell_g}\right)^{-N-2}.
 \label{eq:mixed-coarse-moments}
\end{equation}

We partition the gaps into finitely many classes. Choose integers
$K\ge2$ and $L\ge2$. Within each band
$2^{-k-1}<\ell_g\le2^{-k}$, list the left endpoints in increasing
order and color them cyclically with $2L+1$ colors. Consecutive
endpoints are separated by more than $2^{-k-1}$, so endpoints of
the same color are more than $L2^{-k}$ apart. The final class of a
gap records this color and $k$ modulo $K$.
Thus two gaps in one class either have length ratio at least
$2^{K-1}$, or lie in the same band with separation greater than
$L\max(\ell_g,\ell_h)$.

For any class $\Gamma$, the gap count and
\eqref{eq:mixed-fine-moments} and~\eqref{eq:mixed-coarse-moments} yield
\begin{equation}
 \sup_{h\in\Gamma}
 \sum_{\substack{g\in\Gamma\\g\ne h}}
       |\langle\psi_h,\varphi_g\rangle|^p
 \le C_{N,p}\left[
   2^{-(K-1)[p(N+1+s/2)-s]}
   +2^{-(K-1)ps/2}
   +L^{s-p(N+2)}\right].
 \label{eq:separated-moment-columns}
\end{equation}
For clarity, at the finer scale $\ell_g\asymp2^{-j}\ell_h$,
the annulus $d_{gh}\lesssim2^m\ell_h$ contains at most
$C2^{(j+m)s}$ endpoints. Its contribution is bounded by
\[
 C_{N,p}2^{-j[p(N+1+s/2)-s]}
          2^{-m[p(N+2)-s]}.
\]
At the coarser scale $\ell_g\asymp2^j\ell_h$, the corresponding
bound is
\[
 C_{N,p}2^{-jps/2}2^{-m[p(N+2)-s]}.
\]
For the same band, the annular sum starts at $2^m\gtrsim L$.
Summing these geometric series proves
\eqref{eq:separated-moment-columns}. The exponents are positive
by \eqref{eq:moment-order-choice}. Choose $K,L$ so large that its
right side is at most $c_N^p/2$.

Fix a finite set $F\subset\Gamma$. Let $J_F:\ell^2(F)\to L^2(\mu)$
be the isometry taking the standard vectors to $e_g$, and let
$V_Ff=(\langle f,\varphi_g\rangle)_{g\in F}$.
The matrix of $V_FRJ_F$ is
\[
 A_F=(z_h\langle\psi_h,\varphi_g\rangle)_{g,h\in F}.
\]
Write $A_F=D_F+E_F$, where $D_F$ is diagonal. Then
\[
 \|D_F\|_{S^p}^p
 =\sum_{g\in F}z_g^p|\langle\psi_g,\varphi_g\rangle|^p
 \ge c_N^p\sum_{g\in F}z_g^p.
\]
The $p$th power of the Schatten quasinorm is subadditive for
$0<p\le1$. Applying this first to the operators of rank one corresponding to the entries of $E_F$
and then to $D_F=A_F-E_F$ gives
\[
 \|E_F\|_{S^p}^p
 \le\sum_{h\in F}z_h^p
       \sum_{\substack{g\in F\\g\ne h}}
          |\langle\psi_h,\varphi_g\rangle|^p
 \le\frac{c_N^p}{2}\sum_{h\in F}z_h^p,
\]
and
\[
 c_N^p\sum_{g\in F}z_g^p
 \le\|V_FRJ_F\|_{S^p}^p+\|E_F\|_{S^p}^p
 \le C_N^p\|R\|_{S^p}^p
       +\frac{c_N^p}{2}\sum_{g\in F}z_g^p.
\]
Absorb the last sum. Let $F$ increase to $\Gamma$, and sum over
the finitely many classes. This proves
\eqref{eq:local-moment-schatten-lower}.
\end{proof}

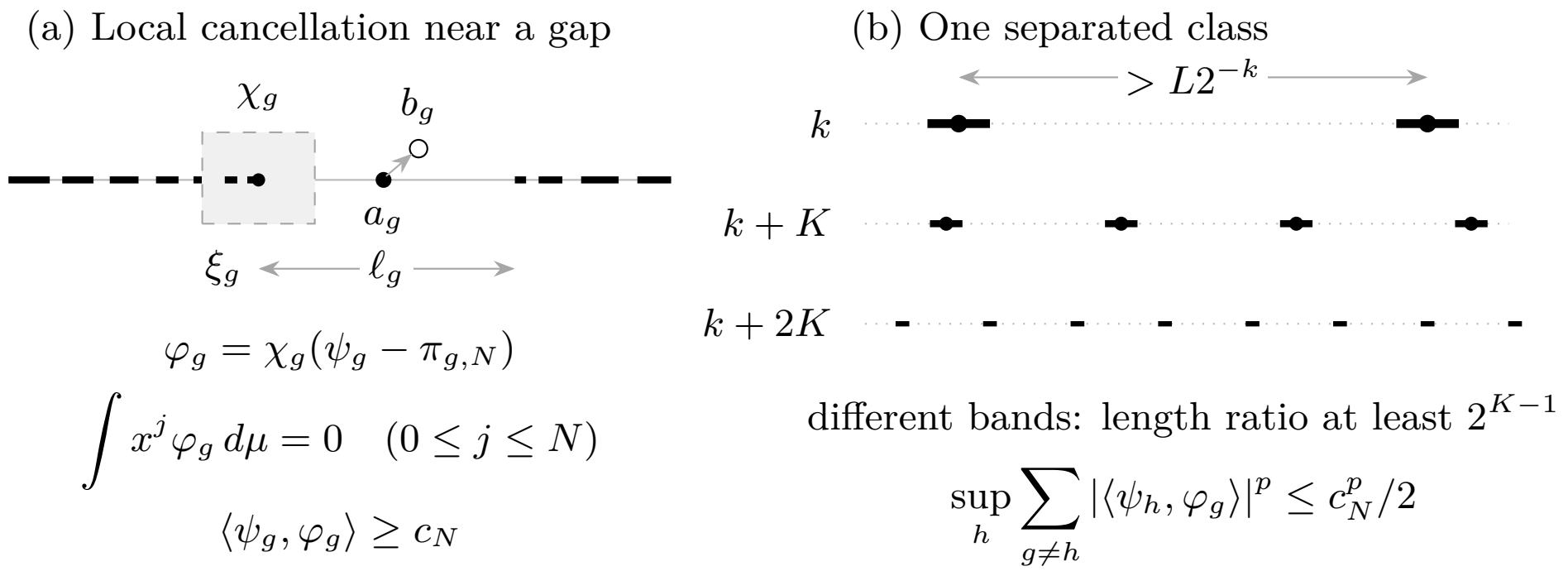
\begin{figure}[htbp]
\centering
\begin{tikzpicture}[x=1cm,y=1cm,font=\footnotesize,>=Stealth]
 % Endpoint localization. Only interval positions are geometric data.
 \node[anchor=west] at (0,3.9) {(a) Local cancellation near a gap};
 \draw[gray!55] (0,2.7)--(5.3,2.7);
 \foreach \a/\b in {0/.33,.43/.68,.81/1.04,1.18/1.36,1.48/1.62,
                    1.73/1.83,1.91/2,4.05/4.14,4.24/4.43,4.58/4.88,5.03/5.3}
   \draw[black,line width=1.6pt] (\a,2.7)--(\b,2.7);
 \fill[gray!12] (1.55,2.35) rectangle (2.45,3.08);
 \draw[dashed,gray!70] (1.55,2.35) rectangle (2.45,3.08);
 \draw[black,line width=1.6pt] (1.73,2.7)--(1.83,2.7) (1.91,2.7)--(2,2.7);
 \fill (2,2.7) circle (1.5pt);
 \fill (3,2.7) circle (1.8pt);
 \fill[white] (3.28,2.95) circle (2pt);
 \draw (3.28,2.95) circle (2pt);
 \draw[->,gray!65] (3.02,2.73)--(3.23,2.91);
 \node[anchor=north] at (1.72,2.28) {$\xi_g$};
 \node[anchor=north] at (3,2.62) {$a_g$};
 \node[anchor=south] at (3.28,3.01) {$b_g$};
 \node[anchor=south] at (2,3.12) {$\chi_g$};
 \draw[<->,gray!70] (2,1.99)--(4.05,1.99);
 \node[fill=white,inner sep=2pt] at (3.025,1.99) {$\ell_g$};
 \node[anchor=north,align=center] at (2.65,1.63)
 {$\varphi_g=\chi_g(\psi_g-\pi_{g,N})$\\[4pt]
  $\displaystyle\int x^j\varphi_g\,d\mu=0\quad(0\le j\le N)$\\[4pt]
  $\langle\psi_g,\varphi_g\rangle\ge c_N$};

 \begin{scope}[xshift=.6cm]
 % Separated classes, schematic finite part.
 \node[anchor=west] at (6,3.9) {(b) One separated class};
 \foreach \y/\lab in {3.15/{k},2.35/{k+K},1.55/{k+2K}}{
   \draw[dotted,gray!55] (6.25,\y)--(11.4,\y);
   \node[anchor=east] at (6.12,\y) {$\lab$};
 }
 \draw[black,line width=2pt] (6.75,3.15)--(7.25,3.15)
                              (10.5,3.15)--(11,3.15);
 \fill (7,3.15) circle (2pt);
 \fill (10.75,3.15) circle (2pt);
 \draw[<->,gray!70] (7,3.52)--(10.75,3.52);
 \node[fill=white,inner sep=1.5pt] at (8.875,3.52) {$>L2^{-k}$};
 \foreach \x in {6.9,8.3,9.7,11.1}{
   \draw[black,line width=1.6pt] (\x-.13,2.35)--(\x+.13,2.35);
   \fill (\x,2.35) circle (1.6pt);
 }
 \foreach \x in {6.55,7.25,7.95,8.65,9.35,10.05,10.75,11.45}{
   \draw[black,line width=1.2pt] (\x-.055,1.55)--(\x+.055,1.55);
 }
 \node[anchor=north,align=center] at (8.8,1.12)
 {different bands: length ratio at least $2^{K-1}$\\[5pt]
  $\displaystyle
   \sup_h\sum_{g\ne h}|\langle\psi_h,\varphi_g\rangle|^p
       \le c_N^p/2$};
\end{scope}
\end{tikzpicture}
\caption{Local moment tests and the separation used in the spectral lower
bound. The shaded interval supports the smooth cutoff
$\chi_g(x)=\chi((x-\xi_g)/\ell_g)$ near one endpoint;
the polynomial $\pi_{g,N}$ is the weighted projection expressed in the
original variable, and its subtraction imposes the vanishing moments. The right panel shows a schematic part of one class after
separating both scales and positions. Segment lengths and the displayed
finite pieces of the support are illustrative; the proof uses the gap
counting estimate on the full regular set.}
\label{fig:moment-spectral-separation}
\end{figure}

Figure~\ref{fig:moment-spectral-separation} describes the local projection
and the separation used in this proof.

\begin{proof}[Proof of Theorem~\ref{thm:geometric-schatten}]
For $\zeta_g>0$, put $\psi_g=r_g/\zeta_g$; if $\zeta_g=0$, put
$\psi_g=e_g$. We first record the bounds
\begin{equation}
 c\le\|\psi_g\|_2\le C,
 \qquad
 |\psi_g(x)|\le
 C\frac{\ell_g^{1-s/2}}{\ell_g+|x-\xi_g|}.
 \label{st:normalized-displacement}
\end{equation}
The upper pointwise bound follows from
\eqref{st:pole-distance} and $w_g\asymp\ell_g^{2-s}$.
Annular summation then gives the upper $L^2$ bound.

For completeness, the lower bound controls possible cancellation
between the two numerator terms. Write $g=(u_g,v_g)$, fix
$0<\varepsilon<1/4$, and put
\[
 J_-=B_E(u_g,\varepsilon\ell_g),\qquad
 J_+=B_E(v_g,\varepsilon\ell_g),\qquad
 t=(x-a_g)/\ell_g.
\]
Both balls have mass comparable to $\ell_g^s$.
The normalized denominators are bounded above there, so
\[
 \|r_g\|_2^2
 \ge c\ell_g^{-s}
       \int_{J_-\cup J_+}|\alpha_gt+\beta_g|^2\,d\mu.
\]
Let $t_-,t_+$ be the averages of $t$ on $J_-,J_+$.
These numbers are uniformly bounded and satisfy $t_+-t_-\ge1$.
Jensen's inequality gives
\[
 \begin{aligned}
 \ell_g^{-s}\int_{J_-\cup J_+}
       |\alpha_gt+\beta_g|^2\,d\mu
 &\ge c\bigl(|\alpha_gt_-+\beta_g|^2
             +|\alpha_gt_++\beta_g|^2\bigr)\\
 &\ge c\bigl(|\alpha_g|^2+|\beta_g|^2\bigr).
 \end{aligned}
\]
The last inequality follows by solving the two equations with
coefficient matrix having rows $(t_-,1)$ and $(t_+,1)$.
This proves the lower bound in \eqref{st:normalized-displacement}.

Use the partitions $\mathcal Q_k$ of Section~\ref{sec:prelim}.
We regard $(k,Q)$ as a tree vertex, so equal underlying cells at
different levels remain distinct. Write $R\preceq Q$ if $R$ is
a descendant of $Q$, including $Q$ itself.
Assign $g$ to the unique level $k$ with
$2^{-k-1}<\ell_g\le2^{-k}$, and then to the cell
$Q(g)\in\mathcal Q_k$ containing its left endpoint $\xi_g$.
By \eqref{var:cell-geometry},
\[
 \mu(Q(g))\asymp\ell_g^s,
 \qquad \operatorname{diam}Q(g)\le C\ell_g.
\]
The local gap count shows that at most $M_0$ gaps are assigned to
any vertex, where $M_0$ depends only on the regularity data.

Annular summation in \eqref{st:normalized-displacement}, now
centered at an arbitrary point of $Q(g)$, gives
\begin{equation}
 |\langle f,\psi_g\rangle|
 \le C\mu(Q(g))^{1/2}\inf_{Q(g)}M_\mu f.
 \label{st:schatten-maximal-pairing}
\end{equation}
Indeed, the annulus of radius $2^j\ell_g$ contributes at most
$C\ell_g^{s/2}2^{-j(1-s)}M_\mu f(x)$ for every $x\in Q(g)$.
The sum converges because $s<1$.

For a finitely supported scalar array $z=(z_g)$, define
\[
 D_zf=\sum_g z_g\langle f,e_g\rangle\psi_g,
 \qquad
 b_z(Q)=\sum_{Q(g)=Q}|z_g|^2,
 \qquad
 \Sigma_z(Q)^2=
 \frac1{\mu(Q)}\sum_{R\preceq Q}b_z(R)\mu(R).
\]
The tree Carleson embedding and
\eqref{st:schatten-maximal-pairing} imply
\begin{equation}
 \|D_z\|_{2\to2}\le C\sup_Q\Sigma_z(Q).
 \label{st:tree-synthesis-bound}
\end{equation}
For completeness, let $L=\sup_Q\Sigma_z(Q)^2$.
For any nonnegative function $F$, select the maximal cells among
those for which $\inf_QF>t$. They are disjoint, and each contains
all its selected descendants. Thus
\[
 \sum_{\inf_QF>t}b_z(Q)\mu(Q)
 \le L\,\mu\{F>t\}.
\]
Integration of the distribution estimate with $F=M_\mu f$ gives
\[
 \sum_g |z_g|^2|\langle f,\psi_g\rangle|^2
 \le CL\|M_\mu f\|_2^2
 \le CL\|f\|_2^2.
\]
Since $(e_g)$ is orthonormal, this bounds $D_z^*$ and proves
\eqref{st:tree-synthesis-bound}.

Fix $t>0$, and retain the coordinates assigned to
\[
 \mathcal F_t=\{Q:\Sigma_z(Q)>t\}.
\]
Their contribution $R_t$ has rank at most $M_0\#\mathcal F_t$.
For the remaining coefficients, the squared tree packing norm is
at most $t^2$. To see this, take any vertex $P$.
If $P\notin\mathcal F_t$, its original packing sum is at most
$t^2\mu(P)$. Otherwise, group the remaining coefficients below
$P$ according to the first vertex $S$ on their
path that lies outside $\mathcal F_t$.
These cells at the first exit are disjoint, and each contributes at most
$t^2\mu(S)$. Their total mass is at most $\mu(P)$.
Consequently,
\begin{equation}
 \operatorname{rank}R_t\le M_0\#\{Q:\Sigma_z(Q)>t\},
 \qquad \|D_z-R_t\|_{2\to2}\le Ct.
 \label{st:adaptive-spectral-truncation}
\end{equation}
This gives the required approximation by operators of finite rank.

We next estimate the distribution of the profile used in this truncation.
For a nonnegative sequence $b$ on the vertices, put
\[
 Hb(Q)=\frac1{\mu(Q)}\sum_{R\preceq Q}b(R)\mu(R)
       =\sum_{j\ge0}H_jb(Q),
\]
where $H_j$ retains descendants exactly $j$ levels below $Q$.
The cell geometry gives
\[
 \frac{\mu(R)}{\mu(Q)}\le C2^{-js}
 \quad\text{if }R\preceq Q
 \text{ is }j\text{ levels below }Q.
\]
For $u\ge1$, Jensen's inequality on the descendant partition yields
\[
 \sum_Q|H_jb(Q)|^u
 \le\sum_Q\sum_{\substack{R\preceq Q\\
                 \operatorname{level}R-\operatorname{level}Q=j}}
       \frac{\mu(R)}{\mu(Q)}|b(R)|^u
 \le C2^{-js}\sum_R|b(R)|^u.
\]
For $0<u<1$, subadditivity gives instead
\[
 \sum_Q|H_jb(Q)|^u
 \le C^u2^{-jsu}\sum_R|b(R)|^u.
\]
Every vertex has at most one ancestor at a prescribed level,
which justifies the last summation in both cases.
Summing the geometric series, using the triangle inequality when
$u\ge1$ and its $u$-power version when $u<1$, proves
\begin{equation}
 \sum_Q|Hb(Q)|^u\le C_u\sum_Q|b(Q)|^u,
 \qquad 0<u<\infty.
 \label{st:tree-hardy-schatten}
\end{equation}

Let $n(t;D_z)$ count the singular values of $D_z$ greater than $t$.
The approximation estimate \eqref{st:adaptive-spectral-truncation}
gives
\[
 n(Ct;D_z)\le M_0\#\{Q:\Sigma_z(Q)>t\}.
\]
Integrating this estimate and using \eqref{st:tree-hardy-schatten}
gives, for every $0<r<\infty$,
\begin{align*}
 \|D_z\|_{S^r}^r
 &\le C_r\sum_Q\Sigma_z(Q)^r\\
 &\le C_r\sum_Qb_z(Q)^{r/2}
 \le C_r\sum_g|z_g|^r.
\end{align*}
The last step uses the uniform multiplicity $M_0$.
All estimates are uniform over finite coordinate sets and over the
admissible complex poles.
Applying this bound to $z_g=\zeta_g$ and then to the tails proves
the upper estimate and convergence in $S^r$ whenever
$\sum_g\zeta_g^r<\infty$.

Suppose conversely that $R\in S^r$ for $1<r<\infty$, and put
$r'=r/(r-1)$. For every finitely supported $z$, orthonormality gives
\[
 \operatorname{tr}(D_z^*R)
 =\sum_g\overline{z_g}\zeta_g\|\psi_g\|_2^2.
\]
H\"older's inequality for Schatten classes and the upper estimate
at exponent $r'$ imply
\[
 \left|\sum_g\overline{z_g}\zeta_g\|\psi_g\|_2^2\right|
 \le\|R\|_{S^r}\|D_z\|_{S^{r'}}
 \le C_{r'}\|R\|_{S^r}\|z\|_{\ell^{r'}}.
\]
Duality over finite coordinate sets and
\eqref{st:normalized-displacement} yield
$\|\zeta\|_{\ell^r}\le C_r\|R\|_{S^r}$.

For $0<r\le1$, apply Lemma~\ref{lem:local-moment-schatten} to the
same numerator coefficients and the same pole array. Its lower bound
combines with the upper estimate already proved. This proves
\eqref{st:rational-schatten-criterion} at every positive exponent.
\end{proof}

\begin{corollary}[Schatten criterion for pole displacements]
\label{prop:schatten-displacements}
For every admissible pole array, put
$\eta_g=|b_g-a_g|/\ell_g$ and define
$T1=0$, $Te_g=q_g-e_g$. Then
\[
 \begin{gathered}
 T\in S^r\quad\Longleftrightarrow\quad
 (\eta_g)\in\ell^r(\mathcal G),\\
 \|T\|_{S^r}\asymp_r
       \left(\sum_g\eta_g^r\right)^{1/r},
 \qquad 0<r<\infty.
 \end{gathered}
\]
\end{corollary}

\begin{proof}
Apply Theorem~\ref{thm:geometric-schatten} with
$\alpha_g=0$ and $\beta_g=(b_g-a_g)/\ell_g$.
\end{proof}

The spectral truncation \eqref{st:adaptive-spectral-truncation}
selects gaps through the accumulated coefficient mass below their
cells. The criteria are independent of an ordering of the gaps and
allow different coefficient scales in different parts of the support.
The local moment tests in Lemma~\ref{lem:local-moment-schatten}
supply the lower estimate for $0<r\le1$, including the trace class.
Thus both the criterion for rational columns and its displacement
specialization hold at every positive exponent.

\begin{figure}[!htbp]
\centering
\begin{tikzpicture}[
 x=1cm,y=1cm,font=\footnotesize,
 selected/.style={circle,draw=DiagramPrimary,fill=DiagramPrimary,
                  text=white,inner sep=2.5pt},
 firstexit/.style={circle,draw=DiagramPrimary,fill=white,
                   double,double distance=.5pt,inner sep=2.5pt},
 ordinary/.style={circle,draw=DiagramSecondary,fill=white,inner sep=2.5pt},
 branch/.style={draw=DiagramSecondary,line width=.45pt},
 region/.style={draw=DiagramSecondary,fill=DiagramSecondary!9,
                line width=.35pt}
]
 \node at (6,5.45) {$\mathcal F_t=\{Q:\Sigma_z(Q)>t\}$};

 % Disjoint subtrees rooted at the first exits from the chosen root.
 \path[region] (1.7,2.38)--(.7,.65)--(2.7,.65)--cycle;
 \path[region] (3.7,1.34)--(3.08,.65)--(4.32,.65)--cycle;
 \path[region] (5.5,1.34)--(4.88,.65)--(6.12,.65)--cycle;
 \path[region] (8.8,3.42)--(6.95,.65)--(10.65,.65)--cycle;

 \coordinate (p) at (6,4.75);
 \coordinate (a) at (3.2,3.65);
 \coordinate (b) at (8.8,3.65);
 \coordinate (c) at (1.7,2.60);
 \coordinate (d) at (4.6,2.60);
 \coordinate (e) at (3.7,1.56);
 \coordinate (f) at (5.5,1.56);
 \coordinate (g) at (8.0,2.60);
 \coordinate (h) at (9.6,2.60);
 \draw[branch] (p)--(a) (p)--(b)
               (a)--(c) (a)--(d) (d)--(e) (d)--(f)
               (b)--(g) (b)--(h);
 \draw[branch,dashed] (g)--(7.62,1.96) (g)--(8.38,1.96)
                      (h)--(9.22,1.96) (h)--(9.98,1.96);

 \node[selected] at (p) {$P$};
 \node[selected] at (a) {$A$};
 \node[firstexit] at (b) {$B$};
 \node[firstexit] at (c) {$C$};
 \node[selected] at (d) {$D$};
 \node[firstexit] at (e) {$E$};
 \node[firstexit] at (f) {$F$};
 \node[selected] at (g) {$G$};
 \node[ordinary] at (h) {$H$};

 \node[anchor=west] at (9.22,3.65) {$\Sigma_z(B)\le t$};
 \node at (8.8,1.36) {$\Sigma_z(G)>t$};
 \node[anchor=north] at (1.7,.57) {$C$-subtree};
 \node[anchor=north] at (3.7,.57) {$E$-subtree};
 \node[anchor=north] at (5.5,.57) {$F$-subtree};
 \node[anchor=north] at (8.8,.57) {$B$-subtree};

 \filldraw[draw=DiagramPrimary,fill=DiagramPrimary]
       (.75,-.27) circle (.085);
 \node[anchor=west] at (.96,-.27)
       {selected vertex: $\Sigma_z(Q)>t$};
 \draw[draw=DiagramPrimary,fill=white,double,double distance=.5pt]
       (6.65,-.27) circle (.085);
 \node[anchor=west] at (6.88,-.27) {first-exit good cell};
 \node at (5.9,-.88)
 {$\operatorname{rank}R_t\le M_0\#\mathcal F_t,
   \qquad\|D_z-R_t\|_{2\to2}\le Ct$};
\end{tikzpicture}
\caption{Adaptive spectral truncation. The shaded subtrees have disjoint
 roots $B,C,E,F$ at the first exit below $P$. The bad vertices
 $\mathcal F_t$ need not form a stopping tree or a sparse family:
 the good cell $B$ contains the bad descendant $G$.
 First exits partition the remaining coefficients below $P$ as in the proof;
 coefficients assigned to bad vertices are still removed.}
\label{fig:carleson-spectral-stopping}
\end{figure}
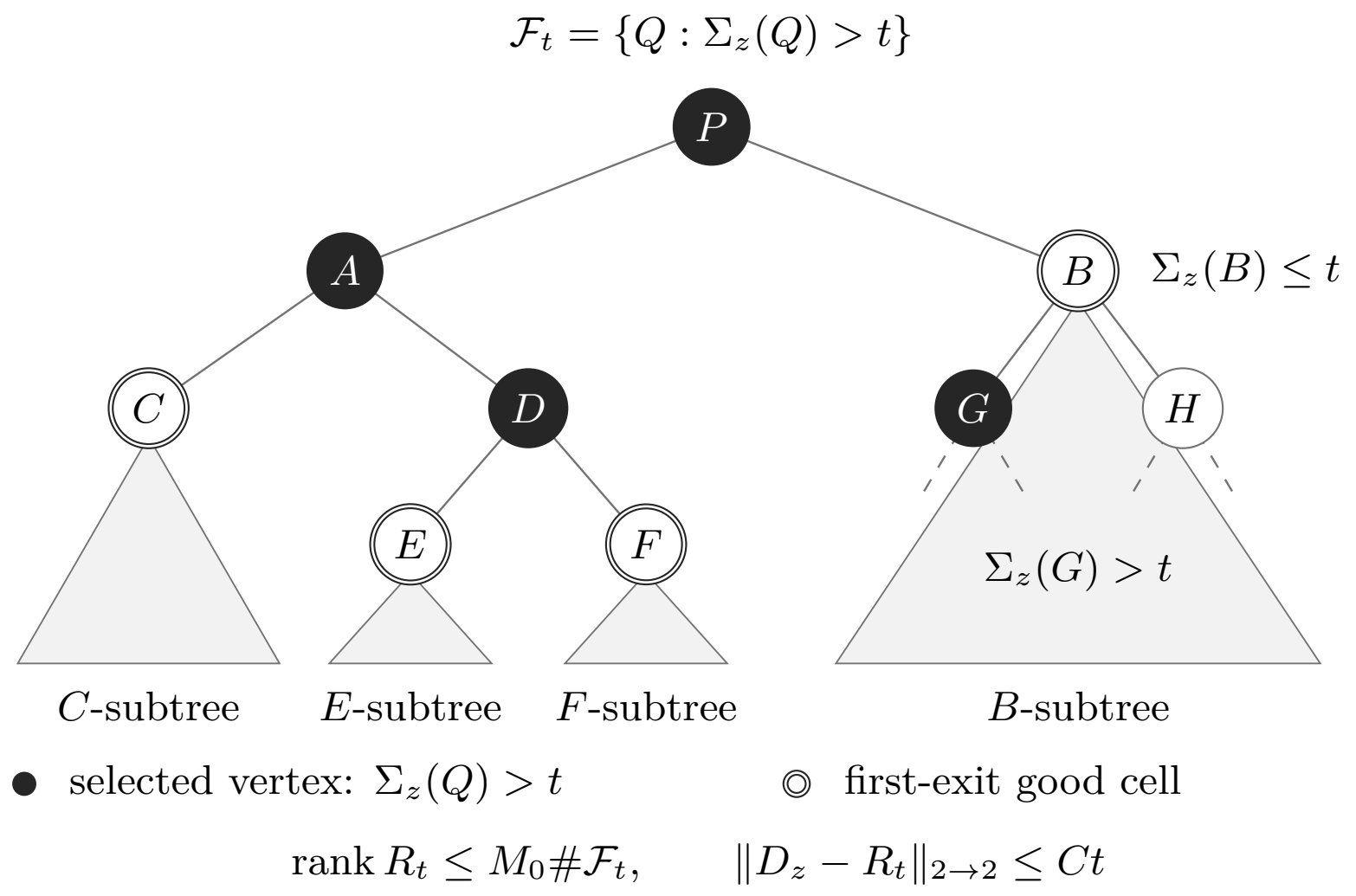

Figure~\ref{fig:carleson-spectral-stopping} shows the decomposition at the first exit. Figure~\ref{fig:moment-spectral-separation} shows the
local moment construction and the separation used for the lower bound.

\FloatBarrier
\subsection{Sparse control on the gap filtration and weighted stability}
\label{subsec:weighted-gap-consequences}

The gap geometry also permits a weighted version of the packing law.
The sparse domination theorem used here is classical
\cite[Theorem~6.1]{Lorist2021}. We record its consequences for the
prescribed Cauchy family, including the lower bound for the operator
norm and the bounds for the essential norm.

For $1<p<\infty$, a positive weight $v$ belongs to $A_p(\mu)$ if
\[
 [v]_{A_p(\mu)}
 =\sup_{B=B_E(x,r)}
 \left(\frac1{\mu(B)}\int_Bv\,d\mu\right)
 \left(\frac1{\mu(B)}\int_Bv^{-1/(p-1)}\,d\mu\right)^{p-1}
 <\infty.
\]
Write $X_{p,v}=L^p(v\,d\mu)$ and
$\kappa_p=\max\{1,(p-1)^{-1}\}$.
Let $\mathcal Q_{\mathrm{gap}}=\bigcup_{k\ge0}\mathcal Q_k$, with repeated cells
identified. A subfamily $\mathcal S\subset\mathcal Q_{\mathrm{gap}}$ is
$\vartheta$-sparse if there are pairwise disjoint measurable sets
$D_Q\subset Q$ with $\mu(D_Q)\ge\vartheta\mu(Q)$.
Put
\[
 \mathcal A_{\mathcal S}f
 =\sum_{Q\in\mathcal S}
       \left(\frac1{\mu(Q)}\int_Q|f|\,d\mu\right)\mathbf1_Q.
\]

\begin{proposition}[Weighted packing and compactness]
\label{prop:weighted-packing}
There is $\vartheta>0$, depending only on the regularity data,
with the following properties.
For every bounded $f$, each canonical multiplier $T_\beta$,
$\|\beta\|_\infty\le1$, admits a $\vartheta$-sparse family of gap
cells such that
\[
 |T_\beta f|\le C\mathcal A_{\mathcal S}f
 \quad\text{almost everywhere}.
\]
If $\mathcal C(\lambda)<\infty$, the pole correction admits such
a family with
\[
 |Tf|\le C\mathcal C(\lambda)^{1/2}\mathcal A_{\mathcal S}f.
\]
The sparse family may depend on the operator and $f$.

For every $v\in A_p(\mu)$, the canonical family is an unconditional
basis of $X_{p,v}$, with coefficient functionals given by integration
against the original $e_g$ with respect to $\mu$. Its multipliers satisfy
\[
 \|T_\beta\|_{X_{p,v}\to X_{p,v}}
 \le C_p[v]_{A_p(\mu)}^{\kappa_p}\|\beta\|_\infty.
\]
The pole correction, initially defined on finite canonical expansions,
extends boundedly to $X_{p,v}$ if and only if
$\mathcal C(\lambda)<\infty$. In this case
\[
 \begin{aligned}
 c_p[v]_{A_p(\mu)}^{-\kappa_p}\mathcal C(\lambda)^{1/2}
 &\le\|T\|_{X_{p,v}\to X_{p,v}}
 \le C_p[v]_{A_p(\mu)}^{\kappa_p}\mathcal C(\lambda)^{1/2},\\
 c_p[v]_{A_p(\mu)}^{-2\kappa_p}\gamma(\lambda)^{1/2}
 &\le\|T\|_{\mathrm{ess},X_{p,v}\to X_{p,v}}
 \le C_p[v]_{A_p(\mu)}^{\kappa_p}\gamma(\lambda)^{1/2}.
 \end{aligned}
\]
In particular, compactness is equivalent to $\gamma(\lambda)=0$,
independently of the chosen $p$ and $A_p$ weight.
Moreover, for a constant $C_p$ depending only on $p$ and the fixed
geometric data,
\[
 C_p[v]_{A_p(\mu)}^{\kappa_p}\gamma(\lambda)^{1/2}<1
\]
implies that $I+T$ is invertible on $X_{p,v}$ and that
$\{1\}\cup(q_g)$ is an unconditional basis there.
\end{proposition}

\begin{proof}
The multiplier kernels and the kernels of the pole corrections have
H\"older exponent $1-s>0$ in both variables.
The $L^2$ norm and the kernel constants are bounded by $C$ for
$\|\beta\|_\infty\le1$, and by $C\mathcal C(\lambda)^{1/2}$
for $T$. The same statements hold for every subfamily of gaps,
with its own packing norm.

We also identify the full correction with its kernel off the
diagonal. Let $f\in L^2(\mu)$ be supported in a closed set $F\subset E$.
Since $\mu(E)=1$, one also has $f\in L^1(\mu)$.
For $x\notin F$, the absolute sums used to prove
\eqref{st:kernel-size} dominate the finite integral kernels by
$C\mathcal C(\lambda)^{1/2}\operatorname{dist}(x,F)^{-s}|f(y)|$
on $F$. Dominated convergence identifies their integral limit.
The finite corrections converge strongly in $L^2$; an almost
everywhere convergent subsequence therefore identifies this limit
with $Tf$ almost everywhere outside $F$.

Regularity converts the kernel bounds
to the Dini kernel hypotheses of
\cite[Theorem~6.1]{Lorist2021}; the modulus is a constant times
$t^{1-s}$. That theorem gives the sparse bounds and the weighted
upper estimates.

The sparse cubes can be replaced by cells of $\mathcal Q_{\mathrm{gap}}$.
Let $Q$ be one of the metric cubes in the cited theorem, with
$0<d=\operatorname{diam}Q<1$. Choose a dyadic scale $r$ with
$d<r\le2d$. No complementary gap of length at least $r$ separates
two points of $Q$, so $Q$ lies in one cell $R(Q)$ at this scale.
The inner ball property of metric cubes and
\eqref{var:cell-geometry} give
\[
 Q\subset R(Q),\qquad \mu(R(Q))\le C\mu(Q).
\]
For cubes with diameter comparable to one, use $R(Q)=E$.
If $D_Q$ are the disjoint major subsets for the original sparse
family, then
\[
 \mu(D_Q)\ge\vartheta_0\mu(Q)
             \ge c\vartheta_0\mu(R(Q)).
\]
For each fixed $R$, these disjoint subsets lie in $R$, so the
number of cubes mapping to $R$ is uniformly bounded.
Retain one $D_Q$ for each distinct $R$. These sets show that the
resulting family of gap cells is sparse, while
\[
 \frac1{\mu(Q)}\int_Q|f|\,d\mu
 \le C\frac1{\mu(R(Q))}\int_{R(Q)}|f|\,d\mu
\]
and the bounded multiplicity give the stated pointwise domination.

The coefficient functionals are continuous on $X_{p,v}$ because
each $e_g$ is bounded and
$\sigma=v^{-1/(p-1)}$ is integrable.
Proposition~\ref{thm:approximation-norming} approximates every member
of $C(E)$ uniformly by finite canonical expansions.
Since $v\,d\mu$ is a finite Borel measure on the compact set $E$,
$C(E)$ is dense in $X_{p,v}$. The canonical span is therefore
dense in $X_{p,v}$. The uniform bound for finite multipliers proves
unconditional convergence in every enumeration, by applying the
bound to the difference from a finite expansion.

For the weighted lower estimate, use
\[
 \sigma\in A_{p'}(\mu),\qquad
 [\sigma]_{A_{p'}(\mu)}
 =[v]_{A_p(\mu)}^{1/(p-1)},\qquad
 \frac{\kappa_{p'}}{p-1}=\kappa_p.
\]
Interpolation with change of measure \cite{SteinWeiss1958} gives,
for every bounded operator $S$ on both spaces,
\[
 \|S\|_{L^2(\mu)\to L^2(\mu)}^2
 \le
 \|S\|_{X_{p,v}\to X_{p,v}}\,
 \|S\|_{X_{p',\sigma}\to X_{p',\sigma}}.
\]
Indeed, interpolation at parameter $1/2$ has exponent $2$ and
weight $v^{1/p}\sigma^{1/p'}=1$.
Apply this to $S=T$ and combine the unweighted lower packing
estimate with the weighted upper estimate on $X_{p',\sigma}$.
If $\mathcal C(\lambda)>0$, cancellation of its square root gives
the asserted lower bound. If $\mathcal C(\lambda)=0$, then $T=0$.
The same reasoning applies to every finite or infinite subfamily
and, in particular, to $T_r=T(I-P_r)$.

Conversely, suppose the initially defined $T$ extends boundedly
to $X_{p,v}$. For a finite gap family $\mathcal F$, write
$T_{\mathcal F}=TP_{\mathcal F}$. The lower estimate for finite families
and the multiplier bound give
\[
 \mathcal C(\lambda\mathbf1_{\mathcal F})^{1/2}
 \le C_p[v]_{A_p(\mu)}^{2\kappa_p}
                 \|T\|_{X_{p,v}\to X_{p,v}}.
\]
Taking the supremum over finite $\mathcal F$ proves
$\mathcal C(\lambda)<\infty$.

The projections $P_r$ converge strongly on both $X_{p,v}$ and
$X_{p',\sigma}$ and have norm at most
$C_p[v]_{A_p(\mu)}^{\kappa_p}$ on the first space.
Under the integral pairing with respect to $\mu$, their adjoints
on the second space are again $P_r$.
For any compact $K$ on $X_{p,v}$,
\[
 \|K(I-P_r)\|_{X_{p,v}\to X_{p,v}}\longrightarrow0.
\]
This follows from the uniform convergence of the adjoints
on the compact closure of the image of the unit ball under $K^*$.
The operators $TP_r$ of finite rank give the upper bound for the
essential norm. For the reverse estimate,
\[
 \begin{aligned}
 C_p[v]_{A_p(\mu)}^{\kappa_p}\|T-K\|_{X_{p,v}\to X_{p,v}}
 &\ge\liminf_{r\downarrow0}
       \|(T-K)(I-P_r)\|_{X_{p,v}\to X_{p,v}}\\
 &\ge c_p[v]_{A_p(\mu)}^{-\kappa_p}
                         \gamma(\lambda)^{1/2}.
 \end{aligned}
\]
Take the infimum over compact $K$.

Finally, choose $r$ so small that
\[
 \zeta=C_p[v]_{A_p(\mu)}^{\kappa_p}
                  \mathcal C_r(\lambda)^{1/2}<1.
\]
The same upper bound $\zeta$ holds for every subfamily of the tail.
The Neumann series makes $I+T_r$ invertible.
The finite replacement argument of
Section~\ref{subsec:tail-stability} applies with the dual space
$X_{p',\sigma}$ in place of $L^{p'}(\mu)$.
To make this substitution explicit, write
$A_{\mathcal F}=I+T_{\mathcal F}$. The original partial fraction identity represents the functional
\[
 \Lambda_{\mathcal F}(f)
 =\int_E A_{\mathcal F}^{-1}f\,d\mu
\]
by the finite product $R_{\mathcal F}$. Hence
\[
 \|R_{\mathcal F}\|_{X_{p',\sigma}}
 \le\frac{\|1\|_{X_{p',\sigma}}}{1-\zeta}.
\]
Reflexivity gives a weakly convergent subsequence.
The test functions $1$ and $(x-z)^{-1}$ belong to $X_{p,v}$
for $z\notin E$, so the same identities for the Cauchy transform pass
to its limit. Multiplication by $(x-a_j)^{-1}$ is bounded on
$X_{p',\sigma}$ for every remaining fixed gap.
Thus the coordinate functionals are bounded, and their finite
replacement matrix is the same invertible Cauchy matrix.
This proves invertibility of $I+T$ and the final basis assertion.
\end{proof}

The powers of the weight characteristic give quantitative bounds;
their optimality for this particular Cauchy family is not asserted.
Sparse domination and the weighted upper bounds are established
results. The additional conclusion is that the explicit
displacement packing and its vanishing tail still govern the
weighted operator and essential norms.

\section{Curve transport}
\label{subsec:curve-proof}
\subsection{Local oscillation and square Dini compactness}
\label{subsec:curve-square-dini}

The Dini condition in Theorem~\ref{thm:curve-transport} gives
estimates at both endpoints. On the interior $L^p$ scale it can be
weakened to a square Dini condition, even when the deformation changes
the support. We first give a local sufficient condition that records
where the derivative oscillates.

Square Dini conditions also occur in estimates for nonlocal operators.
Dong, Jung, and Kim~\cite{DongJungKim2026} prove a sharp such threshold
for operators with variable coefficients on Euclidean Sobolev spaces.
Here the gap geometry and canonical orthogonality give compactness
of the prescribed Cauchy transport.

Let $\phi\in C^1([0,1];\C)$ satisfy
\[
 |\phi(x)-\phi(y)|\ge m_0|x-y|\qquad(x,y\in[0,1]),
 \qquad L=\|\phi'\|_\infty.
\]
The parameter poles $b_g$ are real and satisfy the admissibility
condition of Theorem~\ref{thm:curve-transport}. Write $c_b=c_*-c_0$,
choose $C_*=8+4/c_b$, and define
\[
 \begin{split}
 \omega_g(t)&=\sup\{|\phi'(v)-\phi'(b_g)|:
                   v\in[0,1],\ |v-b_g|\le t\},\\
 \mathfrak a_g&=\sum_{j=0}^{\infty}2^{-j(1-s)}
            \omega_g\bigl(\min(1,C_*2^j\ell_g)\bigr),
 \qquad \Lambda_g=\mathfrak a_g^2\ell_g^s.
 \end{split}
\]
Use the packing norm $\mathcal C(\Lambda)$ and the tail norm
$\mathcal C_r(\Lambda)$ defined as in Theorem~\ref{thm:stability}.
The correction columns are
\[
 h_g(x)=\sqrt{w_g}\left(
      \frac{\phi'(b_g)}{\phi(x)-\phi(b_g)}-\frac1{x-b_g}\right),
 \qquad H1=0,\quad He_g=h_g.
\]

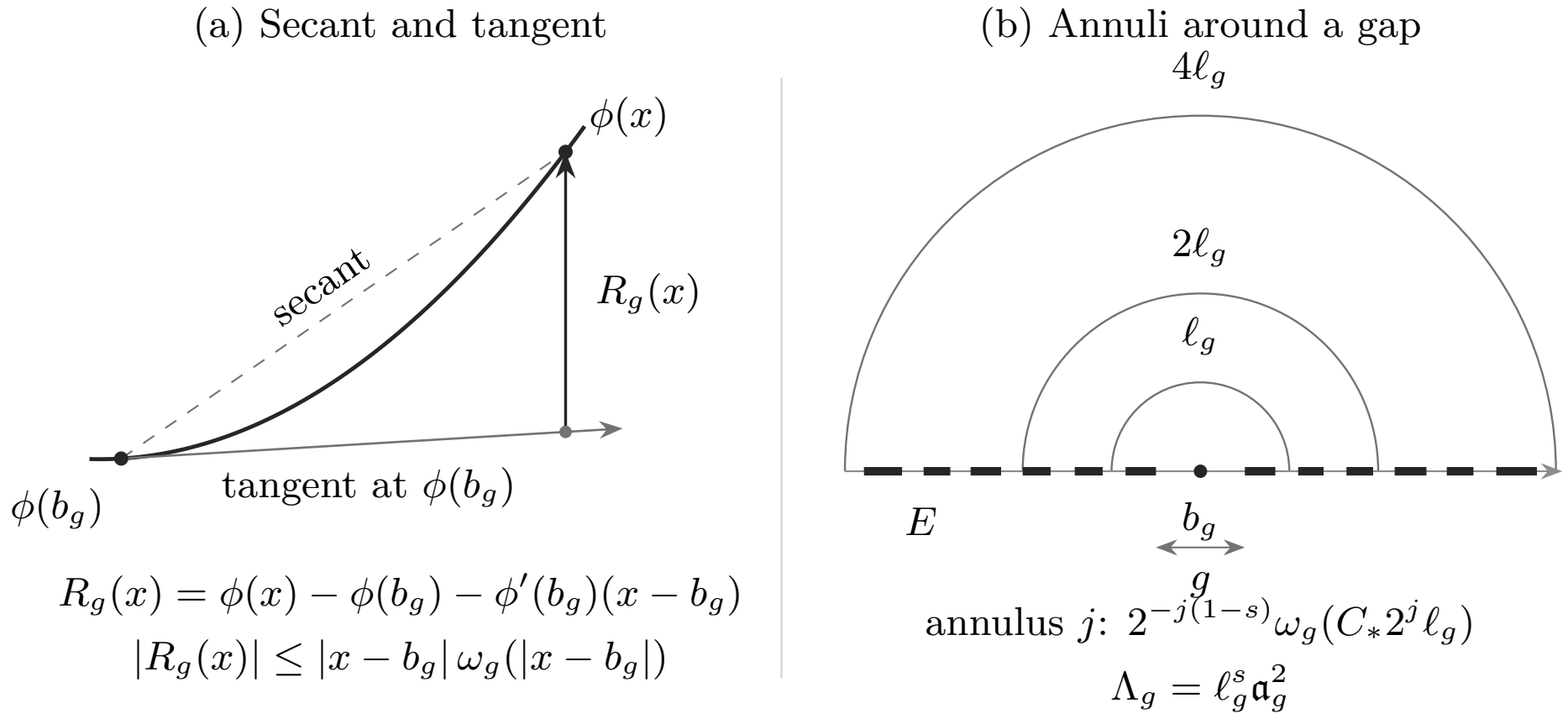
\begin{figure}[htbp]
\centering
\begin{tikzpicture}[x=1cm,y=1cm,font=\footnotesize,>=Stealth]
  \node at (2.55,4.55) {(a) Secant and tangent};
  \begin{scope}[shift={(.35,.15)}]
    \draw[DiagramPrimary,line width=.9pt,domain=-.25:3.65,samples=70]
      plot (\x,{1+.06*\x+.18*\x*\x});
    \draw[->,DiagramSecondary,line width=.5pt] (0,1)--(3.95,1.237);
    \draw[DiagramSecondary,dashed] (0,1)--(3.5,3.415);
    \draw[->,DiagramPrimary,line width=.65pt] (3.5,1.21)--(3.5,3.415);
    \fill[DiagramPrimary] (0,1) circle (1.6pt);
    \fill[DiagramPrimary] (3.5,3.415) circle (1.6pt);
    \fill[DiagramSecondary] (3.5,1.21) circle (1.4pt);
    \node[anchor=north east] at (0,.90) {$\phi(b_g)$};
    \node[anchor=south west] at (3.55,3.40) {$\phi(x)$};
    \node[anchor=west] at (3.60,2.30) {$R_g(x)$};
    \node[anchor=south,rotate=34.6] at (1.70,2.17) {secant};
    \node[anchor=north] at (1.95,1.07) {tangent at $\phi(b_g)$};
  \end{scope}
  \node[align=center] at (2.55,.08)
    {$R_g(x)=\phi(x)-\phi(b_g)-\phi'(b_g)(x-b_g)$};
  \node[align=center] at (2.55,-.43)
    {$|R_g(x)|\le |x-b_g|\,\omega_g(|x-b_g|)$};

  \draw[black!18] (5.55,-.58)--(5.55,4.15);
  \begin{scope}[shift={(6.2,0)}]
    \node at (2.65,4.55) {(b) Annuli around a gap};
    \foreach \rr in {.70,1.40,2.80}{
      \draw[DiagramSecondary,line width=.45pt]
        ({2.65+\rr},1.05) arc[start angle=0,end angle=180,radius=\rr];
    }
    \draw[->,black!45] (-.15,1.05)--(5.50,1.05);
    \foreach \a/\b in {0/.30,.47/.68,.84/1.08,1.25/1.47,1.66/1.88,2.06/2.30,3.0/3.24,3.42/3.62,3.80/4.01,4.18/4.43,4.59/4.80,4.98/5.30}{
      \draw[DiagramPrimary,line width=1.9pt] (\a,1.05)--(\b,1.05);
    }
    \fill[DiagramPrimary] (2.65,1.05) circle (1.5pt);
    \node[anchor=north] at (2.65,.94) {$b_g$};
    \node[anchor=south] at (2.65,1.83) {$\ell_g$};
    \node[anchor=south] at (2.65,2.53) {$2\ell_g$};
    \node[anchor=south] at (2.65,3.93) {$4\ell_g$};
    \node[anchor=north] at (.45,.90) {$E$};
    \draw[<->,DiagramSecondary] (2.30,.45)--(3.00,.45);
    \node[anchor=north] at (2.65,.39) {$g$};
    \node[align=center] at (2.65,-.13)
      {annulus $j$: $2^{-j(1-s)}\omega_g(C_*2^j\ell_g)$};
    \node[align=center] at (2.65,-.66)
      {$\Lambda_g=\ell_g^s\mathfrak a_g^2$};
  \end{scope}
\end{tikzpicture}
\caption{The local oscillation estimate. The secant defect $R_g(x)$
measures the error from the tangent at the transported pole. On the
parameter line, annular integration gives the factor $2^{-j(1-s)}$.
The sum of these contributions is $\mathfrak a_g$, and its square,
weighted by $\ell_g^s$, enters the Carleson packing condition.
The arc and the support at a finite stage are schematic.}
\label{fig:secant-annular-oscillation}
\end{figure}

Figure~\ref{fig:secant-annular-oscillation} illustrates the secant defect
and its annular decomposition.

\begin{proposition}[Local oscillation criterion]
\label{prop:curve-local-oscillation}
If $\mathcal C(\Lambda)<\infty$, the formal operator $H$ extends
boundedly to every $L^p(\mu)$, $1<p<\infty$, and
\begin{equation}
 \|H\|_{2\to2}\le\frac C{m_0}\mathcal C(\Lambda)^{1/2}.
 \label{cv:local-oscillation-L2}
\end{equation}
If $\mathcal C_r(\Lambda)\to0$ as $r\downarrow0$, the finite
correction sums converge to $H$ in operator norm on every such
$L^p$, and $H$ is compact there. The assertions hold uniformly for
arbitrary subfamilies of gaps.
\end{proposition}

\begin{proof}
The fundamental theorem of calculus gives
\[
 |h_g(x)|\le\frac{\sqrt{w_g}}{m_0}
              \frac{\omega_g(|x-b_g|)}{|x-b_g|}.
\]
Fix $z\in B_E(\xi_g,\ell_g)$. Divide $E$ into the ball
$B_E(\xi_g,2\ell_g)$ and the successive annuli about $\xi_g$.
Centrality handles the first ball. On the $j$th subsequent annulus,
$|x-b_g|\asymp2^j\ell_g$. The definition of the maximal function,
upper regularity, and $w_g\asymp\ell_g^{2-s}$ give
\[
 \begin{aligned}
 |\langle f,h_g\rangle|
 &\le\frac C{m_0}\ell_g^{1-s/2}
       \sum_{j\ge0}\frac{\omega_g(\min(1,C_*2^j\ell_g))}
                         {2^j\ell_g}
           (2^j\ell_g)^s M_\mu f(z)\\
 &=\frac C{m_0}\ell_g^{s/2}\mathfrak a_gM_\mu f(z).
 \end{aligned}
\]
For radii above one, use $\mu(E)=1\le(2^j\ell_g)^s$.
Taking the infimum in $z$ yields
\[
 |\langle f,h_g\rangle|
 \le\frac C{m_0}\Lambda_g^{1/2}
              \inf_{B_E(\xi_g,\ell_g)}M_\mu f.
\]
For a finite set $\mathcal F$ of gaps, let
$H_{\mathcal F}f=\sum_{g\in\mathcal F}\langle f,e_g\rangle h_g$.
Parseval, the Carleson embedding \eqref{st:carleson-embedding}, and
the maximal theorem imply
\[
 \begin{aligned}
 \|H_{\mathcal F}^*f\|_2^2
 &=\sum_{g\in\mathcal F}|\langle f,h_g\rangle|^2\\
 &\le\frac C{m_0^2}\mathcal C(\Lambda)
                     \|M_\mu f\|_2^2
 \le\frac C{m_0^2}\mathcal C(\Lambda)\|f\|_2^2.
 \end{aligned}
\]
This proves \eqref{cv:local-oscillation-L2} first on the finite
canonical span and then on $L^2$.

The passage to other exponents uses kernel estimates that require
no Dini summability. Centrality and the lower chord bound imply
\[
 |h_g(x)|\le C(1+L/m_0)|e_g(x)|.
\]
For $x,x'\in E$, exact subtraction gives
\[
 \begin{split}
 h_g(x)-h_g(x')=\sqrt{w_g}\left(
 \frac{\phi'(b_g)(\phi(x')-\phi(x))}
      {(\phi(x)-\phi(b_g))(\phi(x')-\phi(b_g))}
 -\frac{x'-x}{(x-b_g)(x'-b_g)}\right).
 \end{split}
\]
Since $\phi$ is $L$-Lipschitz and the distances to $a_g$ and $b_g$
are comparable on $E$, it follows that
\[
 \begin{aligned}
 |h_g(x)-h_g(x')|
 &\le C(1+L^2/m_0^2)
       \frac{\sqrt{w_g}|x-x'|}{|x-a_g|\,|x'-a_g|}\\
 &= C(1+L^2/m_0^2)|e_g(x)-e_g(x')|.
 \end{aligned}
\]
Thus the kernel
$K_{\mathcal F}(x,y)=\sum_{g\in\mathcal F}h_g(x)e_g(y)$ obeys,
by \eqref{eq:absolute-kernel-size} and~\eqref{eq:absolute-kernel-smoothness},
\[
 \begin{aligned}
 |K_{\mathcal F}(x,y)|&\le C_\phi|x-y|^{-s},\\
 |K_{\mathcal F}(x,y)-K_{\mathcal F}(x',y)|
   &\le C_\phi\frac{|x-x'|^{1-s}}{|x-y|},
 \qquad |x-x'|\le\tfrac12|x-y|.
 \end{aligned}
\]
The same continuity estimate holds in the second variable.
The constant $C_\phi$ depends only on the fixed geometric data,
$L$, and $m_0$. It is independent of $\mathcal F$ and of the
oscillation packing. Lemma~\ref{lem:vector-cz}, with scalar input
and output, now gives uniform bounds on every interior $L^p$.
The density of finite canonical expansions identifies the resulting
extensions.

If $\mathcal F\subset\{g:\ell_g<r\}$, the same $L^2$ proof gives
\[
 \|H_{\mathcal F}\|_{2\to2}
       \le\frac C{m_0}\mathcal C_r(\Lambda)^{1/2}.
\]
Fix $p\ne2$, choose $q\in(1,\infty)$ so that $p$ lies strictly
between $2$ and $q$, and let $0<\theta<1$ satisfy
$1/p=\theta/2+(1-\theta)/q$. The uniform $L^q$ bound and
interpolation give
\[
 \sup_{\mathcal F\subset\{\ell_g<r\}}
       \|H_{\mathcal F}\|_{p\to p}
 \le C_{p,q,\phi,\mathcal C(\Lambda)}
               \mathcal C_r(\Lambda)^{\theta/2}
 \longrightarrow0.
\]
The case $p=2$ is already proved. Only finitely many gaps have
length at least $r$, so the finite sums converge in operator norm
and their limit is compact.
\end{proof}

\begin{corollary}[Square Dini transport on general curves]
\label{cor:general-square-dini}
Under the preceding geometric hypotheses, let $\omega$ be the
modulus of continuity of $\phi'$. If
\[
 \mathfrak D_2(\phi)^2
       :=\sum_{n\ge0}\omega(2^{-n})^2<\infty,
\]
then $H$ is compact on every $L^p(\mu)$, $1<p<\infty$, and
\[
 \|H\|_{2\to2}\le C\mathfrak D_2(\phi)/m_0.
\]
Finiteness of this sum is equivalent to
$\int_0^1\omega(t)^2\,dt/t<\infty$.
The full transported family is Bessel if and only if the original
parameter displacement packing $\mathcal C(\lambda)$ is finite.
When this packing is finite,
\begin{equation}
 \|T+H\|_{\mathrm{ess},L^p\to L^p}
       \asymp_p\gamma(\lambda)^{1/2},\qquad1<p<\infty.
 \label{cv:general-square-dini-essential}
\end{equation}
Consequently $I+T+H$ is Fredholm exactly when $I+T$ is Fredholm,
with the same index. For canonical parameter poles, $I+H$ is
Fredholm of index zero.
\end{corollary}

\begin{proof}
Put $q=1-s$ and $\omega_n=\omega(2^{-n})$. If
$2^{-n-1}<\ell_g\le2^{-n}$, choose a fixed integer $J$ with
$2^J\ge C_*$ and split the defining sum for $\mathfrak a_g$
at $j=n-J$. After increasing a fixed constant to include $n<J$,
\[
 \mathfrak a_g\le C a_n,
 \qquad a_n=\sum_{k=0}^n2^{-(n-k)q}\omega_k.
\]
Indeed, the terms with $j\le n-J$ reindex with $k=n-j-J$;
the remaining geometric tail is bounded by
$C2^{-nq}\omega_0$. Discrete Young's inequality gives
\[
 \sum_{n\ge0}a_n^2
 \le\left(\sum_{j\ge0}2^{-jq}\right)^2
           \sum_{n\ge0}\omega_n^2.
\]
The gap count at a fixed scale shows, in every testing ball of radius
$R$, that the gaps in a single length band with $\ell_g\le R$
contribute at most $CR^sa_n^2$ to $\sum\Lambda_g$. Therefore
\[
 \mathcal C(\Lambda)\le C\sum_{n\ge0}a_n^2,
 \qquad
 \mathcal C_{2^{-N}}(\Lambda)
        \le C\sum_{n\ge\max(0,N-1)}a_n^2\longrightarrow0.
\]
Proposition~\ref{prop:curve-local-oscillation} gives the first
assertions. Monotonicity of $\omega$ gives the equivalence with
the integral condition by integration on dyadic intervals.

On finite canonical expansions, the full transported operator is
$I+T+H$. Since $H$ is bounded on $L^2$, its addition preserves
boundedness of the synthesis operator, proving the Bessel criterion.
It is compact on every interior $L^p$, so addition of $H$ preserves
the essential norm of $T$, Fredholmness, and the Fredholm index.
Apply Theorem~\ref{thm:stability} to obtain
\eqref{cv:general-square-dini-essential}.
\end{proof}

The Fredholm conclusion does not assert injectivity. The relative
smallness condition in Theorem~\ref{thm:curve-transport} still gives
invertibility, while the exact M\"obius identity provides another
class in which invertibility is preserved without smallness.

\subsection{Dini estimates at the endpoints}
\begin{proof}[Dini estimates in Theorem~\ref{thm:curve-transport}]
Throughout this proof,
the parameter poles $b_g$ are real and satisfy
$|b_g-a_g|\leq c_0\ell_g$. Let $\phi:[0,1]\to\C$ be continuously
differentiable and suppose that
\[
 |\phi(x)-\phi(y)|\geq m_0|x-y|
 \qquad(x,y\in[0,1])
\]
for some $m_0>0$. In particular, $\phi$ is injective and
$|\phi'(x)|\geq m_0$. Define the modulus of continuity of its derivative
and the Dini sum by
\[
 \omega_\phi(t)=
 \sup_{\substack{u,v\in[0,1]\\|u-v|\leq t}}
             |\phi'(u)-\phi'(v)|,
 \qquad
 \mathfrak D_\phi=\sum_{j=0}^{\infty}\omega_\phi(2^{-j})<\infty.
\]
Monotonicity gives
\[
 (\log2)\sum_{j\geq1}\omega_\phi(2^{-j})
 \leq\int_0^1\frac{\omega_\phi(t)}t\,dt
 \leq(\log2)\mathfrak D_\phi.
\]
Thus the discrete condition is equivalent to finiteness of the
usual Dini integral. Write $\omega=\omega_\phi$ for the rest of the proof.
Put $\Gamma=\phi(E)$ and $\mu_\phi=\phi_\#\mu$. The pullback
\[
 \mathcal U_\phi:L^p(\mu_\phi)\longrightarrow L^p(\mu),
 \qquad \mathcal U_\phi f=f\circ\phi,
\]
is an isometric isomorphism for $1\leq p\leq\infty$. Define the
rational functions with tangent normalization on $\Gamma$ by
\[
 Q_g(z)=\frac{\sqrt{w_g}\,\phi'(b_g)}{z-\phi(b_g)}.
\]
The denominator is nonzero on $\Gamma$, since $b_g\notin E$ and
$\phi$ is injective.
The correction to the family on the real line is
\[
 \mathcal U_\phi Q_g=q_g+h_g,
 \qquad
 h_g(x)=\sqrt{w_g}\left(
 \frac{\phi'(b_g)}{\phi(x)-\phi(b_g)}-
 \frac1{x-b_g}\right).
\]
Writing the numerator as an integral gives
\[
 \phi(x)-\phi(b_g)-\phi'(b_g)(x-b_g)
 =\int_{b_g}^{x}\bigl(\phi'(u)-\phi'(b_g)\bigr)\,du.
\]
Consequently,
\begin{equation}
 |h_g(x)|\leq
 \frac{\sqrt{w_g}}{m_0}
 \frac{\omega(|x-b_g|)}{|x-b_g|},
 \qquad x\in E.
 \label{cv:remainder}
\end{equation}
The derivative in this estimate is taken along the real interval;
no holomorphic extension of $\phi$ is assumed.
Figure~\ref{fig:curve-deformation} illustrates the transported poles
and the normalization in this identity.

\begin{figure}[htbp]
\centering
\begin{tikzpicture}[x=1cm,y=1cm,font=\small,>=Stealth]
 \draw[->,black!45] (0,0)--(5.1,0);
 \foreach \a/\b in {0/.55,1.1/1.65,3.3/3.85,4.4/4.95}
   \draw[black,line width=2pt] (\a,0)--(\b,0);
 \fill[DiagramSecondary] (2.475,0) circle (2pt);
 \node[anchor=north] at (2.475,-.08) {$b_g$};
 \node[anchor=north] at (.2,-.08) {$E$};
 \draw[->,thick,black!60] (5.4,.5)--(6.35,.5)
   node[midway,above] {$\phi$};
 \begin{scope}[shift={(6.8,0)}]
  \draw[black!35,thin,domain=0:4.95,samples=50] plot (\x,{.055*\x*\x});
  \foreach \a/\b in {0/.55,1.1/1.65,3.3/3.85,4.4/4.95}
    \draw[DiagramPrimary,line width=2pt,domain=\a:\b,samples=15]
      plot (\x,{.055*\x*\x});
  \fill[DiagramSecondary] (2.475,.337) circle (2pt);
  \draw[->,DiagramSecondary,thick] (2.475,.337)--(3.1,.507);
  \node[anchor=north] at (2.475,.20) {$\phi(b_g)$};
  \node[anchor=south] at (3.15,.54) {$\phi'(b_g)$};
  \node[anchor=north] at (.25,-.09) {$\phi(E)$};
 \end{scope}
 \node[anchor=north,align=center] at (5.9,-.75)
 {$\displaystyle Q_g(z)=\frac{\sqrt{w_g}\,\phi'(b_g)}{z-\phi(b_g)}$\\[6pt]
  $\displaystyle Q_g(\phi(x))=q_g(x)+h_g(x),\qquad
    |h_g(x)|\le\frac{\sqrt{w_g}}{m_0}
             \frac{\omega_\phi(|x-b_g|)}{|x-b_g|}$};
\end{tikzpicture}
\caption{Transport of poles and tangent normalization under a smooth
injective deformation. The support is shown at a finite stage of the
Cantor construction. The correction is controlled by the modulus of
continuity of the derivative. Under the Dini condition, the operator
obtained by summing these corrections is compact. The transported poles
are prescribed geometrically and need not be zeros of the Cauchy
transform of the image measure.}
\label{fig:curve-deformation}
\end{figure}

For a finite family $\mathcal F\subset\mathcal G$, let
\[
 H_{\mathcal F}f(x)
 =\sum_{g\in\mathcal F}\langle f,e_g\rangle h_g(x).
\]
Its integral kernel and a nonnegative upper bound are
\[
 L_{\mathcal F}(x,y)
 =\sum_{g\in\mathcal F}h_g(x)\overline{e_g(y)},
 \qquad
 S_{\mathcal F}(x,y)
 =\frac1{m_0}\sum_{g\in\mathcal F}
 \frac{w_g\omega(|x-b_g|)}{|x-b_g|\,|y-a_g|}.
\]
By \eqref{cv:remainder}, $|L_{\mathcal F}|\leq S_{\mathcal F}$.
We prove both Schur estimates directly, uniformly in $\mathcal F$.
Centrality and upper regularity, followed by summation over distance
annuli, give
\[
 \begin{aligned}
 \int_E\frac{d\mu(y)}{|y-a_g|}
 &\leq C\sum_{k\geq0}
       (2^k\ell_g)^{-1}(2^k\ell_g)^s\\
 &\leq\frac{C\ell_g^{s-1}}{1-2^{s-1}}.
 \end{aligned}
\]
The fixed factor $c_*$ in the minimum distance is included in the constants.
Multiplication by $w_g\lesssim\ell_g^{2-s}$ shows that
\begin{equation}
 \int_E S_{\mathcal F}(x,y)\,d\mu(y)
 \leq\frac C{m_0}\sum_{g\in\mathcal F}
           \ell_g\frac{\omega(|x-b_g|)}{|x-b_g|}.
 \label{cv:row-discrete}
\end{equation}
Set $c_b=c_*-c_0>0$. Both distances from $b_g$ to the endpoints
of $g$ are at least $c_b\ell_g$. Fix $x\in E$ and put
$d=|x-b_g|$. Take the subinterval of $g$ from $b_g$ towards the
endpoint farther from $x$, and denote it by $I_g(x)$.
This interval has length at least $c_b\ell_g$. For $u\in I_g(x)$,
\[
 d\leq|x-u|\leq d+\ell_g\leq(1+c_b^{-1})d.
\]
Since $\omega$ is nondecreasing,
\[
 \int_{I_g(x)}\frac{\omega(|x-u|)}{|x-u|}\,du
 \geq\frac{c_b}{1+c_b^{-1}}\,
             \ell_g\frac{\omega(d)}d.
\]
The gaps are disjoint. Summing this estimate in
\eqref{cv:row-discrete} therefore yields
\begin{equation}
 \sup_{x\in E}\int_E S_{\mathcal F}(x,y)\,d\mu(y)
 \leq\frac C{m_0}\int_0^1\frac{\omega(t)}t\,dt
 \leq\frac C{m_0}\mathfrak D_\phi.
 \label{cv:row-schur}
\end{equation}
Here the factor two accounting for the two sides of $x$ is included
in $C$.
For the column estimate, put $q=1-s>0$ and use the bands
\[
 \mathcal G_j=
 \{g:2^{-j-1}\leq\ell_g<2^{-j}\},
 \qquad j\geq0.
\]
If $t=2^{-j-1}$, the endpoint count at a fixed scale gives
\[
 \begin{aligned}
 \sum_{g\in\mathcal G_j}\frac{w_g}{|y-a_g|}
 &\leq Ct^{2-s}\sum_{k\geq0}
                    (2^kt)^{-1}2^{ks}\\
 &\leq Ct^{1-s}\sum_{k\geq0}2^{-k(1-s)}
 \leq C2^{-jq}.
 \end{aligned}
\]
This bound is uniform for $y\in E$ and also holds for every
subfamily of $\mathcal G_j$.
The region $|y-a_g|<t$ contains only a bounded number of such
poles. Centrality bounds its contribution by $Ct^{1-s}$, which
is absorbed in the initial term of the displayed sum.
Choose an integer $L\geq1$ with $2^{-L}<c_b/2$. For
$g\in\mathcal G_j$ and $x\in E$,
\[
 |x-b_g|\geq c_b\ell_g
 \geq c_b2^{-j-1}>2^{-j-L}.
\]
The distance annuli
$2^{-k-1}<|x-b_g|\leq2^{-k}$ are therefore empty for
$k>j+L$. Upper growth for balls with arbitrary centers gives
\[
 \begin{aligned}
 \int_E\frac{\omega(|x-b_g|)}{|x-b_g|}\,d\mu(x)
 &\leq\sum_{k=0}^{j+L}
     2^{k+1}\omega(2^{-k})\mu(B(b_g,2^{-k}))\\
 &\leq C\sum_{k=0}^{j+L}\omega(2^{-k})2^{kq}.
 \end{aligned}
\]
The possible boundary points of these balls have zero measure.
Combining the last two bounds and writing $\omega_k=\omega(2^{-k})$,
we obtain
\[
 \begin{aligned}
 \int_E S_{\mathcal F}(x,y)\,d\mu(x)
 &\leq\frac C{m_0}\sum_{j\geq0}\sum_{k=0}^{j+L}
                           \omega_k2^{-(j-k)q}\\
 &=\frac C{m_0}\sum_{k\geq0}\omega_k
             \sum_{j\geq\max(0,k-L)}2^{-(j-k)q}.
 \end{aligned}
\]
The inner sum is explicitly
\[
 \frac{2^{(k-\max(0,k-L))q}}{1-2^{-q}}
 \leq\frac{2^{Lq}}{1-2^{-q}}.
\]
Thus
\begin{equation}
 \sup_{y\in E}\int_E S_{\mathcal F}(x,y)\,d\mu(x)
 \leq\frac C{m_0}\mathfrak D_\phi.
 \label{cv:column-schur}
\end{equation}
The constants in \eqref{cv:row-schur} and
\eqref{cv:column-schur} depend only on $A,s,c_0$.
The two endpoint estimates and interpolation now give
\[
\|H_{\mathcal F}\|_{L^p(\mu)\to L^p(\mu)}
 \leq\frac C{m_0}\mathfrak D_\phi,
 \qquad1\leq p\leq\infty.
\]
The same $C$ works for every exponent in this range.
We prove convergence in operator norm of the correction series,
including the endpoints. For gaps of length less than $r$,
the global gap count implies
\[
 \begin{aligned}
 \sum_{\ell_g<r}\ell_g
 &\leq\sum_{n\geq0}2^{-n}r\,
        \#\{g:2^{-n-1}r\leq\ell_g<2^{-n}r\}\\
 &\leq Cr^{1-s}\sum_{n\geq0}2^{-n(1-s)}
 \leq Cr^q.
 \end{aligned}
\]
The row argument can be restricted to these same gaps. Splitting
their union at distance $\delta$ from $x$ gives, for $0<\delta<1$,
\[
 \sup_x\int_E S_{\{\ell_g<r\}}(x,y)\,d\mu(y)
 \leq\frac C{m_0}\left(
       \int_0^\delta\frac{\omega(t)}t\,dt
       +\frac{\omega(1)}\delta r^q\right).
\]
The notation denotes a nonnegative infinite sum, obtained by
monotone limits of finite subfamilies. Set $\delta=r^{q/2}$.
The resulting bound tends to zero as $r\downarrow0$.
For the column estimate take $r=2^{-J}$. Its tail is bounded by
\[
 \begin{aligned}
 \sup_y\int_E S_{\{\ell_g<2^{-J}\}}(x,y)\,d\mu(x)
 &\leq\frac C{m_0}\sum_{k\geq0}\omega_k
          \sum_{j\geq\max(J,k-L)}2^{-(j-k)q}\\
 &=\frac C{m_0(1-2^{-q})}\sum_{k\geq0}\omega_k
               2^{-(\max(J,k-L)-k)q}.
 \end{aligned}
\]
For each fixed $k$, the factor after $\omega_k$ tends to zero.
It is bounded by $2^{Lq}$ for all $J,k$. By dominated convergence for
the summable sequence $(\omega_k)$, this bound tends to zero as
$J\to\infty$.
Both tail bounds also hold for arbitrary subfamilies of the omitted
gaps. The Schur estimates and interpolation show that
$(H_{\mathcal F})$ converges in operator norm, along all finite
subsets, to a consistent operator $H$ on every $L^p(\mu)$,
$1\leq p\leq\infty$.
The same row bounds apply in the uniform norm on $E$.
Since the finite sums have continuous ranges, they also converge in
operator norm on $C(E)$ to a compact operator with the same bound.
Every $H_{\mathcal F}$ has finite rank. Consequently, $H$ is compact
on each of these spaces, and
\[
 \|H\|_{L^p(\mu)\to L^p(\mu)}
 \leq C\mathfrak D_\phi/m_0,
 \qquad H1=0,\qquad He_g=h_g.
\]
Compactness requires no smallness assumption on $\mathfrak D_\phi/m_0$.
Let $V_b(c)=\sum_gc_gq_g$ and
$V_\phi(c)=\sum_gc_g\mathcal U_\phi Q_g$, initially on finite
coefficient sequences. If $U(c)=\sum_gc_ge_g$, then
\[
 V_\phi=V_b+HU,
 \qquad \|HU\|_{\ell^2\to L^2}\leq C\mathfrak D_\phi/m_0.
\]
It follows that $V_\phi$ is bounded if and only if $V_b$ is bounded.
The criterion already proved on the real line therefore gives
\[
 (Q_g)\text{ is Bessel in }L^2(\mu_\phi)
 \quad\Longleftrightarrow\quad
 \mathcal C(\lambda)<\infty.
\]
This equivalence holds whenever the Dini quantity is finite, without a
smallness hypothesis on the deformation.
For $1<p<\infty$, the pullback of the perturbed family is the
image of the original basis under $I+T+H$. Indeed,
\[
 (I+T+H)1=1,
 \qquad (I+T+H)e_g=e_g+d_g+h_g=\mathcal U_\phi Q_g.
\]
If $A_b=I+T$ is invertible, factor the pullback operator as
\[
 I+T+H=A_b\bigl(I+A_b^{-1}H\bigr).
\]
Thus $\|A_b^{-1}H\|_{p\to p}<1$ gives the relative basis criterion
in Theorem~\ref{thm:curve-transport}. We also record the direct
sufficient condition in terms of the total packing norm.
Put
\[
 \varepsilon_p=C_p\mathcal C(\lambda)^{1/2}
                         +C\mathfrak D_\phi/m_0.
\]
If $\varepsilon_p<1$, then
\[
 (I+T+H)^{-1}=\sum_{n=0}^{\infty}(-T-H)^n,
 \qquad
 \|(I+T+H)^{-1}\|_{p\to p}\leq(1-\varepsilon_p)^{-1}.
\]
Thus $\{1\}\cup\{Q_g\}$ is an unconditional basis of
$L^p(\mu_\phi)$. At $p=2$ it is a Riesz basis with bounds
$(1-\varepsilon_2)^2$ and $(1+\varepsilon_2)^2$.
The tangent normalization and transported poles are essential to
the formulation. To see why gap zeros cannot simply be transferred,
consider
\[
 \phi_\epsilon(t)=t+i\epsilon t^2,\qquad \epsilon>0,
 \qquad
 F_\epsilon(z)=\int_E\frac{d\mu(t)}{\phi_\epsilon(t)-z}.
\]
This parametrization satisfies the lower chord bound with $m_0=1$,
and $\omega_{\phi_\epsilon}(t)=2\epsilon t$, so that
$\mathfrak D_{\phi_\epsilon}=4\epsilon$.
For every real $a\in[0,1]\setminus E$,
\[
 \phi_\epsilon(t)-\phi_\epsilon(a)
 =(t-a)\bigl(1+i\epsilon(t+a)\bigr),
 \qquad \phi_\epsilon'(a)=1+2i\epsilon a.
\]
Multiplication by the conjugate of the second factor gives
\[
 \Im\frac{\phi_\epsilon'(a)}
             {\phi_\epsilon(t)-\phi_\epsilon(a)}
 =-\frac{\epsilon}{1+\epsilon^2(t+a)^2}.
\]
After integration,
\[
 \Im\bigl(\phi_\epsilon'(a)
                 F_\epsilon(\phi_\epsilon(a))\bigr)
 =-\epsilon\int_E
              \frac{d\mu(t)}{1+\epsilon^2(t+a)^2}<0.
\]
Hence $F_\epsilon$ has no zero on any transported gap arc
$\phi_\epsilon(g)$, however small $\epsilon$ is. The construction therefore
gives a stable rational basis on the curved support under the stated
smallness condition. It does not preserve the reciprocal Herglotz
representation or orthogonality on that support.
\end{proof}

\subsection{M\"obius covariance modulo rank one}
\label{subsec:mobius-covariance}

Tangent normalization gives an exact transformation law under
M\"obius maps. The correction has finite rank, with no smallness
assumption.

\begin{proposition}[M\"obius covariance]
\label{prop:mobius-covariance}
Let $\phi$ be $C^1$ and satisfy the lower chord hypothesis of
Theorem~\ref{thm:curve-transport}, with real parameter poles, and let
\[
 M(z)=\frac{az+b}{cz+d},\qquad ad-bc\ne0,
 \qquad c\phi(t)+d\ne0\quad(0\le t\le1).
\]
Then $M\circ\phi$ satisfies the same qualitative hypotheses.
Dini or square Dini continuity of the derivative is preserved whenever
it holds for $\phi$. On finite canonical expansions,
\begin{equation}
 H_{M\circ\phi}=H_\phi+N,\qquad
 Nf=\mathbf1\sum_g\langle f,e_g\rangle\sqrt{w_g}
                  \frac{c\phi'(b_g)}{c\phi(b_g)+d}.
 \label{cv:mobius-rank-one}
\end{equation}
The series defining $N$ converges in operator norm on every
$L^p(\mu)$, $1\le p\le\infty$, and on $C(E)$.
It has rank at most one and satisfies $N\mathbf1=0$ and $N^2=0$.

For $X=L^\infty(\mu)$, the assertions below concern realizations
obtained as limits in operator norm of finite correction sums. For
$A_\phi$, this requirement is imposed on the sums for $T+H_\phi$.
The Dini transport estimate above and the endpoint pole estimate in
Proposition~\ref{prop:endpoint-pole-stability} provide such realizations
under their respective hypotheses.

On any of these spaces on which $A_\phi=I+T+H_\phi$ is bounded,
\begin{equation}
 A_{M\circ\phi}=A_\phi(I+N).
 \label{cv:mobius-factorization}
\end{equation}
Thus $A_{M\circ\phi}$ is invertible if and only if $A_\phi$ is
invertible, and in that case
\[
 A_{M\circ\phi}^{-1}=(I-N)A_\phi^{-1}.
\]
On each of the stated spaces on which $H_\phi$ is bounded,
the approximation numbers satisfy, for every $n\ge1$,
\[
 a_{n+1}^X(H_{M\circ\phi})\le a_n^X(H_\phi),\qquad
 a_{n+1}^X(H_\phi)\le a_n^X(H_{M\circ\phi}),
\]
where $X$ is any of the stated spaces.  In particular, M\"obius
postcomposition preserves polynomial approximation rates and,
on $L^2(\mu)$, membership in every Schatten class $S^q$, $q>0$.
\end{proposition}

\begin{proof}
The identity
\[
 M(z)-M(w)=\frac{(ad-bc)(z-w)}{(cz+d)(cw+d)}
\]
shows that the lower chord bound is preserved.  The derivative of the composition
is $M'(\phi)\phi'$.  On a compact neighborhood of the arc avoiding
the pole of $M$, both $M'$ and $M''$ are bounded.  Since $\phi$
is Lipschitz, the modulus of this derivative is bounded by
$C(\omega_\phi(t)+t)$, so Dini or square Dini summability is preserved whenever it
holds for $\phi$.

Direct division gives the exact identity
\[
 \frac{M'(w)}{M(z)-M(w)}
 =\frac1{z-w}+\frac{c}{cw+d}.
\]
Apply it with $z=\phi(x)$ and $w=\phi(b_g)$ and multiply by
$\sqrt{w_g}\phi'(b_g)$.  The difference of the two correction
functions is constant in $x$ and equals
$\sqrt{w_g}\kappa_g$, where
\[
 \kappa_g=\frac{c\phi'(b_g)}{c\phi(b_g)+d},\qquad
 \sup_g|\kappa_g|<\infty.
\]
The gap count at a fixed scale and centrality give
\[
 \sup_{y\in E}\sum_{u/2\le\ell_g<u}\frac{w_g}{|y-a_g|}
 \le Cu^{1-s}.
\]
For example, the poles in the annulus of radius $2^j u$ about
$y$ number at most $C2^{js}$, and their contribution is at most
$Cu^{1-s}2^{-j(1-s)}$.  Summing annuli and then dyadic $u$
proves uniform absolute convergence of
\[
 K_N(y)=\sum_g\kappa_g\frac{w_g}{y-a_g}.
\]
Its sum is continuous and bounded on $E$, and the tail over small gaps
is bounded by $Cr^{1-s}$.  Since $\mu(E)=1$, the operators
$f\mapsto\mathbf1\int fK_N\,d\mu$ are bounded on all the
stated spaces, and the finite sums converge in operator norm.
The identity \eqref{cv:mobius-rank-one} holds on finite canonical
expansions. Whenever $H_\phi$ has a bounded extension, it passes
to $L^p(\mu)$, $1\le p<\infty$, and to $C(E)$ by density,
using Proposition~\ref{thm:approximation-norming} at the endpoints.
On $L^\infty(\mu)$ it passes to the specified limits in operator norm. Existence of such a limit for one correction implies its
existence for the other, since the finite sums defining $N$
converge in operator norm.

The canonical functions have mean zero, so $N\mathbf1=0$.
The range of $N$ consists of constants, and hence $N^2=0$ and
$(I+N)^{-1}=I-N$.  Also $A_\phi\mathbf1=\mathbf1$;
therefore $A_\phi N=N$ wherever $A_\phi$ is bounded.
Equation~\eqref{cv:mobius-factorization} follows on the canonical
span and then on each stated realization, together with the
inverse formula. Finally, adding or subtracting an operator
of rank at most one increases the approximation index by at most
one.  The two inequalities and their stated consequences follow.
\end{proof}

\begin{corollary}[Transport to circular arcs]
\label{cor:circular-arc-transport}
For admissible real parameter poles $b_g$, let $M$ be a
M\"obius map with no pole on $[0,1]$.
If $A_b=I+T$ is a bounded isomorphism of $L^p(\mu)$,
$1<p<\infty$, the family $\{1\}\cup\{Q_g\}$ with tangent normalization is an
unconditional basis of $L^p(M_\#\mu)$.  If $A_b$ is a bounded
isomorphism of $L^1(\mu)$ or $C(E)$, the transported family is a
Schauder basis of the corresponding image space when ordered by gap length.
In particular, all these conclusions
hold for the canonical poles, without a smallness condition on
$M$.  For those poles the transport operator after pullback is
exactly $I+N$, with inverse $I-N$.
\end{corollary}

\begin{proof}
Take $\phi(t)=t$ in the proposition, so that $H_\phi=0$.
The factorization is a bounded isomorphism whenever $A_b$ is.
Pullback under $M$ is an isometry for the pushforward measure
and for the uniform norm on the compact support.  The canonical
endpoint bases are supplied by
Proposition~\ref{thm:approximation-norming}.
\end{proof}

\subsection{Approximation at positive regularity}
\label{subsec:higher-curve-approximation}
For dyadic $r$, write $\mathcal Q_r$ for the gap partition at
scale $r$ and $\mathbb E_r$ for its conditional expectation.

\begin{lemma}[Polynomial projections on gap cells]
\label{lem:cell-polynomial-projection}
Fix an integer $k\ge1$. For every gap cell $Q$ at scale $r$, let
$\Pi_Q$ be the orthogonal projection in $L^2(\mu|_Q)$ onto the
restrictions of polynomials of degree less than $k$.
Its integral kernel satisfies
\[
 |\Pi_Q(x,z)|\le\frac{C_k}{\mu(Q)},\qquad x,z\in Q.
\]
Consequently, the operator $\Pi_{k,r}$ obtained by applying $\Pi_Q$
on each cell obeys
\[
 |\Pi_{k,r}f|\le C_k\mathbb E_r|f|,
 \qquad
 \|\Pi_{k,r}\|_{p\to p}\le C_k\quad(1\le p\le\infty),
\]
has rank at most $k\#\mathcal Q_r$, and maps $C(E)$ into $C(E)$.
The constants depend only on $k$ and the regularity data.
\end{lemma}

\begin{proof}
Choose $z_Q\in Q$ and push the probability measure
$\mu|_Q/\mu(Q)$ forward by $x\mapsto(x-z_Q)/r$.
The resulting measure $\nu_Q$ is supported in $[-C_0,C_0]$ and
satisfies
\[
 \nu_Q(B(t,h))\le C h^s\qquad(t\in\R,\ h>0).
\]
This follows from $\mu(Q)\asymp r^s$ and the upper growth bound
for balls with arbitrary centers.
The Gram matrices of $1,t,\ldots,t^{k-1}$ with respect to these
probability measures have a uniform positive lower bound. Otherwise,
after passing to subsequences, the measures would converge weakly and coefficient
vectors of Euclidean norm one would converge to a nonzero polynomial
$P$ with
\[
 \int|P|^2\,d\nu=0.
\]
The limit probability retains the displayed upper growth bound:
apply the Portmanteau inequality to open balls. It is therefore
nonatomic. A nonzero polynomial has finitely many real zeros,
which contradicts the last equality. The Gram matrices also have
a uniform upper bound because all the measures are supported in a
fixed interval.

The Gram matrix formula for the projection kernel, with
the rescaled monomial vector, now gives $C_k/\mu(Q)$.
Integrating this bound gives the pointwise domination by
$C_k\mathbb E_r|f|$ and both endpoint operator bounds.
Interpolation gives every intermediate exponent with the same
constant. The cells are clopen in $E$, so the piecewise polynomial
range is contained in $C(E)$.
\end{proof}

\begin{lemma}[Approximation from bounds on the correction functions]
\label{lem:gap-envelope-approximation}
Let $\sigma>0$ and let $k\ge1$ be an integer satisfying
$k>\sigma+s-1$. For each real parameter pole $b_g$, suppose
$h_g\in C^k([0,1]\setminus\{b_g\};\C)$ and
\begin{equation}
 |h_g(x)|\le B\sqrt{w_g}|x-b_g|^{\sigma-1},
 \qquad
 |h_g^{(k)}(x)|\le B\sqrt{w_g}|x-b_g|^{\sigma-1-k},
 \qquad x\ne b_g.
 \label{cv:gap-envelopes}
\end{equation}
Then the finite sums
\[
 S_{\mathcal F}f=\sum_{g\in\mathcal F}\langle f,e_g\rangle h_g
\]
converge in operator norm on every $L^p(\mu)$, $1\le p\le\infty$,
and on $C(E)$, to a compact operator $S$. Moreover,
\[
 \|S\|_{p\to p}\le CB,\qquad
 a_n^{(p)}(S)\le CBn^{-\sigma/s},
 \qquad n\ge1,\quad1\le p\le\infty.
\]
The same approximation estimate holds on $C(E)$. The constant
depends on $\sigma,k$, the regularity data, and the fixed pole
separation, and is independent of $p$.
\end{lemma}

\begin{proof}
The size estimate gives the absolute majorant in the Dini proof,
with the nondecreasing function $Bt^\sigma$ in place of
$\omega(t)/m_0$. The row and column calculations and their tail
arguments use only monotonicity and
$\sum_{j\ge0}2^{-j\sigma}<\infty$; they therefore apply for
every $\sigma>0$. They prove the norm bound and convergence in operator norm over
arbitrary finite subsets, including on $C(E)$,
where the finite sums have continuous ranges.

For a dyadic $r$, put
\[
 S_rf=\sum_{\ell_g\ge r}\langle f,e_g\rangle h_g,
 \qquad
 R_r=S_r+\Pi_{k,r}(S-S_r).
\]
If $N_r=\#\{g:\ell_g\ge r\}$, the preceding lemma gives
\begin{equation}
 \operatorname{rank}R_r
 \le N_r+k\#\mathcal Q_r
 \le(k+1)N_r+k\le C_k r^{-s}.
 \label{cv:approximation-rank}
\end{equation}

Write $d_g(x)=|x-b_g|$. Choose $A_0>2C_0+2$ and put
$B_0=A_0+C_0$. For $\ell_g<r$, define
\[
 U_{g,r}(x)=d_g(x)^{\sigma-1}
                  \mathbf1_{\{d_g(x)\le B_0r\}},
 \qquad
 V_{g,r}(x)=r^k d_g(x)^{\sigma-1-k}
                  \mathbf1_{\{d_g(x)>A_0r\}}.
\]
We claim that
\[
 |(I-\Pi_{k,r})h_g(x)|
 \le CB\sqrt{w_g}
       \bigl(U_{g,r}(x)+\mathbb E_rU_{g,r}(x)+V_{g,r}(x)\bigr).
\]
If $d_g(x)\le A_0r$, the size estimate applies at every point
of the cell of $x$, and the kernel bound for $\Pi_{k,r}$ proves
the claim. If $d_g(x)>A_0r$, every segment from $x$ to that cell
stays at distance at least $d_g(x)/2$ from $b_g$. Subtract the
Taylor polynomial of $h_g$ at $x$ of degree $k-1$.
The projection reproduces that polynomial, and Taylor's integral
remainder and \eqref{cv:gap-envelopes} give the term $V_{g,r}$.

For the row Schur bound, use
\[
 \sqrt{w_g}\|e_g\|_1\le C\ell_g.
\]
The portions of the gaps from $b_g$ towards the endpoint farther
from $x$ are disjoint, have length comparable to $\ell_g$, and
have distance from $x$ comparable to $d_g(x)$. Hence
\[
 \sum_{\ell_g<r}\ell_gU_{g,r}(x)
 \le C\int_{|z-x|\le Cr}|x-z|^{\sigma-1}\,dz
 \le C_\sigma r^\sigma.
\]
The same bound holds after applying $\mathbb E_r$. The localized gap count
also gives, for $t\ge r$,
\[
 \sum_{\substack{\ell_g<r\\d_g(x)\le2t}}\ell_g
 \le C\sum_{j\ge0}(2^{-j}r)
                  \left(\frac{t}{2^{-j}r}\right)^s
 \le Ct^s r^{1-s}.
\]
Summing over dyadic $t\ge A_0r$ therefore yields
\[
 \sum_{\ell_g<r}\ell_gV_{g,r}(x)
 \le Cr^{k+1-s}\sum_{t\ge A_0r}t^{\sigma+s-1-k}
 \le Cr^\sigma,
\]
where the last exponent is negative by hypothesis.

For the column bound, the pole count at a fixed scale gives
\[
 \sum_{u/2\le\ell_g<u}\frac{w_g}{|y-a_g|}
 \le Cu^{1-s}.
\]
For every gap in this band, annular integration gives
\[
 \int_EU_{g,r}\,d\mu
 \le C\sum_{cu\le t\le Cr}t^{\sigma+s-1},
 \qquad
 \int_EV_{g,r}\,d\mu
 \le Cr^k\sum_{t\ge A_0r}t^{\sigma+s-1-k}
 \le Cr^{\sigma+s-1}.
\]
The integral of $\mathbb E_rU_{g,r}$ is the same as that of
$U_{g,r}$. Interchanging the nonnegative sums over nearby gaps gives
\[
 \begin{aligned}
 \sum_{u\le r}u^{1-s}\sum_{cu\le t\le Cr}t^{\sigma+s-1}
 &\le C\sum_{t\le Cr}t^{\sigma+s-1}
                \sum_{u\le\min(r,Ct)}u^{1-s}\\
 &\le C\sum_{t\le Cr}t^\sigma
 \le Cr^\sigma.
 \end{aligned}
\]
The sum over distant gaps is at most
$Cr^{\sigma+s-1}\sum_{u\le r}u^{1-s}\le Cr^\sigma$.
All scale variables in these displays are dyadic.
Thus both Schur bounds for $S-R_r$ are at most $CBr^\sigma$,
uniformly for finite subfamilies of the omitted gaps. The
limit in operator norm gives
\[
 \sup_{1\le p\le\infty}\|S-R_r\|_{p\to p}
 \le CBr^\sigma,
\]
also in the uniform norm on $C(E)$.
Choosing a dyadic $r\asymp n^{-1/s}$ with rank in
\eqref{cv:approximation-rank} less than $n$ proves the result;
the finitely many smaller indices are controlled by $\|S\|\le CB$.
\end{proof}

For an integer $m\ge0$ and $0<\alpha\le1$, use the norm
\[
 \|r\|_{C^{m,\alpha}}
 =\sum_{j=0}^m\|r^{(j)}\|_\infty+[r^{(m)}]_{C^\alpha}.
\]
The next lemma identifies the finite-dimensional part that must be
removed before applying the preceding approximation lemma.

\begin{lemma}[Jets of divided differences]
\label{lem:divided-difference-jets}
Let $\phi\in C^{m+1,\alpha}([0,1];\C)$ have lower chord constant
$c>0$, and suppose $\|\phi'\|_{C^{m,\alpha}}\le K$.
Set $\sigma=m+\alpha$ and $k=m+1$.
For every $b\in[0,1]$ there is a polynomial $J_{\phi,b}$ of
degree at most $m-1$ (zero when $m=0$) such that
\[
 \frac{\phi'(b)}{\phi(x)-\phi(b)}-\frac1{x-b}
 =J_{\phi,b}(x-b)+\mathcal R_{\phi,b}(x-b),
\]
where, for $0\le j\le k$ and $t\ne0$ with $b+t\in[0,1]$,
\begin{equation}
 |\mathcal R_{\phi,b}^{(j)}(t)|
 \le C|t|^{\sigma-1-j}.
 \label{cv:higher-jet-bounds}
\end{equation}
The coefficients of $J_{\phi,b}$ are uniformly bounded.
If $\phi,\psi$ satisfy the same bounds and
$\Delta=\|\phi'-\psi'\|_{C^{m,\alpha}}$, their jet coefficients
differ by at most $C\Delta$, and
\[
 |(\mathcal R_{\phi,b}-\mathcal R_{\psi,b})^{(j)}(t)|
 \le C\Delta|t|^{\sigma-1-j},
 \qquad0\le j\le k.
\]
All constants depend only on $c,K,m,\alpha$.
\end{lemma}

\begin{proof}
Write $r=\phi'$ and
\[
 D_b(t)=\frac{\phi(b+t)-\phi(b)}t
       =\int_0^1r(b+ut)\,du,
 \qquad |D_b(t)|\ge c.
\]
The value at zero is $D_b(0)=r(b)$, with $|r(b)|\ge c$.
Put
\[
 P_b(t)=\sum_{j=0}^m d_jt^j,
 \qquad d_j=\frac{r^{(j)}(b)}{(j+1)!},
 \qquad e_b=D_b-P_b.
\]
Taylor's theorem for $\phi\in C^{k,\alpha}$ gives remainder
bounds $C|t|^{k+\alpha-j}$ through derivative order $k$.
Dividing that remainder by $t$ and differentiating yields
\begin{equation}
 |e_b^{(j)}(t)|\le C|t|^{\sigma-j},
 \qquad0\le j\le k.
 \label{cv:secant-jet-remainder}
\end{equation}
This uses only derivatives of $\phi$ through order $k=m+1$.
For $j\le m$, the integral formula bounds $D_b^{(j)}$ uniformly;
for $j=k$, \eqref{cv:secant-jet-remainder} bounds it by
$C|t|^{\alpha-1}$. Differentiating the reciprocal therefore gives
\[
 |(D_b^{-1})^{(j)}(t)|\le C\quad(0\le j\le m),
 \qquad |(D_b^{-1})^{(k)}(t)|\le C|t|^{\alpha-1}.
\]

Define $\kappa_0=1$ and, for $1\le j\le m$, recursively set
\[
 \kappa_j=-\frac1{d_0}\sum_{\ell=1}^j d_\ell \kappa_{j-\ell}.
\]
Then $W_b(t)=\sum_{j=0}^m \kappa_jt^j$ satisfies
$P_bW_b-r(b)=O(t^{m+1})$ as a polynomial identity.
The $\kappa_j$ are uniformly bounded rational functions of the finite
jet of $r$, with denominators bounded away from zero.
Let
\[
 U_b(t)=r(b)-D_b(t)W_b(t).
\]
The polynomial identity and \eqref{cv:secant-jet-remainder} imply
\[
 |U_b^{(j)}(t)|\le C|t|^{\sigma-j},
 \qquad0\le j\le k,
\]
since $\sigma\le m+1$ and $|t|\le1$.
The product rule and the reciprocal bounds then give the same
estimate for
\[
 A_b(t)=\frac{U_b(t)}{D_b(t)}
       =\frac{r(b)}{D_b(t)}-W_b(t).
\]
In the only term using the $k$th reciprocal derivative, the factor
$U_b=O(|t|^\sigma)$ gives a stronger power than required.
Set
\[
 J_{\phi,b}(t)=\sum_{j=1}^m \kappa_jt^{j-1},
 \qquad \mathcal R_{\phi,b}(t)=\frac{A_b(t)}t.
\]
The asserted expansion and \eqref{cv:higher-jet-bounds} follow
by differentiating $A_b/t$ through order $k$.

For the difference assertion, Taylor's remainder is linear in
$\phi$, so the bounds for $e_{\phi,b}-e_{\psi,b}$ include the
factor $\Delta$. The coefficient recursion is Lipschitz on the
set of bounded jets where $|d_0|\ge c$. The same is true of reciprocal
derivatives: the orders through $m$ are bounded by $C\Delta$,
and the order $k$ by $C\Delta|t|^{\alpha-1}$.
Apply the preceding product identities to differences, keeping
one difference factor in each term. They give
$|A_{\phi,b}^{(j)}-A_{\psi,b}^{(j)}|
\le C\Delta|t|^{\sigma-j}$ and hence the stated bounds after
division by $t$.
\end{proof}

\begin{proposition}[Approximation at finite regularity]
\label{prop:higher-curve-approximation}
Fix an integer $m\ge0$, $0<\alpha\le1$, and put
$\sigma=m+\alpha$. If $\phi\in C^{m+1,\alpha}$ satisfies
the lower chord bound, then
\[
 a_n^{(p)}(H_\phi)\le Cn^{-\sigma/s},
 \qquad n\ge1,\quad1\le p\le\infty.
\]
The same estimate holds on $C(E)$.
For fixed real parameter poles and parametrizations $\phi,\psi$
with common lower chord and $C^{m,\alpha}$ derivative bounds,
\begin{equation}
 a_n^{(p)}(H_\phi-H_\psi)
 \le C\|\phi'-\psi'\|_{C^{m,\alpha}}n^{-\sigma/s},
 \qquad n\ge1,\quad1\le p\le\infty,
 \label{cv:higher-continuity}
\end{equation}
again including $C(E)$. The constants are independent of $p$.
\end{proposition}

\begin{proof}
Apply Lemma~\ref{lem:divided-difference-jets} at each $b_g$.
The correction splits into
\[
 h_g^\phi(x)=\sqrt{w_g}J_{\phi,b_g}(x-b_g)
                  +\sqrt{w_g}\mathcal R_{\phi,b_g}(x-b_g).
\]
The sums of the polynomial terms converge to an operator $F_\phi$
whose range consists of polynomials of degree at most $m-1$.
To justify convergence at both endpoints, expand these polynomials
in powers of $x$ and use the gap potential estimate
\begin{equation}
 \sup_{y\in E}\sum_{u/2\le\ell_g<u}
       \frac{w_g}{|y-a_g|}\le Cu^{1-s},
 \qquad
 \sup_{y\in E}\sum_g\frac{w_g}{|y-a_g|}\le C.
 \label{cv:global-gap-potential}
\end{equation}
The first inequality follows from the gap count at a fixed scale,
as in the Dini column estimate. Summing it gives the second because $s<1$.
Every polynomial coefficient kernel is therefore a uniformly
convergent sum of continuous functions of $y$, with bounded
supremum. Thus $F_\phi$ has rank at most $m$ and is bounded on
all the stated spaces. The jet difference estimate likewise gives
\[
 \operatorname{rank}(F_\phi-F_\psi)\le m,
 \qquad \|F_\phi-F_\psi\|_{p\to p}\le C\Delta,
 \qquad \Delta=\|\phi'-\psi'\|_{C^{m,\alpha}}.
\]

The remainder terms, and their differences, satisfy the bounds required
in the approximation lemma with exponent $\sigma$ and derivative order
$k=m+1$. The required strict inequality is automatic:
\[
 k-(\sigma+s-1)=2-\alpha-s>0.
\]
Lemma~\ref{lem:gap-envelope-approximation} consequently gives
\[
 a_n^{(p)}(H_\phi-F_\phi)\le Cn^{-\sigma/s},
\]
and
\[
 a_n^{(p)}\bigl(H_\phi-H_\psi-(F_\phi-F_\psi)\bigr)
 \le C\Delta n^{-\sigma/s},
\]
together with the corresponding norm bounds and estimates on $C(E)$.
Adding an operator of rank at most $m$ shifts the approximation
index by at most $m$. For $n>2m+2$ this shift changes the displayed
power by a fixed factor; for the remaining indices use the norm
bounds. This proves both assertions.

In particular, comparison with an affine parametrization shows
that the constant is linear in
$\|\phi'-v\|_{C^{m,\alpha}}$ near any fixed nonzero constant
derivative $v$, while the common geometric bounds are kept fixed.
\end{proof}

\begin{proposition}[Gap dependence and kernel regularity]
\label{prop:curve-kernel-obstruction}
Fix an integer $m\ge0$, $0<\alpha\le1$, and put
$\sigma=m+\alpha$. There are $C,\varepsilon_0>0$ such that,
for every $0<\varepsilon\le\varepsilon_0$, a real-valued map
$\phi\in C^{m+1,\alpha}([0,1])$ satisfies
\[
 \phi|_E=\operatorname{id}|_E,\qquad \phi'|_E=1,\qquad
 \|\phi'-1\|_{C^{m,\alpha}}\le C\varepsilon,\qquad
 |\phi(x)-\phi(y)|\ge\tfrac12|x-y|.
\]
Moreover, $\phi(a_g)=a_g$ and
$\phi'(a_g)=1+\varepsilon\ell_g^\sigma$ for every gap $g$.
For canonical parameter poles $b_g=a_g$, its correction is
\begin{equation}
 H_\phi1=0,\qquad H_\phi e_g=\varepsilon\ell_g^\sigma e_g.
 \label{cv:fixed-support-diagonal}
\end{equation}

If $\sigma<s$, the operator $H_\phi$ has no bounded representing
kernel. If $0<\beta\le1$ and $\beta>\sigma-s$, it has no
representing kernel uniformly $\beta$-H\"older continuous in
either variable. Here a representing kernel means
$K\in L^1(\mu\times\mu)$ such that, for every $f\in C(E)$,
\[
 H_\phi f(x)=\int_E K(x,y)f(y)\,d\mu(y)
 \quad\text{for $\mu$-almost every $x$}.
\]
Uniform H\"older continuity in the first variable means that a
representative satisfies
\[
 |K(x,y)-K(x',y)|\le M|x-x'|^\beta
 \quad(x,x'\in E),
\]
for $\mu$-almost every $y$, with a common constant $M$;
the condition in the second variable is defined analogously.
\end{proposition}

\begin{proof}
Choose a real function
$\eta\in C_c^\infty((-c_*/2,c_*/2))$ with
$\eta(0)=0$ and $\eta'(0)=1$, and set
\[
 \Psi(t)=\sum_g\ell_g^{\sigma+1}
                 \eta\bigl((t-a_g)/\ell_g\bigr),\qquad
 \phi(t)=t+\varepsilon \Psi(t).
\]
Centrality places each summand strictly inside its gap, so their
supports are disjoint. For $0\le j\le m+1$, define
\[
 \Psi_j(t)=\sum_g\ell_g^{\sigma+1-j}
                 \eta^{(j)}\bigl((t-a_g)/\ell_g\bigr).
\]
These functions vanish on $E$. Wherever a summand is nonzero,
$\operatorname{dist}(t,E)\ge(c_*/2)\ell_g$; hence
\begin{equation}
 |\Psi_j(t)|\le C\operatorname{dist}(t,E)^{\sigma+1-j}.
 \label{cv:fixed-support-jet-bound}
\end{equation}
The exponents are positive for $j\le m+1$, so $\Psi_j$ is continuous
at every point of $E$. Off $E$ the derivative identities are
immediate. For $x\in E$ and $1\le j\le m+1$,
\[
 \frac{|\Psi_{j-1}(x+h)-\Psi_{j-1}(x)|}{|h|}
 \le C|h|^{\sigma+1-j}\longrightarrow0.
\]
Thus $\Psi\in C^{m+1}$ and $\Psi^{(j)}=\Psi_j$.

On any one gap, rescaling gives
$[\Psi_{m+1}]_{C^\alpha}\le C$, since its amplitude is
$\ell_g^\alpha$. If $t,z$ lie in different gaps, or one lies
in $E$, choose $x\in E$ between them. By
\eqref{cv:fixed-support-jet-bound},
\[
 |\Psi_{m+1}(t)-\Psi_{m+1}(z)|
 \le C\bigl(|t-x|^\alpha+|z-x|^\alpha\bigr)
 \le C'|t-z|^\alpha.
\]
This proves a global $C^{m+1,\alpha}$ bound for $\Psi$, including
$\alpha=1$. In particular,
$\|\phi'-1\|_{C^{m,\alpha}}\le C\varepsilon$.
Taking $\varepsilon_0$ small enough gives $\phi'\ge1/2$ and
the lower chord bound. The stated values on $E$ and at $a_g$
follow directly from the construction.

For $x\in E$, these identities give
\[
 \sqrt{w_g}\left(
 \frac{\phi'(a_g)}{\phi(x)-\phi(a_g)}-\frac1{x-a_g}\right)
 =\varepsilon\ell_g^\sigma e_g(x).
\]
The convergence in operator norm already established for curve
corrections, together with orthonormality of
$\{1\}\cup\{e_g\}$, proves \eqref{cv:fixed-support-diagonal}.
In particular, $H_\phi$ is self-adjoint on $L^2(\mu)$.

The gap estimates and annular integration give
\begin{equation}
 \|e_g\|_1\le C\ell_g^{s/2},\qquad
 \|e_g\|_\infty\asymp\ell_g^{-s/2},\qquad
 |e_g(u_g)-e_g(v_g)|\asymp\ell_g^{-s/2}.
 \label{cv:fixed-support-kernel-tests}
\end{equation}
For the last comparison, use
\[
 |e_g(u_g)-e_g(v_g)|
 =\frac{\sqrt{w_g}\ell_g}{(a_g-u_g)(v_g-a_g)}
 \asymp\frac{\sqrt{w_g}}{\ell_g}.
\]
If a representing kernel is essentially bounded by $M$, then
\[
 \varepsilon\ell_g^\sigma\|e_g\|_\infty
 =\|H_\phi e_g\|_\infty\le M\|e_g\|_1.
\]
Here continuity and full support of $\mu$ identify the supremum
of $e_g$ with its essential supremum. Therefore
$M\ge c\varepsilon\ell_g^{\sigma-s}$, which is impossible
as $\ell_g\to0$ if $\sigma<s$.

Suppose next that $K$ satisfies the H\"older condition in the
first variable. Fubini's theorem provides $x_0\in E$ with
$K(x_0,\cdot)\in L^1(\mu)$. The uniform difference bound
then makes every section $K(x,\cdot)$ integrable and shows that
\[
 x\longmapsto\int_EK(x,y)e_g(y)\,d\mu(y)
\]
is continuous, with H\"older seminorm at most $M\|e_g\|_1$.
It agrees almost everywhere with the continuous function
$\varepsilon\ell_g^\sigma e_g$ and hence agrees on all of $E$,
by full support of $\mu$. Evaluation at $u_g,v_g$ and
\eqref{cv:fixed-support-kernel-tests} yield
\[
 c\varepsilon\ell_g^{\sigma-s/2}
 \le M\ell_g^\beta\|e_g\|_1
 \le CM\ell_g^{\beta+s/2}.
\]
Thus $M\ge c\varepsilon\ell_g^{\sigma-s-\beta}$, excluding
$\beta>\sigma-s$. Finally, Fubini's theorem and self-adjointness
show that $K^*(x,y)=\overline{K(y,x)}$ also represents $H_\phi$.
A H\"older bound in the second variable for $K$ is a bound in the
first variable for $K^*$, so the same obstruction applies.
\end{proof}

This construction shows that the operator depends on the extension
through the complementary gaps: the support map, its derivative on $E$, and
the image pole positions are fixed, while the tangent factors at
the poles change the coefficient operator.

\subsection{Smooth and analytic deformations}
Exponential approximation from uniform analyticity in one kernel
variable is classical; see \cite[Theorem~11.3]{BirmanSolomyak1977}.
We include a direct argument that gives simultaneous estimates for this
correction at both endpoints and on $C(E)$.

\begin{corollary}[Smooth and analytic deformations]
\label{cor:smooth-analytic-curves}
If $\phi$ is $C^\infty$ and has a lower chord bound, then, for
every $N>0$,
\[
 a_n^{(p)}(H_\phi)\le C_Nn^{-N},
 \qquad n\ge1,\quad1\le p\le\infty.
\]
If $\phi$ extends holomorphically and injectively to a complex
neighborhood of $[0,1]$, there are $C,c>0$ such that
\[
 a_n^{(p)}(H_\phi)\le Ce^{-cn},
 \qquad n\ge1,\quad1\le p\le\infty.
\]
Both conclusions hold on $C(E)$, with constants independent of $p$.
In either case $H_\phi\in S^q$ for every $q>0$ on $L^2(\mu)$.
\end{corollary}

\begin{proof}
The smooth assertion follows from
Proposition~\ref{prop:higher-curve-approximation} by choosing
$m+\alpha\ge Ns$.
For the analytic assertion choose a closed Bernstein ellipse
about $[0,1]$, with parameter greater than one, inside the univalent
neighborhood. For real $b\in[0,1]$ the function
\[
 R_b(z)=\frac{\phi'(b)}{\phi(z)-\phi(b)}-\frac1{z-b}
\]
extends holomorphically across $z=b$, with value
$-\phi''(b)/(2\phi'(b))$. The divided difference
$(\phi(z)-\phi(b))/(z-b)$ is jointly continuous and never
vanishes on the relevant compact set, including its diagonal.
The singularity is therefore removable uniformly in $b$. Univalence excludes other
zeros of the denominator, and compactness gives a uniform bound
for $R_b$ on a slightly smaller fixed ellipse, whose
parameter we denote by $\varrho>1$. In particular, the analytic
hypothesis also supplies a lower chord bound on $[0,1]$.
Its Chebyshev polynomial approximation $P_{N,b}$ of degree less
than $N$ satisfies
\[
 \sup_{b\in[0,1]}\sup_{x\in[0,1]}
       |R_b(x)-P_{N,b}(x)|\le C\varrho^{-N}.
\]
For completeness, compose with
$z=(2+w+w^{-1})/4$. The resulting Laurent series is symmetric
under $w\mapsto w^{-1}$, and Cauchy's coefficient bound on the
fixed annulus makes its tail a geometric series. Pairing opposite
powers of $w$ gives precisely the stated polynomial estimate.

Replace $R_{b_g}$ by $P_{N,b_g}$ in each correction function and
sum the resulting kernels. The polynomial coefficients are uniformly
bounded for each $N$, and \eqref{cv:global-gap-potential} gives
uniform convergence. The resulting operator $A_N$ has polynomial
range of dimension at most $N$, while its error kernel satisfies
\[
 |K_{H_\phi-A_N}(x,y)|
 \le C\varrho^{-N}\sum_g\frac{w_g}{|y-a_g|}
 \le C\varrho^{-N}.
\]
Both Schur bounds, and the bound in the uniform norm on $C(E)$, are
therefore at most $C\varrho^{-N}$. Since
$\operatorname{rank}A_N\le N$, this bounds $a_{N+1}$ and proves
the exponential estimate after adjusting the constant.
The Schatten class assertions follow by summing the singular values.
\end{proof}

\subsection{Sharpness on curves parametrized by arc length}
\label{subsec:curve-sharpness}
\begin{proposition}[Sharpness throughout the regularity scale]
\label{prop:arclength-sharpness}
Fix $0<s<1$, an integer $m\ge0$, and $0<\alpha\le1$, and put
$\sigma=m+\alpha$. There are an $s$-Ahlfors regular Cantor
probability measure $\mu$ and constants $\varepsilon_0,c,C>0$
such that, for every $0<\varepsilon\le\varepsilon_0$, a map
$\phi_\varepsilon\in C^{m+1,\alpha}([0,1];\C)$ satisfies
\[
 |\phi_\varepsilon'|=1,\qquad
 |\phi_\varepsilon(x)-\phi_\varepsilon(y)|\ge\tfrac12|x-y|,
 \qquad
 \|\phi_\varepsilon'-1\|_{C^{m,\alpha}}\le C\varepsilon.
\]
The image $\phi_\varepsilon(E)$ is not contained in any affine line.
For the canonical parameter poles $b_g=a_g$, its correction satisfies
\begin{equation}
 c\varepsilon n^{-\sigma/s}
 \le a_n^{(p)}(H_{\phi_\varepsilon})
 \le C\varepsilon n^{-\sigma/s},
 \qquad n\ge1,\quad1\le p\le\infty.
 \label{cv:arclength-sharp-rate}
\end{equation}
The constants and $\varepsilon_0$ are independent of $p$.
The same upper and lower estimates hold on $C(E)$.
\end{proposition}

\begin{proof}
Let $\rho=2^{-1/s}<1/2$ and use the symmetric Cantor set generated
by $x\mapsto\rho x$ and $x\mapsto1-\rho+\rho x$.
Its natural probability measure assigns every generation-$j$
cylinder mass $2^{-j}=\rho^{js}$ and is $s$-Ahlfors regular.
The generation-$j$ gaps have length and multiplicity
\[
 \ell_j=(1-2\rho)\rho^{j-1},\qquad M_j=2^{j-1}.
\]

Choose an even $\zeta\in C_c^\infty((-1,1))$ with
$0\le\zeta\le1$ and $\zeta(0)=1$. Fix $0<d<c_*/5$ and put
\[
 \chi(t)=\zeta(t/d)
       -\zeta\bigl(2(t-2d)/d\bigr)
       -\zeta\bigl(2(t+2d)/d\bigr).
\]
The three supports are disjoint and lie in $(-c_*/2,c_*/2)$.
Thus $\chi(0)=1$ and $|\chi|\le1$.
Evenness and a change of variables give, for every real $A$,
\[
 \begin{aligned}
 \int_0^\infty\sin(A\chi(t))\,dt
 &=d\int_0^1\sin(A\zeta(v))\,dv
   -\frac d2\int_{-1}^1\sin(A\zeta(v))\,dv=0.
 \end{aligned}
\]
The same calculation on the negative half-line gives
\begin{equation}
 \int_{-\infty}^0\sin(A\chi(t))\,dt
 =\int_0^\infty\sin(A\chi(t))\,dt=0.
 \label{cv:balanced-angle}
\end{equation}

Set $A_g=\varepsilon\ell_g^\sigma$ and define
\[
 \theta(t)=
 \begin{cases}
 A_g\chi((t-a_g)/\ell_g),&t\in g,\\
 0,&t\in E.
 \end{cases}
\]
We give the gluing argument needed for higher regularity.
If smooth functions supported in the same fixed central portions
of the gaps satisfy
\[
 \|f_g^{(j)}\|_\infty\le M\ell_g^{\sigma-j}
 \quad(0\le j\le m),\qquad
 [f_g^{(m)}]_{C^\alpha(g)}\le M,
\]
their extension by zero belongs to $C^{m,\alpha}([0,1])$ with
norm at most $CM$. Indeed, central support gives
$|f_g^{(j)}(x)|\le CM\dist(x,E)^{\sigma-j}$.
For $j<m$, difference quotients at a point of $E$ tend to zero,
because $\sigma-j-1\ge\alpha>0$.
Induction therefore identifies all derivatives through order $m$
with the piecewise derivatives extended by zero. For the last
derivative, points in different gaps satisfy
$\dist(x,E)+\dist(y,E)\le|x-y|$, so the two boundary estimates
give its global $\alpha$-H\"older bound. The cases with one point
in $E$ and with both points in one gap follow directly.

Scaling the fixed profile $\chi$ now proves
\[
 \|\theta\|_{C^{m,\alpha}}\le C\varepsilon,
 \qquad \|\theta\|_\infty\le\varepsilon.
\]
For $|A|\le1$, the fixed profiles
$e^{iA\chi}-1$ have $C^{m,\alpha}$ norm at most $C|A|$,
and the profiles $\cos(A\chi)-1$ have that norm at most $CA^2$.
These bounds follow by differentiating through order $m+1$
and using $|\sin(A\chi)|\le|A|$ and
$|\cos(A\chi)-1|\le A^2/2$; every differentiated term in
the second profile either contains at least two powers of $A$
or a factor $A\sin(A\chi)$. Scaling and the same gluing argument give
\begin{equation}
 \|e^{i\theta}-1\|_{C^{m,\alpha}}\le C\varepsilon,
 \qquad
 \|\cos\theta-1\|_{C^{m,\alpha}}\le C\varepsilon^2.
 \label{cv:angle-regularity}
\end{equation}
For the second estimate, the scaled derivative bounds contain
$\varepsilon^2\ell_g^{2\sigma-j}$, which are bounded by
$\varepsilon^2\ell_g^{\sigma-j}$, and the last H\"older
seminorm is at most $C\varepsilon^2\ell_g^\sigma$.

Define the auxiliary parametrization and its real part by
\[
 \phi_0(t)=\int_0^t e^{i\theta(v)}\,dv,
 \qquad u(t)=\operatorname{Re}\phi_0(t).
\]
For sufficiently small $\varepsilon$, both have lower chord
constant $1/2$, since $u'=\cos\theta\ge1/2$.
The curve $\phi_0$ has unit speed and belongs to $C^{m+1,\alpha}$.
Equation~\eqref{cv:balanced-angle} says that the integral of
$\sin\theta$ vanishes on each side of $a_g$ within every gap.
Absolute integrability and summation over complete gaps imply
\[
 \operatorname{Im}\phi_0(x)=0\quad(x\in E),\qquad
 \operatorname{Im}\phi_0(a_g)=0\quad(g\in\mathcal G).
\]
Consequently,
\[
 \phi_0(x)=u(x)\quad(x\in E),\qquad
 \phi_0(a_g)=u(a_g),\qquad
 \phi_0'(a_g)=e^{iA_g}.
\]
Comparison of $u$ with the identity parametrization in
\eqref{cv:higher-continuity}, together with
\eqref{cv:angle-regularity}, gives
\begin{equation}
 a_n^{(p)}(H_u)\le C\varepsilon^2n^{-\sigma/s},
 \qquad n\ge1,\quad1\le p\le\infty,
 \label{cv:real-part-correction}
\end{equation}
and the same estimate on $C(E)$.

Define the finite diagonal operators
\[
 D_{\mathcal F}f
 =i\sum_{g\in\mathcal F}\tan(A_g)
                    \langle f,e_g\rangle e_g.
\]
For small $\varepsilon$, $|\tan A_g|\le2\varepsilon\ell_g^\sigma$.
The row integral of their absolute kernels is at most
\[
 C\varepsilon\sum_{g\in\mathcal F}
                   \frac{\ell_g^{\sigma+1}}{|x-a_g|}.
\]
For a band $t/2\le\ell_g<t$, centrality and the localized gap
count give
\[
 \sum_{t/2\le\ell_g<t}\frac1{|x-a_g|}
 \le Ct^{-1}\sum_{j\ge0}2^{-j(1-s)}\le Ct^{-1}.
\]
Thus this band contributes at most $C\varepsilon t^\sigma$.
The absolute kernel is symmetric, so the column estimate is the
same. Summing the bands proves
\[
 \sup_{\mathcal F}\|D_{\mathcal F}\|_{p\to p}\le C\varepsilon,
 \qquad1\le p\le\infty.
\]
The tail with $\ell_g<r$ is bounded by $C\varepsilon r^\sigma$
uniformly over subfamilies. Hence the finite sums converge in
operator norm on all these spaces and on $C(E)$ to a compact
operator $D$ with
\[
 D1=0,\qquad De_g=i\tan(A_g)e_g,\qquad
 \|D\|_{p\to p}\le C\varepsilon.
\]

On $E$, the correction functions satisfy
\[
 \begin{aligned}
 h_g^{\phi_0}(x)
 &=\frac{\sqrt{w_g}e^{iA_g}}{u(x)-u(a_g)}-e_g(x),\\
 h_g^u(x)
 &=\frac{\sqrt{w_g}\cos A_g}{u(x)-u(a_g)}-e_g(x).
 \end{aligned}
\]
Since $e^{iA}/\cos A=1+i\tan A$,
\[
 h_g^{\phi_0}=i\tan(A_g)e_g+(1+i\tan(A_g))h_g^u.
\]
Each coefficient functional is bounded on every space under
consideration. Passing from finite sums to their limits in operator
norm therefore gives
\begin{equation}
 H_{\phi_0}=D+H_u(I+D)
 \quad\text{on }L^p(\mu),\ 1\le p\le\infty,
 \text{ and on }C(E).
 \label{cv:arclength-diagonal-factorization}
\end{equation}
This argument uses convergence in operator norm, including at the
endpoints; it does not require finite expansions to be dense in
$L^\infty$.

On the $M_j$-dimensional span of the generation-$j$ functions,
$D$ is multiplication by $i\tan(\varepsilon\ell_j^\sigma)$.
For any integer $n\ge1$, choose $j$ with $n\le M_j<2n$.
An operator of rank less than $n$ has a nonzero kernel vector
in that span. Normalizing it in the space under consideration gives
\[
 a_n^{(p)}(D)
 \ge\tan(\varepsilon\ell_j^\sigma)
 \ge\varepsilon\ell_j^\sigma
 \ge c\varepsilon n^{-\sigma/s}.
\]
The vectors are continuous and bounded, so this argument includes
both endpoints and $C(E)$.
The ideal inequality and \eqref{cv:real-part-correction} give
\[
 a_n^{(p)}\bigl(H_u(I+D)\bigr)
 \le C\varepsilon^2n^{-\sigma/s}.
\]
The addition inequality $a_{2n-1}(A+B)\le a_n(A)+a_n(B)$,
applied to \eqref{cv:arclength-diagonal-factorization}, yields
\[
 a_n^{(p)}(H_{\phi_0})
 \ge a_{2n-1}^{(p)}(D)
      -a_n^{(p)}\bigl(H_u(I+D)\bigr)
 \ge(c_1\varepsilon-C_1\varepsilon^2)n^{-\sigma/s}.
\]
The upper bound follows by comparing $\phi_0$ with the identity
in \eqref{cv:higher-continuity}. Decreasing $\varepsilon_0$ gives
the following bounds, also valid on $C(E)$:
\begin{equation}
 a_n^{(p)}(H_{\phi_0})\asymp\varepsilon n^{-\sigma/s}
 \qquad(n\ge1,\quad1\le p\le\infty),
 \label{cv:auxiliary-arclength-rate}
\end{equation}

Finally, put $\delta=\varepsilon^2$ and bend the curve by
\[
 \phi_\varepsilon(t)=\int_0^t
                   e^{i(\theta(v)+\delta v)}\,dv.
\]
It has unit speed and lower chord constant $1/2$ after decreasing
$\varepsilon_0$ once more. The product rule for $C^{m,\alpha}$,
obtained from Leibniz' rule and the H\"older product inequality,
gives
\[
 \|\phi_\varepsilon'-\phi_0'\|_{C^{m,\alpha}}
 =\|e^{i\theta}(e^{i\delta t}-1)\|_{C^{m,\alpha}}
 \le C\delta.
\]
Indeed, the first factor has bounded norm by
\eqref{cv:angle-regularity}, and the second has norm at most
$C_{m,\alpha}\delta$ for $0<\delta\le1$.
The continuity estimate consequently gives
\[
 a_n^{(p)}(H_{\phi_\varepsilon}-H_{\phi_0})
 \le C\varepsilon^2n^{-\sigma/s}.
\]
Using the addition inequality again with
\eqref{cv:auxiliary-arclength-rate} yields
\[
 a_n^{(p)}(H_{\phi_\varepsilon})
 \ge a_{2n-1}^{(p)}(H_{\phi_0})
      -a_n^{(p)}(H_{\phi_\varepsilon}-H_{\phi_0})
 \ge(c_2\varepsilon-C_2\varepsilon^2)n^{-\sigma/s}.
\]
The bound $\|\phi_\varepsilon'-1\|_{C^{m,\alpha}}\le C\varepsilon$
and \eqref{cv:higher-continuity} give the upper estimate.
Choosing $\varepsilon_0$ sufficiently small for both estimates proves
\eqref{cv:arclength-sharp-rate}.

Both endpoints are accumulation points of $E$, and
\[
 \phi_\varepsilon'(0)=1,\qquad
 \phi_\varepsilon'(1)=e^{i\delta}.
\]
If the image support lay in an affine line, difference quotients
along $E$ approaching each endpoint would place both derivatives
in the one-dimensional real vector space parallel to that line. This
is impossible
for $0<\delta<\pi$.
\end{proof}

The leading term in the lower bound is the diagonal operator given by
the tangent phases in \eqref{cv:arclength-diagonal-factorization};
the smooth bend
preserves its order. The result concerns the specified transport
operator and the full parametrized arc, including its extension
through the gaps.
For these examples, with $\tau=\min(\sigma,1)$,
\[
 \mathfrak D_{\phi_\varepsilon}
 \le\frac{C\varepsilon}{1-2^{-\tau}}.
\]
Thus the examples can be chosen in the range where the basis is stable,
with $b_g=a_g$ and a common smallness threshold for all $p$.
For canonical poles, the Dini bound also gives
\[
 \max\{\|H\|_{C(E)\to C(E)},\|H\|_{L^1\to L^1}\}
 \le C\mathfrak D_\phi/m_0,
\]
where $m_0$ is the lower chord constant.
If this is less than one, $I+H$ transports the endpoint Schauder
bases of Proposition~\ref{thm:approximation-norming} to
$C(\phi(E))$ and $L^1(\mu_\phi)$.
At $p=2$, the sharp examples satisfy
\[
 H_{\phi_\varepsilon}\in S^q
 \quad\Longleftrightarrow\quad q>s/\sigma.
\]

\subsection{A geometric criterion for \texorpdfstring{$C^1$}{C1} deformations that fix the support}
\label{subsec:c1-gap-geometry}

Displacement packing gives a necessary and sufficient condition
when the deformation fixes the support pointwise. In this class, the positions of the transported poles determine compactness;
the change in tangent normalization is described by a compact diagonal multiplier.

\begin{proposition}[$C^1$ transport with fixed support]
\label{prop:c1-support-fixing}
Let $\phi\in C^1([0,1];\C)$ satisfy
\[
 \phi|_E=\operatorname{id}|_E,\qquad
 \|\phi'-1\|_\infty\le\varepsilon_*,
\]
where $0<\varepsilon_*<\min(c_0,1/2)$ is fixed.
Use the canonical parameter poles $a_g$, and put
\[
 \widetilde b_g=\phi(a_g),\qquad
 \eta_g=\frac{|\widetilde b_g-a_g|}{\ell_g},\qquad
 \lambda_g=\eta_g^2\ell_g^s.
\]
Define $H_\phi$ initially on the finite canonical span by
\[
 H_\phi1=0,\qquad
 H_\phi e_g=
 \frac{\sqrt{w_g}\,\phi'(a_g)}{x-\widetilde b_g}-e_g.
\]
Then the transported family is Bessel in $L^2(\mu)$ if and only
if $\mathcal C(\lambda)<\infty$. This is also equivalent to
the existence of a bounded extension of $H_\phi$ to one, and hence every,
interior
$L^p(\mu)$. In that case,
\begin{equation}
 \|H_\phi\|_{\mathrm{ess},L^p\to L^p}
 \asymp_p\gamma(\lambda)^{1/2},\qquad1<p<\infty.
 \label{cv:c1-essential-packing}
\end{equation}
In particular, $H_\phi$ is compact exactly when the displacement
packing has vanishing tails.

If $\omega$ is the modulus of continuity of $\phi'$ and
\begin{equation}
 \int_0^1\frac{\omega(t)^2}{t}\,dt<\infty,
 \label{cv:square-dini}
\end{equation}
then $H_\phi$ is compact on every interior $L^p$, and the
transported family together with $1$ is an unconditional basis
there. No smallness of the square Dini integral is required.
\end{proposition}

\begin{proof}
Every point of $E$ is an accumulation point of $E$. Difference
quotients along $E$ therefore give $\phi'(x)=1$ for $x\in E$.
The same conclusion uses one-sided quotients at $0$ and $1$.
For either endpoint $u_g$ of a gap,
\[
 \phi(a_g)-a_g
 =\int_{u_g}^{a_g}(\phi'(t)-1)\,dt.
\]
Thus $|\widetilde b_g-a_g|\le\varepsilon_*\ell_g$, so these
complex poles are admissible. Also
$|\phi(x)-\phi(y)|\ge(1-\varepsilon_*)|x-y|$.

Write $\theta_g=\phi'(a_g)-1$ and define the diagonal multiplier
$D_\theta1=0$, $D_\theta e_g=\theta_ge_g$.
Since $|\theta_g|\le\omega(\ell_g)$, one has
\[
 \sup_{\ell_g<r}|\theta_g|\longrightarrow0.
\]
The multiplier theorem shows that $D_\theta$ is compact on every
interior $L^p$. The operator $I+D_\theta$ is an isomorphism:
its inverse is the canonical multiplier with coefficients
$1$ on the constant and $(1+\theta_g)^{-1}$ on $e_g$.

Let $T1=0$ and
$Te_g=\sqrt{w_g}/(x-\widetilde b_g)-e_g$.
On the finite canonical span,
\begin{equation}
 I+H_\phi=(I+T)(I+D_\theta),\qquad
 H_\phi-T=D_\theta+TD_\theta.
 \label{cv:c1-diagonal-factorization}
\end{equation}
The scalar factors $1+\theta_g$ are bounded above and away from
zero, so the Bessel assertion follows from Theorem~\ref{thm:stability}.
If the packing is finite, that theorem and
\eqref{cv:c1-diagonal-factorization} give boundedness on every
interior $L^p$.

Conversely, suppose $H_\phi$ has a bounded extension on one such
$L^p$. The first identity gives a bounded extension of $T$.
For a finite set $\mathcal F$ of gaps, let $P_{\mathcal F}$ be
the corresponding canonical coordinate projection.
Its norm is uniformly bounded, while the estimate for finite arrays
in Theorem~\ref{thm:stability} gives
\[
 \mathcal C((\lambda_g\mathbf1_{\mathcal F})_g)^{1/2}
 \le C_p\|TP_{\mathcal F}\|_{p\to p}
 \le C_p'\|T\|_{p\to p}.
\]
Letting $\mathcal F$ increase proves that the packing is finite.
When $T$ is bounded, the difference in the second identity in
\eqref{cv:c1-diagonal-factorization} is compact. Consequently
$H_\phi$ and $T$ have the same essential norm, proving
\eqref{cv:c1-essential-packing}.

Finally, $\eta_g\le\omega(\ell_g)$. For a dyadic $r=2^{-N}$,
the gap count at each fixed scale, applied in every testing ball of
radius $R$, gives
\[
 \mathcal C_r(\lambda)
 \le C\sum_{j\ge N}\omega(2^{-j})^2.
\]
Indeed, the gaps with $2^{-j-1}<\ell_g\le2^{-j}$ and
$\ell_g\le R$ contribute at most
$C R^s\omega(2^{-j})^2$; the bands with $2^{-j}>R$ have
only a bounded number of admissible initial terms and obey the
same estimate after shifting $j$ by one. The discrete sum is
finite and has vanishing tails by \eqref{cv:square-dini}.
The compactness assertion follows. The basis criterion for vanishing
tails makes $I+T$ invertible on each interior $L^p$;
\eqref{cv:c1-diagonal-factorization} then proves the last assertion.
\end{proof}

For this class, displacement packing is the exact condition. It
records the distribution of oscillation among the gaps, which a global
modulus does not determine.

The same boundedness and compactness equivalences hold on
$L^p(v\,d\mu)$ for $v\in A_p(\mu)$. Indeed,
Proposition~\ref{prop:weighted-packing} makes $D_\theta$ compact
on these spaces as well. The factorization
\eqref{cv:c1-diagonal-factorization} transfers the weighted
packing criterion and the weighted bounds for the essential norm to
$H_\phi$. Vanishing tail packing also gives weighted basis transport.

\begin{corollary}[Schatten criterion for pole position and tangent]
\label{cor:c1-joint-schatten}
Under the hypotheses of Proposition~\ref{prop:c1-support-fixing}, put
\[
 \zeta_g=
 \left(
   |\phi'(a_g)-1|^2+
   \frac{|\phi(a_g)-a_g|^2}{\ell_g^2}
 \right)^{1/2}.
\]
For every $0<r<\infty$, the formal correction $H_\phi$ extends
to $S^r(L^2(\mu))$ if and only if $(\zeta_g)\in\ell^r(\G)$,
and
\[
 \|H_\phi\|_{S^r}\asymp_r
       \left(\sum_g\zeta_g^r\right)^{1/r}.
\]
The constants are uniform under the stated geometric and
admissibility bounds.
\end{corollary}

\begin{proof}
Write
\[
 \delta_g=\phi(a_g)-a_g,\qquad
 \theta_g=\phi'(a_g)-1,\qquad
 b_g=\phi(a_g).
\]
These complex poles are admissible by
Proposition~\ref{prop:c1-support-fixing}. Since $\phi$ fixes $E$
pointwise,
\[
 H_\phi e_g(x)=
 \frac{\sqrt{w_g}\,[\theta_g(x-a_g)+\delta_g]}
      {(x-a_g)(x-b_g)}.
\]
Apply Theorem~\ref{thm:geometric-schatten} with
\[
 \alpha_g=\theta_g,\qquad \beta_g=\delta_g/\ell_g.
\]
Its coefficient size is exactly $\zeta_g$. The equivalence and
the two norm estimates follow for every $0<r<\infty$.
\end{proof}

\begin{proposition}[$C^1$ realization of prescribed gap displacements]
\label{prop:c1-gap-realization}
There are $c_1,C>0$, depending only on the fixed geometric data, such
that every real array $(\varepsilon_g)$ satisfying
\[
 \sup_g|\varepsilon_g|\le c_1,\qquad
 \lim_{r\downarrow0}\sup_{\ell_g<r}|\varepsilon_g|=0
\]
is realized by an increasing $C^1$ diffeomorphism
$\phi:[0,1]\to[0,1]$ with
\[
 \begin{gathered}
 \phi|_E=\operatorname{id}|_E,\qquad \phi'|_E=1,
 \qquad \|\phi-\operatorname{id}\|_{C^1}
       \le C\sup_g|\varepsilon_g|,\\
 \tfrac12\le\phi'\le\tfrac32,
 \qquad \phi(a_g)=a_g+\varepsilon_g\ell_g,
 \qquad \phi'(a_g)=1.
 \end{gathered}
\]
Its formal transport correction is exactly the operator for pole
perturbations with $b_g=a_g+\varepsilon_g\ell_g$.
\end{proposition}

\begin{proof}
Choose a real $\chi\in C_c^\infty((-c_*/2,c_*/2))$ such that
$\chi(0)=1$ and $\chi'(0)=0$. Set
\[
 u(t)=
 \begin{cases}
 \varepsilon_g\ell_g\chi((t-a_g)/\ell_g),&t\in g,\\
 0,&t\in E,
 \end{cases}
 \qquad \phi(t)=t+u(t).
\]
Uniform centrality places every bump strictly inside its gap.
On its support, $\operatorname{dist}(t,E)\ge(c_*/2)\ell_g$.
The piecewise candidate for the derivative is
\[
 v(t)=
 \begin{cases}
 \varepsilon_g\chi'((t-a_g)/\ell_g),&t\in g,\\
 0,&t\in E.
 \end{cases}
\]
It is continuous. To check this at $E$, the finitely many gaps
of length at least $r$ have bump supports at positive distance
from $E$, while all remaining derivative values are bounded by
$\|\chi'\|_\infty\sup_{\ell_g<r}|\varepsilon_g|$, which tends to zero.

For $x\in E$ and $t\to x$ in a bump support,
\[
 \frac{|u(t)-u(x)|}{|t-x|}
 \le C|\varepsilon_g|\longrightarrow0.
\]
The limit follows by excluding the finitely many large gaps as above;
outside the bump supports the quotient is zero. Hence $u'=v$ on all of
$[0,1]$, with one-sided derivatives at the endpoints.
This proves $u\in C^1$ and the claimed norm bound.
Choose $c_1$ so that
$c_1\|\chi'\|_\infty\le1/2$ and $Cc_1\le\varepsilon_*$.
Then $\phi'$ lies between $1/2$ and $3/2$; since $\phi$ fixes
$0$ and $1$, it is the asserted diffeomorphism. Its values and
derivatives at $a_g$ follow from the two normalizations of $\chi$.
For $x\in E$,
\[
 \frac{\sqrt{w_g}\,\phi'(a_g)}{\phi(x)-\phi(a_g)}
 =\frac{\sqrt{w_g}}{x-a_g-\varepsilon_g\ell_g},
\]
which proves the exact operator identity on the canonical span.
\end{proof}

\begin{corollary}[The $C^1$ compactness trichotomy]
\label{cor:c1-transport-trichotomy}
On the middle-thirds Cantor probability measure, each of the
following occurs for increasing diffeomorphisms arbitrarily close
to the identity in $C^1$, fixing $E$ pointwise and satisfying
$\phi'|_E=\phi'(a_g)=1$:
\begin{enumerate}
 \item The correction is compact on every interior $L^p$ and
       the transported family is an unconditional basis there.
 \item The correction is bounded and noncompact on every
       interior $L^p$, while the transported family is a Riesz
       basis on $L^2$.
 \item The transported family is not Bessel on $L^2$.
\end{enumerate}
For displacements of the form
$\varepsilon_g=a n^{-\beta}$ in generation $n$, with $a>0$ sufficiently
small and $\beta>0$, the family is Bessel exactly when
$\beta>1/2$. In that range the first alternative holds.
The derivative is Dini continuous exactly when $\beta>1$.
For these examples, boundedness of the correction
on $L^1(\mu)$ holds exactly when $\beta>1$. More generally, suppose
$0\le\eta_n\le c_1$ and $\eta_n\to0$, where $c_1$ is the constant in
Proposition~\ref{prop:c1-gap-realization}. For the deformation constructed
there with $\varepsilon_g=\eta_n$ on every gap of generation $n$,
\begin{equation}
 \begin{gathered}
 H_\phi:L^1\to L^1\text{ is bounded}
 \quad\Longleftrightarrow\quad
 \sum_{n\ge1}\eta_n<\infty,\\
 \|H_\phi\|_{1\to1}\asymp\sum_{n\ge1}\eta_n
 \quad\text{when finite}.
 \end{gathered}
 \label{cv:c1-endpoint-threshold}
\end{equation}
\end{corollary}

\begin{proof}
Proposition~\ref{prop:c1-gap-realization} realizes each of the
generation sequences above. Corollary~\ref{cor:generation-classification}
gives the exact condition
\[
 \sum_{n\ge1}a^2n^{-2\beta}<\infty
 \quad\Longleftrightarrow\quad\beta>1/2.
\]
It also gives compactness and every interior basis assertion
when this sum is finite. Taking $0<\beta\le1/2$ proves the
third alternative. The $C^1$ distance tends to zero with $a$.

For the Dini assertion, let $\omega$ be the modulus of $\phi'$.
The fixed bump profile gives
\[
 c_\beta aN^{-\beta}\le\omega(3^{-N})
 \le C_\beta aN^{-\beta},\qquad N\ge1.
\]
For the lower bound, compare a point where $\chi'$ is nonzero
with a point of $E$ on the boundary of the same gap of generation $N$;
their distance is at most $3^{-N}$.
For the upper bound, a derivative bump in generation $n$ has
amplitude at most $Ca n^{-\beta}$ and Lipschitz constant at
most $Ca n^{-\beta}3^n$.
If the two points are in different supports, split the difference
at the intervening zero values. Thus
\[
 \omega(3^{-N})
 \le Ca\sup_{n\ge1}
     n^{-\beta}\min(1,3^{n-N})
 \le C_\beta aN^{-\beta}.
\]
The triadic Dini sum now converges exactly when $\beta>1$.
In particular, the interval $1/2<\beta\le1$ gives compact
transport beyond the Dini hypothesis.

We prove \eqref{cv:c1-endpoint-threshold}. Let $P_N$ retain the
constant and all gaps of generations at most $N$, and put
$T_N=H_\phi P_N$ on the canonical span. Proposition~\ref{thm:approximation-norming}
gives $\sup_N\|P_N\|_{1\to1}<\infty$. The adjoint satisfies
\[
 T_N^*1(x)=\sum_{\operatorname{gen}(g)\le N}
                   \frac{w_gF(b_g)}{x-a_g},
 \qquad b_g=a_g+\eta_n\ell_g.
\]
The numerator is nonnegative since $F(a_g)=0$ and $F$ is
increasing in each real gap. On the leftmost generation-$N$
cylinder $Q_N=E\cap[0,3^{-N}]$, every denominator is negative.
For the leftmost generation-$n$ gap
$g_n=(3^{-n},2\cdot3^{-n})$, the derivative estimate in the
central part of the gap gives
\[
 w_{g_n}F(b_{g_n})\ge c\eta_n\ell_{g_n},\qquad
 0<a_{g_n}-x\le2\ell_{g_n}\quad(x\in Q_N,\ n\le N).
\]
Consequently
\[
 T_N^*1(x)\le-c\sum_{n\le N}\eta_n\quad(x\in Q_N),
 \qquad \|T_N\|_{1\to1}\ge c\sum_{n\le N}\eta_n.
\]
The cylinder has positive measure. If $H_\phi$ were bounded on
$L^1$, the uniform projection bound would force the partial sums
on the right to remain bounded. This proves necessity and the
lower norm estimate. Conversely, the absolute endpoint Schur
estimate in Proposition~\ref{prop:endpoint-pole-stability} bounds
$\|H_\phi\|_{1\to1}$ by $C\sum_n\eta_n$; dyadic length bands
and Cantor generations give the same summability condition.
That proposition also gives endpoint invertibility when the sum
is finite. The stated $\beta$ threshold follows.

For the second alternative, label the two Cantor children by
$0$ and $1$, and let $J_k$ be the cylinder with address
$1^{k-1}0$. These cylinders are pairwise disjoint and have
length $L_k=3^{-k}$ and measure $2^{-k}=L_k^s$.
Choose integers $M_k\to\infty$. In each $J_k$, take the
first $M_k$ generations of descendant gaps as the active gaps, and set
\[
 \varepsilon_g=\frac{a}{\sqrt{M_k}}
 \quad\text{on those gaps},\qquad
 \varepsilon_g=0\quad\text{elsewhere}.
\]
Only finitely many active gaps occur in any fixed finite set of
these blocks. It follows that
$\sup_{\ell_g<r}|\varepsilon_g|\to0$, so the realization proposition
applies.

At each active generation inside $J_k$, the packing mass is
$a^2\mu(J_k)/(2M_k)$. Hence
\begin{equation}
 \sum_{g\text{ active in }J_k}\lambda_g
 =\frac{a^2}{2}\mu(J_k).
 \label{cv:c1-block-energy}
\end{equation}
More generally, for every Cantor cylinder $I$,
\[
 \sum_{g\subset I}\lambda_g\le\frac{a^2}{2}\mu(I).
\]
If $I\subset J_k$, at most $M_k$ active generations relative to that cylinder
contribute, each with mass at most
$a^2\mu(I)/(2M_k)$. Otherwise the active roots contained in
$I$ are disjoint, and \eqref{cv:c1-block-energy} gives the
same bound after summing their measures.

This cylinder estimate implies
$\mathcal C(\lambda)\le Ca^2$.
Indeed, choose $n$ with $3^{-n}\le R<3^{1-n}$.
A testing ball of radius $R$ meets only a bounded number of
generation-$n$ cylinders. Gaps of smaller length with an
endpoint in the ball lie inside these cylinders; the generation-$n$
gaps contribute at most $Ca^2R^s$ by the count at a fixed scale.
Applying the cylinder estimate proves the asserted ball bound.

For every $r>0$, choose $k$ so large that $L_k<r$.
A ball of radius $L_k$ centered at an endpoint of $J_k$ contains
all active gap endpoints in $J_k$, and all the gaps in
\eqref{cv:c1-block-energy} have length less than $r$.
Consequently
\[
 \mathcal C_r(\lambda)\ge a^2/2,
 \qquad a^2/2\le\gamma(\lambda)\le\mathcal C(\lambda)
 \le Ca^2.
\]
Theorem~\ref{thm:stability} gives boundedness and noncompactness
on every interior $L^p$. It also gives $\|H_\phi\|_{2\to2}\le Ca$.
For sufficiently small $a$, the Neumann series makes $I+H_\phi$
invertible on $L^2$, proving the Riesz basis assertion.
\end{proof}

The block construction separates pointwise decay from vanishing
Carleson packing: each individual displacement tends to zero,
but arbitrarily small cylinders retain a fixed normalized
quadratic mass. It answers the $C^1$ compactness question
negatively even when the transported family remains a Riesz basis.

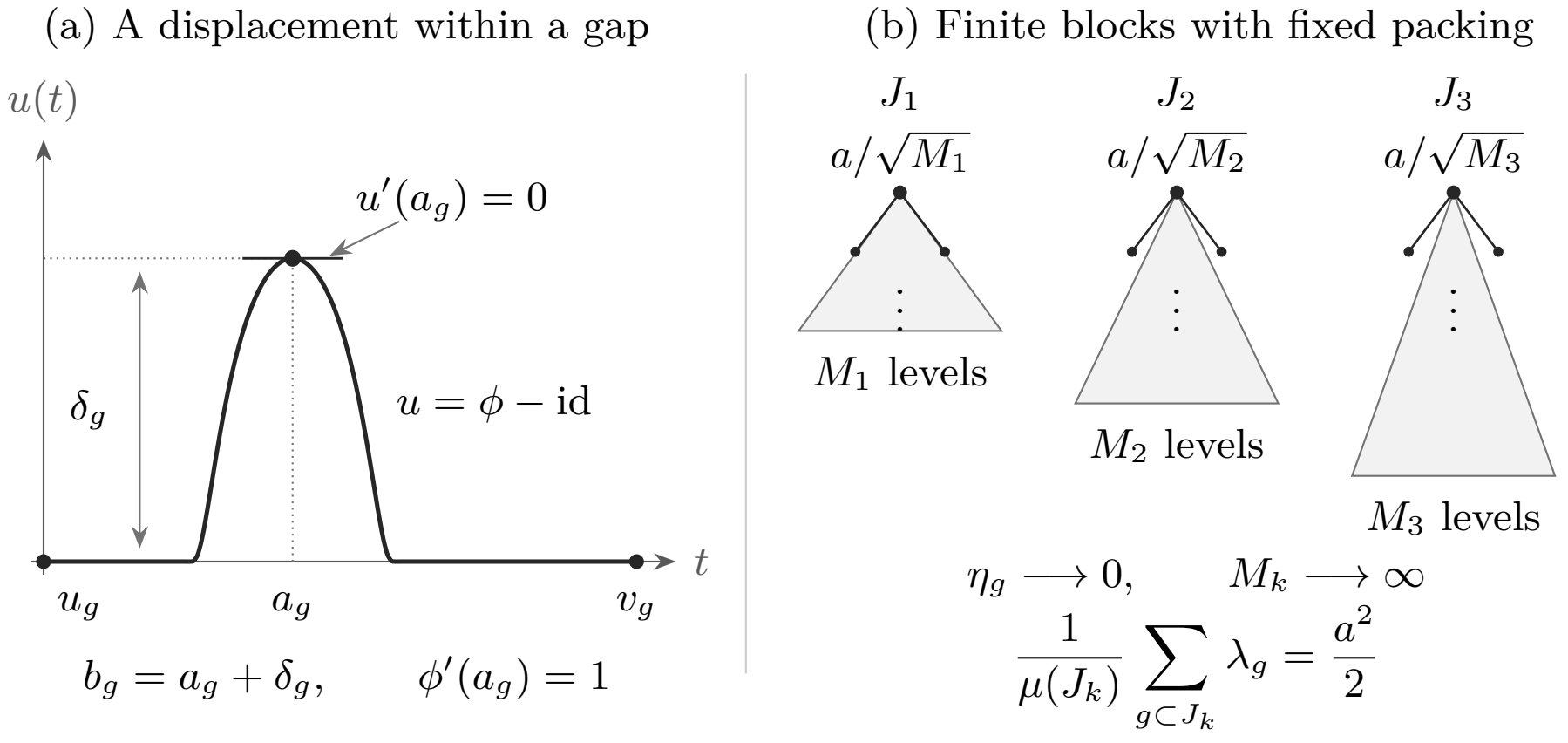
\begin{figure}[htbp]
\centering
\begin{tikzpicture}[x=1cm,y=1cm,font=\footnotesize,>=Stealth]
 \node at (2.75,4.95) {(a) A displacement within a gap};
 \draw[->,black!65] (.35,.90)--(5.25,.90) node[right] {$t$};
 \draw[->,black!65] (.45,.80)--(.45,4.10) node[above] {$u(t)$};
 \draw[DiagramPrimary,line width=.95pt] plot[smooth] coordinates {(0.45000,0.90000) (0.48000,0.90000) (0.51000,0.90000) (0.54000,0.90000) (0.57000,0.90000) (0.60000,0.90000) (0.63000,0.90000) (0.66000,0.90000) (0.69000,0.90000) (0.72000,0.90000) (0.75000,0.90000) (0.78000,0.90000) (0.81000,0.90000) (0.84000,0.90000) (0.87000,0.90000) (0.90000,0.90000) (0.93000,0.90000) (0.96000,0.90000) (0.99000,0.90000) (1.02000,0.90000) (1.05000,0.90000) (1.08000,0.90000) (1.11000,0.90000) (1.14000,0.90000) (1.17000,0.90000) (1.20000,0.90000) (1.23000,0.90000) (1.26000,0.90000) (1.29000,0.90000) (1.32000,0.90000) (1.35000,0.90000) (1.38000,0.90000) (1.41000,0.90000) (1.44000,0.90000) (1.47000,0.90000) (1.50000,0.90000) (1.53000,0.90000) (1.56000,0.90001) (1.59000,0.90565) (1.62000,0.95330) (1.65000,1.06331) (1.68000,1.21899) (1.71000,1.39741) (1.74000,1.58188) (1.77000,1.76237) (1.80000,1.93346) (1.83000,2.09255) (1.86000,2.23868) (1.89000,2.37179) (1.92000,2.49230) (1.95000,2.60084) (1.98000,2.69814) (2.01000,2.78494) (2.04000,2.86193) (2.07000,2.92974) (2.10000,2.98896) (2.13000,3.04010) (2.16000,3.08357) (2.19000,3.11976) (2.22000,3.14896) (2.25000,3.17143) (2.28000,3.18735) (2.31000,3.19684) (2.34000,3.20000) (2.37000,3.19684) (2.40000,3.18735) (2.43000,3.17143) (2.46000,3.14896) (2.49000,3.11976) (2.52000,3.08357) (2.55000,3.04010) (2.58000,2.98896) (2.61000,2.92974) (2.64000,2.86193) (2.67000,2.78494) (2.70000,2.69814) (2.73000,2.60084) (2.76000,2.49230) (2.79000,2.37179) (2.82000,2.23868) (2.85000,2.09255) (2.88000,1.93346) (2.91000,1.76237) (2.94000,1.58188) (2.97000,1.39741) (3.00000,1.21899) (3.03000,1.06331) (3.06000,0.95330) (3.09000,0.90565) (3.12000,0.90001) (3.15000,0.90000) (3.18000,0.90000) (3.21000,0.90000) (3.24000,0.90000) (3.27000,0.90000) (3.30000,0.90000) (3.33000,0.90000) (3.36000,0.90000) (3.39000,0.90000) (3.42000,0.90000) (3.45000,0.90000) (3.48000,0.90000) (3.51000,0.90000) (3.54000,0.90000) (3.57000,0.90000) (3.60000,0.90000) (3.63000,0.90000) (3.66000,0.90000) (3.69000,0.90000) (3.72000,0.90000) (3.75000,0.90000) (3.78000,0.90000) (3.81000,0.90000) (3.84000,0.90000) (3.87000,0.90000) (3.90000,0.90000) (3.93000,0.90000) (3.96000,0.90000) (3.99000,0.90000) (4.02000,0.90000) (4.05000,0.90000) (4.08000,0.90000) (4.11000,0.90000) (4.14000,0.90000) (4.17000,0.90000) (4.20000,0.90000) (4.23000,0.90000) (4.26000,0.90000) (4.29000,0.90000) (4.32000,0.90000) (4.35000,0.90000) (4.38000,0.90000) (4.41000,0.90000) (4.44000,0.90000) (4.47000,0.90000) (4.50000,0.90000) (4.53000,0.90000) (4.56000,0.90000) (4.59000,0.90000) (4.62000,0.90000) (4.65000,0.90000) (4.68000,0.90000) (4.71000,0.90000) (4.74000,0.90000) (4.77000,0.90000) (4.80000,0.90000) (4.83000,0.90000) (4.86000,0.90000) (4.89000,0.90000) (4.92000,0.90000) (4.95000,0.90000)};
 \fill[DiagramPrimary] (.45,.90) circle (1.6pt);
 \fill[DiagramPrimary] (4.95,.90) circle (1.6pt);
 \draw[DiagramSecondary,densely dotted] (2.34,.90)--(2.34,3.20);
 \draw[DiagramSecondary,densely dotted] (.45,3.20)--(2.34,3.20);
 \fill[DiagramPrimary] (2.34,3.20) circle (1.8pt);
 \draw[DiagramPrimary,line width=.55pt] (1.96,3.20)--(2.72,3.20);
 \draw[<->,DiagramSecondary] (1.18,1.00)--(1.18,3.10);
 \node[anchor=east] at (1.08,2.05) {$\delta_g$};
 \node[anchor=north] at (2.34,.80) {$a_g$};
 \node[anchor=north west] at (.42,.80) {$u_g$};
 \node[anchor=north] at (4.95,.80) {$v_g$};
 \node at (3.55,3.65) {$u'(a_g)=0$};
 \draw[->,DiagramSecondary] (3.15,3.48)--(2.62,3.22);
 \node at (3.87,2.10) {$u=\phi-\operatorname{id}$};
 \node at (2.75,.02) {$b_g=a_g+\delta_g,\qquad \phi'(a_g)=1$};
 \draw[black!20] (5.78,.05)--(5.78,4.6);
 \node at (9.22,4.95) {(b) Finite blocks with fixed packing};
 \foreach \x/\dep/\idx in {6.95/1.05/1,9.05/1.6/2,11.15/2.15/3}{
   \pgfmathsetmacro{\basey}{3.7-\dep}
   \filldraw[draw=DiagramSecondary,fill=DiagramPrimary!6,line width=.4pt]
       (\x,3.7)--(\x-.77,\basey)--(\x+.77,\basey)--cycle;
   \draw[DiagramPrimary,line width=.5pt]
       (\x,3.7)--(\x-.34,3.25)
       (\x,3.7)--(\x+.34,3.25);
   \fill[DiagramPrimary] (\x,3.7) circle (1.5pt);
   \fill[DiagramPrimary] (\x-.34,3.25) circle (1.1pt);
   \fill[DiagramPrimary] (\x+.34,3.25) circle (1.1pt);
   \node at (\x,4.45) {$J_{\idx}$};
   \node at (\x,3.99) {$a/\sqrt{M_{\idx}}$};
   \node at (\x,2.91) {$\vdots$};
   \node[anchor=north] at (\x,\basey-.06) {$M_{\idx}$ levels};
 }
 \node at (9.2,.78) {$\eta_g\longrightarrow0,\qquad M_k\longrightarrow\infty$};
 \node at (9.2,.08)
 {$\displaystyle\frac{1}{\mu(J_k)}
       \sum_{g\subset J_k}\lambda_g=\frac{a^2}{2}$};
\end{tikzpicture}
\caption{Two parts of the construction that fixes the support.
Left: the displacement $u=\phi-\operatorname{id}$ is a smooth central
bump with $u(a_g)=\delta_g$ and $u'(a_g)=0$. It moves the canonical
pole to $b_g$ while preserving its tangent; the deformation
fixes the entire support $E$.
Right: a schematic arrangement of disjoint Cantor cylinders, each
carrying a finite active tree. Its relative displacement is $a/\sqrt{M_k}$ on
every active gap and its normalized quadratic packing is $a^2/2$.
The depths increase and the individual displacements vanish, but the
packing persists in arbitrarily small cylinders. Tree widths and bump
height are schematic.}
\label{fig:c1-gap-blocks}
\end{figure}

Figure~\ref{fig:c1-gap-blocks} illustrates the gap realization and the
construction using finite trees.

\begin{remark}[Noncollinear support]
The three alternatives also occur on a fixed noncollinear circular
image of $E$. Put $M(z)=z/(1+iz)$ and postcompose the preceding
maps with $M$. Although these maps need not have Dini continuous
derivatives, the algebraic
identity in the proof of Proposition~\ref{prop:mobius-covariance}
still applies to every finite canonical sum and gives
\[
 H_{M\circ\phi}=H_\phi+N,\qquad
 Nf=1\sum_g\langle f,e_g\rangle
                  \sqrt{w_g}\frac{i}{1+i\phi(a_g)}.
\]
The estimate at a fixed scale used in that proposition gives
\[
 \sup_{y\in E}\sum_{u/2<\ell_g\le u}
                  \frac{w_g}{|y-a_g|}\le Cu^{1-s}.
\]
Since $|i/(1+i\phi(a_g))|\le1$, the kernel of this functional converges uniformly on $E$. Thus $N$ is bounded on
every $L^p$, including both endpoints, and on $C(E)$; it has
$N1=0$ and $N^2=0$. Boundedness and compactness on the interior
$L^p$ spaces, and failure of boundedness on $L^2$, are unchanged
by this addition of rank one.
In the bounded case,
$I+H_{M\circ\phi}=(I+H_\phi)(I+N)$, so the Riesz basis
property is preserved as well. The image measure is always
$M_\#\mu$, because $\phi$ fixes $E$. The image $M(E)$ is not
contained in a line, and the lower chord bound follows directly
from the difference formula for M\"obius maps. These maps are
arbitrarily $C^1$-close to $M$. Choosing $M_\delta(z)=z/(1+i\delta z)$
with $\delta>0$ small also gives noncollinear examples arbitrarily
close to the identity.
\end{remark}

\subsection{Sharpness for a prescribed modulus}
\label{subsec:prescribed-modulus-sharpness}

\begin{proposition}[Sharpness for a prescribed modulus]
\label{prop:sharp-modulus-transport}
Let $\omega:[0,1]\to[0,1]$ be continuous, nondecreasing, and
concave, with $\omega(0)=0$ and $\omega(1)=1$.
There is an absolute $\varepsilon_0>0$ with the following property.
For every $0<\varepsilon\le\varepsilon_0$, the middle-thirds
Cantor measure admits an increasing $C^1$ diffeomorphism
$\phi:[0,1]\to[0,1]$ satisfying
\[
 \begin{gathered}
 \phi|_E=\operatorname{id}|_E,\qquad
 \tfrac12\le\phi'\le\tfrac32,\qquad
 \|\phi-\operatorname{id}\|_{C^1}\le C\varepsilon,\\
 |\phi'(x)-\phi'(y)|\le C\varepsilon\omega(|x-y|),
 \qquad \phi'(a_g)=1,
 \end{gathered}
\]
whose correction with canonical parameter poles has the exact
criteria
\[
 \begin{split}
 H_\phi\text{ is bounded on an interior }L^p
 &\quad\Longleftrightarrow\quad
      \int_0^1\frac{\omega(t)^2}{t}\,dt<\infty,\\
 H_\phi\text{ is bounded on }L^1
 &\quad\Longleftrightarrow\quad
      \int_0^1\frac{\omega(t)}{t}\,dt<\infty.
 \end{split}
\]
In the first case $H_\phi$ is compact on every interior $L^p$,
and the transported family is an unconditional basis on each of
these spaces. In the second case it is also compact on $L^1$.
The same two boundedness criteria can be realized by arbitrarily
small $C^1$ deformations whose image support is not contained in
an affine line.
\end{proposition}

\begin{proof}
Choose the fixed real bump $\chi$ from
Proposition~\ref{prop:c1-gap-realization}, with $\chi(0)=1$ and
$\chi'(0)=0$. On a gap $g$, put
\[
 u(t)=\varepsilon\omega(\ell_g)\ell_g
                 \chi((t-a_g)/\ell_g),
 \qquad u|_E=0,\qquad \phi(t)=t+u(t).
\]
Since $\omega(\ell_g)\to0$ as $\ell_g\to0$, the gluing proof in
that proposition gives $u\in C^1$, $u'|_E=0$, and
$\|u\|_{C^1}\le C\varepsilon$. Decrease
$\varepsilon_0$ to obtain the derivative bounds and admissibility.
The pole shifts and tangents are
\[
 \phi(a_g)=a_g+\varepsilon\omega(\ell_g)\ell_g,
 \qquad \phi'(a_g)=1.
\]
Thus the transport correction is exactly the pole correction for
these positive shifts.

For $0<h\le\ell\le1$, concavity gives
\[
 \frac h\ell\omega(\ell)\le\omega(h).
\]
The derivative bump on a gap of length $\ell$ has size at most
$C\varepsilon\omega(\ell)$ and Lipschitz constant at most
$C\varepsilon\omega(\ell)/\ell$. Hence its oscillation between
two points a distance $h$ apart is bounded by
\[
 C\varepsilon\omega(\ell)\min(1,h/\ell)
       \le C\varepsilon\omega(h).
\]
This proves the asserted modulus within one gap. For points in
different gaps, split the difference at intervening endpoints,
where $u'=0$; both distances to those endpoints are at most
$|x-y|$. The same estimate, with at most twice the constant,
handles those points and points of $E$.

A generation-$n$ gap has length $3^{-n}$, and there are
$2^{n-1}$ such gaps. Its relative shift is
$\eta_n=\varepsilon\omega(3^{-n})$. For
$s=\log2/\log3$, the packing mass of that generation is exactly
\[
 \sum_{\operatorname{gen}(g)=n}\lambda_g
   =2^{n-1}\varepsilon^2\omega(3^{-n})^2\,3^{-ns}
   =\tfrac12\varepsilon^2\omega(3^{-n})^2.
\]
If the square Dini integral diverges, the total packing is infinite,
so Theorem~\ref{thm:stability} excludes boundedness on every
interior $L^p$. If it converges, the gap count at a fixed scale in every testing ball
of radius $R$ bounds the generation-$n$ packing mass there by
$C\varepsilon^2R^s\omega(3^{-n})^2$ whenever $3^{-n}\le R$.
Consequently
\[
 \mathcal C_r(\lambda)
 \le C\varepsilon^2\sum_{3^{-n}<r}\omega(3^{-n})^2
 \longrightarrow0.
\]
Theorem~\ref{thm:stability} gives compactness and the interior
basis assertions.
Monotonicity of $\omega$ identifies the triadic square sum with
the square Dini integral.

The criterion for positive displacements by generation
\eqref{cv:c1-endpoint-threshold} applies to the sequence
$\eta_n=\varepsilon\omega(3^{-n})$. It gives boundedness on $L^1$
exactly when $\sum_n\omega(3^{-n})<\infty$, which is equivalent
to the Dini integral. The endpoint Schur estimate also gives
convergence in operator norm of the tails in that case, and hence
compactness on $L^1$.

For a noncollinear image support, let
\[
 M_\varepsilon(z)=\frac{z}{1-i\varepsilon z},
 \qquad \psi=M_\varepsilon\circ\phi.
\]
The image $M_\varepsilon(E)$ lies on a nondegenerate circular arc
and is not contained in any affine line. Direct differentiation,
$\omega(h)\ge h$, and the bounds for $\phi$ give
\[
 \|\psi-\operatorname{id}\|_{C^1}\le C\varepsilon,
 \qquad
 |\psi'(x)-\psi'(y)|\le C\varepsilon\omega(|x-y|).
\]
The lower chord bound is also uniform. The algebraic identity
in \eqref{cv:mobius-rank-one} gives, on finite canonical expansions,
\[
 H_\psi=H_\phi+N,\qquad
 Nf=1\sum_g\langle f,e_g\rangle\sqrt{w_g}
                    \frac{-i\varepsilon\phi'(a_g)}
                         {1-i\varepsilon\phi(a_g)}.
\]
This identity needs no Dini assumption. Its coefficient factor is
uniformly bounded. The global potential estimate
\eqref{cv:global-gap-potential} makes the defining integral kernel
of $N$ bounded, and its finite sums converge in operator norm
on every $L^p$, including $p=1$ and $p=\infty$.
Thus $N$ has rank at most one, and adding it preserves all
boundedness and compactness conclusions above.
\end{proof}

For a prescribed concave modulus, the square Dini integral is
therefore exactly the threshold for a uniform interior compactness
guarantee: finiteness gives compactness for every deformation whose derivative modulus is bounded by a constant times
$\omega$ and which has a lower chord bound; divergence allows
an arbitrarily small deformation with no bounded interior extension.
For $L^1$, the corresponding threshold is the Dini integral.
The global modulus conditions are sufficient conditions for
individual maps. The local oscillation packing can be finite when
a chosen global modulus fails either summability condition.

\section{Reconstruction, variation, and approximation}
\label{subsec:reconstruction-proof}
\begin{proof}[Proof of Proposition~\ref{thm:approximation-norming}]
Write \(g=(u_g,v_g)\) and \(n_r=\#\{g:\ell_g\ge r\}\).
The gap count \eqref{eq:prelim-gapcount} and
\eqref{sp:gap-bounds} will be used throughout.
They first give the following potential estimate for a single band:
\begin{equation}
 \sup_{x\in E}\sum_{t\le\ell_g<2t}
       \frac{w_g}{|x-a_g|}
 \le Ct^{1-s}.
 \label{eq:reconstruction-band-potential}
\end{equation}
Indeed, \(|x-a_g|\ge ct\), and the annulus
\(|x-a_g|\asymp2^mt\) contains at most \(C2^{ms}\) poles.
Its contribution is at most
\[
 C2^{ms}\frac{t^{2-s}}{2^mt}
 =Ct^{1-s}2^{-m(1-s)}.
\]
The sum converges because \(s<1\). The same estimates hold for
half-open bands with the opposite endpoint convention.

Let \(0<\alpha\le1\) and \(f\in C^\alpha(E)\).
Cancellation and the location of \(a_g\) give
\[
 \begin{aligned}
 \langle f,e_g\rangle
 &=\sqrt{w_g}\int_E\frac{f(t)-f(u_g)}{t-a_g}\,d\mu(t),\\
 |\langle f,e_g\rangle|
 &\le C\sqrt{w_g}[f]_{C^\alpha(E)}
               \int_E|t-a_g|^{\alpha-1}\,d\mu(t).
 \end{aligned}
\]
There is no mass at distances smaller than \(c_*\ell_g\).
Upper regularity on the annuli
\(|t-a_g|\asymp2^j\ell_g\) yields
\[
 \int_E|t-a_g|^{\alpha-1}\,d\mu(t)
 \le C_{\alpha,s}
 \begin{cases}
  \ell_g^{\alpha+s-1},&\alpha<1-s,\\
  \log(e/\ell_g),&\alpha=1-s,\\
  1,&\alpha>1-s.
 \end{cases}
\]
The three cases are the geometric sum of
\((2^j\ell_g)^{\alpha+s-1}\) from scale \(\ell_g\) to scale one.
Multiplication by \(|e_g(x)|\), followed by
\eqref{eq:reconstruction-band-potential}, therefore gives
\[
 \sup_{x\in E}\sum_{t\le\ell_g<2t}
       |\langle f,e_g\rangle e_g(x)|
 \le C_{\alpha,s}[f]_{C^\alpha(E)}
 \begin{cases}
  t^\alpha,&\alpha<1-s,\\
  t^{1-s}\log(e/t),&\alpha=1-s,\\
  t^{1-s},&\alpha>1-s.
 \end{cases}
\]
Summing bands below \(r\) gives
\begin{equation}
 \sup_{x\in E}\sum_{\ell_g<r}
       |\langle f,e_g\rangle e_g(x)|
 \le C_{\alpha,s}[f]_{C^\alpha(E)}\Psi_{\alpha,s}(r).
 \label{rec:holder-tail}
\end{equation}
For the critical sum, with \(q=2^{-(1-s)}\),
\[
 \sum_{j\ge0}(2^{-j}r)^{1-s}\log\frac{e}{2^{-j}r}
 =r^{1-s}\left(
       \frac{\log(e/r)}{1-q}+\frac{q\log2}{(1-q)^2}\right).
\]
The other two sums are geometric.

Thus the orthonormal expansion of \(f\) converges absolutely and
uniformly. Its sum equals \(f\) in \(L^2(\mu)\) by
Proposition~\ref{thm:basis}; continuity and full support of \(\mu\)
give equality everywhere on \(E\).
For the coordinate function, cancellation gives
\(\langle t,e_g\rangle=\sqrt{w_g}\), so
\[
 x=m_1+\sum_g\frac{w_g}{x-a_g},\qquad x\in E,
 \qquad m_1=\int_Et\,d\mu(t),
\]
with absolute uniform convergence.
Subtracting this identity at \(x\ne y\) gives
\[
 1=-\sum_g\frac{w_g}{(x-a_g)(y-a_g)}.
\]
The series of divided differences is absolutely convergent, since its
absolute sum is at most
\[
 \frac1{|x-y|}\sum_gw_g
       \left(\frac1{|x-a_g|}+\frac1{|y-a_g|}\right)<\infty.
\]
Consequently, the finite projection kernel
\[
 J_r(x,y)=1+\sum_{\ell_g\ge r}
              \frac{w_g}{(x-a_g)(y-a_g)}
\]
has the alternative representation
\begin{equation}
 J_r(x,y)=-\sum_{\ell_g<r}
              \frac{w_g}{(x-a_g)(y-a_g)},\qquad x\ne y.
 \label{eq:projection-tail-identity}
\end{equation}

We record an absolute estimate for a single band:
\begin{equation}
 \sum_{t\le\ell_g<2t}|e_g(x)e_g(y)|
 \le Ct^{-s}\left(1+\frac{|x-y|}{t}\right)^{-1}.
 \label{rec:absolute-band-kernel}
\end{equation}
Put \(d=|x-y|\).
If \(d\le t\), centrality gives
\(|x-a_g|\asymp|y-a_g|\) in the band. The annular count gives
\[
 \sum_{t\le\ell_g<2t}|e_g(x)e_g(y)|
 \le Ct^{2-s}\sum_{m\ge0}
           (2^mt)^{-2}2^{ms}
 \le Ct^{-s}.
\]
If \(d>t\), split the poles into the regions
\(|x-a_g|<d/3\), \(|y-a_g|<d/3\), and their complement.
The first region contributes at most
\[
 \frac C d\sum_{t\le\ell_g<2t}\frac{w_g}{|x-a_g|}
 \le C\frac{t^{1-s}}d;
\]
the second has the same bound. On the complement both distances
are at least \(d/3\) and are comparable. Its contribution is at most
\[
 Ct^{2-s}\sum_{m\ge0}(2^md)^{-2}(2^md/t)^s
 \le Ct^{2-2s}d^{s-2}
 =C\frac{t^{1-s}}d\left(\frac td\right)^{1-s}
 \le C\frac{t^{1-s}}d.
\]
This proves \eqref{rec:absolute-band-kernel}, also for every
subfamily of a band.

For \(d<r\), sum the retained bands using their bound near the diagonal.
For \(d\ge r\), use \eqref{eq:projection-tail-identity} and sum
the omitted bands using their bound away from the diagonal. The result is
\begin{equation}
 |J_r(x,y)|\le C\Phi_r(x,y),\qquad
 \Phi_r(x,y)=r^{-s}
       \left(1+\frac{|x-y|}{r}\right)^{-1}.
 \label{eq:projection-kernel-envelope}
\end{equation}
The estimate includes \(x=y\). Upper regularity gives
\begin{equation}
 \sup_x\int_E\Phi_r(x,y)\,d\mu(y)
 +\sup_y\int_E\Phi_r(x,y)\,d\mu(x)
 \le C\left(1+\sum_{m\ge0}2^{-m(1-s)}\right)\le C.
 \label{var:envelope-integrals}
\end{equation}
The row and column bounds imply
\[
 \sup_{0<r\le1}\|P_r\|_{L^p(\mu)\to L^p(\mu)}\le C,
 \qquad1\le p\le\infty.
\]
Here and below \(P_1f=\langle f,1\rangle1\); the largest gap
has length less than one.
Keeping the local averages of \(|f|\) in the same annular sum gives
\begin{equation}
 P_*f(x):=\sup_{0<r\le1}|P_rf(x)|\le CM_\mu f(x).
 \label{eq:projection-maximal-domination}
\end{equation}
In particular, \(P_*\) is of weak type \((1,1)\) and is bounded
on \(L^p\), \(1<p\le\infty\).
Only finitely many gap lengths exceed each positive number, so
the family of distinct projections is countable and the supremum
is measurable.

The estimate \eqref{rec:holder-tail} proves uniform convergence
for Lipschitz functions. Polynomial restrictions are Lipschitz and
dense in \(C(E)\). For such an approximation \(h\) to \(f\),
\[
 \|P_rf-f\|_\infty
 \le (C+1)\|f-h\|_\infty+\|P_rh-h\|_\infty.
\]
Thus \(P_rf\to f\) uniformly for every \(f\in C(E)\).
Density and the uniform operator bounds give convergence in
\(L^p(\mu)\) for \(1\le p<\infty\).
For \(f\in L^1\), \eqref{eq:projection-maximal-domination}
and the maximal weak estimate yield
\[
 \mu\left\{\limsup_{r\downarrow0}|P_rf-f|>\tau\right\}
 \le \frac C\tau\|f-h\|_1
\]
for every Lipschitz \(h\): use
\(\limsup|P_rf-f|\le P_*(f-h)+|f-h|\).
Letting \(\|f-h\|_1\to0\) proves reconstruction almost everywhere.

Now enumerate the gaps as \(g_1,g_2,\ldots\) in nonincreasing
order of length, with any ordering among equal lengths, and set
\[
 Q_0f=\langle f,1\rangle1,\qquad
 Q_Nf=Q_0f+\sum_{j=1}^N\langle f,e_{g_j}\rangle e_{g_j}.
\]
Such an enumeration exists by the gap count. If
\(t/2\le\ell_{g_N}<t\), with \(t\) dyadic, then
\[
 Q_N=P_t+A_{\mathcal F},\qquad
 A_{\mathcal F}f=\sum_{g\in\mathcal F}\langle f,e_g\rangle e_g,
 \quad
 \mathcal F\subset\{g:t/2\le\ell_g<t\}.
\]
All gaps of length at least \(t\) occur before \(g_N\).
The absolute band estimate \eqref{rec:absolute-band-kernel},
with scale \(t/2\), gives both
\[
 \sup_N\|Q_N\|_{L^p\to L^p}\le C
 \quad(1\le p\le\infty),\qquad
 \sup_N|Q_Nf(x)|\le CM_\mu f(x).
\]
Every finite basis span is eventually fixed by \(Q_N\), and
these spans are dense in \(C(E)\) and \(L^1(\mu)\), as follows
from convergence of \(P_r\).
Thus \(Q_N\) converges to the identity in both endpoint norms.
Each coefficient functional is continuous on these spaces because
each \(e_g\) is bounded and continuous. Applying it to a convergent
expansion proves uniqueness. Hence the ordered family is a
Schauder basis of \(C(E)\) and \(L^1(\mu)\).
The same maximal approximation argument also gives \(Q_Nf\to f\)
almost everywhere for every \(f\in L^1\).

We prove the variation estimate for all ordinary partial sums:
\[
 V_\rho(Qf)
 =\sup_{0\le N_0<\cdots<N_M}
       \left(\sum_{i=1}^M
             |Q_{N_i}f-Q_{N_{i-1}}f|^\rho\right)^{1/\rho},
 \qquad \rho>2.
\]
The cutoff family is a subsequence, since \(P_r=Q_{n_r}\).
It is enough first to prove the \(L^2\) bound.
Put \(r_k=2^{-k}\), \(P_k=P_{r_k}\),
\(D_k=P_k-\mathbb E_k\), and
\(\Delta_j=\mathbb E_j-\mathbb E_{j-1}\) for \(j\ge1\).
The gap filtration gives
\begin{equation}
 f=\mathbb E_0f+\sum_{j\ge1}\Delta_jf\quad\hbox{in }L^2,
 \qquad
 \sum_{j\ge1}\|\Delta_jf\|_2^2
       =\|f-\mathbb E_0f\|_2^2.
 \label{var:martingale-resolution}
\end{equation}

We need the following continuity estimate:
\begin{equation}
 |J_r(x,y)-J_r(x,y')|
 \le C\frac{|y-y'|}{r}\Phi_r(x,y),
 \qquad |y-y'|\le r/2.
 \label{var:retained-lipschitz}
\end{equation}
Let \(\nu_r=\sum_{\ell_g\ge r}w_g\delta_{a_g}\) and
\(h=|y-y'|\le r/2\).
Centrality makes \(|y-a|\asymp|y'-a|\) on
\(\operatorname{supp}\nu_r\). Hence
\[
 |J_r(x,y)-J_r(x,y')|
 \le Ch\int\frac{d\nu_r(a)}{|x-a|\,|y-a|^2}.
\]
If \(d=|x-y|\le r\), both distances are comparable and
\eqref{eq:dual-far-three} bounds the integral by \(Cr^{-1-s}\).
If \(d>r\), split as in the band estimate. The regions near
\(x\), near \(y\), and away from both contribute respectively
\[
 \begin{aligned}
 Cd^{-2}\int_{|a-x|<d/3}|x-a|^{-1}\,d\nu(a)
      &\le Cd^{-1-s},\\
 Cd^{-1}\int_{|a-y|\ge c_*r}|y-a|^{-2}\,d\nu(a)
      &\le Cd^{-1}r^{-s},\\
 C\int_{|a-y|\ge d/3}|y-a|^{-3}\,d\nu(a)
      &\le Cd^{-1-s}.
 \end{aligned}
\]
Here we used \eqref{eq:dual-near-one},
\eqref{eq:dual-far-two}, and \eqref{eq:dual-far-three}.
The sum is at most \(Cr^{-s}/d\), proving
\eqref{var:retained-lipschitz}.

For \(j>k\), the function \(\Delta_jf\) has zero integral on
each \(Q\in\mathcal Q_{j-1}\), and
\(\mathbb E_k\Delta_jf=0\).
Choose \(z_Q\in Q\). Then
\[
 D_k\Delta_jf(x)
 =\sum_{Q\in\mathcal Q_{j-1}}
    \int_Q\bigl(J_{r_k}(x,y)-J_{r_k}(x,z_Q)\bigr)
                         \Delta_jf(y)\,d\mu(y).
\]
For \(j-k\) larger than a fixed geometric constant,
\(\diam Q\le2C_0r_j\le r_k/2\). By
\eqref{var:retained-lipschitz} and
\eqref{var:envelope-integrals}, both Schur norms of the
difference kernel are at most \(C2^{k-j}\).
The remaining finitely many differences use
\(\|P_k\|_{2\to2}=\|\mathbb E_k\|_{2\to2}=1\).
Thus
\[
 \|D_k\Delta_jf\|_2\le C2^{k-j}\|\Delta_jf\|_2,
 \qquad j>k.
\]
For \(j\le k\), the function \(v=\Delta_jf\) is constant on
the level-\(j\) cells and \(\mathbb E_kv=v\).
Since \(P_k1=1\),
\[
 D_kv(x)=\int_{\{y:Q_j(y)\ne Q_j(x)\}}
        J_{r_k}(x,y)\bigl(v(y)-v(x)\bigr)\,d\mu(y).
\]
Distinct level-\(j\) cells are separated by \(r_j\).
The envelope \eqref{eq:projection-kernel-envelope} gives
\[
 \sup_x\int_{|x-y|\ge r_j}|J_{r_k}(x,y)|\,d\mu(y)
 \le Cr_k^{1-s}\sum_{m\ge0}(2^mr_j)^{s-1}
 \le C2^{-(k-j)(1-s)},
\]
and the same column bound follows by symmetry.
The integral term and the multiplication term in \(D_kv\)
therefore obey
\[
 \|D_k\Delta_jf\|_2
 \le C2^{-(k-j)(1-s)}\|\Delta_jf\|_2,\qquad j\le k.
\]
Both cases are summarized by
\[
 \|D_k\Delta_jf\|_2\le Ca_{k-j}\|\Delta_jf\|_2,\qquad
 a_n=
 \begin{cases}
  2^{-n(1-s)},&n\ge0,\\
  2^n,&n<0,
 \end{cases}
 \quad a\in\ell^1(\mathbb Z).
\]
Since \(D_k1=0\), Minkowski's inequality, convolution on
\(\ell^2\), and \eqref{var:martingale-resolution} yield
\begin{equation}
 \left(\sum_{k\ge0}\|D_kf\|_2^2\right)^{1/2}
 \le C\|a\|_{\ell^1}
       \left(\sum_{j\ge1}\|\Delta_jf\|_2^2\right)^{1/2}
 \le C\|f\|_2.
 \label{var:square-comparison}
\end{equation}
First take finite martingale sums, then use \(L^2\) approximation
and Fatou.
For any scalar sequence, its \(\rho\)-variation is at most twice
its \(\ell^2\) norm when \(\rho\ge2\).
L\'epingle's inequality and \eqref{var:square-comparison}
therefore give
\begin{equation}
 \|V_\rho(P_kf:k\ge0)\|_2\le C_\rho\|f\|_2.
 \label{var:dyadic-variation}
\end{equation}

For \(t=r_k\), let
\[
 \mathcal B_k=\{g:t/2\le\ell_g<t\},\qquad
 c_g=\langle f,e_g\rangle,\qquad
 S_kf(x)=\sum_{g\in\mathcal B_k}|c_g|\,|e_g(x)|.
\]
The localization estimate \eqref{sp:size} and the gap count give
\[
 \sum_{g\in\mathcal B_k}t^{s/2}|e_g(x)|
 \le C\sum_{m\ge0}2^{-m(1-s)}\le C,
 \qquad
 \|e_g\|_1\le Ct^{s/2}\quad(g\in\mathcal B_k).
\]
Thus Cauchy--Schwarz gives
\[
 |S_kf(x)|^2
 \le
 \left(\sum_{g\in\mathcal B_k}
       |c_g|^2t^{-s/2}|e_g(x)|\right)
 \left(\sum_{g\in\mathcal B_k}t^{s/2}|e_g(x)|\right)
 \le C\sum_{g\in\mathcal B_k}
       |c_g|^2t^{-s/2}|e_g(x)|.
\]
After integration and summation,
\begin{equation}
 \left\|\left(\sum_{k\ge0}|S_kf|^2\right)^{1/2}\right\|_2^2
 \le C\sum_g|c_g|^2\le C\|f\|_2^2.
 \label{var:band-square-l2}
\end{equation}

Put \(N_k=n_{r_k}\), so \(Q_{N_k}=P_k\).
Between \(N_k\) and \(N_{k+1}\), only gaps of \(\mathcal B_k\)
are added, one at a time. Any successive increments use disjoint
subfamilies of that band; hence
\[
 V_\rho(Q_Nf:N_k\le N\le N_{k+1})\le S_kf.
\]
This statement holds for every ordering among gaps of equal length.
The elementary long and short variation decomposition gives
\[
 V_\rho(Qf)\le C_\rho\left(
       V_\rho(P_kf:k\ge0)
       +\left(\sum_{k\ge0}|S_kf|^2\right)^{1/2}\right).
\]
Indeed, insert the relevant dyadic block boundaries into each
increment of a finite increasing sequence of indices.
The resulting middle increments have ordered, disjoint interiors
and are controlled by dyadic variation. Inside any block, the
remaining end pieces and original internal increments are controlled
by its short variation. Apply
\[
 |a+b+c|^\rho
 \le3^{\rho-1}(|a|^\rho+|b|^\rho+|c|^\rho)
\]
and use \(\ell^2\subset\ell^\rho\).
Empty blocks cause no change. Combining
\eqref{var:dyadic-variation} and \eqref{var:band-square-l2} proves
\begin{equation}
 \|V_\rho(Qf)\|_2\le C_\rho\|f\|_2.
 \label{rec:ordered-variation-l2}
\end{equation}
This is the comparison with martingales used in
\cite[Lemma~3.2 and Section~4]{JonesSeegerWright2008},
with the gap filtration and the absolute band bounds as its
geometric estimates.

The passage to all exponents uses the absolute kernel estimates
already established in \eqref{eq:absolute-kernel-size} and
\eqref{eq:absolute-kernel-smoothness}.
For \(M\ge1\), let \(B_M\) be the finite-dimensional space of
complex sequences \(v=(v_0,\ldots,v_M)\) with \(v_0=0\), equipped
with the norm
\[
 \|v\|_{B_M}
 =\sup_{0\le N_0<\cdots<N_j\le M}
       \left(\sum_{i=1}^j
            |v_{N_i}-v_{N_{i-1}}|^\rho\right)^{1/\rho}.
\]
Requiring the first coordinate to be zero makes this a norm.
Define the linear operator
\[
 T_Mf=(Q_Nf-Q_0f)_{N=0}^M,\qquad
 K_M(x,y)=
       \left(\sum_{j=1}^Ne_{g_j}(x)e_{g_j}(y)\right)_{N=0}^M.
\]
By \eqref{rec:ordered-variation-l2},
\(\|T_M\|_{L^2\to L^2(B_M)}\le C_\rho\), uniformly in \(M\).
For fixed \(j\), the coordinate sequence
\((\mathbf1_{\{j\le N\}})_{N=0}^M\) has variation one.
The triangle inequality in \(B_M\) therefore gives, for
\(h=|x-x'|\le|x-y|/2\),
\[
 \begin{aligned}
 \|K_M(x,y)\|_{B_M}
 &\le\sum_g|e_g(x)e_g(y)|
 \le C|x-y|^{-s},\\
 \|K_M(x,y)-K_M(x',y)\|_{B_M}
 &\le\sum_g|e_g(x)-e_g(x')|\,|e_g(y)|
 \le C\frac{h^{1-s}}{|x-y|}.
 \end{aligned}
\]
The same continuity estimate holds in the second variable.
These constants are independent of \(M\) and of the chosen
ordering within a band.

The finite-dimensional Banach version of
\eqref{eq:prelim-cz-conclusion} now gives weak type \((1,1)\)
and the strong estimates for \(1<p\le2\), uniformly in \(M\).
For the duality step, finite dimensionality identifies
\(L^2(\mu;B_M)^*\) with \(L^2(\mu;B_M^*)\).
Thus
\[
 T_M^*:L^2(\mu;B_M^*)\longrightarrow L^2(\mu)
\]
has the same \(L^2\) bound. Its kernel is \(K_M(y,x)^*\),
which satisfies the same estimates in operator norm.
The same singular integral bound for the adjoint and duality give
\(2<p<\infty\). The constants are independent of the dimensions
of these sequence spaces.

Since each \(Q_N\) has a bounded finite kernel, approximation
in \(L^1\) makes every fixed coordinate converge pointwise.
Fatou consequently extends the finite-\(M\) weak estimate from
bounded functions to all \(L^1\).
As \(M\to\infty\), the norms \(\|T_Mf(x)\|_{B_M}\) increase
to \(V_\rho(Qf)(x)\). Monotone convergence gives
\begin{equation}
 \begin{gathered}
  \sup_{\tau>0}\tau\mu\{V_\rho(Qf)>\tau\}
      \le C_\rho\|f\|_1,\\
  \|V_\rho(Qf)\|_p\le C_{p,\rho}\|f\|_p,
       \qquad1<p<\infty,\quad \rho>2.
 \end{gathered}
 \label{var:full-variation}
\end{equation}
The estimates \eqref{var:full-variation} also hold for
\(V_\rho(Pf)\), since \(P_r=Q_{n_r}\).
Both variations are suprema over countable families.

We next identify the sharp approximation rates.
The gap count and orthonormality give
\[
 \dim\mathcal V_r=1+n_r\le1+Cr^{-s}.
\]
For \(v\in\mathcal V_r\), one has \(P_rv=v\); hence
\[
 E_r(f)\le\|f-P_rf\|_\infty
       =\|(I-P_r)(f-v)\|_\infty
       \le(1+C)\|f-v\|_\infty.
\]
Taking the infimum in \(v\) proves
\[
 E_r(f)\asymp\|f-P_rf\|_\infty,\qquad f\in C(E).
\]
The upper H\"older rate is therefore \eqref{rec:holder-tail}.

To obtain matching lower bounds, we find an omitted gap near zero
at every sufficiently small scale.
Endpoint regularity gives, after enlarging \(A\) if necessary,
\[
 A^{-1}t^s\le\mu([0,t])\le At^s,\qquad0<t\le1.
\]
Choose \(\vartheta\in(0,1)\) with
\(A\vartheta^s\le(2A)^{-1}\). Then
\(\mu([\vartheta r,r])\ge(2A)^{-1}r^s>0\).
Choose \(t\in E\cap[\vartheta r,r]\) and apply the porosity
estimate from the construction of the gap filtration to \([0,t]\).
It provides a gap \(g_r\subset(0,t)\) with
\[
 cr\le\ell_{g_r}<r,\qquad
 0<a_{g_r}<r,\qquad
 \frac{w_{g_r}}{a_{g_r}}\ge cr^{1-s}.
\]
Equality \(\ell_{g_r}=r\) would force \(g_r=(0,r)\),
isolating \(0\), which is impossible.

Choose \(L>2\) so that \(L^s\ge2A^2\).
Whenever \(L\rho\le1\),
\begin{equation}
 \mu([\rho,L\rho])
 \ge(A^{-1}L^s-A)\rho^s
 \ge(2A)^{-1}(L\rho)^s.
 \label{rec:endpoint-annulus}
\end{equation}
For increasing \(f\) on \([0,1]\), cancellation gives
\[
 \langle f,e_g\rangle
 =\sqrt{w_g}\int_E
       \frac{f(t)-f(a_g)}{t-a_g}\,d\mu(t)\ge0,
 \qquad e_g(0)=-\frac{\sqrt{w_g}}{a_g}<0.
\]
For the H\"older functions below, absolute convergence yields
\[
 |f(0)-P_rf(0)|
 =\sum_{\ell_g<r}\frac{w_g}{a_g}
       \int_E\frac{f(t)-f(a_g)}{t-a_g}\,d\mu(t).
\]
Every summand is nonnegative.

If \(0<\alpha<1-s\), take \(f(t)=t^\alpha\), whose
\(C^\alpha\) seminorm is at most one.
Writing \(a=a_{g_r}<r\), the inequality
\[
 \frac{t^\alpha-a^\alpha}{t-a}
 \ge(1-2^{-\alpha})t^{\alpha-1},\qquad t\ge2a,
\]
and \eqref{rec:endpoint-annulus} on \([2r,2Lr]\) give
\[
 \int_E\frac{t^\alpha-a^\alpha}{t-a}\,d\mu(t)
 \ge c_\alpha(2Lr)^{\alpha-1}\mu([2r,2Lr])
 \ge c_{\alpha,s,A}r^{\alpha+s-1}.
\]
The term \(g_r\) alone is therefore at least \(cr^\alpha\).

At \(\alpha=1-s\), use the same function and the disjoint
annuli
\[
 [2L^jr,2L^{j+1}r],\qquad j\ge0,\quad2L^{j+1}r\le1.
\]
Their endpoints have zero measure, and each contributes at least
a positive constant to the integral, because
\[
 \int_{[\rho,L\rho]}t^{-s}\,d\mu(t)
 \ge(L\rho)^{-s}\mu([\rho,L\rho])\ge(2A)^{-1}.
\]
There are at least \(c_L\log(1/r)\) such annuli for small \(r\).
Multiplication by \(w_{g_r}/a_{g_r}\) gives
\(cr^{1-s}\log(e/r)\).

For \(1-s<\alpha\le1\), take \(f(t)=t\).
Its \(C^\alpha\) seminorm on \([0,1]\) is at most one, and
\[
 \langle t,e_g\rangle
 =\sqrt{w_g}\int_E\frac{t-a_g}{t-a_g}\,d\mu(t)
 =\sqrt{w_g}.
\]
The selected tail term is \(w_{g_r}/a_{g_r}\ge cr^{1-s}\).
These fixed functions prove all three lower rates.
Together with bounded projection, this establishes
\[
 \sup_{[f]_{C^\alpha(E)}\le1}E_r(f)
 \asymp_\alpha\Psi_{\alpha,s}(r)
\]
for every sufficiently small \(r\).

This sharpness concerns the prescribed rational spaces.
Indeed, for \(r=2^{-k}\), the space
\(\mathcal W_k=\operatorname{span}\{\mathbf1_Q:Q\in\mathcal Q_k\}\)
lies in \(C(E)\), since the cells are relatively open and closed.
Its dimension is at most \(1+Cr^{-s}\), and
\[
 |f(x)-\mathbb E_kf(x)|
 \le\frac1{\mu(Q)}\int_Q|f(x)-f(y)|\,d\mu(y)
 \le C_0^\alpha r^\alpha[f]_{C^\alpha(E)},\qquad x\in Q.
\]
Thus spaces of comparable dimension can avoid the saturation
and the critical logarithm.

Finally, let \(S\subset E\) be Borel with \(\mu(S)\ge\theta>0\),
and let \(v\in\mathcal V_r\).
Write \(M=\esssup_S|v|\), and choose a compact
\[
 K\subset\{x\in S:|v(x)|\le M\},\qquad \mu(K)\ge\theta/2.
\]
For \(\sigma=\mu|_K/\mu(K)\), upper regularity gives
\[
 \sigma(B(u,\rho))\le\min(1,C\theta^{-1}\rho^s),
 \qquad u\in[0,1].
\]
Tonelli's theorem therefore yields
\[
 \begin{aligned}
 \int_K\log\frac1{|t-u|}\,d\sigma(t)
 &=\int_0^\infty\sigma\{|t-u|<e^{-z}\}\,dz\\
 &\le\int_0^\infty\min(1,C\theta^{-1}e^{-sz})\,dz\\
 &=\frac{\log(C/\theta)+1}{s}=:L_\theta .
 \end{aligned}
\]
All logarithmic singularities are integrable, uniformly in \(u\).
For a polynomial \(p\) of degree at most \(n\ge1\), choose
\(t_0,\ldots,t_n\in K\) maximizing the Vandermonde product.
The points are distinct because \(K\) has positive nonatomic mass.
Varying \(t_i\), taking logarithms, and integrating against \(\sigma\)
gives
\[
 \log\prod_{j\ne i}|t_i-t_j|
 \ge\sum_{j\ne i}\int_K\log|t-t_j|\,d\sigma(t)
 \ge-nL_\theta.
\]
Lagrange interpolation, with numerator factors bounded by one on
\([0,1]\), implies
\begin{equation}
 \|p\|_{[0,1]}\le(n+1)e^{nL_\theta}\sup_K|p|.
 \label{rec:polynomial-norming}
\end{equation}
This is the polynomial norming estimate; compare
\cite{Brudnyi2002,BrudnyiBrudnyi2006}.

Put
\[
 B_r(x)=\prod_{\ell_g\ge r}(x-a_g).
\]
Then \(p=vB_r\) is a polynomial of degree at most \(n_r\),
and \(|B_r|\le1\) on \([0,1]\).
For fixed \(x\in E\), let
\[
 N_x(t)=\#\{g:\ell_g\ge r,\ |x-a_g|<t\}.
\]
This count vanishes for \(t<c_*r\).
For larger \(t\), every counted gap has
\(\ell_g\le t/c_*\), so its left endpoint lies in \(B(x,Ct)\).
Those endpoints are \(r\)-separated. Disjoint balls of radius
\(r/8\) and regularity consequently give
\[
 N_x(t)\le C(t/r)^s,\qquad c_*r\le t\le1.
\]
Integration of this count yields
\[
 -\log|B_r(x)|
 =\int_0^1\frac{N_x(t)}t\,dt
 \le Cr^{-s}\int_{c_*r}^1t^{s-1}\,dt
 \le Cr^{-s}.
\]
There is therefore no additional \(\log(1/r)\) in the
denominator estimate.
If \(n_r=0\), the norming assertion is immediate.
Otherwise, \eqref{rec:polynomial-norming}, the bound
\(n_r\le Cr^{-s}\), and \(\sup_K|p|\le M\) give
\[
 \begin{aligned}
 \|v\|_{C(E)}
 &\le e^{Cr^{-s}}\|p\|_{[0,1]}\\
 &\le e^{Cr^{-s}}(n_r+1)e^{n_rL_\theta}M\\
 &\le\exp\!\left(Cr^{-s}\log\frac C\theta\right)
             \esssup_S|v|.
 \end{aligned}
\]
This proves the Remez consequence and completes the argument.
\end{proof}

\section{Further questions}
\label{sec:future}

\subsection{Geometric conditions for the lower basis bound}
For admissible complex poles, consider
\[
 \mathcal A_b:\C\oplus\ell^2(\mathcal G)\longrightarrow L^2(\mu),
 \qquad
 \mathcal A_b(u,c)=u1+\sum_gc_gq_g.
\]
Theorem~\ref{thm:stability} characterizes boundedness by finite
packing and proves invertibility for sufficiently small tail
packing. Corollary~\ref{cor:dyadic-prescribed-excess} shows that
finite packing can also give any finite excess or infinite excess.
Theorem~\ref{prop:exact-chain-index} determines the exact basis,
frame, and Fredholm regimes for finitely many convergent endpoint
chains at each fixed contraction ratio of the thin Cantor set.

The corresponding classification on a general regular set remains
open. A geometric criterion must distinguish completeness from a
positive lower synthesis bound and must account for persistent
interactions across scales. Even on a fixed Cantor set, the exact
criterion for an arbitrary bounded Carleson array is not supplied
by the theorem for finitely many chains. The uniform paraproduct comparison
and its reverse H\"older sufficient condition give partial information.
The description in terms of model spaces in Appendix~\ref{app:modelspaces}
connects the remaining questions with interpolation and uniqueness
for reproducing kernels.

\subsection{Geometric criteria for general deformations}
For general $C^1$ maps satisfying a lower chord bound,
Proposition~\ref{prop:curve-local-oscillation} gives a local
Carleson sufficient condition for boundedness and compactness.
The square Dini condition is sharp as a uniform hypothesis on a
concave derivative modulus, by
Proposition~\ref{prop:sharp-modulus-transport}. An exact local
criterion for an individual deformation remains open. It would
need to account for cancellations that the positive oscillation
profile does not record.

For maps sufficiently close to the identity in $C^1$ that fix the
support, Proposition~\ref{prop:c1-support-fixing}
already identifies the exact displacement packing. When the
support changes, a necessary and sufficient condition must retain
both the measured support and the positions and tangents at the
parameter poles. The theorem for two weights in
\cite{LiShenSawyerWick2026} gives a related comparison for fractional Riesz transforms on Dini curves, with a different operator and
testing structure. Injectivity beyond the smallness conditions and
the exact M\"obius invariance is a separate problem: compactness
alone yields Fredholm information, not a basis theorem.

\subsection{Trace ideals and spectral asymptotics}
Theorem~\ref{thm:geometric-schatten} settles the Schatten criterion
for the rational columns considered here at every positive
exponent. Corollary~\ref{cor:c1-joint-schatten} gives its exact form in terms of position and tangent for deformations fixing the support.
For general deformations, a local Schatten criterion must also
include the change of support. Embeddings of model spaces
\cite{Baranov2009} and metric operator ideals
\cite{Hytonen2026} provide comparisons for this question.

The noncollinear examples parametrized by arc length in
Proposition~\ref{prop:arclength-sharpness} show that the uniform
approximation power $(m+\alpha)/s$ cannot be improved. More
structured measures or curves may admit spectral asymptotics with
a leading constant. The upper and lower bounds proved here do not give
such an asymptotic formula. Any proposed criterion must respect
M\"obius covariance modulo rank one and the dependence on
$\phi(b_g)$ and $\phi'(b_g)$, which are evaluated outside the support.

\appendix
\section{Comparison with model spaces}
\label{app:modelspaces}
Assume throughout this appendix that the parameter poles $b_g$ are real.
We give the normalizations for the comparison with model spaces.
Put
\[
 \Theta(z)=\frac{F(z)-i}{F(z)+i},\qquad z\in\C_+.
\]
Since $\Im F>0$, one has $|\Theta|<1$ in $\C_+$. On
$\R\setminus E$, the function $F$ is real and hence $|\Theta|=1$.
The identity $|E|=0$ therefore makes $\Theta$ an inner function.
Use the Hardy norm
$\|G\|_{H^2}^2=\sup_{y>0}\int_\R|G(x+iy)|^2\,dx$, and let
$K_\Theta=H^2(\C_+)\ominus\Theta H^2(\C_+)$.
The reproducing kernel obtained by projecting the Hardy kernel is
\[
 k_\zeta(z)=\frac{i}{2\pi}
       \frac{1-\overline{\Theta(\zeta)}\Theta(z)}{z-\overline\zeta},
 \qquad \zeta,z\in\C_+.
\]
In this normalization, the Clark unitary has the concrete form
\begin{equation}
 (\mathcal W f)(z)
 =\frac{\displaystyle\int_E f(t)(t-z)^{-1}\,d\mu(t)}
        {\sqrt\pi\,(F(z)+i)}.
 \label{ms:clark-unitary}
\end{equation}
This is the usual Clark representation; see \cite{Clark1972}.
We verify the constants directly.
For $\zeta\in\C_+$ define
\[
 v_\zeta(t)=
 \frac1{\sqrt\pi\,(\overline{F(\zeta)}-i)(t-\overline\zeta)}.
\]
The inner product convention gives
$\langle f,v_\zeta\rangle=(\mathcal Wf)(\zeta)$.
Partial fractions give
\[
 \int_E\frac{d\mu(t)}{(t-\overline\zeta)(t-z)}
 =\frac{F(z)-\overline{F(\zeta)}}{z-\overline\zeta}.
\]
Moreover,
\[
 1-\overline{\Theta(\zeta)}\Theta(z)
 =\frac{2i(\overline{F(\zeta)}-F(z))}
        {(\overline{F(\zeta)}-i)(F(z)+i)}.
\]
Substitution proves
\[
 \mathcal Wv_\zeta=k_\zeta,
 \qquad
 \langle v_\zeta,v_\omega\rangle
 =k_\zeta(\omega)=\langle k_\zeta,k_\omega\rangle.
\]
The kernels span $K_\Theta$. The functions $v_\zeta$ span
$L^2(\mu)$: a vector orthogonal to all of them has identically zero
Cauchy transform, hence all its moments vanish. Polynomial density
then makes that vector zero.
Thus the displayed Gram identity proves that $\mathcal W$ is unitary.
For a real $b$ inside a gap, passage to the boundary is permitted
because $F$ is analytic there. We obtain
\begin{equation}
 \mathcal W^{-1}k_b(t)
 =\frac1{\sqrt\pi\,(F(b)-i)(t-b)},
 \qquad
 \|k_b\|_{H^2}^2
 =\frac{F'(b)}{\pi(1+F(b)^2)}
 =\frac{|\Theta'(b)|}{2\pi}.
 \label{ms:boundary-kernel}
\end{equation}
For $\kappa_b=k_b/\|k_b\|_{H^2}$ this becomes
\[
 \mathcal W^{-1}\kappa_b(t)
 =\frac{F(b)+i}{\sqrt{1+F(b)^2}}
     \frac1{\sqrt{F'(b)}(t-b)}.
\]
The first factor has modulus one. Consequently,
\[
 \mathcal Wq_g
 =\sqrt\pi\sqrt{w_g}(F(b_g)-i)k_{b_g},
 \qquad
 \|q_g\|_2^2=w_gF'(b_g)\asymp1.
\]
The perturbed functions are therefore normalized reproducing kernels
up to factors of modulus one and positive factors bounded above and
away from zero.
The constant vector must be retained in this identification:
\[
 \mathcal W1=\frac{F}{\sqrt\pi(F+i)}
            =\frac{1+\Theta}{2\sqrt\pi},
 \qquad
 \mathcal We_g=-i\kappa_{a_g}.
\]
By Proposition~\ref{thm:basis}, these vectors together form an
orthonormal basis. The additional vector is the Clark vector at
infinity; omitting it would omit a one-dimensional summand.
For $G=\mathcal W f$, evaluation of \eqref{ms:clark-unitary} at
$b_g$ yields the exact identity
\[
 \sum_g|\langle f,q_g\rangle|^2
 =\pi\sum_gw_g(1+F(b_g)^2)|G(b_g)|^2.
\]
Thus the Bessel condition is the embedding
$K_\Theta\hookrightarrow L^2(\rho)$ for
\begin{equation}
 \rho=\pi\sum_gw_g(1+F(b_g)^2)\delta_{b_g}.
 \label{ms:embedding-measure}
\end{equation}
Write $\delta_g=b_g-a_g$. Every point of the segment joining
$a_g$ to $b_g$ remains a fixed multiple of $\ell_g$ from $E$.
The denominator comparison therefore gives
\[
 F'(a_g+u\delta_g)\asymp F'(a_g)=w_g^{-1},
 \qquad 0\le u\le1.
\]
Since $F(a_g)=0$ and $F'$ is positive on the gap,
\[
 |F(b_g)|
 =|\delta_g|\int_0^1F'(a_g+u\delta_g)\,du
 \asymp\frac{|\delta_g|}{w_g}.
\]
The atomic weights in \eqref{ms:embedding-measure} consequently satisfy
\[
\pi w_g(1+F(b_g)^2)
 \asymp w_g+\frac{\delta_g^2}{w_g}
 \asymp\ell_g^{2-s}+\eta_g^2\ell_g^s.
\]
This identifies the geometric perturbation weight inside the
embedding problem for the model space.
The embedding norm is also determined quantitatively. If
$V_b(c_g)=\sum_gc_gq_g$ and $U(c_g)=\sum_gc_ge_g$, then
$V_b=U+D$. Since $\|U\|=1$ and $\|q_g\|_2\asymp1$,
\[
 \|V_b\|\le1+\|D\|,\qquad
 \|D\|\le1+\|V_b\|,\qquad \|V_b\|\ge c.
\]
Together with Theorem~\ref{thm:stability}, this proves
\[
 \|I_\rho:K_\Theta\to L^2(\rho)\|^2
 =\|V_b\|^2\asymp1+\mathcal C(\lambda).
\]
This is an explicit evaluation for the measure
\eqref{ms:embedding-measure}, within the general embedding problem
treated in \cite{LaceySawyerShenUriarteWick2018}.
The statement about the essential norm in Theorem~\ref{thm:stability} concerns
the difference $D=V_b-U$. The operator $V_b$ itself is never
compact: the coordinate vectors tend weakly to zero in $\ell^2$,
whereas their images $q_g$ have norms bounded below. By taking
adjoints, the corresponding embedding $I_\rho$ is not compact
either. This distinction persists when the perturbation $D$ is
compact.

The following computation shows that the hypothesis of comparable
kernel norms in \cite[Theorem~1.1]{Baranov2005} and
\cite[Theorem~7.2]{BaranovFricainMashreghi2010} need not hold here.
Equation~\eqref{ms:boundary-kernel} gives
\[
 \frac{\|k_{b_g}\|_{H^2}}{\|k_{a_g}\|_{H^2}}
 =\left(\frac{w_gF'(b_g)}{1+F(b_g)^2}\right)^{1/2}
 \asymp
 \left(1+\frac{\delta_g^2}{w_g^2}\right)^{-1/2}.
\]
On the same sparse family of gaps, now choose real displacements
$|\delta_g|=\eta\ell_g$ with fixed $\eta>0$. Then
\[
 \frac{|\delta_g|}{w_g}\asymp\eta\ell_g^{s-1}\longrightarrow\infty,
 \qquad
 \frac{\|k_{b_g}\|_{H^2}}{\|k_{a_g}\|_{H^2}}
 \longrightarrow0.
\]
Hence no uniform positive lower bound in that comparability
hypothesis can hold for these displacements.
The vector at infinity does not affect this obstruction.
The automorphism $\phi(z)=-1/(z-2)$ sends infinity to zero and
satisfies
\[
 \tfrac14\le|\phi'(b)|\le1,\qquad 0\le b\le1.
\]
For $\widehat\Theta=\Theta\circ\phi^{-1}$, the norm formula for the
boundary kernel gives
\[
 \|\widehat k_{\phi(b)}\|_{H^2}^2
 =\frac{|\Theta'(b)|}{2\pi|\phi'(b)|}.
\]
The norm ratios at the finite gap poles change only by bounded factors, and
still tend to zero after moving the additional node to the real line.

A fixed change of Clark parameter has a smaller displacement
scale than the examples of Corollary~\ref{cor:dyadic-prescribed-excess}.
For a fixed real number $t$, let $b_g(t)$ solve $F(b_g(t))=t$
in the gap. On a fixed central portion of every sufficiently
small gap, $F'\asymp w_g^{-1}$, and the values of $F$ at the
two ends of this portion tend to opposite infinities. The mean
value theorem therefore gives
\[
 |b_g(t)-a_g|\asymp |t|w_g
       \asymp |t|\ell_g^{2-s}\quad(t\ne0),
 \qquad
 \frac{|b_g(t)-a_g|}{\ell_g}\lesssim_t\ell_g^{1-s}.
\]
The local gap count \eqref{eq:prelim-gapcount}, summed over
dyadic scales below $r$, gives
$\mathcal C_r(\lambda)\lesssim_t r^{2(1-s)}$.
Thus these tails have vanishing packing. A constant nonzero
fractional displacement along an infinite chain does not arise
from one fixed finite Clark parameter.

The inner function is outside the one-component class. Indeed, its
Clark measure corresponding to $\mu$ is nonatomic, whereas the
one-component condition forces discrete Clark measures; see
\cite{Bessonov2015}. Related constructions using Ahlfors regular
measures and Cantor gaps appear in
\cite{BellavitaDellepianeLuo2026}. This qualitative property is
established background for the quantitative comparisons above.

\subsection*{Clark dilations and the original gap family}
Redundant kernel frames obtained by enlarging a function space and
projecting an orthonormal kernel basis are classical; see
Ortega-Cerd\`a and Seip~\cite[pp.~797--798]{OrtegaCerdaSeip2002} and
al-Sa'di and Weber~\cite[Theorem~3]{AlSadiWeber2017}.
The decomposition
\[
 K_{\Theta B}=K_\Theta\oplus\Theta K_B
\]
shows that projecting a complete orthonormal Clark family from
$K_{\Theta B}$ onto $K_\Theta$ gives a Parseval frame with excess
$\dim K_B$. The original gap geometry imposes an additional restriction.

\begin{proposition}[Obstruction to a full Clark projection]
\label{prop:clark-projection-obstruction}
Let $E,\mu$, and the admissible real poles $b_g$ be as above.
Suppose that $\{1\}\cup(q_g)$ is a frame in $L^2(\mu)$.
Let $B$ be an inner function in $\C_+$. Assume that a complete
orthonormal Clark family for $K_{\Theta B}$ consists of the node
at infinity and exactly one node $b_g$ in each original gap.
Under orthogonal projection onto $K_\Theta$ and the original
unitary $\mathcal W^{-1}$, assume the infinity coordinate is a
nonzero multiple of $1$. Then the following assertions hold.
\begin{enumerate}
\item If $B$ extends holomorphically through every point of $E$,
then $B$ is constant and the frame has zero excess.
\item If the frame has finite excess, then $B$ is constant and
its excess is zero.
\end{enumerate}
In particular, the families with positive finite excess of
Corollary~\ref{cor:dyadic-prescribed-excess} admit no such full
Clark projection for any inner factor $B$.
A realization of this form with infinite excess would require $B$ to
fail to extend holomorphically through at least one point of $E$.
\end{proposition}

\begin{proof}
Let $|\gamma|=1$ be the common Clark parameter. At every gap node,
the boundary equation gives
\[
 (\Theta B)(b_g)=\gamma,
 \qquad
 F(b_g)=Q(b_g),\qquad
 Q=i\frac{B+\gamma}{B-\gamma}.
\]
The exceptional constant $B\equiv\gamma$ has no such finite
nodes outside the support and is excluded by the hypothesis.

For unnormalized boundary kernels, reproducing evaluation on
$K_\Theta$ gives
\[
 P_{K_\Theta}k_{b_g}^{\Theta B}=k_{b_g}^\Theta.
\]
This identity also holds for angular boundary kernels. Thus the
projected normalized Clark vectors, after applying $\mathcal W^{-1}$,
have the form $v_0=c_0'1$ and $v_g=c_g'q_g$, with all scalars
nonzero, by \eqref{ms:boundary-kernel} and the infinity hypothesis.

First suppose that $B$ extends holomorphically through every point
of $E$. Compactness gives a neighbourhood on which this extension
exists. The function $Q$ is meromorphic there and real on the
real line away from its poles. If $Q$ has no pole on $E$, it is
bounded in a smaller neighbourhood. The gap derivative estimate
then gives, for every sufficiently small gap,
\[
 |b_g-a_g|\lesssim\ell_g^{2-s},\qquad
 \frac{|b_g-a_g|}{\ell_g}\lesssim\ell_g^{1-s}.
\]
The local gap count \eqref{eq:prelim-gapcount}, summed over
dyadic scales, yields
\[
 \mathcal C_r(\lambda)\lesssim r^{2(1-s)}\longrightarrow0.
\]
Thus the synthesis of $\{1\}\cup(q_g)$ is a compact perturbation
of the canonical unitary synthesis, by Theorem~\ref{thm:stability}.
It has Fredholm index zero. Since it is surjective by the frame
hypothesis, its kernel is zero.

On the tail where $Q$ is bounded, the kernel norm formulas are
\[
 \|k_{b_g}^\Theta\|^2=\frac{|\Theta'(b_g)|}{2\pi},\qquad
 \|k_{b_g}^{\Theta B}\|^2
 =\frac{|\Theta'(b_g)|+|B'(b_g)|}{2\pi}.
\]
Here $B'$ is bounded near $E$, whereas
\[
 |\Theta'(b_g)|\asymp F'(b_g)\asymp\ell_g^{s-2},
 \qquad w_gF'(b_g)\asymp1.
\]
Consequently $|c_g'|$ is bounded above and below. The remaining
finitely many coordinates, including infinity, have the same
property after adjusting the constants. This bounded invertible
diagonal rescaling preserves the kernel dimension. The projected
Parseval frame therefore has zero excess, so $\dim K_B=0$ and
$B$ is constant.

It remains to exclude a pole $x_0\in E$ of $Q$.
Choose $\kappa\in(0,1)$ with $A\kappa^s<(2A)^{-1}$.
For every sufficiently small $R>0$, regularity gives
\[
 \mu\bigl(B(x_0,R)\setminus B(x_0,\kappa R)\bigr)
 \ge (2A)^{-1}R^s>0.
\]
Choose $y\in E$ in this annulus. The porosity estimate used in
the construction of the gap filtration supplies a gap $g$ between
$x_0$ and $y$ with
\[
 c_{\rm por}\kappa R\le\ell_g\le R.
\]
Since $b_g$ is uniformly central in $g$, one has
$|b_g-x_0|\asymp\ell_g$. As $R\downarrow0$, these gaps shrink
to $x_0$. Annular summation and the pole of $Q$ give
\[
 |F(b_g)|\le\int_E\frac{d\mu(t)}{|t-b_g|}
       \lesssim\ell_g^{s-1},
 \qquad
 |Q(b_g)|\gtrsim\ell_g^{-1}.
\]
These estimates contradict $F(b_g)=Q(b_g)$ because $s>0$.
This proves the first assertion for every regular support in
the manuscript, without a self-similarity assumption.

For the second assertion, $B$ need not extend through $E$.
Let $A_q$ be the synthesis of the full family consisting of the constant and the gap functions
and let $A_v$ be that of the projected normalized Clark family,
after applying $\mathcal W^{-1}$. Completeness in the larger model space gives
\[
 \dim\ker A_v=\dim K_B.
\]
The projection is contractive, while $\|q_g\|_2\asymp1$ and
$\|1\|_2=1$. Hence the nonzero scalar factors $c_j'$ are
uniformly bounded above. Their diagonal operator $D_c$ on
$\ell^2(\{0\}\cup\mathcal G)$ is bounded and injective, and
\[
 A_v=A_qD_c,\qquad
 D_c(\ker A_v)\subseteq\ker A_q.
\]
If $\dim\ker A_q<\infty$, this injection implies
$\dim K_B<\infty$. A finite-dimensional model space corresponds
to a finite Blaschke product, including a constant of degree zero.
Such a factor extends through $E$, so the first assertion applies.
Notice that this argument does not require a lower bound on
$|c_j'|$ before the finite-dimensional conclusion.
\end{proof}

The proposition concerns the complete projected Clark family with
exactly the original gap coordinates and the infinity coordinate.
Deleting additional Clark nodes changes the synthesis kernel and
falls outside its hypotheses. The proposition also makes no claim
about other dilations, changes to an equivalent norm, or a construction of infinite excess whose factor does not extend through $E$.

\section*{AI use declaration}
ChatGPT assisted with language editing.
The author takes full responsibility for the mathematical content,
correctness, and accuracy of the references.

\end{document}